\documentclass[a4paper, oneside]{amsart}

\usepackage{amsthm, amsmath, amssymb, amstext, amsfonts, bm}
\usepackage{a4wide}
\usepackage{hyperref}
\usepackage{enumitem}
\usepackage[svgnames]{xcolor}
\usepackage{romannum}
\AtBeginDocument{\pagenumbering{arabic}}
\usepackage{makecell}
\usepackage{tabularray}
\usepackage{multirow}
\usepackage{subfiles}
\usepackage{float}
\usepackage{booktabs} 
\usepackage{arydshln}
\definecolor{smoked}{RGB}{216, 212, 204}
\definecolor{mauve}{RGB}{200, 55, 171}
\definecolor{apricot}{RGB}{250, 144, 4}
\definecolor{sky}{RGB}{66, 169, 244}
\definecolor{plum}{RGB}{76, 0, 102}
\definecolor{darkgreen}{RGB}{0,128,0}

\definecolor{lightmauve}{RGB}{232, 173, 220}
\definecolor{lightapricot}{RGB}{253, 211, 155}
\definecolor{lightsky}{RGB}{178, 221, 251}
\definecolor{lightplum}{RGB}{184, 153, 192}

\hypersetup{
   colorlinks,
    linkcolor={mauve},
    citecolor={plum},
    urlcolor={plum}
}
\usepackage{graphicx}
\usepackage{tikz,tqft}
\usetikzlibrary{matrix,arrows,bending,cd,decorations.markings,calc,patterns,angles,quotes,backgrounds,shapes,math,hobby}
\tikzset{>=latex}
\usepackage{pgfplots}
    \pgfplotsset{compat=1.6}
\usepackage{caption}
\usepackage{subfig}

\theoremstyle{definition}
\newtheorem{defn}{Definition}[section]
\newtheorem{rem}[defn]{Remark}

\newtheorem{ex}[defn]{Example}
\newtheorem{notation}[defn]{Notation}

\theoremstyle{plain}
\newtheorem{thm}[defn]{Theorem}
\newtheorem{lem}[defn]{Lemma}
\newtheorem{fact}[defn]{Fact}
\newtheorem{claim}[defn]{Claim}
\newenvironment{proofclaim}
 {\proof}
 {\endproof}
\newtheorem{prop}[defn]{Proposition}
\newtheorem{thmx}{Theorem}

\newtheorem{cor}[defn]{Corollary}
\newtheorem{conj}[defn]{Conjecture}

\numberwithin{equation}{section}

\DeclareMathOperator{\id}{id}

\DeclareMathOperator{\Ok}{Ok}
\newcommand{\sslash}{\mathbin{/\mkern-3mu/}}

\DeclareMathOperator{\trace}{Tr}

\DeclareMathOperator{\R}{\mathbb{R}}
\DeclareMathOperator{\Q}{\mathbb{Q}}
\DeclareMathOperator{\C}{\mathbb{C}}
\DeclareMathOperator{\Z}{\mathbb{Z}}
\DeclareMathOperator{\N}{\mathbb{N}}

\DeclareMathOperator{\CP}{\mathbb{CP}}
\DeclareMathOperator{\HH}{\mathbb{H}}

\newcommand{\SL}{\operatorname{SL}}
\newcommand{\GL}{\operatorname{GL}}
\newcommand{\SU}{\operatorname{SU}}
\newcommand{\psl}{\operatorname{PSL}_2 \R}
\newcommand{\pgl}{\operatorname{PGL}_2 \R}

\DeclareMathOperator{\Hom}{Hom}
\DeclareMathOperator{\Out}{Out}
\DeclareMathOperator{\Aut}{Aut}
\DeclareMathOperator{\Rep}{Rep}

\DeclareMathOperator{\Teich}{Teich}
\newcommand{\RepDT}[1]{\Rep^\mathrm{DT}_{#1}}

\DeclareMathOperator{\PMod}{PMod}

\newcommand{\SpherePk}{\Sigma}

\newcommand{\cm}[1]{#1}
\newcommand{\sam}[1]{#1}

\title{Tykhyy's conjecture on finite mapping class group orbits}

\author{Samuel Bronstein}
\address[S.~Bronstein]{DMA, \'Ecole Normale Supérieure, PSL, 75005 Paris, France.}
\email{samuel.bronstein@ens.fr}

\author{Arnaud Maret}
\address[A.~Maret]{Sorbonne Université and Université Paris Cité, CNRS, IMJ-PRG, F-75005 Paris, France.}
\email{maret@imj-prg.fr}

\makeatletter
\newcommand{\addresseshere}{%
  \enddoc@text\let\enddoc@text\relax
}
\makeatother

\date{\today}

\begin{document}
\begin{abstract}
We classify the finite orbits of the mapping class group action on the character variety of
Deroin--Tholozan representations of punctured spheres. In particular, we prove that the action
has no finite orbits if the underlying sphere has 7 punctures or more. When the sphere has six
punctures, we show that there is a unique 1-parameter family of finite orbits. Our methods also
recover Tykhyy's classification of finite orbits for 5-punctured spheres. The proof is inductive
and uses Lisovyy--Tykhyy's classification of finite mapping class group orbits for 4-punctured
spheres as the base case for the induction.

Our results on Deroin--Tholozan representations cover the last missing cases to complete the
proof of Tykhyy's Conjecture on finite mapping class group orbits for $\SL_2\C$ representations
of punctured spheres, after the recent work by Lam--Landesman--Litt.
\end{abstract}


\maketitle

\section{Introduction}
\subsection{Results and overview}
Our goal is to pursue a long series of works to understand all the finite mapping class group
orbits of conjugacy classes of representations $\pi_1\Sigma\to\SL_2\C$ where $\Sigma$ is a sphere
with at least three punctures. This work initially focuses on \emph{Deroin--Tholozan} (DT)
representations (Section~\ref{sec:DT-representations}), a special kind of representations
introduced in~\cite{deroin-tholozan} that take values in $\psl$ and can be parametrized
by \emph{chains of hyperbolic triangles} (Definition~\ref{defn:triangle-chain}). Unlike Fuchsian representations, DT representations
are \emph{totally elliptic} \cm{(Lemma~\ref{lem:totally-elliptic})}, i.e.~all simple closed curves, \cm{including peripheral loops around punctures}, are sent to \cm{non-trivial} elliptic elements of $\psl$.
There are rare examples of DT representations with discrete image. These examples give rise to
finite orbits of the mapping class group action (Corollary~\ref{cor:DT-discrete-implies-finite-orbit}).
\cm{DT representations are discrete when the corresponding triangle chain fits well within
a Coxeter tessellation of the hyperbolic plane (in the sense of Example~\ref{ex:discrete-DT}), like on the illustration below.}
\begin{center}
\begin{tikzpicture}[font=\sffamily,decoration={markings, mark=at position 1 with {\arrow{>}}}]
\node[anchor=east] at (-0.5,0) {\includegraphics[width=6cm]{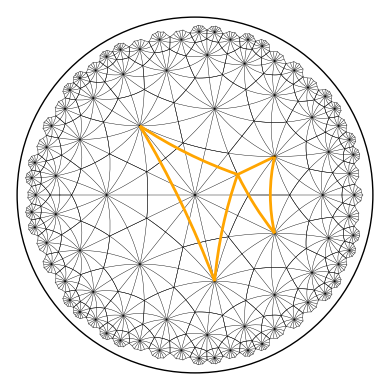}};
\node[anchor=west] at (0.5,0) {\includegraphics[width=6cm]{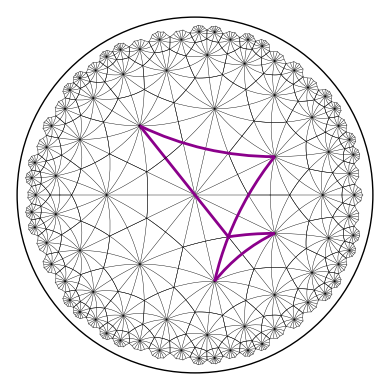}};
\end{tikzpicture}
\end{center}
Even though not all finite mapping class group orbits come from discrete representations, finite orbits remain fairly uncommon. The examples coming from DT representations are particularly interesting because the image of a DT representation is always infinite and Zariski dense.

\begin{thmx}[Theorem~\ref{thm:no-finite-orbit-for-n-geq-7}]\label{intro:tykhyy-DT}
For a punctured sphere~$\Sigma$ with $7$ punctures or more, all the mapping class group orbits of conjugacy classes of DT representations are infinite.
\end{thmx}

When $\Sigma$ is a sphere with $6$ punctures, there are examples of finite mapping class group orbits coming from DT representations. They occur for representations that map all six peripheral loops of $\Sigma$ into the same conjugacy class of elliptic elements inside $\psl$. Each such choice of elliptic conjugacy class leads to a unique finite orbit. These orbits are of \emph{pullback} type (Section~\ref{sec:pullback-orbits}), in the sense that they can be pulled back from a representation of a pair of pants via a family of ramified coverings.

\begin{thmx}[Theorems~\ref{thm:angle-vector-alpha-with-finite-orbits-n=6} \&~\ref{thm:existence-finite-orbit-n=6}]\label{intro:finite-orbit-n=6}
Let $\Sigma$ be a sphere with $6$ punctures and $\rho\colon\pi_1\Sigma\to\psl$ be a DT representation
whose conjugacy class belongs to a finite mapping class group orbit.
Then $\rho$ maps all the peripheral loops of $\Sigma$ to the same elliptic conjugacy class of $\psl$.

Conversely, for every \cm{$\theta\in (5\pi/3,2\pi)$}, there is a unique finite mapping class group orbit among the
DT representations mapping the peripheral loops of $\Sigma$ to elliptic elements of $\psl$ with the
same rotation angle $\theta$. This orbit consists of 40 points which are listed in
Table~\ref{tab:UFO-orbit}, Appendix~\ref{app:tables}.
\end{thmx}

We express orbit points in terms of their action-angle coordinates (Section~\ref{sec:action-angle-coordinates}), which were developed by the second author in~\cite{action-angle}. \cm{Equivalently, one may represent them via triangle chains in the hyperbolic plane (Section \ref{sec:triangle-chains}). For a DT representation $\rho\colon \pi_1\Sigma \to \mathrm{PSL}_2(\mathbb{R})$, the vertices and interior angles of the associated triangle chain encode the fixed points and rotation angles of the images under $\rho$ of a chosen generating set of $\pi_1\Sigma$, making it possible to reconstruct $\rho$ efficiently. The triangle chains associated to the 40 orbit points from Theorem \ref{intro:finite-orbit-n=6} each consist of four triangles; in one case, the configuration takes the shape of a ``jester's hat'' when drawn in the upper half-plane like on the illustration below. We shaded the four triangles of the chain in blue to aid visualization.}
\begin{center}
\begin{tikzpicture}    
\node[anchor=south west,inner sep=0] at (0,0) {\includegraphics[width=7cm]{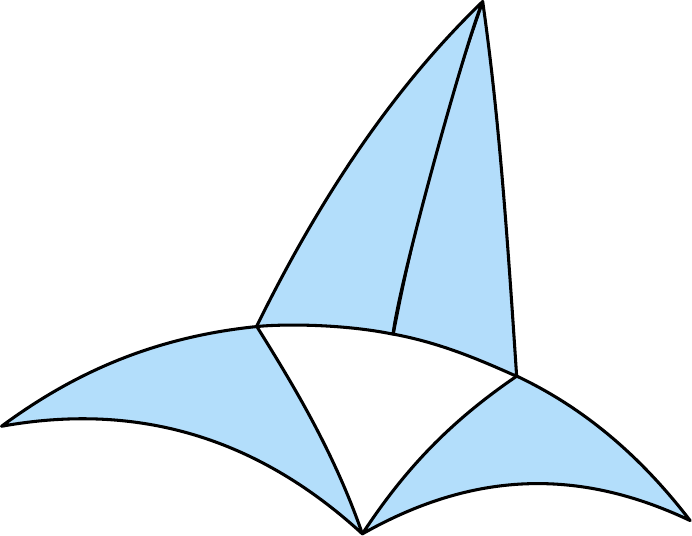}};
\end{tikzpicture}
\end{center}
For this reason, we'll refer to the finite mapping class group orbits of Theorem~\ref{intro:finite-orbit-n=6} as \emph{jester's hat orbits}. Jester's hat orbits already appeared in several places in the literature such as in Diarra's work on pullback orbits~\cite{diarra-pull-back}. The new contribution of Theorem~\ref{intro:finite-orbit-n=6} is that jester's hat orbits are the only possible kind of finite orbits for DT representations of $6$-punctured spheres.

\cm{Moving beyond DT representations, our next and main result completes the classification of all finite mapping class group orbits for $\SL_2\C$ representations of genus-0 surface groups. This problem has a long and rich history (Section~\ref{sec:history-and-motivations}), going back to the early $20^\textrm{th}$ century and the study of ``special'' (or algebraic) solutions to isomonodromy differential equations such as Painlevé~\Romannum{6}.} A classification of all finite mapping class group orbits \cm{for $\SL_2\C$ representations} has been completed in the case of 4-punctured spheres by Lisovyy--Tykhyy~\cite{LT} and \cm{a computer-aided classification was achieved by Tykhyy for 5-punctured spheres~\cite{tykhyy}.} \cm{Our methods use Lisovyy--Tykhyy's results, but are independent of Tykhyy's later work on 5-punctured spheres.} In particular, we recover Tykhyy's list of finite orbits for 5-punctured spheres in the special case of DT representations (Theorem~\ref{thm:angle-vector-alpha-with-finite-orbits-n=5}). Those finite orbits are of three different kinds, all of pullback type. The first two come as 1-parameter families. \cm{Motivated by the shape of triangle chains representing particular orbit points, all of which consist of three triangles,} we refer to them as \emph{hang-glider orbits} (Section~\ref{sec:proof-thm-finite-orbits-n=5-1}) and \emph{sand clock orbits} (Section~\ref{sec:proof-thm-finite-orbits-n=5-2}). They are of respective length $9$ and $12$ and their orbit points are listed in Tables~\ref{tab:hang-glider-orbit} \&~\ref{tab:sand-clock-orbit} in Appendix~\ref{app:tables}. The last orbit is rather exceptional as it only exists for a very particular choice of peripheral elliptic classes. We call it the \emph{bat orbit} (Section~\ref{sec:proof-thm-finite-orbits-n=5-3}). It is made of 105 orbit points, all listed in Table~\ref{tab:bat-orbit} in Appendix~\ref{app:tables}.
\begin{center}
\vspace{3mm}
\begin{tikzpicture}    
\node[anchor=east,inner sep=0] at (-3,0) {\includegraphics[width=2.5cm]{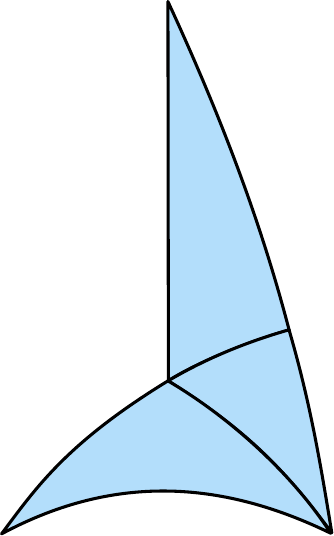}};
\node[anchor=center,inner sep=0] at (0,0) {\includegraphics[width=4cm]{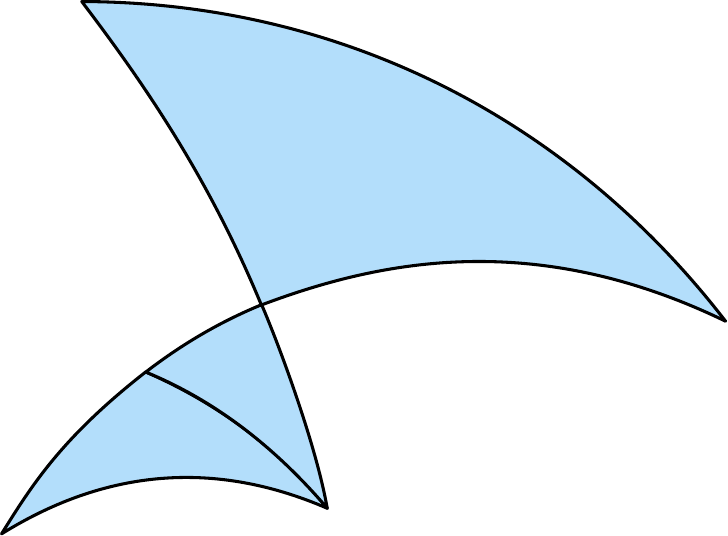}};
\node[anchor=west,inner sep=0] at (3,0) {\includegraphics[width=4cm]{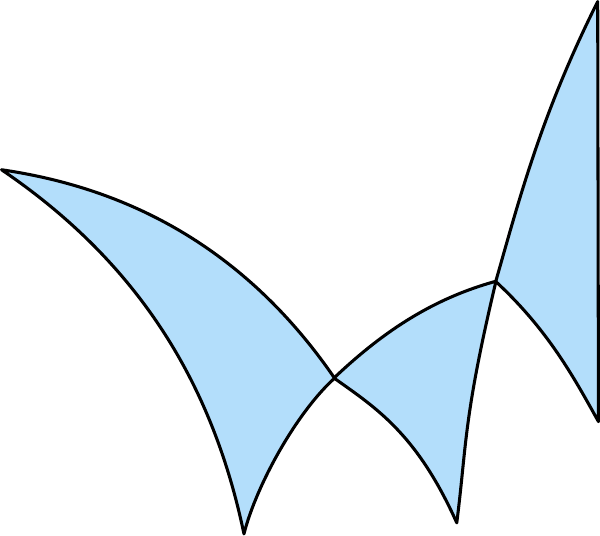}};

\node at (-4,-2.2) {\textrm{hang-glider}};
\node at (0,-2.2) {\textrm{sand clock}};
\node at (5.8,-2.2) {\textrm{bat}};
\end{tikzpicture}
\vspace{3mm}
\end{center}

Tykhyy conjectured that if the sphere $\Sigma$ has 7 punctures or more, then the conjugacy class of a Zariski dense representation $\rho\colon\pi_1\Sigma\to\SL_2\C$ has a finite mapping class group orbit only when $\rho$ maps all but at most 6 peripheral loops of $\Sigma$ to plus or minus the identity matrix~\cite[Section~11]{tykhyy}. \cm{The full statement of Tykhyy's Conjecture, which we state as Conjecture~\ref{conj:tykhyy}, actually enumerates all possible kinds of finite mapping class group orbits for spheres with an arbitrary number of punctures. While the conjecture contains a detailed list of all finite orbits for Zariski dense representations, it doesn't provide explicit lists in the upper triangular and finite image cases. We'll discuss these two cases in Sections~\ref{sec:non-zariski-dense-orbits} and~\ref{sec:finite-image} (see also Remark~\ref{rem:complete-list-of-orbits}). Our main result is the completion of the proof of Tykhyy's Conjecture as stated in Conjecture~\ref{conj:tykhyy}.}

\begin{thmx}[Theorem~\ref{thm:tykhyy}]\label{intro:tykhyy}
Tykhyy's Conjecture \cm{(Conjecture~\ref{conj:tykhyy})} is true.
\end{thmx}

\cm{Some particular cases of the conjecture have been established already. We mentioned the classification of finite orbits by Lisovyy--Tykhyy for $4$-punctured spheres~\cite{LT}. Cousin--Moussard covered the case of representations with image in the group of upper triangular matrices~\cite[Theorem~2.3.4]{cousin-moussard}. Diarra proved that there cannot be any finite orbit of pullback type if $\Sigma$ has 7 punctures or more~\cite[Théorème~5.1]{diarra-pull-back}. More recently, Lam--Landesman--Litt proved Tykhyy's Conjecture for Zariski dense representations, under the assumption that at least one peripheral loop of $\Sigma$ is sent to an element of $\SL_2\C$ with infinite order~\cite[Corollary~1.1.8]{lll}. The new contributions of this note are the following.}
\begin{enumerate}
\item \cm{We settle the last remaining cases needed to prove Tykhyy's Conjecture in full generality; namely, Zariski dense representations $\pi_1\Sigma\to\SL_2\C$ with finite order monodromy at each puncture of $\Sigma$. It turns out that such representations are Galois conjugate to DT representations whenever the associated mapping class group orbits are finite (Proposition~\ref{prop:Zd+non-pullback-implies-DT}). This allows us to conclude the proof of Theorem~\ref{intro:tykhyy} from Theorems~\ref{intro:tykhyy-DT} \&~\ref{intro:finite-orbit-n=6}.}

\item \cm{By reducing the remaining cases of Tykhyy's conjecture to DT representations and classifying all corresponding finite orbits, we also answer a recent question of Litt~\cite[Question~2.4.1]{litt} asking whether the classification obtained in~\cite{lll} can be completed without the infinite order assumption on peripheral monodromies. This yields an alternative formulation of the classification (Corollary~\ref{cor:litt-version}).}
\end{enumerate}

\cm{We learned through personal communications that Deroin--Landesman--Litt--Tholozan have independently obtained a proof that Tykhyy's Conjecture holds for a sufficiently large number of punctures on $\Sigma$ using different techniques.}

\subsection{Ideas of the proofs}
Our approach to prove Theorem~\ref{intro:tykhyy-DT} is inductive. Lisovyy--Tykhyy's classification of finite mapping class group orbits for 4-punctured spheres (\cite{LT}) will serve as base case for the induction. The induction step works as follows. When we work with a sphere $\Sigma$ punctured at $n$ points, then we can consider a special kind of pants decompositions of $\Sigma$ which are called \emph{chained pants decompositions} (Section~\ref{sec:triangle-chains}). They partition $\Sigma$ into $n-2$ pair of pants, only two of which contain two of the original punctures \cm{and the others only contain one original puncture}. When we glue back two neighbouring pairs of pants along their common pants curve, we obtain a sub-surface $\Sigma'\subset \Sigma$ homeomorphic to a 4-punctured sphere (Section~\ref{sec:restrecting-beta_i-to-a-finite-set}). A chained pants decomposition of $\Sigma$ gives rise to $n-3$ such sub-surfaces $\Sigma'$, one for each pants curve. A representation $\rho\colon\pi_1\Sigma\to\SL_2\C$ restricts to a representation $\rho'\colon\Sigma'\to \SL_2\C$ for each sub-surface $\Sigma'$. If the mapping class group orbit of the conjugacy class of $\rho$, denoted by $[\rho]$, is finite, then the mapping class group orbit of each restriction $[\rho']$ is also finite. Because we picked the sub-surfaces $\Sigma'$ to be 4-punctured spheres, the restrictions $\rho'$ therefore all belong to the list of finite orbits for 4-punctured spheres. 

We'll use this observation to infer that the rotation angles of the image by $\rho$ of all pants
curves must be realized by a finite orbit in the 4-punctured case. A classical invariant that helps understanding the possible values of the rotation angles is the trace field of a representation. We use a variant of it which we call the \emph{non-peripheral trace field} (Definition~\ref{def:non-peripheral-trace-field}). By carefully choosing the initial chained pants decomposition (Lemma~\ref{lem:existence-regular-triangle-chain}), we can constrain the list of possible values for these rotation angles to a finite set of $15$ values (Lemma~\ref{lem:finite-orbit-implies-beta_i-in-some-finite-list}). This alone already shows that $n\leq 18$ (Remark~\ref{rem:n_leq_18}).  In order to lower the upper bound on $n$ down to $6$, we need to analyse which combinations of restrictions $\rho'$ are possible (Lemma~\ref{lem:possible-partial-sequences-beta-angles}).

The analysis that we carry out to prove Theorem~\ref{intro:tykhyy-DT} also applies when studying
finite mapping class group orbits for 6-punctured spheres. We use it to prove the first part of the statement of Theorem~\ref{intro:finite-orbit-n=6}: only representations with the same elliptic peripheral behaviour give rise to finite orbits. The existence of jester's hat orbits has been known since the work of Diarra~\cite{diarra-pull-back}. It can also be established in terms of triangle chains. One advantage of the triangle chain model of DT representations is that it reduces the computations of orbit points to elementary hyperbolic geometry. We use our methods to give a new argument that jester's hat orbits are finite. We work with an explicit (minimal) generating family of the mapping class group of $\Sigma$ made of Dehn twists (Lemma~\ref{lem:mcg-generators} from Appendix~\ref{apx:generators-of-pmod}). This generating family is particularly pleasant to work with because the action of each Dehn twist it contains can be described geometrically in terms of triangle chains (Section~\ref{sec:action-of-Dehn-twists} and Example~\ref{ex:example-Dehn-twist-action}). To identify every point in a finite orbit, we'll use the algorithm described in Section~\ref{sec:algorithm-to-compute-orbit-points} in Appendix~\ref{apx:algo-orbit}. Given the length of jester's hat orbits, we only run some of these computations by hand and rely on a computer otherwise. We describe a routine that helps us approximate orbit points in Appendix~\ref{apx:algo-orbit}.\footnote{All pieces of code that we wrote have been made public as Jupyter notebooks on GitHub (\url{https://github.com/shmulik377/FiniteOrbits}).}

The tricky part of the proof of Theorem~\ref{intro:finite-orbit-n=6} is the uniqueness statement. In order to prove that there are no finite mapping class group orbits for $6$-punctured spheres other than the jester's hat orbits, we have to make sure that no finite orbits can be made of \emph{singular} triangle chains only \cm{(Definition~\ref{def:regular-singular-triangle-chains}). (A triangle chain is singular if all three vertices of at least one triangle in the chain coincide.)} We get rid of that possibility by constructing explicit Dehn twists that map a point whose triangle chain is singular to a point with a regular triangle chain.

In order to complete the proof of Tykhyy's Conjecture (Theorem~\ref{intro:tykhyy}), we invoke an alternative by Corlette--Simpson~\cite{corlette-simpson} and Loray--Pereira--Touzet~\cite{loraypereiratouzet} (Theorems~\ref{thm:corlette-simpson} \&~\ref{thm:loray-pereira-touzet}). Its relation to finite mapping class group orbits was already exploited in~\cite{lll}. In broad words, the alternative says that a Zariski dense representation $\rho\colon\pi_1\Sigma\to\SL_2\C$ whose conjugacy class belongs to a finite mapping class group orbit is either of pullback type (classified by Diarra~\cite{diarra-pull-back}), or the associated local system supports a variation of Hodge structures on every punctured Riemann sphere. In the latter case, we'll say in short that $\rho$ is a \emph{universal} variation of Hodge structures \cm{(Definition~\ref{defn:universal-CVHS})}. Such representations $\rho$ are rigid and valued in the integers of a number field. They also preserve a Hermitian form on that
number field. When we consider the Galois conjugates of $\rho$, we obtain three further possibilities.
\begin{enumerate}
\item All Galois conjugates of $\rho$ preserve a Hermitian metric of signature $(2,0)$ or $(0,2)$. In that case, $\rho$ has finite image. The mapping class group orbit of its conjugacy class is automatically finite (Section~\ref{sec:finite-image}).
\item Some Galois conjugate of $\rho$ preserves a Hermitian metric of signature $(1,1)$ and $\rho$ has at least one infinite order peripheral monodromy. As it turns out, in such a scenario, all Galois conjugates of $\rho$ preserve a Hermitian metric of signature $(1,1)$. The classification of finite mapping class group orbits in this case was done by Lam--Landesman--Litt~\cite{lll}. They used Katz's middle convolution to reduce the study of finite mapping class group orbits to complex reflection groups~\cite[Corollary~1.1.7]{lll}, see also Corollary~\ref{cor:litt-version} and the discussion beforehand. All the finite orbits obtained in this way have been listed by Vayalinkal~\cite{amal}. 
\item Some Galois conjugate of $\rho$ preserves a Hermitian metric of signature $(1,1)$ and all peripheral monodromies of $\rho$ are elliptic of finite order. \sam{Since $\rho$ is a universal variation of Hodge structures, this Galois conjugate must be a DT representation (Proposition~\ref{prop:Zd+non-pullback-implies-DT}). This is a consequence of} the characterization of DT representations as universal variations of Hodge structures that was already observed by Deroin--Tholozan~\cite[Theorem~5 and discussion thereafter]{deroin-tholozan} \cm{(Proposition~\ref{prop:universal-CVHS-implies-DT}).} The finite mapping class group orbits that arise from such representations are therefore classified by Theorems~\ref{intro:tykhyy-DT} \&~\ref{intro:finite-orbit-n=6}, \cm{and by \cite{LT} for $4$-punctured spheres.}
\end{enumerate}

\subsection{Motivations and related works}\label{sec:history-and-motivations}
\subsubsection{Dynamics on character varieties}
Our original motivation comes from questions related to dynamics on character varieties (Section~\ref{sec:character-varieties}). The goal is to understand the mapping class group action of a surface $\Sigma$ on the space of conjugacy classes of group homomorphisms from $\pi_1\Sigma$ into a Lie group $G$. This action preserves the Liouville measure associated to the Goldman symplectic form~\cite{goldman-symplectic}. We call it the \emph{Goldman measure} in short. The flavours of the dynamics vary substantially depending on the nature of $G$ and on the component of the character variety that's being acted on. For instance, the action is known to be ergodic when $G$ is a compact Lie group. The original case when $G=\SU(2)$ was treated by Goldman~\cite{goldman-ergodicity} (see also~\cite{goldman-xia}). The result was later generalized to arbitrary compact Lie groups by Pickrell--Xia~\cite{pickrell-xia-1, pickrell-xia-2}. Less is known for non-compact Lie groups. Goldman conjectured that the action should be ergodic on all the non-Teichmüller components of the $\psl$-character varieties for closed hyperbolic surfaces~\cite{goldman-conjecture}. The case of surfaces of genus two was treated by Marché--Wolff~\cite{marche-wolff, marche-wolff-2} and the conjecture remains open for higher \cm{genera}. The second author established ergodicity of the mapping class group action on DT components~\cite{maret-ergodicity} adapting the methods developed in~\cite{goldman-xia, marche-wolff}. Lam--Landesman--Litt later observed that DT components are isomorphic to certain components of relative character varieties of representations in unitary groups via Katz's middle convolution~\cite[Remark~4.3.4]{lll}. Combined with the results of~\cite{pickrell-xia-2}, this gives a new proof of the ergodicity of the mapping class group action on DT components.

Two possible directions to push the study of the mapping class group action a step further than ergodicity are concerned with minimality properties (see Goldman's Problem~\cite[Problem~2.7]{goldman-conjecture}), as well as understanding invariant measures. Previte--Xia proved that for all surfaces $\Sigma$ except punctured spheres, if a representation $\rho\colon\pi_1\Sigma\to\SU(2)$ has dense image, then the mapping class group orbit of its conjugacy class is dense~\cite{previte-xia-1, previte-xia-2}. It's interesting to notice that Golsefidy--Tamam recently identified counter-examples to Previte-Xia's results in \cm{genera} 1 and 2 and proposed a revised statement~\cite[Corollary~93]{golsefidy-tamam}. \cm{Previte--Xia~\cite[Theorem~3.3]{previte-xia-minimality}, and later Cantat--Loray~\cite[Theorem~C]{cantat-loray}, proved an analogous statement for 4-punctured spheres: any infinite mapping class group orbit contained in a relative $\SU(2)$-character variety or in a DT component is dense. This result was recently generalized to DT components for spheres with an arbitrary number of punctures~\cite[Theorem~A]{gianluca-yohann-arnaud}. It shows that the mapping class group action on DT components that don't contain any finite orbit is minimal. This is for instance the case when the underlying sphere has $7$ punctures or more by Theorem~\ref{intro:tykhyy-DT}.}

In a recent work, Cantat--Dupont--Martin-Baillon proved that any stationary ergodic measure on a relative $\SL_2\C$-character variety of representations of a 4-punctured sphere is either supported on a finite orbit, or it is supported on a compact component and coincides with the Goldman measure~\cite{stationary-measures}. \cm{In particular, the mapping class group action is uniquely ergodic on all DT components of 4-punctured spheres that don't contain any finite orbit. The question was raised in~\cite{gianluca-yohann-arnaud} whether this result can be generalized to all DT components: namely, whether an ergodic measure on a DT component for a sphere with an arbitrary number of punctures is either supported on a finite orbit or is the Goldman measure.}

\subsubsection{Isomonodromy differential equations}
A different context in which finite mapping class group orbits of conjugacy classes of representations $\pi_1\Sigma\to\SL_2\C$ arise naturally is the study of algebraic solutions to isomonodromy differential equations such as Painlevé~\Romannum{6}, or Garnier \cm{and Schlesinger} systems. \cm{According to Iwasaki, classifying algebraic solutions to Painlevé~\Romannum{6} and finite mapping class group orbits in the $\SL_2\C$-character variety of a 4-punctured sphere are equivalent problems~\cite{iwasaki}.} The quest for algebraic solutions to \cm{the Painlevé~\Romannum{6}} equation produced many examples of finite mapping class group orbits over the years, including contributions by Andreev--Kitaev~\cite{andreev-kitaev}, Boalch~\cite{boalch-3, boalch-1, boalch-4, boalch-2}, Dubrovin~\cite{dubrovin}, Dubrovin--Mazzocco~\cite{dubrovin-mazzocco}, Hitchin~\cite{hitchin-1, hitchin-2}, and Kitaev~\cite{kitaev, kitaev-2}. Boalch presented a list of 45 known exceptional algebraic solutions to Painlevé~\Romannum{6} in~\cite{boalch-talk} which he then summarized in~\cite{boalch}. This list turned out to be complete, as later confirmed by Lisovyy--Tykhyy~\cite{LT}. 

\cm{The Painlevé~\Romannum{6} equation generalizes to Garnier systems. Algebraic solutions of these systems give rise to monodromy representations into $\SL_2\C$ with finite mapping class group orbits for spheres with more than four punctures, and all such finite orbits arise in this way by a result of Cousin~\cite[Theorem~C]{cousin-isomonodromy}.} A partial classification of finite orbits for 5-punctured spheres was achieved by Calligaris--Mazzocco~\cite{calligaris-mazzocco} and later extended by Tykhyy~\cite{tykhyy}. Diarra classified a special kind of algebraic solutions to Garnier systems for arbitrary spheres: those whose monodromies are of pullback type~\cite{diarra-pull-back} (see also the discussion in Section~\ref{sec:pullback-orbits}). In particular, Diarra proved that algebraic solutions of this kind cease to exist if the underlying sphere has 7 punctures or more in accordance with Tykhyy's Conjecture (Conjecture~\ref{conj:tykhyy}).

\cm{We briefly mention that passing from a finite mapping class group orbit to an explicit algebraic isomonodromic family in the de Rham moduli space of flat connections is a delicate problem. Even for 4-punctured spheres, constructing the corresponding families of Fuchsian systems requires substantial analytic and computational effort. This was carried out for all finite orbits prior to Lisovyy--Tykhyy's work~\cite{LT} by several authors and notably Boalch~\cite{boalch-1}. In this note, we focus on classifying finite orbits and will not attempt to compute the corresponding algebraic isomonodromy families.}

\subsubsection{Higher genus}
Leaving the world of punctured spheres, all the finite mapping class group orbits of conjugacy classes of \emph{reductive representations} $\pi_1\Sigma\to\SL_2\C$ \cm{(Definition~\ref{defn:reductive-and-irreducible-representations})} are known in the case where $\Sigma$ has positive genus and a non-negative number of punctures. The classification is due to Biswas--Gupta--Mj--Whang~\cite{positive-genus}. The variety of finite orbits is less rich than for punctured spheres. For instance, they proved that if the genus of $\Sigma$ is at least two, then all the finite orbits come from representations with finite image. The case of non-reductive representations was covered by Cousin--Heu~\cite{cousin-heu} who reached similar conclusions on the type of finite orbits.

So far, we've only considered finite mapping class group orbits for $\SL_2\C$ surface group representations. \cm{The systematic study of finite mapping class group orbits for representations into $\GL_r\C$ remains largely undeveloped. An important result in this direction is due to Landesman--Litt, who answered a question of Whang: for a fixed rank $r$, does there exist a constant $C(r)$ such that any representation $\pi_1\Sigma\to \GL_r\C$ of a surface $\Sigma$ of genus $g\ge C(r)$ inducing a finite mapping class group orbit must have finite image? Landesman--Litt proved that such a bound does exist, and that one may take $C(r)=r^2$~\cite{LL}.} A more detailed account on the history of finite orbits, and their connection to algebraic geometry, low dimensional topology and differential equations can be found in Litt's survey paper~\cite{litt}.

\subsection{Structure of the paper}
Chapter~\ref{chap:background} reviews some background knowledge. After fixing some notation about character varieties and mapping class group dynamics (Section~\ref{sec:character-varieties}), we give a condensed introduction to DT representations (Section~\ref{sec:DT-representations}). We mention their parametrization by triangle chains, as well as how to extract action-angle coordinates from that model. The chapter ends with two brief recaps on ``trigonometric number fields'' (Section~\ref{sec:number-fields}) and on (discrete) triangle groups (Section~\ref{sec:triangle-groups}).

We spend some time in Chapter~\ref{chap:examples} presenting various examples of finite mapping class group orbits with the hope to train the reader's intuition on the topic. We cover the cases of representations with abelian or finite image (Sections~\ref{sec:abelian-image} \&~\ref{sec:finite-image}), discrete representations (Section~\ref{sec:discrete-finite-orbits}), and finite orbits of pullback types (Section~\ref{sec:pullback-orbits}). The list of examples is expanded in Chapter~\ref{sec:4-punctured} where we review Lisovyy--Tykhyy's classification. We also study some particular finite orbits in details and provide the action-angle coordinates of their orbit points (Section~\ref{sec:some-finite-orbits}).

The proof of Theorem~\ref{intro:tykhyy-DT} is explained in Chapter~\ref{chap:tykhyy}. The classification of finite orbits for 6-punctured spheres (Theorem~\ref{intro:finite-orbit-n=6}) is established in Chapter~\ref{sec:classification}. We also recover Tykhyy's classification for 5-punctured spheres later in Chapter~\ref{sec:classification} (Theorem~\ref{thm:angle-vector-alpha-with-finite-orbits-n=5}). Along the way, we introduce the names of jester's hat, hang-glider, sand clock, and bat orbits. Finally, we explain how our results complete the proof of Tykhyy's Conjecture in Chapter~\ref{chap:tykhyy-conjecture-general}.

Three appendices are provided at the end of the paper. Appendix~\ref{apx:generators-of-pmod} is about explicit generating families for mapping class groups of punctured spheres. The algorithm we use to compute all orbits points in a finite mapping class group orbit, as well as a description of how we implement it on a computer, is presented in Appendix~\ref{apx:algo-orbit}. The last appendix (Appendix~\ref{app:tables}) contains several tables, among which the ones listing the action-angle coordinates of orbit points for finite mapping class group orbits of 5-punctured and 6-punctured spheres.

\subsection{Acknowledgements}
We are extremely grateful to Nicolas Tholozan for introducing us to the topic of finite mapping class group orbits and for clarifying both the significance of the Corlette--Simpson alternative and the crucial connection between DT representations and variations of Hodge structures, as detailed in Chapter~\ref{chap:tykhyy-conjecture-general}. We extend our appreciation to Philip Boalch, Bertrand Deroin, Elisha Falbel, Josh Lam, Aaron Landesman, Daniel Litt, Frank Loray, and Julien Marché for fruitful conversations and their guidance. Thanks to Anna Felikson and Detchat Samart for sharing some insights on their personal work. \cm{We are also grateful to Simon André, Renat Gontsov, and the anonymous referees for numerous comments and suggestions.} 

\section{Background}\label{chap:background}
\subsection{Character varieties and mapping class group action}\label{sec:character-varieties}
A character variety is, by definition, the space of conjugacy classes of group homomorphisms from the fundamental group of an oriented surface $\Sigma$ (possibly with punctures) into a Lie group $G$. More precisely, we first consider the set $\Hom(\pi_1\Sigma,G)$ of group homomorphisms $\pi_1\Sigma\to G$ and we equip it with the compact-open topology (where $\pi_1\Sigma$ is given the discrete topology). The elements $\rho\in \Hom(\pi_1\Sigma,G)$ will also be called \emph{representations} of the fundamental group of $\Sigma$, sometimes simply representations of $\Sigma$. There is an action of $G$ on $\Hom(\pi_1\Sigma,G)$ by post-composition by inner automorphisms of $G$. 
\begin{defn}\label{defn:character-varieties}
\cm{The quotient $\Hom(\pi_1\Sigma,G)/G$, equipped with the quotient topology, will be called the \emph{character variety} of the pair $(\Sigma,G)$. We'll denote it by}
\[
\Rep(\Sigma, G).
\]
\end{defn}
There are alternative definitions of character varieties that are more sophisticated and ensure a better topological structure. For instance, when $G=\SL_2\C$, it's sometimes convenient (as we'll do in Section~\ref{sec:fricke-relation}) to define the character variety of the pair $(\Sigma,\SL_2\C)$ to be the algebraic \cm{GIT} quotient $\Hom(\pi_1\Sigma,\SL_2\C)\sslash\SL_2\C$, \cm{or, equivalently, as the topological quotient of the space of all reductive representations.}
\begin{defn}\label{defn:reductive-and-irreducible-representations}
\cm{A representation $\rho\colon\pi_1\Sigma\to\SL_2\C$ is \emph{reductive} if the linear action of $\rho(\pi_1\Sigma)$ on $\mathbb{C}^2$ is semisimple, i.e.~every invariant subspace has an invariant complement. The representation is \emph{irreducible} if the action is irreducible, i.e.~there are no proper invariant subspaces.}
\end{defn}

\cm{When $\Sigma$ has genus $g\geq 0$ and a positive number $n\geq 1$ of punctures forming a set $\mathcal{P}$, its fundamental group $\pi_1\Sigma$ is a free group and $\Hom(\pi_1\Sigma,G)$ is isomorphic to $G^{2g+n-1}$.} In order to get a more interesting space of representations, we fix peripheral data. We'll use the name of \emph{peripheral loops} for counter-clockwise loops on $\Sigma$ that isolate one of the punctures from the others. We denote by $G/G$ the set of conjugacy classes in $G$. If $\mathcal{C}\in (G/G)^{\mathcal{P}}$ is a tuple of conjugacy classes indexed on $\mathcal{P}$, then the set of all the representations $\pi_1\Sigma\to G$ that map every peripheral loop around $p\in\mathcal{P}$ to some element of $\mathcal{C}_p$ is denoted by $\Hom_\mathcal{C}(\pi_1\Sigma, G)$. 
\begin{defn}\label{defn:relative-character-variety}
\cm{The topological quotient $\Hom_\mathcal{C}(\pi_1\Sigma, G)/G$ is called the \emph{$\mathcal{C}$-relative character variety} of the pair $(\Sigma,G)$ and we'll write it}
\[
\Rep_\mathcal{C}(\Sigma, G).
\]
\end{defn}
\cm{This note is concerned with mapping class group dynamics on (relative) character varieties.}
\begin{defn}\label{defn:pure-mcg}
\cm{The \emph{pure mapping class group} of $\Sigma$ is defined as the group of isotopy classes of
orientation-preserving homeomorphisms of $\SpherePk$ that fix each puncture in $\mathcal{P}$ individually. We'll denote it by~$\PMod(\SpherePk)$.}
\end{defn}
The pure mapping class group is finitely generated by a special kind of elements called \emph{Dehn twists}. We'll describe a generating family of Dehn twists in Appendix~\ref{apx:generators-of-pmod}. For precise statements and more considerations on mapping class groups of surfaces, the reader is referred to Farb--Margalit's book~\cite{mcg-primer}. There is a natural morphism from $\PMod(\SpherePk)$ to the group $\Out(\pi_1\SpherePk)$ of outer automorphisms of $\pi_1\SpherePk$. Since we are focusing on the pure mapping class group, the image of the morphism $\PMod(\SpherePk)\to \Out(\pi_1\SpherePk)$ is actually contained inside the subgroup $\Out^\star(\pi_1\SpherePk)$ of $\Out(\pi_1\SpherePk)$ which consists of outer automorphisms of $\pi_1\Sigma$ that preserve conjugacy classes of peripheral loops around each puncture in $\mathcal{P}$. The Dehn--Nielsen--Baer Theorem says that the morphism $\PMod(\SpherePk)\to \Out^\star(\pi_1\SpherePk)$ is an isomorphism, see~\cite[Theorem~8.8]{mcg-primer} for more details. 

Since $\PMod(\SpherePk)$ is isomorphic to $\Out^\star(\pi_1\SpherePk)$,
it acts by pre-composition on~$\Rep_\mathcal{C}(\Sigma, G)$.
More precisely, if $[\rho]\in \Rep_\mathcal{C}(\Sigma, G)$ is the conjugacy class of the representation $\rho\colon\pi_1\Sigma\to G$ and $\tau\in \PMod(\SpherePk)$ is the equivalence class of an automorphism $\varphi\colon \pi_1\Sigma\to\pi_1\Sigma$, then $\tau.[\rho]$ is the conjugacy class of $\rho\circ \varphi^{-1}$. The peripheral structure is of course preserved by the action of $\PMod(\SpherePk)$. Most of the time, we'll refer to the action of $\PMod(\SpherePk)$ on $\Rep_\mathcal{C}(\Sigma, G)$ as the \emph{mapping class group action} and to its orbits as \emph{mapping class group orbits}, dropping the adjective ``pure'' for simplicity.

\subsection{Recap on DT representations}\label{sec:DT-representations}
\subsubsection{Brief history}
Through computer simulations, Benedetto--Goldman were the first to spot a compact component among the real points of the relative character variety of representations of a 4-punctured sphere into $\SL_2\C$~\cite{benedetto-goldman}. Their discovery was particularly interesting as they also proved that, for some choices of peripheral parameters, the representations whose conjugacy class belong to this compact component are valued inside $\SL_2\R$, and not into $\SU(2)$. They explained that all these compact components are smooth spheres and already hinted at the connection between the corresponding representations and non-embedded \emph{sand clock-shaped quadrilaterals} in the hyperbolic plane. These representations have the property of mapping every simple closed curve on the underlying 4-punctured sphere to elliptic elements of $\SL_2\R$. In today's terminology, we say that they're \emph{totally elliptic}. Cantat--Loray later explained that the compact components of $\SL_2\R$ representations found by Benedetto--Goldman can be equivariantly identified with relative $\SU(2)$-character varieties of 4-punctured spheres using \emph{Okamoto transformations}~\cite{cantat-loray} (we'll expand on this in Section~\ref{sec:finite-orbit-n=4}). 

About two decades after Benedetto--Goldman, Deroin--Tholozan discovered
analogous compact components of totally elliptic representations of
arbitrary punctured spheres into $\SL_2\R$~\cite{deroin-tholozan}. They
proved that these compact components are isomorphic to complex projective
spaces whose real dimension is twice the number of punctures minus six. At
the time, they called them \emph{supra-maximal} representations; a name that
was later abandoned to avoid confusion with the maximal representations from
higher Teichmüller theory. They also suggested that these representations
can be constructed from \emph{necklaces} of hyperbolic triangles. Around the
same time, Mondello studied the topology of relative $\psl$-character
varieties and used Higgs bundles methods to enumerate their connected
components~\cite{mondello}. One implication of his work is that the compact
component of Deroin--Tholozan is unique in the corresponding relative
character variety. 

Some years later, the second author built \cm{upon} the ideas of non-embedded hyperbolic quadrilaterals from Benedetto--Goldman and necklaces of hyperbolic triangles from Deroin--Tholozan to construct action-angle coordinates to parametrize these compact components~\cite{action-angle}. At the same time, the name of \emph{Deroin--Tholozan representations} was introduced. 

In recent works, Lam--Landesman--Litt described how the compact character varieties of Deroin--Tholozan representations can be equivariantly identified with (some component of) relative character varieties of representations into a compact Lie group~\cite{lll}. The isomorphism comes from \emph{Katz's middle convolution}---a generalization of Okamoto transformations.

\subsubsection{Definition}\label{sec:definition-DT}
For every integer $n\geq 3$, we fix an oriented sphere $\SpherePk$ punctured at $n$ points.
The set of punctures is denoted by $\mathcal{P}$. The fundamental group of $\SpherePk$ can be presented as
\begin{equation}\label{eq:geometric-presentation}
\pi_1\SpherePk =\langle c_1,\dots,c_n\, |\, c_1\cdots c_n =1\rangle
\end{equation}
by carefully choosing each $c_i$ as the homotopy class of a peripheral loop around one of
the punctures. Such a presentation of $\pi_1\SpherePk$ is called \emph{geometric}.
For every angle vector $\alpha\in(0,2\pi)^\mathcal{P}$, we introduce
the~\emph{$\alpha$-relative character variety}\footnote{\cm{The notions of relative $\psl$ character varieties for hyperbolic, parabolic, or mixed peripheral monodromies may be defined as in Definition~\ref{defn:relative-character-variety}. However, those won't be relevant for this note because DT representations will be required to have elliptic peripheral monodromy (Definition~\ref{defn:DT-representation}).}} 
\begin{equation*}
\Rep_{\alpha}(\SpherePk,\psl)
\end{equation*}
as the space of conjugacy classes of representations $\rho\colon\pi_1\SpherePk\rightarrow \psl$
that map counter-clockwise loops around each puncture $p\in\mathcal{P}$ to an elliptic element
of~$\psl$ of rotation angle $\alpha_p$. 
\begin{defn}\label{defn:rotation-angle}
An elliptic element of~$\psl$
has \emph{rotation angle} $\vartheta\in (0,2\pi)$ if it is conjugate to
\begin{equation*}
	\pm\begin{pmatrix}
		 \cos(\vartheta/2) & \sin(\vartheta/2)\\
		-\sin(\vartheta/2) & \cos(\vartheta/2)
	\end{pmatrix}.
\end{equation*}
\end{defn}
We'll often refer to $\alpha$ as the vector of \emph{peripheral angles} and denote by $[\rho]$ the conjugacy class of the representation $\rho$. 
\begin{fact}\label{fact:rep_alpha-smooth}
\cm{If $\sum_{p\in\mathcal{P}}\alpha_p\notin 2\pi\Z$, then every representation whose conjugacy class belongs to the $\alpha$-relative character variety $\Rep_{\alpha}(\SpherePk,\psl)$ has Zariski dense image, hence is irreducible. In particular, $\Rep_{\alpha}(\SpherePk,\psl)$ is smooth.}
\end{fact}
\begin{proof}
\cm{If $[\rho]\in \Rep_{\alpha}(\SpherePk,\psl)$, then by definition $\rho(c_1),\ldots, \rho(c_n)$ are non-trivial elliptic elements of $\psl$. Their fixed points may coincide only when $\sum_{p\in\mathcal{P}}\alpha_p\in 2\pi\Z$. So, if $\sum_{p\in\mathcal{P}}\alpha_p\notin 2\pi\Z$, then $\rho(\pi_1\SpherePk)$ contains non-trivial elliptic elements with different fixed points. This implies that $\rho$ has Zariski dense image in $\psl$. Smoothness follows from Goldman's work~\cite{goldman-symplectic}.}
\end{proof}

Whenever $\alpha$ satisfies
\begin{equation}\label{eq:angle-condition}
\sum_{p\in\mathcal P}\alpha_p> 2\pi(n-1),
\end{equation}
Deroin--Tholozan proved that $\Rep_{\alpha}(\SpherePk,\mathrm{PSL}_2\mathbb{R})$ \cm{(which is smooth by Fact~\ref{fact:rep_alpha-smooth})} contains a
compact component which is isomorphic (as a symplectic toric manifold) to $\CP^{n-3}$~\cite{deroin-tholozan}. Conjugation by the non-trivial outer automorphism of $\psl$ produces an analogous compact component when $\sum \alpha_p<2\pi$. 
It's a consequence of Mondello's work on the topology of relative $\psl$-character varieties that each of these compact components is unique in its respective relative character variety~\cm{\cite[Corollary~4.17]{mondello}. Moreover, Mondello's result also shows that $\Rep_{\alpha}(\SpherePk,\mathrm{PSL}_2\mathbb{R})$ doesn't contain any compact component when $2\pi\leq\sum\alpha_p\leq 2\pi(n-1)$, except isolated points when $\sum\alpha_p\in 2\pi\Z$. Those isolated points are conjugacy classes of abelian representations.}
\begin{defn}\label{defn:DT-representation}
\cm{We'll refer to the compact components of Deroin--Tholozan in the regimes $\sum\alpha_p<2\pi$ or $\sum\alpha_p>2\pi (n-1)$ as \emph{DT components} and denote them by
\[
\RepDT{\alpha}\subset \Rep_{\alpha}(\SpherePk,\mathrm{PSL}_2\mathbb{R}).
\]
The representations whose conjugacy classes belong to $\RepDT{\alpha}$ will be called \emph{DT representations}.}
\end{defn}
\cm{By definition of the $\alpha$-relative character variety $\Rep_{\alpha}(\SpherePk,\mathrm{PSL}_2\mathbb{R})$, DT representations map all peripheral curves on $\SpherePk$ to non-trivial elliptic elements of $\psl$.} 
\begin{lem}[{\cite[Lemma~3.5]{deroin-tholozan}}]\label{lem:totally-elliptic}
\cm{DT representations are actually \emph{totally elliptic} in the sense that they map every simple closed curve of $\SpherePk$, not only peripheral loops, to a \cm{non-trivial} elliptic element of~$\psl$.}
\end{lem}

\cm{DT representations can always be lifted to representations into $\SL_2\R$ because $\pi_1\Sigma$ is a free group. There are exactly $2^{n-1}$ possible lifts obtained by choosing the lifts of $\rho(c_1),\ldots,\rho(c_{n-1})$ freely, which forces the sign of the lift of $\rho(c_n)$. If $\overline{\rho}\colon\pi_1\Sigma\to\SL_2\R$ denotes one of those lifts, then all the other lifts are obtained by multiplying an even number of $\overline{\rho}(c_i)$ by $-1$. In particular, the sign of the product of traces $\prod_i \trace\overline{\rho}(c_i)$ is independent of the choice of lift. The following result computes that sign when $\Sigma$ is a 4-punctured sphere.}
\begin{lem}\label{lem:sign-product-of-traces-n=4}
\cm{Let $\Sigma$ be a 4-punctured sphere $(n=4)$ and let $\alpha = (\alpha_1,\alpha_2,\alpha_3,\alpha_4)$ be a vector of peripheral angles satisfying $\alpha_1\geq \alpha_2\geq \alpha_3\geq \alpha_4$ and $\sum_{i}\alpha_i>6\pi$---this is Condition~\eqref{eq:angle-condition}. If $\rho\colon\pi_1\Sigma\to\psl$ is a DT representation with peripheral angles $\alpha$ and $\overline{\rho}\colon\pi_1\Sigma\to\SL_2\R$ is any lift of $\rho$ with peripheral traces $t_i=\trace\overline{\rho}(c_i)$, then the following may happen:
\[
\begin{cases}
\alpha_4=\pi \text{ and } \prod_i t_i =0\\
\alpha_4<\pi \text{ and } \prod_i t_i >0\\
\alpha_4>\pi \text{ and } \prod_i t_i <0.
\end{cases}
\]}
\end{lem}
\begin{proof}
\cm{Consider the simple closed curve on $\Sigma$ represented by the fundamental group element $b=(c_1c_2)^{-1}$. Lemma~\ref{lem:totally-elliptic} says that $\rho(b)$ is a non-trivial elliptic element of $\psl$. We denote its rotation angle (Definition~\ref{defn:rotation-angle}) by $\beta\in (0,2\pi)$. Let $R_b$ and $R_i$ for $i=1,2,3,4$ be the unique elements of $\SL_2\R$ that respectively lift $\rho(b)$ and $\rho(c_i)$, and are conjugate in $\SL_2\R$ to the following matrices:
\[
R_b\sim \begin{pmatrix}
		 \cos(\beta/2) & \sin(\beta/2)\\
		-\sin(\beta/2) & \cos(\beta/2)
	\end{pmatrix}\quad\text{and}\quad 
R_i\sim \begin{pmatrix}
		 \cos(\alpha_i/2) & \sin(\alpha_i/2)\\
		-\sin(\alpha_i/2) & \cos(\alpha_i/2)
	\end{pmatrix}.
\]
The products $R_1R_2R_b$ and $R_b^{-1}R_3R_4$ are each equal to $\id$ or to $-\id$, where $\id\in\SL_2\R$ is the identity matrix. In order to determine the correct sign, we use a result of Orevkov that describes all possible conjugacy classes for products of three elliptic elements in $\SL_2\R$~\cite[Table~3]{orevkov}. In particular, it says whether such a product may be equal to $\id$ or $-\id$. Following the notation of~\cite{orevkov}, we denote by $\mathfrak{c}_3^\theta$ the conjugacy class in $\SL_2\R$ of
\[
\begin{pmatrix}
		 \cos(\theta) & -\sin(\theta)\\
		\sin(\theta) & \cos(\theta)
	\end{pmatrix}.
\]
With this notation, we have $R_i\in \mathfrak{c}_3^{2\pi-\alpha_i/2}=-\mathfrak{c}_3^{\pi-\alpha_i/2}$ and $R_b\in \mathfrak{c}_3^{2\pi-\beta/2}=-\mathfrak{c}_3^{\pi-\beta/2}$. The inequalities~\eqref{eq:inequalities-polytope-beta} below applied to the case of a 4-punctured sphere read
\begin{align*}
\pi &> \pi-\frac{\beta}{2} + \pi-\frac{\alpha_1}{2} + \pi-\frac{\alpha_2}{2} >0\\
\pi &> \frac{\beta}{2} + \pi-\frac{\alpha_3}{2} + \pi-\frac{\alpha_4}{2} >0.
\end{align*}
According to the first line of~\cite[Table~3]{orevkov}, we deduce that $\id\notin\mathfrak{c}_3^{\pi-\alpha_1/2}\mathfrak{c}_3^{\pi-\alpha_2/2}\mathfrak{c}_3^{\pi-\beta/2}$ and $\id\notin \mathfrak{c}_3^{\beta/2}\mathfrak{c}_3^{\pi-\alpha_3/2}\mathfrak{c}_3^{\pi-\alpha_4/2}$. This means that $R_1R_2R_b=\id$ and $R_b^{-1}R_3R_4=-\id$, and thus $R_1R_2R_3R_4=-\id$. In other words, we may assume that the lift $\overline{\rho}$ is given by $\overline{\rho}(c_i)=R_i$ for $i=1,2,3$ and $\overline{\rho}(c_4)=-R_4$. Therefore, since $\trace R_i=2\cos(\alpha_i/2)$, we conclude that 
\[
\prod_{i=1}^4 t_i= -2^4\prod_{i=1}^4 \cos\left(\frac{\alpha_i}{2}\right).
\]
Since we're assuming that $\alpha_1\geq \alpha_2\geq \alpha_3\geq \alpha_4$ and $\sum_{i}\alpha_i>6\pi$, we must have $\alpha_1\geq \alpha_2\geq \alpha_3>\pi$ which implies $\cos(\alpha_i/2)<0$ for $i=1,2,3$. In other words, the sign of $\prod_i t_i$ is the same as the sign of $\cos(\alpha_4/2)$, which gives the desired conclusion.}
\end{proof}

\subsubsection{Triangle chains}\label{sec:triangle-chains}
DT representations are parametrized by certain polygonal objects in the hyperbolic plane
called \emph{triangle chains}. The correspondence was made explicit by the second author
in~\cite{action-angle}. To associate a triangle chain to a DT representation, one starts by picking 
a pants decomposition $\mathcal{B}$ of $\SpherePk$. We'll always work
with~\emph{chained pants decompositions}, meaning that every pair of pants contains at least
one of the punctures of $\SpherePk$. The next step consists in finding a geometric presentation
of $\pi_1\SpherePk$ which is \emph{compatible} with $\mathcal{B}$. \cm{Before defining what we mean by compatible, recall that a geometric
presentation of $\pi_1\SpherePk$ consists of $n$ generators $c_1,\ldots, c_n$ which are homotopy classes
of peripheral loops around the punctures of $\SpherePk$ and satisfy the unique relation $c_1\cdots c_n=1$.}
\begin{defn}\label{defn:compatible-generators}
\cm{The geometric generators $c_1,\ldots, c_n$ are said to be \emph{compatible} with the chained pants decomposition $\mathcal{B}$, if the~$n-3$
pants curves of $\mathcal{B}$ are represented by the fundamental group
elements}
\begin{equation}\label{eq:pants-curves}
b_i=(c_1\cdots c_{i+1})^{-1} \text{ for } i=1,\ldots,n-3.
\end{equation}
\end{defn}
\cm{In a slight abuse of terminology, we'll sometimes say that the pants decomposition $\mathcal{B}$ is given by the curves $b_1,\ldots,b_{n-3}$, even though those are fundamental group elements.}

It's always possible to find a geometric presentation of $\pi_1\Sigma$ which is compatible
with a given chained pants decomposition, as explained in~\cite[Appendix~B]{action-angle}.
Conversely, if a geometric presentation of $\pi_1\SpherePk$ with generators $c_1,\ldots,c_n$ is
given, then the pants decomposition of $\SpherePk$ given by $b_1,\ldots, b_{n-3}$ is called
the \emph{standard} chained pants decomposition of $\Sigma$ associated to the given geometric presentation.  
\begin{center}
\vspace{2mm}
\begin{tikzpicture}[scale=1.4, decoration={
    markings,
    mark=at position 0.6 with {\arrow{>}}}]
  \draw[postaction={decorate}] (0,-.5) arc(-90:-270: .25 and .5) node[midway, left]{$c_1$};
  \draw[black!40] (0,.5) arc(90:-90: .25 and .5);
  \draw[apricot, postaction={decorate}] (2,.5) arc(90:270: .25 and .5) node[midway, left]{$b_1$};
  \draw[lightapricot] (2,.5) arc(90:-90: .25 and .5);
  \draw[apricot, postaction={decorate}] (4,.5) arc(90:270: .25 and .5) node[midway, left]{$b_2$};
  \draw[lightapricot] (4,.5) arc(90:-90: .25 and .5);
  \draw[apricot, postaction={decorate}] (6,.5) arc(90:270: .25 and .5) node[midway, left]{$b_3$};
  \draw[lightapricot] (6,.5) arc(90:-90: .25 and .5);
  \draw[postaction={decorate}] (8,.5) arc(90:270: .25 and .5) node[midway, left]{$c_6$};
  \draw (8,.5) arc(90:-90: .25 and .5);
  
  \draw (.5,1) arc(180:0: .5 and .25) node[midway, above]{$c_2$};
  \draw[postaction={decorate}] (.5,1) arc(-180:0: .5 and .25);
  \draw (2.5,1) arc(180:0: .5 and .25)node[midway, above]{$c_3$};
  \draw[postaction={decorate}] (2.5,1) arc(-180:0: .5 and .25);
  \draw (4.5,1) arc(180:0: .5 and .25)node[midway, above]{$c_4$};
  \draw[postaction={decorate}] (4.5,1) arc(-180:0: .5 and .25);
  \draw (6.5,1) arc(180:0: .5 and .25)node[midway, above]{$c_5$};
  \draw[postaction={decorate}] (6.5,1) arc(-180:0: .5 and .25);
   
  \draw (0,.5) to[out=0,in=-90] (.5,1);
  \draw (1.5,1) to[out=-90,in=180] (2,.5);
  \draw (0,-.5) to[out=0,in=180] (2,-.5);
  
  \draw (2,.5) to[out=0,in=-90] (2.5,1);
  \draw (3.5,1) to[out=-90,in=180] (4,.5);
  \draw (2,-.5) to[out=0,in=180] (4,-.5);
  
  \draw (4,.5) to[out=0,in=-90] (4.5,1);
  \draw (5.5,1) to[out=-90,in=180] (6,.5);
  \draw (4,-.5) to[out=0,in=180] (6,-.5);
  
  \draw (6,.5) to[out=0,in=-90] (6.5,1);
  \draw (7.5,1) to[out=-90,in=180] (8,.5);
  \draw (6,-.5) to[out=0,in=180] (8,-.5);
\end{tikzpicture}
\vspace{2mm}
\end{center}
We wish to emphasize that a geometric presentation of $\pi_1\SpherePk$ induces a
bijection~$\mathcal{P}\to \{1,\ldots,n\}$ which gives a labelling of the punctures. In practice,
we'll use this labelling to index variables, such as the entries of the peripheral angles vector~$\alpha$,
on~$\{1,\ldots,n\}$ rather than on $\mathcal{P}$. 
Whenever we'll be working with some fixed labelling of the punctures, we'll favour the \cm{parentheses} notation to list the angles in $\alpha$, and stick to the curly bracket notations when it's not the case.

Let us now explain how to associate a triangle chain to a DT
representation $\rho\colon\pi_1\SpherePk\to\psl$ from a chained pants decomposition $\mathcal{B}$ of $\SpherePk$ and a compatible geometric
presentation of $\pi_1\SpherePk$. 
\begin{defn}\label{defn:triangle-chain}
The \emph{$\mathcal{B}$-triangle chain of $\rho$} is defined as the assembly of the $n-2$ hyperbolic
triangles constructed as follows.
\begin{itemize}
    \item \cm{In the hyperbolic plane,} draw the fixed points $C_1,\ldots,C_n$ of $\rho(c_1),\ldots,\rho(c_n)$ and the fixed
points $B_1,\ldots,B_{n-3}$ of $\rho(b_1),\ldots,\rho(b_{n-3})$. We're using here that~$\rho$
is totally elliptic \cm{(Lemma~\ref{lem:totally-elliptic})} in order to say that $\rho(b_1),\ldots,\rho(b_{n-3})$ are elliptic.
    \item Draw a geodesic segment between two of these points if the corresponding curves
on~$\SpherePk$ belong to the same pair of pants. We end up with a chain of $n-2$ triangles whose
vertices are $(C_1,C_2,B_1), (B_1,C_3,B_2),\ldots, (B_{n-3},C_{n-1},C_n)$.
\end{itemize}
\begin{center}
\resizebox{9.9cm}{!}{
\begin{tikzpicture}[scale=1.1]
\node[anchor=south west,inner sep=0] at (0,0) {\includegraphics[width=9.9cm]{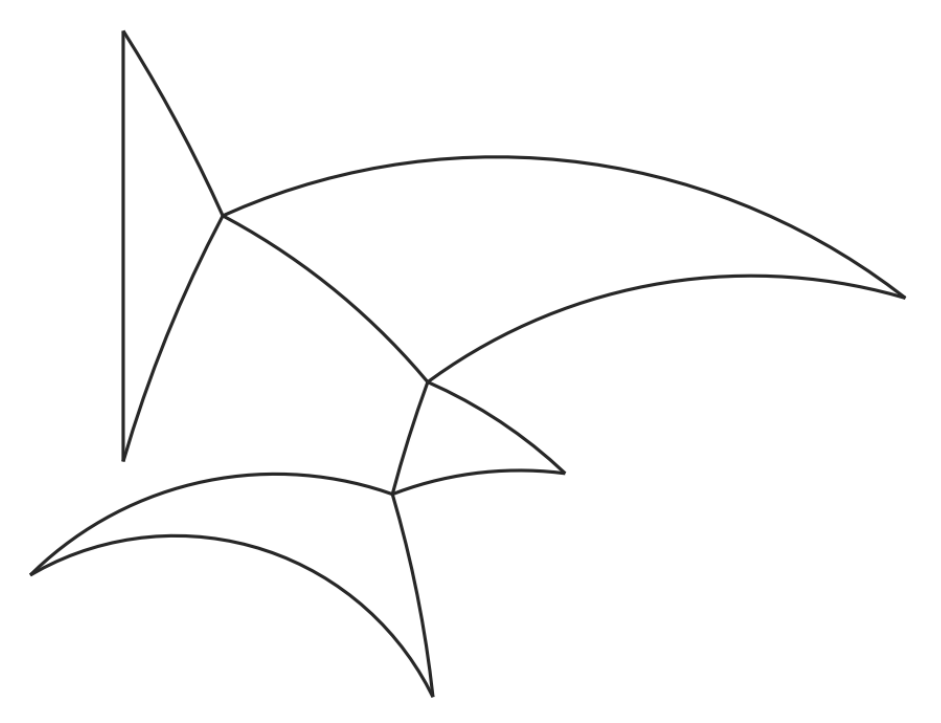}};

\begin{scope}
\fill (1.19,2.53) circle (0.07) node[left]{$C_1$};
\fill (1.19,6.62) circle (0.07) node[left]{$C_2$};
\fill (8.65,4.08) circle (0.07) node[below right]{$C_3$};
\fill (5.4,2.4) circle (0.07) node[below right]{$C_4$};
\fill (4.14,0.28) circle (0.07) node[below right]{$C_5$};
\fill (.31,1.43) circle (0.07) node[left]{$C_6$};
\end{scope}

\begin{scope}[apricot]
\fill (2.15,4.86) circle (0.07);
\fill (2.25,5.2) node{$B_1$};
\fill (4.1,3.27) circle (0.07) node[left]{$B_2$};
\fill (3.76,2.2) circle (0.07);
\fill (4.05,2.05) node{$B_3$};
\end{scope}
\end{tikzpicture}
}
\end{center}
The vertices $C_1,\ldots,C_n$ are called \emph{exterior vertices} and $B_1,\ldots,B_{n-3}$ are called \emph{shared vertices}.
\end{defn}
\cm{When we conjugate $\rho$ by an element of $\psl$, its triangle chain changes by a global orientation-preserving isometry of the hyperbolic plane. We will sometimes talk about the \emph{$\mathcal{B}$-triangle chain of $[\rho]$} which is thus only really defined up to orientation-preserving isometries of the hyperbolic plane.}
\begin{lem}[{\cite[Corollary~3.6]{action-angle}}]\label{lem:properties-triangle-chains}\cm{Assume that we are in the regime $\alpha_1+\cdots+\alpha_n> 2\pi(n-1)$, i.e.~that the peripheral angle vector $\alpha$ satisfies~\eqref{eq:angle-condition}. Then the important geometric features of a triangle chain are the following:
\begin{itemize}
    \item The interior angle at the exterior vertex $C_i$ is $\pi-\alpha_i/2$.
    \item The two interior angles on each side of a shared vertex add up to $\pi$.
    \item The vertices triples  $(C_1,C_2,B_1), (B_1,C_3,B_2),\ldots, (B_{n-3},C_{n-1},C_n)$ are clockwise oriented.
\end{itemize}
Moreover, any chain of $n-2$ triangles satisfying the above properties comes from a DT representation whose conjugacy class belongs to $\RepDT{\alpha}$.}
\end{lem}

\cm{A triangle inside a chain may be degenerate. This happens only when all of its three vertices coincide.} It's however not possible that all the triangles in a chain are degenerate since the total area of the triangles is $\lambda/2$, where $\lambda=\alpha_1+\cdots+\alpha_n-2\pi(n-1)$ is a positive constant. 
\begin{defn}\label{def:regular-singular-triangle-chains}
\cm{When a triangle in a chain is degenerate, the chain is said to be \emph{singular}. On the other hand, a chain is \emph{regular} if all $n-2$ triangles have positive area.}
\end{defn}
It turns out, however, that for any given DT representation $\rho$, it's always possible to find a chained pants decomposition $\mathcal{B}$ of $\SpherePk$ and a compatible presentation of $\SpherePk$ such that the $\mathcal{B}$-triangle chain of $\rho$ is regular.
\begin{lem}[{\cite[Proposition~2]{aaron-arnaud}}]\label{lem:existence-regular-triangle-chain}
For any DT representation $\rho$, there exists a geometric presentation of $\pi_1\SpherePk$ and standard chained pants decomposition $\mathcal{B}$ such that the $\mathcal{B}$-triangle chain of $\rho$ is regular.    
\end{lem}

\subsubsection{Action-angle coordinates}\label{sec:action-angle-coordinates}
The realization of DT representations as triangle chains can be used to produce action-angle coordinates for $\RepDT{\alpha}$, as explained in~\cite{action-angle}. \cm{As for triangle chains, the action-angle coordinates depend on a choice of chained pants decomposition $\mathcal{B}$ of $\SpherePk$ and a compatible geometric presentation of $\pi_1\SpherePk$ with generators $c_1,\ldots,c_n$ (Definition~\ref{defn:compatible-generators}), so that the pants curves of $\mathcal{B}$ are represented by the fundamental group elements $b_i=(c_1\cdots c_{i+1})^{-1}$ from~\eqref{eq:pants-curves}. For simplicity, we will assume that we are in the regime $\alpha_1+\cdots+\alpha_n>2\pi (n-1)$, i.e.~that the peripheral angle vector $\alpha$ satisfies~\eqref{eq:angle-condition}.}

The action coordinates $\beta_1,\ldots,\beta_{n-3}\in (0,2\pi)$ of $[\rho]$ are defined to be the rotation angles \cm{(Definition~\ref{defn:rotation-angle})} of $\rho(b_1),\ldots,\rho(b_{n-3})\in \psl$. Recall that the elements~$\rho(b_1),\ldots,\rho(b_{n-3})$ are all elliptic because DT representations are totally elliptic \cm{(Lemma~\ref{lem:totally-elliptic})}. It's another important feature of triangle chains that the interior angles on both sides of the shared vertex $B_i$ are equal to $\pi-\beta_i/2$ and $\beta_i/2$ \cm{(Lemma~\ref{lem:properties-triangle-chains})}, as shown on the picture below. The angle coordinates $\gamma_1,\ldots,\gamma_{n-3}\in \R/2\pi\Z$ of $[\rho]$ are then defined as the angles ``between'' consecutive triangles in the chain.
\begin{center}
\resizebox{9.9cm}{!}{
\begin{tikzpicture}[scale = 1.1, font=\sffamily,decoration={markings, mark=at position 1 with {\arrow{>}}}]
\node[anchor=south west,inner sep=0] at (0,0) {\includegraphics[width=9.9cm]{fig/fig-triangles-black}};

\begin{scope}
\fill (1.19,2.53) circle (0.07) node[left]{$C_1$};
\fill (1.19,6.62) circle (0.07) node[left]{$C_2$};
\fill (8.65,4.08) circle (0.07) node[below right]{$C_3$};
\fill (5.4,2.4) circle (0.07) node[below right]{$C_4$};
\fill (4.14,0.28) circle (0.07) node[below right]{$C_5$};
\fill (.31,1.43) circle (0.07) node[left]{$C_6$};
\end{scope}

\begin{scope}[apricot]
\fill (2.15,4.86) circle (0.07);
\fill (2.25,5.2) node{$B_1$};
\fill (4.1,3.27) circle (0.07) node[left]{$B_2$};
\fill (3.76,2.2) circle (0.07);
\fill (4.05,2.07) node{$B_3$};
\end{scope}
\draw[thick, postaction={decorate}, sky] (2.8,5.1) arc (10:115:.7) node[near end, above right]{$\gamma_1$};
\draw[thick, postaction={decorate}, sky] (4.7,2.97) arc (-20:32:.7) node[near end, below right]{$\gamma_2$};
\draw[thick, postaction={decorate}, sky] (3.95,1.6) arc (-80:15:.6) node[near end, below right]{$\gamma_3$};

\draw[thick, apricot] (1.9,5.35) arc (109:250:.5);
\draw[apricot] (2.4,4.2) node{\tiny $\pi-\beta_1/2$};
\draw[thick, plum] (2.6,4.62) arc (-30:21:.5);
\draw[plum] (3,4.8) node{\tiny $\beta_1/2$};
\draw[thick, apricot] (4.5,3.55) arc (25:130:.5) node[midway, above]{\tiny $\pi-\beta_2/2$};
\draw[thick, plum] (4.35,3.13) arc (-20:-103:.3);
\draw[plum] (4.4,2.83) node{\tiny $\beta_2/2$};
\draw[thick, apricot] (4.25,2.35) arc (5:78:.4) node[at end, left]{\tiny $\pi-\beta_3/2$};
\draw[thick, plum] (3.85,1.8) arc (-80:-192:.4);
\draw[plum] (3.2,1.8) node{\tiny $\beta_3/2$};
\end{tikzpicture}
}
\end{center}
It's important to remember that the angle coordinates $\gamma_1,\ldots,\gamma_{n-3}$ are only entirely defined when the triangle chain of $[\rho]$ is regular \cm{(Definition~\ref{def:regular-singular-triangle-chains})}. If some triangles are degenerate to a single point, then some of the $\gamma_1,\ldots,\gamma_{n-3}$ are not defined any more. For each degenerate triangle, we ``lose'' an angle coordinate.  The range of attainable values for the action coordinates $\beta_1,\ldots,\beta_{n-3}$ over $\RepDT{\alpha}$ is the (moment) polytope inside $\R^{n-3}$ defined by the $n-2$ inequalities
\begin{equation}\label{eq:inequalities-polytope-beta}
\begin{cases}
     \beta_1 \geq  4\pi-\alpha_1-\alpha_2, &  \\
     \beta_{i+1}-\beta_i \geq 2\pi-\alpha_{i+2}, & i=1,\ldots n-4,\\
     \beta_{n-3} \leq \alpha_{n-1}+\alpha_n-2\pi. 
\end{cases}
\end{equation}
The inequalities are the algebraic translation of the non-negativity of the area of each triangle in a chain. The fibres of the map $(\beta_1,\ldots,\beta_{n-3})\colon\RepDT{\alpha}\to \R^{n-3}$ over the above polytope define a Lagrangian toric fibration of $\RepDT{\alpha}$. 

\subsubsection{The action of Dehn twists}\label{sec:action-of-Dehn-twists}
\cm{The Dehn twists $\tau_{b_1},\ldots,\tau_{b_{n-3}}$ along the curves $b_1,\ldots,b_{n-3}$ of the pants decomposition $\mathcal{B}$ chosen to define the action-angle coordinates on $\RepDT{\alpha}$ are some of the few elements of the pure mapping class group $\PMod(\SpherePk)$ whose action on $\RepDT{\alpha}$ can be described geometrically using $\mathcal{B}$-triangle chains. Recall from~\eqref{eq:pants-curves} that $b_i=(c_1\cdots c_{i+1})^{-1}$, where $c_1,\ldots,c_n$ are the generators of $\pi_1\Sigma$ compatible with $\mathcal{B}$.}

\cm{The action of the Dehn twists $\tau_{b_i}$ can be described as follows. Since the simple closed curve on $\Sigma$ represented by $b_i$ is separating, the Dehn twist $\tau_{b_i}$ maps $[\rho]$ to the conjugacy class of the representation}
\begin{equation}\label{eq:Dehn-twist-action}
c_j\mapsto\begin{cases}
\rho(c_j) &\text{if } j\leq \cm{i+1},\\
\rho(b_i)\rho(c_j)\rho(b_i)^{-1} &\text{if } j>\cm{i+1}.
\end{cases}
\end{equation}
\cm{Formula~\eqref{eq:Dehn-twist-action} originates from Goldman's work~\cite[Section~4]{goldman-invariant}, and is stated in~\cite[Equation~4.8]{goldman-xia}. We'll write $\tau_{b_i}.[\rho]$ for the image of $[\rho]$ by $\tau_{b_i}$. }
\cm{\begin{lem}\label{lem:Dehn-twist-action-on-action-angle-coordinates}
	Denote by $(\beta_1,\ldots,\beta_{n-3},\gamma_1,\ldots,\gamma_{n-3})$ the action-angle coordinates of $[\rho]\in\RepDT{\alpha}$ obtained from the pants decomposition $\mathcal{B}$ and those of $\tau_{b_i}.[\rho]$ by~$(\beta_1',\ldots,\beta_{n-3}',\gamma_1',\ldots,\gamma_{n-3}')$. It then holds that
\[
\beta_j'=\beta_j\quad\text{and}\quad \gamma_j'=\begin{cases}\gamma_j&\text{if } j\neq i\\
\gamma_i-\beta_i&\text{if } j= i.
\end{cases}
\]
In particular, $\tau_{b_i}$ acts linearly in the action-angle coordinates. 
\end{lem}
\begin{proof}
We denote by $\rho'$ the representation~\eqref{eq:Dehn-twist-action} so that $\tau_{b_i}.[\rho]=[\rho']$. Let $C_1,\ldots, C_n$ be the exterior vertices in the $\mathcal{B}$-triangle chain of $\rho$.
By the definition of $\rho'$, the fixed points of $\rho'(c_j)$ is $C_j$ if $j\leq i+1$ and $\rho(b_i)C_j$ if $j>i+1$. In other words, the $\mathcal{B}$-triangle chain of $\rho'$ is obtained from the triangle chain of $\rho$ by an anti-clockwise
rotation of the sub-chain made of the triangles $(B_i,C_{i+2},B_{i+1}),\ldots, (B_{n-3},C_{n-1},C_n)$ by an angle $\beta_i$ around the vertex $B_i$, while leaving the triangles $(C_1,C_2,B_1), \ldots, (B_{i-1},C_{i+1},B_i)$ unaffected. This
geometric description of the action of the Dehn twist $\tau_{b_i}$ implies the desired conclusion about action-angle coordinates.
\end{proof}}

There is another set of Dehn twists whose action on $\RepDT{\alpha}$ can \cm{be computed using $\mathcal{B}$-triangle chains, but not as directly as in Lemma~\ref{lem:Dehn-twist-action-on-action-angle-coordinates}.} These are the Dehn twists along the simple closed curves represented by the fundamental group elements $c_ic_{i+1}\cdots c_{j-1}c_j$ for $1\leq i<j\leq n$. The interesting cases are when $j-i\leq n-2$, for otherwise the Dehn twists are trivial. The Dehn twist along the curve $c_ic_{i+1}\cdots c_{j-1}c_j$ will be denoted by 
\begin{equation}\label{eq:Dehn-twist-tau_i,j}
\cm{\tau_{i,j}\in \PMod(\SpherePk).}
\end{equation}
Observe that in this notation $\tau_{b_i}=\tau_{1,i+1}$. 

\cm{For a fixed index $i> 1$, in order to describe the action of the Dehn twists $\tau_{i,j}$ for $i<j\leq n$, we consider a new pants decomposition $\mathcal{B}[i]$ of $\Sigma$ given by the curves 
\begin{equation}\label{eq:pants-curve-b'}
\cm{b_k'=(c_i\cdots c_{i+k})^{-1} \quad\text{for}\quad k=1,\ldots, n-3,}
\end{equation}
where indices are taken modulo $n$. The curves $b_k'$ and the pants decomposition $\mathcal{B}[i]$ depend on the index $i$. To avoid having too many indices, we won't make the dependence on $i$ explicit in the notation $b_k'$. Note that $\tau_{i,j}$ is the Dehn twist along the curve $b_k'$ for $k=j-i$. The Dehn twists $\tau_{i,j}$ act linearly in the action-angle coordinates obtained from $\mathcal{B}[i]$ by Lemma~\ref{lem:Dehn-twist-action-on-action-angle-coordinates}, but their action is harder to express in terms of the action-angle coordinates obtained from $\mathcal{B}$ when $i>1$. In other words, $\mathcal{B}[i]$ provides the best coordinate system on $\RepDT{\alpha}$ to compute the action of the Dehn twists $\tau_{i,j}$ when $i>1$, while $\mathcal{B}=\mathcal{B}[1]$ is best-suited in the case $i=1$.}

 \cm{A geometric presentation of $\pi_1\Sigma$ which is compatible with $\mathcal{B}[i]$ is given by the cyclic permutation $c_i,c_{i+1},\ldots,c_n,c_1,\ldots, c_{i-1}$ of the geometric generators compatible with $\mathcal{B}$. This means that the $\mathcal{B}$-triangle chain and the $\mathcal{B}[i]$-triangle chain of a given $[\rho]\in\RepDT{\alpha}$ share the same exterior vertices $\{C_1,\ldots, C_n\}$ but have different shared vertices (see Definition~\ref{defn:triangle-chain}). }

\cm{In practice, if we are given the action-angle coordinates of some point $[\rho]\in\RepDT{\alpha}$ with respect to the pants decomposition $\mathcal{B}$ and we wish to compute the action-angle coordinates of $\tau_{i,j}.[\rho]$ (also with respect to the pants decomposition $\mathcal{B}$), then we can apply the following procedure. The idea is to change our coordinate system on $\RepDT{\alpha}$ from the one given by the pants decomposition $\mathcal{B}$ to the one given by $\mathcal{B}[i]$ because the action of $\tau_{i,j}$ is linear in the coordinates given by $\mathcal{B}[i]$ by Lemma~\ref{lem:Dehn-twist-action-on-action-angle-coordinates}. More precisely, we apply the following steps, which we illustrate in Example~\ref{ex:example-Dehn-twist-action} below.}
\begin{enumerate}
\item The action-angle coordinates of $[\rho]$ completely determine its $\mathcal{B}$-triangle chain which we call $T$. \cm{This triangle chain is well-suited to visualize the action of the Dehn twists $\tau_{b_1},\ldots,\tau_{b_{n-3}}$ (Lemma~\ref{lem:Dehn-twist-action-on-action-angle-coordinates}), but not the action of $\tau_{i,j}$ when $i\geq 2$.}
\item The next step consists \cm{of} constructing the $\mathcal{B}[i]$-triangle chain $T'$ of $[\rho]$. We know that $T'$ has the same exterior vertices as $T$, so we only need to find the shared vertices $B_1',\ldots,B_{n-3}'$ of $T'$. We also know \cm{from the properties of triangles chains stated in Lemma~\ref{lem:properties-triangle-chains}} that the first triangle in the chain $T'$ has clockwise oriented vertices $(C_i, C_{i+1}, B_1')$ with respective interior angles $(\pi-\alpha_i/2, \pi-\alpha_{i+1}/2, \pi-\beta_1'/2)$. The angle $\beta_1'$---the first action coordinate of $[\rho]$ with respect to $\mathcal{B}[i]$---is still unknown at this point. However, since \cm{$T$ determines the locations of the vertices $C_i$ and $C_{i+1}$, and since the two angles $\pi-\alpha_i/2$ and $\pi-\alpha_{i+1}/2$, as well as the orientation of the triangle $(C_i, C_{i+1}, B_1')$ are also given}, the vertex $B_1'$ is uniquely determined. Once $B_1'$ has been located, we deduce the value of $\beta_1'$. The second triangle in the chain $T'$ has clockwise oriented vertices $(B_1', C_{i+2}, B_2')$ and respective interior angles $(\beta_1'/2, \pi-\alpha_{i+2}/2, \pi-\beta_2'/2)$. As before, this data determines the vertex $B_2'$ uniquely, \cm{and thus the angle $\beta_2'$}. We can iterate the procedure until we have constructed all of $T'$.
\item Next, we construct the $\mathcal{B}[i]$-triangle chain of $\tau_{i,j}.[\rho]$ from $T'$ by rotating the \cm{sub-chain of} triangles to the right of $B_j'$ in an anti-clockwise fashion and by an angle $\beta_j'$, \cm{as we explained in the proof of Lemma~\ref{lem:Dehn-twist-action-on-action-angle-coordinates}.}
\item It's time to infer the $\mathcal{B}$-triangle chain of $\tau_{i,j}.[\rho]$ from its $\mathcal{B}[i]$-triangle chain by performing the same steps as in (2) once again. The action-angle coordinates of $\tau_{i,j}.[\rho]$ obtained from $\mathcal{B}$ can then be measured directly from its $\mathcal{B}$-triangle chain.
\end{enumerate}

\begin{ex}\label{ex:example-Dehn-twist-action}
We illustrate \cm{each of the four steps in the above procedure} to compute the image of a triangle chain by a Dehn twist on a concrete example in the case where $\Sigma$ is a 4-punctured sphere (another example can be found in the proof of Lemma~\ref{lem:beta_i-for-obits-of-type-33}). \cm{We fix some} geometric generators $c_1,c_2,c_3,c_4$ of $\pi_1\Sigma$ and let $\mathcal{B}$ be the associated standard pants decomposition of $\Sigma$. \cm{According to~\eqref{eq:pants-curves}, $\mathcal{B}$ consists of the only curve $b=b_1=(c_2c_1)^{-1}$ (we're dropping the index on $b$ to simplify the notation).} Let $\alpha=(12\pi/7, 12\pi/7, 10\pi/7, 12\pi/7)$. \cm{The DT component $\RepDT{\alpha}$ contains a point $[\rho]$ with action-angle coordinates $(\beta,\gamma)=(\beta_1,\gamma_1)=(\pi,3\pi/4)$.\footnote{As we'll see below in Section~\ref{sec:finite-orbits-of-type-8-33}, the point $[\rho]$ belongs to a finite mapping class group orbit ``of Type~8'' in Lisovyy--Tykhyy's nomenclature, also referred to as the ``Klein solution'' by Boalch. It's made of 7 orbit points.} The image of $[\rho]$ by $\tau_{1,2}=\tau_b$ is easy to compute: it's the point with coordinate $(\beta,\gamma)=(\pi,\pi/4)$ according to Lemma~\ref{lem:Dehn-twist-action-on-action-angle-coordinates}. A less straight-forward exercise is to compute the $(\beta,\gamma)$ action-angle coordinates of $\tau_{2,3}.[\rho]$. We'll do it step by step.}

\begin{enumerate}
\item \cm{First, we draw the $\mathcal{B}$-triangle chain $T$ of $[\rho]$.} It consists of two triangles $(C_1,C_2,B)$ and $(B,C_3,C_4)$. The triangle $(C_1,C_2,B)$ has respective interior angles $(\pi/7, \pi/7, \pi/2)$ and the triangle $(B,C_3,C_4)$ has interior angles $(\pi/2,2\pi/7,\pi/7)$. \cm{It turns out that $T$ fits well on the grid given by a $(2,3,7)$-tessellation of the hyperbolic plane in the sense of Example~\ref{ex:discrete-DT}: the vertices of $T$ are tessellation vertices and the interior angles of $T$ are integer multiples of the tessellation angle at the vertex. The configuration is depicted by the following illustration.} For convenience, we switch to the Poincaré disk model. 
\begin{center}
\begin{tikzpicture}[font=\sffamily,decoration={markings, mark=at position 1 with {\arrow{>}}}]
\node[anchor=south west,inner sep=0] at (0,0) {\includegraphics[width=8cm]{fig/example-B-triangle-chain-1}};

\fill (4.4,2.2) circle (0.07) node[left]{$C_1$};
\fill (2.9,5.4) circle (0.07) node[left]{$C_2$};
\fill (4.85,4.45) circle (0.07) node[above]{$B$};
\fill (5.65,4.8) circle (0.07) node[right]{$C_3$};
\fill (5.6,3.25) circle (0.07) node[left]{$C_4$};
\end{tikzpicture}
\end{center}

\item \cm{Since we are interested in the Dehn twist $\tau_{2,3}$, we have $i=2$ and $j=3$ in the notation~\eqref{eq:Dehn-twist-tau_i,j}, and so the adapted pants decomposition $\mathcal{B}[2]$ we want to consider consists of one curve: $b'=b'_1=(c_2c_3)^{-1}$ (take $i=2$ and $k=1$ in~\eqref{eq:pants-curve-b'}).} We shall now describe the $\mathcal{B}[2]$-triangle chain $T'$ of $[\rho]$. Since both $T$ and $T'$ have the same exterior vertices, all we need to do is to find the shared vertex $B'$ of $T'$. The two triangles of $T'$ have vertices $(C_2,C_3,B')$ and $(B',C_4,C_1)$, with respective interior angles $(\pi/7, 2\pi/7, \pi-\beta'/2)$ and $(\beta'/2, \pi/7, \pi/7)$. Using the grid, we can find the unique point $B'$ that satisfies those angle conditions. Once we've found $B'$, we can deduce that $\beta'=4\pi/3$.
\begin{center}
\begin{tikzpicture}[font=\sffamily,decoration={markings, mark=at position 1 with {\arrow{>}}}]
\node[anchor=south west,inner sep=0] at (0,0) {\includegraphics[width=8cm]{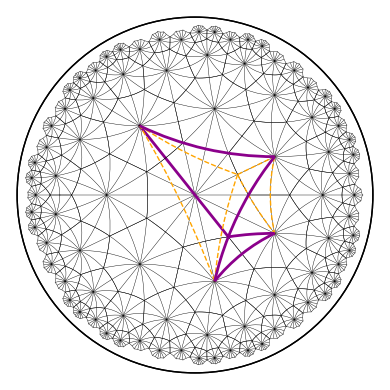}};

\fill (4.4,2.2) circle (0.07) node[left]{$C_1$};
\fill (2.9,5.4) circle (0.07) node[left]{$C_2$};
\fill (4.67,3.15) circle (0.07) node[left]{$B'$};
\fill (5.65,4.8) circle (0.07) node[right]{$C_3$};
\fill (5.6,3.25) circle (0.07) node[right]{$C_4$};
\end{tikzpicture}
\end{center}

\item The action of $\tau_{2,3}$ on $[\rho]$ is easy to describe using $T'$, \cm{as we explained in the proof of Lemma~\ref{lem:Dehn-twist-action-on-action-angle-coordinates}.} It simply rotates the triangle $(B',C_4,C_1)$ anti-clockwise by an angle $\beta'=4\pi/3$ around $B'$. The resulting triangle chain is the $\mathcal{B}[2]$-triangle chain of $\tau_{2,3}.[\rho]$ and it consists of two triangles with vertices $(C_2,C_3,B')$ and $(B', C_4^{\textrm{new}}, C_1^{\textrm{new}})$.
\begin{center}
\begin{tikzpicture}[font=\sffamily,decoration={markings, mark=at position 1 with {\arrow{>}}}]
\node[anchor=south west,inner sep=0] at (0,0) {\includegraphics[width=8cm]{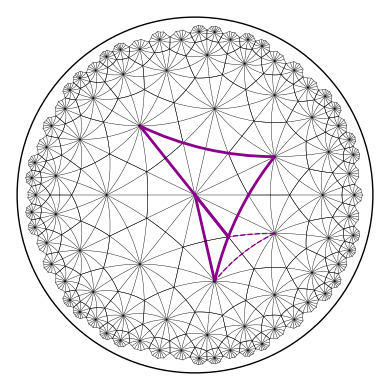}};

\fill (4.4,2.2) circle (0.07) node[left]{$C_4^{\textrm{new}}$};
\fill (2.9,5.4) circle (0.07) node[left]{$C_2$};
\fill (4.67,3.15) circle (0.07);
\draw (4.85,2.9) node{$B'$};
\fill (5.65,4.8) circle (0.07) node[right]{$C_3$};
\fill (4,4) circle (0.07) node[left]{$C_1^{\textrm{new}}$};

\draw[thick, postaction={decorate}, mauve] (5.05,3) arc (0:145:.4);
\end{tikzpicture}
\end{center}

\item The last step consists in \cm{determining} the $\mathcal{B}$-triangle chain of $\tau_{2,3}.[\rho]$ from its $\mathcal{B}[2]$-triangle chain. We know that both chains share the same exterior vertices, so that the $\mathcal{B}$-triangle chain of $\tau_{2,3}.[\rho]$ is made of the two triangles $(C_1^{\textrm{new}},C_2,B^{\textrm{new}})$ and $(B^{\textrm{new}}, C_3, C_4^{\textrm{new}})$. Their respective angles are $(\pi/7,\pi/7,\pi-\beta^{\textrm{new}}/2)$ and $(\beta^{\textrm{new}}/2,2\pi/7, \pi/7)$. These angle relations determine the location of $B^{\textrm{new}}$ uniquely, \cm{and it can be determined with the help of the grid.} \cm{From the resulting $\mathcal{B}$-triangle chain of $\tau_{2,3}.[\rho]$,} we deduce that the action-angle coordinates of $\tau_{2,3}.[\rho]$ \cm{with respect to the pants decomposition} $\mathcal{B}$ are $(\beta^{\textrm{new}},\gamma^{\textrm{new}})=(2\pi/3,\pi)$.
\begin{center}
\begin{tikzpicture}[font=\sffamily,decoration={markings, mark=at position 1 with {\arrow{>}}}]
\node[anchor=south west,inner sep=0] at (0,0) {\includegraphics[width=8cm]{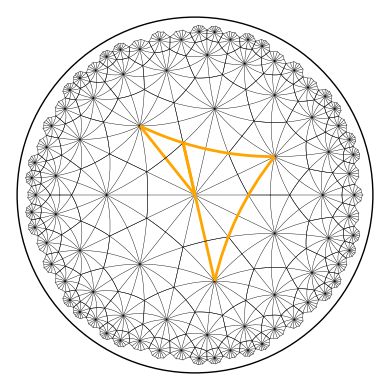}};

\fill (4.4,2.2) circle (0.07) node[left]{$C_4^{\textrm{new}}$};
\fill (2.9,5.4) circle (0.07) node[left]{$C_2$};
\fill (3.8,5.1) circle (0.07) node[above]{$B^{\textrm{new}}$};
\fill (5.65,4.8) circle (0.07) node[right]{$C_3$};
\fill (4,4) circle (0.07) node[left]{$C_1^{\textrm{new}}$};

\draw[thick, postaction={decorate}, sky] (4.6,4.9) arc (-5:150:.75);
\node[sky] (N) at (8,8) {$\gamma^{\textrm{new}}$};
\draw (N.west) edge[out=180,in=65,->] (4.1, 5.8);

\draw[thick, apricot] (4.4,4.9) arc (-5:-80:.6);
\node[apricot] (M) at (8,0) {$\beta^{\textrm{new}}/2$};
\draw (M.west) edge[out=180,in=-65,->] (4.3,4.5);
\end{tikzpicture}
\end{center}
\end{enumerate}
\ensuremath{\lozenge}
\end{ex}

We've just explained how to compute the action-angle coordinates of the image of $[\rho]$ by any of the Dehn twists $\tau_{i,j}$ \cm{for $1\leq i<j\leq n$} in a purely geometric manner using triangle chains. The advantage of this method is that it makes it possible to work out some trigonometric formulae to express the coordinates of $\tau_{i,j}.[\rho]$ in terms of those of $[\rho]$, determining in this way the exact values of the action-angles coordinates of $\tau_{i,j}.[\rho]$.

There is \cm{another} remarkable fact about the family of Dehn twists $\{\tau_{i,j}:1\leq i<j\leq n\}$. It turns out that this set is also a generating family for $\PMod(\Sigma)$. \cm{This will be shown later in Lemma~\ref{lem:mcg-generators} in Appendix~\ref{apx:generators-of-pmod} as a consequence of a result of Ghaswala--Winarski~\cite{mcg-generators}. Consequently, the routine above is all we need to compute entire mapping class group orbits using $\mathcal{B}$-triangle chains.} In practice, computing a whole orbit with this method can be tedious for a human being. For this reason, we'll sometimes favour a slightly different method which we can easily run on a computer and which will help us identify all orbit points. A description of this method can be found in Appendix~\ref{apx:algo-orbit}.

Before ending this recap on DT representations and Dehn twists actions, we point out one last useful fact. \cm{It's a characterization of the fixed points of $\tau_{i,j}$ in $\RepDT{\alpha}$ in terms of triangle chains.}
\cm{\begin{lem}\label{lem:fixed-points-Dehn-twists}
Let $C_1,\ldots,C_n$ denote the exterior vertices of the $\mathcal{B}$-triangle chain of $[\rho]\in\RepDT{\alpha}$. Given two indices $1\leq i<j\leq n$, the Dehn twist $\tau_{i,j}$ fixes $[\rho]$ if and only if $C_i=\cdots = C_j$ or $C_{j+1}=\cdots=C_n=C_1=\cdots=C_{i-1}$. In particular, if the $\mathcal{B}[i]$-triangle chain of $[\rho]$ is regular (Definition~\ref{def:regular-singular-triangle-chains}), then $\tau_{i,j}$ doesn't fix $[\rho]$.
\end{lem}}
\begin{proof}
\cm{We will consider the fundamental group element $b=(c_i\cdots c_{j})^{-1}$. Since $b$ represents a simple closed curve on $\Sigma$, its image $\rho(b)$ is a non-trivial elliptic element of $\psl$ by total ellipticity of DT representations (Lemma~\ref{lem:totally-elliptic}). We can adapt~\eqref{eq:Dehn-twist-action} to write the image of $[\rho]$ by $\tau_{i,j}$ as the conjugacy class of the representation 
\begin{equation}\label{eq:Dehn-twist-action-general}
\rho'(c_k)=\begin{cases}
\rho(c_k)&\text{if } i\leq k\leq j,\\
\rho(b)\rho(c_k)\rho(b)^{-1}&\text{if } k<i \text{ or } j<k.
\end{cases}
\end{equation}
}

\cm{If the points $C_i,\ldots,C_j$ coincide, or if the points $C_{j+1},\ldots,C_n,C_1,\ldots,C_{i-1}$ coincide, say with a point $X$, then $\rho(b)$ fixes $X$. In that case, $\rho$ and $\rho'$ are conjugate by $\rho(b)$. This implies that $\tau_{i,j}$ fixes $[\rho]$.}

\cm{Conversely, if $\tau_{i,j}.[\rho]=[\rho]$, then $\rho$ and $\rho'$ are conjugate. In other words, there exists $A\in\psl$ such that $\rho=A\rho' A^{-1}$. This happens if and only if $A\in Z(\rho(c_i))\cap\cdots\cap Z(\rho(c_{j}))$ and $A\rho(b)\in Z(\rho(c_{1}))\cap\cdots\cap Z(\rho(c_{i-1}))\cap Z(\rho(c_{j+1}))\cap\cdots\cap Z(\rho(c_n))$, where $Z(M)\subset\psl$ denotes the commutator of $M$. Now, recall that if $M$ and $M'$ are two non-trivial elliptic elements of $\psl$ such that $Z(M)\cap Z(M')\neq \{\id\}$, then $M$ and $M'$ have the same unique fixed point in the hyperbolic plane. So, if $A\neq \id$, then $\rho(c_i),\ldots,\rho(c_{j})$ have the same fixed point, i.e.~$C_i=\cdots=C_{j}$. On the other hand, if $A=\id$, then $A\rho(b)\neq \id$ because $\rho(b)$ is a non-trivial elliptic element and thus $C_{j+1}=\cdots=C_n=C_1=\cdots=C_{i-1}$.}

\cm{Finally, if the $\mathcal{B}[i]$-triangle chain of $[\rho]$ is regular, then in particular the first and last triangles in the chain are non-degenerate. This means that $C_i\neq C_{i+1}$ and $C_{i-2}\neq C_{i-1}$, where indices are taken modulo $n$. By the previous characterization, we conclude that $\tau_{i,j}.[\rho]\neq [\rho]$.}
\end{proof}

\subsection{Trigonometric number fields}\label{sec:number-fields}
In several places throughout the paper, we'll need to answer questions of the form: what are all
the numbers $r\in (0,2\pi)$ that are rational multiples of $\pi$ and such that $2\cos(r/2)$
belongs to some fixed number field $K=\Q(\zeta)$ ($\zeta$ will always be a real algebraic number).
For instance, in the simplest case where $K=\Q$, the answer is given by Niven's Theorem~\cite[Corollary~3.12]{niven}. It
asserts that if $2\cos(r/2)$ is rational and $r\in \pi\Q\cap (0,2\pi)$, then $2\cos(r/2)$ is
actually one of the numbers $\{-\frac{1}{2},0,\frac{1}{2}\}$ and $r\in \{2\pi/3, \pi, 4\pi/3\}$. 
For more general number fields, an answer can be found in a note of Lehmer~\cite{lehmer}
which boils down to tabulating the integers $n$ such that $\varphi(n)\leq 2\deg(K)$, where $\varphi$ is the Euler $\varphi$-function.

Panraksa--Samart--Sriwongsa~\cite{niven-generalization} recently introduced an algorithm to find all $r\in\pi\Q$ such that $2\cos(r/2)$ is in a number field of degree~$D$. Here's how their algorithm works. Consider the real function $f\colon [-2,2]\rightarrow[-2,2]$ defined by $f(x)=x^2-2$. The function $f$ being polynomial, it restricts to a self-function of any number field $K=\Q(\zeta)$ where $\zeta$ is a real algebraic number.
\begin{thm}[\cite{niven-generalization}]\label{thm:rational-cosines}
	Let $K$ be a number field of degree~$D$. Then we have
	\begin{equation*}
		2\cos(\pi\Q)\cap K=\mathrm{PrePer} (f,K),
	\end{equation*}
	where $\mathrm{PrePer}(f,K)$ stands for the set of pre-periodic points\footnote{A number in $K$ is said to be \emph{pre-periodic} is some iterate of it under $f$ is periodic.} of $f$ in $K$.
	Moreover, any $x\in\mathrm{Per}(f,K)$ satisfies $f^{D}(x)=0$, where $f^D$ refers to the $D^\textrm{th}$ iterate of $f$.
\end{thm}
Obviously, if $y\in\mathrm{PrePer}(f,K)$, then so do all iterates $f^{n}(y)$ of $y$ for
any~$n\geq 0$. 
\cm{Note that $\mathrm{PrePer} (f,K)$ is a finite set by a result of Northcott~\cite{northcott}}. 
In order to identify all $r\in\pi\Q\cap (0,2\pi)$ such that $2\cos(r/2)\in K$, the following procedure can be applied. 
\begin{enumerate}
\item First, factorize the polynomial $f^{D}(x)-x$ in irreducible factors and compute the roots of all factors of degree at most $D$. Discard all the roots that do not belong to $K$. By Theorem~\ref{thm:rational-cosines}, the resulting collection of roots correspond to the periodic points of $f$ in $K$. 
\item For each root selected in (1), compute the tree of its successive preimages by $f$. Stop computing preimages along a branch as soon as one preimage doesn't belong to $K$. The resulting set of numbers, along with the corresponding roots from (1), are all the pre-periodic points of $f$. By Theorem~\ref{thm:rational-cosines}, these are all the possible values of $2\cos(r/2)$ inside $K$ when $r\in \pi\Q\cap (0,2\pi)$. 
\item Finally, for all possible value of $2\cos(r/2)$ that belongs to $K$, compute the corresponding value of $r\in \pi\Q\cap (0,2\pi)$.
\end{enumerate} 
We provide an example of the above procedure in Example~\ref{ex:trigonometric-number-fields} below.

For the purpose of this paper, $\zeta$ will always be of the form $\cos(\pi/N)$ for some $N\geq 3$. In that case, the degree of $\zeta$ is related to Euler $\varphi$-function. 
\begin{fact}\label{fact:degree-cos-pi-over-N}
For any integer~$N\geq 3$, the number $\cos(2\pi/N)$
is algebraic of degree $\varphi(N)/2$. So, if $N$ is odd, then the degree of $\cos(\pi/N)$ is $\varphi(N)/2$, and if $N$ is even, the degree of $\cos(\pi/N)$ is $\varphi(N)$. In particular, if $N\geq 3$ is odd, then $\cos(\pi/2N)$ does not belong to $\Q(\cos(\pi/N))$.
\end{fact}
\begin{proof}
\cm{This fact is standard (see e.g.~\cite{lang}); we provide a short argument for the sake of completeness. Let $\omega_N=e^{2\pi i/N}$ and note that $2\cos(2\pi/N)=\omega_N+\omega_N^{-1}$. The cyclotomic field $\Q(\omega_N)$ has degree $[\Q(\omega_N):\Q]=\varphi(N)$~\cite[Chapter~\Romannum{4}, Theorem~2]{lang}. If we pick any embedding of $\Q(e^{2\pi i/N})$ in $\C$, then the fixed point set of the complex conjugation is the maximal real subfield $\Q(\omega_N + \omega_N^{-1})=\Q(\cos(2\pi/N))\subset \Q(\omega_N)$. In particular, $[\Q(\omega_N):\Q(\cos(2\pi/N))]=2$ and thus $[\Q(\cos(2\pi/N)):\Q]=\varphi(N)/2$. We conclude that the degree of $\cos(2\pi/N)$ is $\varphi(N)/2$ and the degree of $\cos(\pi/N)$ is $\varphi(2N)/2$.}
\end{proof}
The last statement in Fact~\ref{fact:degree-cos-pi-over-N} is useful in the second step of the above algorithm since the two preimages by $f$ of $2\cos(\pi/N)$ are $\pm 2\cos(\pi/2N)$. We'll also need a criterion to decide whether $\cos(k\pi/a)$ belongs to $\Q(\cos(\pi/b))$. Using Chebyshev polynomials, we observe that when $\gcd(a,k)=1$, $\cos(k\pi/a) \in\Q(\cos(\pi/b))$ if and only if $\cos(\pi/a)\in\Q(\cos(\pi/b))$. Furthermore, if $\gcd(a,b)=1$, then we also have the following fact.
\begin{fact}\label{fact:intersection-number-fields}
If $a,b\geq 3$ are two integers, then 
\[
\Q(\cos(\pi/a))\cap\Q(\cos(\pi/b))=\Q\big(\cos(\pi/\gcd(a,b))\big).
\]
In particular, if $\gcd(a,b)=1$, then $\cos(\pi/a)$ belongs to $\Q(\cos(\pi/b))$ if and only if $\cos(\pi/a)\in\Q$, which is equivalent to $a\in \{1,2,3\}$ by Niven's Theorem.
\end{fact}
\begin{proof}
\cm{Again, let $\omega_{N}=e^{2\pi i/N}$. The intersection formula for cyclotomic fields reads $\Q(\omega_{N})\cap \Q(\omega_{M})=\Q(\omega_{\gcd(N,M)})$. Letting $N=2a$ and $M=2b$, and taking maximal real subfields leads to the desired identity.}
\end{proof}

We anticipate some of the upcoming computations and apply Panraksa--Samart--Sriwongsa's algorithm to a selection of five number fields $K$ of degree $D$. Here's what we obtain.
\begin{table}[ht]
    \centering
    \begin{tblr}{c|c|c|c}
        $K$ & $D$ & $2\cos(\pi\Q)\cap K$ & $\{r\in\pi\Q\cap(0,2\pi):2\cos(r/2)\in K\}$ \\
        \hline\hline
        $\Q(\sqrt{2})$ & 2 & $\{\pm 2, \pm 1,0,\pm \sqrt{2}\}$ & $\{\frac{\pi}{2},\frac{2\pi}{3},\pi,\frac{4\pi}{3},\frac{3\pi}{2}\}$ \\
        \hline
        $\Q(\sqrt{5})$ & 2 & $\{\pm 2,\pm 1,0,\frac{\pm 1\pm\sqrt{5}}{2}\}$ & $\{\frac{2\pi}{5},\frac{2\pi}{3},\frac{4\pi}{5},\pi,\frac{6\pi}{5},\frac{4\pi}{3},\frac{8\pi}{5}\}$ \\
        \hline
        $\Q(\cos(\pi/7))$ & 3 & $\{\pm 1, 2\cos(\frac{n\pi}{7}) : n\in\Z\}$ & $\{\frac{2\pi}{7},\frac{4\pi}{7},\frac{2\pi}{3},\frac{6\pi}{7},
		\pi,\frac{8\pi}{7},\frac{4\pi}{3},\frac{10\pi}{7},\frac{12\pi}{7}\}$ \\
        \hline
        $\Q(\cos(\pi/9))$ & 3 & $\{0, 2\cos(\frac{n\pi}{9}):n\in\Z\}$ & $\{\frac{2\pi}{9},\frac{4\pi}{9},\frac{2\pi}{3},\frac{8\pi}{9},
		\pi,\frac{10\pi}{9},\frac{4\pi}{3},\frac{14\pi}{9},\frac{16\pi}{9}\}$ \\
        \hline
        $\Q(\cos(\pi/18))$ & 6 & $\{2\cos(\frac{n\pi}{18}):n\in\Z\}$ & $\{\frac{n\pi}{9}:n=1,\ldots,17\}$ \\
        \end{tblr}
        \caption{The list of all numbers $r\in(0,2\pi)$ that are rational multiples of $\pi$ and such that $2\cos(r/2)$ belongs to some number field $K$ of degree $D$.}
        \label{tab:rational-angles}
\end{table}

\begin{ex}\label{ex:trigonometric-number-fields}
We sketch the details of the computations leading to Table~\ref{tab:rational-angles} for $K=\Q(\cos(\pi/7))$ (see also~\cite[Examples~15--19]{niven-generalization} for more examples).
\cm{Let $f(x)=x^2-2$.} Note that $K$ is a cubic field by Fact~\ref{fact:degree-cos-pi-over-N}. So, we start by factoring the polynomial $f^{3}(x)-x$ which gives
\begin{equation*}
	f^{3}(x)-x=(x-2)(x+1)(x^3-3x+1)(x^3+x^2-2x-1).
\end{equation*}
The eight roots are $2$, $-1$, $2\cos(2\pi/9)$, $2\cos(4\pi/9)$, $2\cos(8\pi/9)$, $2\cos(2\pi/7)$, 
$2\cos(4\pi/7)$, and $2\cos(6\pi/7)$. Among these roots, the only ones in $K$
are $2$, $-1$, $2\cos(2\pi/7)$, $2\cos(4\pi/7)$, and $2\cos(6\pi/7)$ because $\gcd(7,9)=1$ (Fact~\ref{fact:intersection-number-fields}). For each root, we compute its succesive preimages by $f$ and stop as soon as a preimage does not belong to~$K$. For instance, the preimages of $2\cos(2\pi/7)$ are $2\cos(\pi/7)$ and $2\cos(6\pi/7)$ (which we already accounted for). Furthermore, the preimages of $2\cos(\pi/7)$ are $\pm 2\cos(\pi/14)$ which don't belong to $K$ by Fact~\ref{fact:degree-cos-pi-over-N}, so we stop. Similarly, by considering preimages of $2\cos(4\pi/7)$ and $2\cos(6\pi/7)$, we obtain the new values $2\cos(5\pi/7)$, respectively $2\cos(3\pi/7)$. The resulting numbers are listed in the third column of Table~\ref{tab:rational-angles} and the corresponding rational multiple of $r$ in the fourth column. 
\ensuremath{\lozenge}
\end{ex}

\subsection{Triangle groups}\label{sec:triangle-groups}
As we'll see, the image of a DT representation coming from a finite mapping class group orbit will often turn out to be a triangle group. To fix some notation, we include a short recap about triangle groups. Let $(p,q,r)$ be a triple of positive real numbers such that $1/p+1/q+1/r<1$. Up to hyperbolic isometries, there is a unique triangle $T_{p,q,r}$ in the hyperbolic plane with interior angles $(\pi/p,\pi/q,\pi/r)$. There are two groups of hyperbolic isometries associated with the triangle $T_{p,q,r}$. 
\begin{enumerate}
    \item The first one is the subgroup of $\pgl$ generated by the (orientation-reversing) reflections through the sides of $T_{p,q,r}$. It has various name in the literature, one being the \emph{reflection triangle group} of $T_{p,q,r}$. It is commonly denoted by $\Delta(p,q,r)$. 
    \item The second one is the subgroup of $\psl$ generated by the three rotations around the vertices of the triangle whose rotation angles are twice the interior angle at the corresponding vertex. It's called the \emph{rotation triangle group} of $T_{p,q,r}$. We'll denote it by $D(p,q,r)$.
\end{enumerate}
The rotation triangle group $D(p,q,r)$ can also be defined as the index-2 subgroup of the reflection triangle group $\Delta(p,q,r)$ generated by products of two reflections. 

If the image of a representation into $\psl$ is a finite-index subgroup of a rotation triangle group
(typically for DT representations coming from a finite mapping class orbit), then the question
of its discreteness boils down to the discreteness of the triangle group. Felikson~\cite{felikson}
established the exhaustive list of discrete triangle groups, so we will be able
to determine when a DT representation is discrete. Conversely, any discrete DT representation
belongs to a finite orbit of the mapping class group action, as we'll explain in Section~\ref{sec:discrete-finite-orbits}.
\begin{thm}[\cite{felikson}]\label{thm:felikson}
The reflection triangle group $\Delta(p,q,r)$ is a discrete subgroup of $\pgl$, or equivalently the rotation triangle group $D(p,q,r)$ is a discrete subgroup of $\psl$, if and only if the triple $(p,q,r)$ can be found in the first column of Table~\ref{tab:discrete-reflection-groups}. The numbers $a$, $b$, and $c$ always denote positive integers.
\begin{table}[ht]
    \centering
    \begin{tblr}{c|c|c}
        Triangle & \makecell{Hyperbolicity condition\\ on $a,b,c\in \N$} & Rotation triangle group\\
        \hline\hline
        $(a,b,c)$ & $1/a+1/b+1/c<1$ & $D(a,b,c)$\\
        \hline
        $(\frac{a}{2},b,b)$ & $1/a+1/b<1/2$ & $D(\frac{a}{2},b,b)=D(2,a,b)$, if $\gcd(a,2)=1$ \\
        \hline
        $(2,\frac{a}{2},a)$ & $a\geq 7$ & $D(2,\frac{a}{2},a)=D(2,3,a)$, if $\gcd(a,2)=1$\\
        \hline
        $(\frac{a}{2},a,a)$ & $a\geq 7$ & $D(\frac{a}{2},a,a)=D(2,4,a)$, if $\gcd(a,2)=1$\\
        \hline
        $(3,\frac{a}{3},a)$ & $a\geq 7$ & $D(3,\frac{a}{3},a)=D(2,3,a)$, if $\gcd(a,3)=1$\\
        \hline
        $(\frac{a}{4},a,a)$ & $a\geq 7$ & $D(\frac{a}{4},a,a)=D(2,3,a)$, if $\gcd(a,4)=1$\\
        \hline
        $(\frac{a}{2},\frac{a}{2},\frac{a}{2})$ & $a\geq 7$ & $D(\frac{a}{2},\frac{a}{2},\frac{a}{2})=D(2,3,a)$, if $\gcd(a,2)=1$\\
        \hline
        $(3,\frac{7}{2},7)$ & none & $D(3,\frac{7}{2},7)=D(2,3,7)$\\
    \end{tblr}
    \caption{Felikson's list of hyperbolic triangles with discrete reflection groups.}
    \label{tab:discrete-reflection-groups}
\end{table}
\end{thm}

\section{Examples of finite mapping class group orbits}\label{chap:examples}

\subsection{Overview}
We give a short recap about several different types of representations $\rho\colon \pi_1\Sigma\to \SL_2\C$ whose conjugacy class $[\rho]$ is known to belong to a finite mapping class group orbit. The material presented in this section is classical; we include it for the sake of completeness and to help the reader develop an intuition about finite mapping class group orbits. \cm{It is not required for the proof of the main theorems and may be skipped.} We cover representations with abelian image (Section~\ref{sec:abelian-image}), finite image (Section~\ref{sec:finite-image}), and discrete image (Section~\ref{sec:discrete-finite-orbits}). \cm{The case of totally non-hyperbolic representations is discussed in Section~\ref{sec:totally-non-hyperbolic}.} We introduce finite orbits of ``pullback'' types in Section~\ref{sec:pullback-orbits} and we provide examples of Zariski dense representations giving rise to arbitrarily long finite orbits (Section~\ref{sec:arbitrarily-long-orbits}), as well as examples with non-Zariski dense image (Section~\ref{sec:non-zariski-dense-orbits}).

\subsection{Representations with abelian image}\label{sec:abelian-image}
The most basic example of finite mapping class group orbits are arguably induced by representations with abelian image. 
\begin{lem}\label{lem:abelian-implies-fixed-point}
\cm{Let $\Sigma$ be a punctured sphere.} The conjugacy class of a representation $\rho\colon\pi_1\Sigma\to\SL_2\C$ whose image lies in an abelian subgroup of $\SL_2\C$ is always fixed globally by the mapping class group action. (In this statement, $\SL_2\C$ can be replaced by any other Lie group, but it's important that $\Sigma$ is a punctured sphere.)
\end{lem}
\begin{proof}
It's convenient to fix an auxiliary geometric presentation of $\pi_1\Sigma$ with generators $c_1,\ldots,c_n$. As we explained in Section~\ref{sec:character-varieties}, the pure mapping class group of $\Sigma$ acts on $[\rho]$ by pre-composing $\rho$ with an automorphism of $\pi_1\Sigma$ that conjugates each generator $c_i$ individually. If $\varphi$ denotes that automorphism of $\pi_1\Sigma$, then $\rho(\varphi^{-1}(c_i))$ is conjugate to $\rho(c_i)$ by an element in the image of $\rho$. Since we're assuming that the image of $\rho$ is abelian, we have $\rho(\varphi^{-1}(c_i))=\rho(c_i)$. We conclude that $[\rho]$ is globally fixed by the mapping class group action.
\end{proof}

Not all the fixed points of the mapping class group action come from representations with abelian image \cm{as we'll see in Lemma~\ref{lem:fixed-point-implies-abelian}.} For instance, any representation of a 4-punctured sphere that sends a peripheral loop to $\pm\id$ gives rise to a finite orbit. However, the converse of Lemma~\ref{lem:abelian-implies-fixed-point} holds when \cm{sufficiently many peripheral monodromies are non-scalar}. (For the converse to be true, it's important that the target group is $\SL_2\C$.) The proof relies on the following fact about centralizers in $\SL_2\C$.

\begin{fact}\label{fact:centralizers-SL2C}
If $A,B\in\SL_2\C\setminus\{\pm\id\}$ are two elements with distinct centralizers $Z(A)\neq Z(B)$, then their intersection $Z(A)\cap Z(B)$ is equal to $Z(\SL_2\C)=\{\pm \id\}$.
\end{fact}
\begin{proof}
\sam{Let $A,B\in\SL_2\C\setminus\{\pm\id\}$. The centralizer of $A$ is always the product of the one-parameter subgroup of $\SL_2\C$ containing $A$ with $\pm\id$.
	This is clear if $A$ is diagonalizable, and still true when $A$ is not as then it is conjugate, up to multiplication by $-\id$ to the standard \cm{parabolic} element of $\SL_2\C$.
	Now, the intersection of $Z(A)$ and $Z(B)$ is the product of $\pm\id$ with the intersection of two one-parameter subgroups of $\SL_2\C$, so it is equal to $Z(A)$ or to $\pm\id$,
	as claimed.}
\end{proof}

\begin{lem}\label{lem:fixed-point-implies-abelian}
	\sam{Let $\Sigma$ be a punctured sphere with at least four punctures and $\rho\colon \pi_1\Sigma\to\SL_2\C$ be a representation. The conjugacy class of $\rho$ is a fixed point of the mapping class group action if and only if one of the following properties is verified:
	\begin{itemize}
		\item All but at most three peripheral loops of $\Sigma$ are mapped by $\rho$ to $\pm\id$.
		\item The image of $\rho$ is abelian.
	\end{itemize}}
\end{lem}
\begin{proof}
	\sam{As before, it's convenient to fix an auxiliary geometric presentation of $\pi_1\Sigma$ with generators $c_1,\ldots,c_n$.
	If $\rho$ has abelian image, then $[\rho]$ is a fixed point of the mapping class group action by Lemma~\ref{lem:abelian-implies-fixed-point}. Now, assume that all but at most three of the generators $c_1,\ldots,c_n$ are mapped to $\pm\id$.}
	\begin{claim}
	\sam{The Dehn twist $\sigma_{i,j}$ along the simple closed curve represented by the fundamental group element $c_ic_j$ fixes $[\rho]$.}
	\end{claim}
	\begin{proofclaim}
	\sam{The simple closed curve represented by $c_ic_j$ separates $\Sigma$ into two sub-spheres $\Sigma_1\cup\Sigma_2$. There is one sub-sphere, say $\Sigma_1$, such that $\rho$ sends all but at most two peripheral loops on $\Sigma_1$ to $\pm\id$. This means that the restriction of $\rho$ to $\pi_1\Sigma_1$ has abelian image, and thus the Dehn twist $\sigma_{i,j}$ fixes $[\rho]$. Since the set of Dehn twists $\sigma_{i,j}$ for $1\leq i<j\leq n-1$ generates $\PMod(\Sigma)$ by Lemma~\ref{lem:GW-mcg-generators} from Appendix~\ref{apx:generators-of-pmod}, we deduce that $[\rho]$ is a fixed point of the mapping class group action.}
	\end{proofclaim}
	
	\sam{Conversely, let us now consider a fixed point $[\rho]$ of the mapping class group action, and assume that it sends at least four peripheral curves to non-scalar elements of $\SL_2\C$.}
	We'll proceed by contradiction and suppose that $\rho$ has non-abelian image. This means that we can find two of the generators $c_1,\ldots,c_n$ whose images under $\rho$ don't commute. Up to relabelling the generators, we can even assume
	that $\rho(c_1)\notin Z(\rho(c_2))$. We consider the Dehn twist $\tau$ along the curve $c_1c_2$ \cm{(this is the Dehn twist $\tau_{1,2}=\sigma_{1,2}$ in the notation of Appendix~\ref{apx:generators-of-pmod}).} Since we're assuming that $\Sigma$ has at least
	four punctures, $\tau$ is a non-trivial element of $\PMod(\Sigma)$. \cm{Using~\eqref{eq:Dehn-twist-action}, we see that the image of $[\rho]$ by $\tau$, which we denote by $\tau.[\rho]$, is the conjugacy class of the representation}
\[
c_i\mapsto \begin{cases} \rho(c_i) &\text{if } i=1,2,\\
\rho(c_1c_2)\rho(c_i)\rho(c_1c_2)^{-1} &\text{if } i\geq 3.
\end{cases}
\]
	Since $[\rho]$ is a global fixed point of the mapping class group action by assumption, we must have $[\rho]=\tau.[\rho]$. This means that there exists some $A\in \SL_2\C$, such that $A\in Z(\rho(c_1))\cap Z(\rho(c_2))$ and $A\rho(c_1c_2)\in Z(\rho(c_3))\cap\cdots\cap Z(\rho(c_n))$. Our assumption that $\rho(c_1)\notin Z(\rho(c_2))$ implies $Z(\rho(c_1))\neq Z(\rho(c_2))$. Because we're also assuming that $\rho(c_1)$ and $\rho(c_2)$ are different from $\pm\id$, we can apply Fact~\ref{fact:centralizers-SL2C} to conclude that $A=\pm \id$. The assumption $\rho(c_1)\notin Z(\rho(c_2))$ also implies that $\rho(c_1c_2)\neq \pm\id$, which means that $Z(\rho(c_3))\cap\cdots\cap Z(\rho(c_n))$ contains a \sam{non-scalar} element. Since all $\rho(c_3),\ldots,\rho(c_n)$ are different from $\pm\id$, Fact~\ref{fact:centralizers-SL2C} implies $Z(\rho(c_3))=\cdots=Z(\rho(c_n))$.

To finish the proof, we consider $\rho(c_2)$ and $\rho(c_3)$. If $\rho(c_2)\notin Z(\rho(c_3))$, then the same argument as above will give $Z(\rho(c_1))=Z(\rho(c_4))=\cdots=Z(\rho(c_n))$ and so $Z(\rho(c_1))=Z(\rho(c_3))=\cdots=Z(\rho(c_n))$. Since $\rho(c_2)=\rho(c_1)^{-1}\rho(c_n)^{-1}\cdots\rho(c_3)^{-1}$, we obtain $\rho(c_1)\in Z(\rho(c_2))$; a contradiction. On the contrary, if $\rho(c_2)\in Z(\rho(c_3))$, then $Z(\rho(c_2))=Z(\rho(c_3))$ by Fact~\ref{fact:centralizers-SL2C} and so $Z(\rho(c_2))=\cdots=Z(\rho(c_n))$. As before, this implies $\rho(c_2)\in Z(\rho(c_1))$, or equivalently $\rho(c_1)\in Z(\rho(c_2))$; a contradiction again.
\end{proof}

\subsection{Representations with finite image}\label{sec:finite-image}
\sam{Other examples} of finite mapping class group orbits are given by representations with finite image.
\begin{lem}\label{lem:finite-image-implies-finite-orbit}
The conjugacy class of a representation $\rho\colon\pi_1\Sigma\to\SL_2\C$ whose image is a finite subgroup of $\SL_2\C$ belongs to a finite mapping class group orbit. (This statement remains true for general surfaces and any target Lie group.)
\end{lem}
\begin{proof}
If $\Gamma$ denotes a finite subgroup of $\SL_2\C$, then there are at most finitely many representations $\pi_1\Sigma\to \SL_2\C$ whose image is contained in $\Gamma$ because $\pi_1\Sigma$ is finitely generated. So, if the image of $\rho$ is contained in $\Gamma$, then any point in the mapping class group orbit of $[\rho]$ is the conjugacy class of another representation whose image is also contained in $\Gamma$. This implies that the mapping class group orbit of $[\rho]$ is finite.
\end{proof}
\begin{rem}
	\sam{Since we did not use any conjugation in the proof of Lemma~\ref{lem:finite-image-implies-finite-orbit}, representations with finite image are more generally finite orbits of the $\rm{Aut}(\pi_1\Sigma)$-action on
	$\Hom(\pi_1\Sigma,\SL_2\C)$ by pre-composition.}
\end{rem}

The converse of Lemma~\ref{lem:finite-image-implies-finite-orbit} is not true. We'll encounter many examples of finite mapping class group orbits coming from $\SL_2\C$-representations of punctured spheres with infinite image later in the paper (this is, for instance, always the case for finite orbits coming from DT representations).

The finite subgroups of $\SL_2\C$ are all conjugate to a finite subgroup of $\SU(2)$ since they are compact. The list of finite subgroups of $\SU(2)$ was computed by Klein from the classification of Platonic solids~\cite{klein}. He proved that the only finite subgroups of $\SU(2)$ are cyclic groups (which are abelian), along with dihedral, tetrahedral, octahedral, and icosahedral subgroups. Many examples of representations with image contained in these ``Platonic subgroups'' of $\SU(2)$ can be obtained as monodromies of algebraic solutions to the Painlevé \Romannum{6} equation, such as the ones constructed by Boalch in~\cite{boalch-finite-image}.

\cm{To our knowledge, the exact enumeration of all finite orbits coming from representations with finite image when $\Sigma$ has six punctures or more is yet to be accomplished. For $5$-punctured spheres, a list can be found in~\cite{tykhyy}. For $4$-punctured spheres, in the icosahedral case for instance, Boalch~\cite{boalch-1} enumerated all 52 finite orbits among the naive list of 26,668 conjugacy classes of representations computed by Hall~\cite{hall}. Those orbits make up for most of the 45 exceptional orbits (see Theorem~\ref{thm:LT-classification}).}

\subsection{Representations with discrete image}\label{sec:discrete-finite-orbits}
A larger family of representations that give rise to finite orbits are representations with discrete image whose conjugacy class belongs to a compact component of the character variety. If the compact component consists of $\SU(2)$-representations, then a representation with discrete image actually has finite image. The situation is more interesting in the case of DT representations whose images can be both discrete and infinite. \cm{We start with two standard facts.}
\begin{fact}\label{fact:discrete-image-discrete-Aut-orbit}
	\sam{Let $\Gamma$ be a finitely generated group and $G$ be a Lie group. The automorphism group $\Aut(\Gamma)$ acts on $\Hom(\Gamma,G)$ by pre-composition. If $\rho\colon\Gamma\rightarrow G$ is a representation with discrete image, then the  $\Aut(\Gamma)$-orbit of $\rho$ is a discrete and closed subset of $\Hom(\Gamma,G)$.}
\end{fact}
\begin{proof}
	\sam{
	Let $(\varphi_n)\subset\Aut(\Gamma)$ be a sequence of automorphisms such that $\rho\circ\varphi_n$ converges towards some $\rho_\infty\in\Hom(\Gamma,G)$. 
	In particular, for every $\gamma\in\Gamma$, $\rho(\varphi_n(\gamma))$ converges towards $\rho_\infty(\gamma)$. Since $\rho(\Gamma)$ is a discrete subgroup of $G$, it's a closed subgroup and therefore $\rho_\infty(\gamma)\in \rho(\Gamma)$. Moreover, the sequence $\rho(\varphi_n(\gamma))$ must be stationary, so there is $N=N(\gamma)>0$ such that for $n\geq N$, $\rho(\varphi_n(\gamma))=\rho_\infty(\gamma)$.
	Applying this observation to a finite generating family for $\Gamma$, we obtain that $\rho\circ\varphi_n$ is a stationary sequence of representations and thus $\rho_\infty = \rho\circ\varphi_n$ for large enough $n$. We conclude that the $\Aut(\Gamma)$-orbit of $\rho$
	is discrete and closed in $\Hom(\Gamma,\rho)$, as desired.}
\end{proof}
When $G=\psl$, an analogous statement to Fact~\ref{fact:discrete-image-discrete-Aut-orbit} holds for the mapping class group action on the character varieties. This is due to the following fact.
\begin{fact}\label{fact:discrete-subgroup-discrete-conj-classes}
	If $\Gamma\subset\SL_2\R$ is a discrete subgroup, then the set of traces $\trace(\Gamma)=\{\trace(\gamma) : \gamma\in\Gamma\}$ is a discrete and closed subset of $\R$.
\end{fact}
\begin{proof}
	\sam{By Selberg's Lemma, $\Gamma$ has a torsion-free finite-index subgroup $\Gamma'$. By discreteness of $\Gamma$, every elliptic element $\gamma\in\Gamma$
	must therefore be of finite order, and its order is bounded by the index of $\Gamma'$ in $\Gamma$. In particular, the intersection $\trace(\Gamma)\cap (-2,2)$ is finite, and thus so is $\trace(\Gamma)\cap [-2,2]$.
	Now, $\Gamma'$ being torsion-free, the quotient $\mathbb H/\Gamma'$ is a finite-type hyperbolic surface and its length spectrum is a discrete and closed subset of $\R$ (more precisely: there are only finitely many closed geodesics whose length is at most $L$ for every $L>0$, see e.g.~\cite[Proposition~2.28]{borthwick}).
	This implies that $\trace(\Gamma')\cap \R\setminus [-2,2]\subset \R$ is discrete and closed too.
	Since $\Gamma'$ has finite index in $\Gamma$, we can show that $\trace(\Gamma)\cap \R\setminus [-2,2]$ is also a discrete and closed subset of $\R$. Indeed,
	since $\Gamma'$ has finite index, there is some $N>0$ such that for every $\gamma\in\Gamma$, $\gamma^N\in\Gamma'$.
	From the Cayley--Hamilton Theorem applied to $\rm{SL}_2\R$, there is a non-constant polynomial $P$ such that $\trace(A^N)=P(\trace(A))$ for any $A\in\rm{SL}_2\R$.
	This shows that $\trace(\Gamma)\subset P^{-1}(\trace(\Gamma'))$. As $\trace(\Gamma')$ is discrete and $P$ is proper, we deduce that $\trace(\Gamma)$ is discrete and closed,
	as desired.}
\end{proof}
We are now ready to prove our main statement regarding discrete representations.
\begin{cor}\label{cor:discrete-rep-discrete-orbit}
	\sam{Let $\alpha$ be a vector of peripheral angles satisfying $\sum_{p\in\mathcal{P}}\alpha_p\notin 2\pi\Z$ and let $[\rho]\in\Rep_\alpha(\SpherePk,\psl)$. If $\rho$ has discrete image, then the mapping class group orbit of $[\rho]$ is a discrete and closed subset of $\Rep_\alpha(\SpherePk,\psl)$.}
\end{cor}
\begin{proof}
	\sam{Assume there is a sequence $(\sigma_n)$ of automorphisms of $\pi_1\Sigma$ whose conjugacy classes give elements of $\PMod(\Sigma)$ and such that $[\rho\circ\sigma_n]$ converges to $[\rho_\infty]\in\Rep_\alpha(\SpherePk,\psl)$. The conclusion is equivalent to proving that the sequence $[\rho\circ\sigma_n]$ is stationary.}
	
	\sam{Since $\pi_1\SpherePk$ is a free group, the representation $\rho$ can be lifted to a (non-unique) representation $\overline{\rho}\colon \pi_1\Sigma\to\SL_2\R$. Its conjugacy class $[\overline{\rho}]$ lies in some relative character variety $\Rep_{\mathcal{C}}(\Sigma,\SL_2\R)$ (Definition~\ref{defn:relative-character-variety}). The projection $\SL_2\R\to\psl$ defines an equivariant homeomorphism between $\Rep_{\mathcal{C}}(\Sigma,\SL_2\R)$ and $\Rep_\alpha(\SpherePk,\psl)$. In particular, the sequence $[\overline{\rho}\circ\sigma_n]$ converges to some $[\overline{\rho_\infty}]$ in $\Rep_{\mathcal{C}}(\Sigma,\SL_2\R)$. This means that for any $\gamma\in\pi_1\SpherePk$, $\trace(\overline{\rho}\circ\sigma_n(\gamma))$ converges to $\trace(\overline{\rho_\infty}(\gamma))$. Since $\mathrm{Im}(\overline{\rho})$ is a discrete subgroup of $\SL_2\R$ by assumption, we know from Fact~\ref{fact:discrete-subgroup-discrete-conj-classes} that $\trace(\mathrm{Im}(\overline{\rho}))$ is a discrete and closed subset of $\R$. Therefore, the sequence $\trace(\overline{\rho}\circ\sigma_n(\gamma))$ must be stationary.}
	
	\sam{As we're assuming $\sum_{p\in\mathcal{P}}\alpha_p\notin 2\pi\Z$, it follows that all the representations $\overline{\rho}\circ\sigma_n$ are irreducible (Definition~\ref{defn:reductive-and-irreducible-representations}). Now, it follows from the work of Procesi~\cite{procesi} that there is a finite collection $\gamma_1,\ldots,\gamma_k\in\pi_1\Sigma$ such that the conjugacy class in $\Rep_{\mathcal{C}}(\Sigma,\SL_2\R)$ of an irreducible representation $\phi$ is uniquely determined by the values $\trace(\phi(\gamma_i))$ for $i=1,\ldots,k$. This means that the sequence $[\overline\rho\circ\sigma_n]$ is stationary, and thus so is the sequence $[\rho\circ\sigma_n]$.}
\end{proof}
\begin{rem}
	\sam{A more general version of Corollary~\ref{cor:discrete-rep-discrete-orbit} holds. If $\rho$ is a discrete and irreducible representation of a surface group into a linear group of real or complex matrices such that $\trace(\mathrm{Im}(\rho))$ is a discrete and closed subset of $\R$ or $\C$ (recall that this is automatic for representations into $\SL_2\R$ by Fact~\ref{fact:discrete-subgroup-discrete-conj-classes}), then the mapping class group of $[\rho]$ will be discrete. This is because Procesi's theorem holds for any matrix rank. Examples of such $\rho$ include any representation whose image is a lattice, or more generally Anosov representations, but well outside the scope of this paper. However, it's not sufficient that $\rho$ is discrete and irreducible for the orbit of $[\rho]$ to be discrete. Counter-examples with image in $\SL_2\C$ have been described by Souto--Storm~\cite{souto-storm}.}
\end{rem}

Corollary~\ref{cor:discrete-rep-discrete-orbit} relates the discreteness of a representation with the discreteness of the corresponding mapping class group orbit. This result can be seen as an analogue of a theorem by Previte--Xia and Golsefidy--Tamam which says that (under certain conditions) dense $\SU(2)$-representations have dense mapping class group orbits~\cite{previte-xia-1, previte-xia-2, golsefidy-tamam}.

\begin{cor}\label{cor:DT-discrete-implies-finite-orbit}
The conjugacy class of a discrete DT representation always belongs to a finite mapping class group orbit inside $\RepDT{\alpha}$.
\end{cor}
\begin{proof}
We know from Corollary~\ref{cor:discrete-rep-discrete-orbit} that the orbit of a discrete DT representation is a discrete and closed subset of $\RepDT{\alpha}$. Since $\RepDT{\alpha}$ is compact, the orbit is finite.
\end{proof}

It's possible to obtain examples of discrete DT representations (with infinite image) from certain representations with finite image inside $\SU(2)$.
The procedure is known as an Okamoto transformation (Section~\ref{sec:finite-orbit-n=4}). An example of this process applied to representations with image contained in a finite
octahedral subgroup of $\SU(2)$ is described by Boalch in~\cite[Section 5]{boalch-finite-image}. \cm{Another way to construct discrete DT representations
is from a Coxeter tessellation of the hyperbolic plane which consists of a tiling generated by reflections in the sides of a geodesic triangle with interior angles $(\pi/p, \pi/q, \pi/r)$ for some positive integers $p,q,r$ satisfying $1/p + 1/q + 1/r<1$. We'll call it a \emph{$(p,q,r)$-tessellation}.}
\begin{ex}\label{ex:discrete-DT}
\cm{Given such a tessellation, assume that we can construct a triangle chain $T$ that satisfies the angle and orientation properties of Lemma~\ref{lem:properties-triangle-chains} and which ``fits well'' within the tessellation, meaning that it satisfies the following properties:
\begin{itemize}
\item The exterior vertices and the shared vertices of $T$ are vertices of the tessellation;
\item The interior angles at each vertex of $T$ are integer multiples of the tessellation angle at the vertex.
\end{itemize}
Note that we don't require the sides of the triangles in $T$ to align with the tessellation edges. When such a $T$ exists, then the image of the associated DT representation is a subgroup of the rotation triangle group associated to the tessellation (Section~\ref{sec:triangle-groups}), hence a discrete subgroup of $\psl$. In particular, the mapping class group orbit of the conjugacy class of this representation is finite by Corollary~\ref{cor:DT-discrete-implies-finite-orbit}. Below are two examples of chains of two triangles fitting well on a $(2,3,7)$-tessellation.
\begin{center}
\begin{tikzpicture}[font=\sffamily,decoration={markings, mark=at position 1 with {\arrow{>}}}]
\node[anchor=east] at (-0.5,0) {\includegraphics[width=6cm]{fig/example-B-triangle-chain-1}};
\node[anchor=west] at (0.5,0) {\includegraphics[width=6cm]{fig/example-Bprime-triangle-chain-1-no-shadow}};
\end{tikzpicture}
\end{center}
Those two examples are borrowed from the finite orbit described in Example~\ref{ex:example-Dehn-twist-action}.\ensuremath{\lozenge}}
\end{ex}

\subsection{Pullback representations}\label{sec:pullback-orbits}
The next source of examples of finite mapping class groups orbit is obtained by pulling back representations of a sphere with three punctures via a family of ramified coverings.
\sam{Those pullback representations were first studied by Doran, who called them \emph{geometric isomonodromic deformations}~\cite{doran}. Their classification was completed
by Diarra~\cite{diarra-pull-back}.}
It turns out to be quite crucial to consider a family of coverings and not just a single one. Let us explain the procedure in more details, starting with an example that we learned
from Loray (many more examples can be found in Landesman's notes~\cite{landesman-notes}). 

\begin{ex}\label{ex:pullback}
Consider the degree-2 family of ramified coverings $f_c\colon \CP^1\to \CP^1$ given by $f_c(z)=(c+1)^2z/(c+z)^2$, where $c$ belongs to the parameter space $C=\CP^1\setminus\{0,\pm 1,\infty\}$. The covering $f_c$ ramifies over $\infty$ and $c'=(c+1)^2/4c$. The preimages of $0$ are $\{0,\infty\}$ and the preimages of $1$ are $\{1,\lambda(c)\}$. A rapid computation shows that $\lambda(c)=c^2$. 

\begin{center}
\begin{tikzpicture}

\draw[-] (0,0) to (6,0);

\draw[-] (2,1) to (4,1);
\draw[-] (2,2) to (4,2);
\draw[-] (2,1) to[out=180, in=0] (0,2);
\draw[-] (2,2) to[out=180, in=0] (0,1);
\draw[-] (4,1) to[out=0, in=180] (6,2);
\draw[-] (4,2) to[out=0, in=180] (6,1);

\draw[dashed] (1,0) to (1, 1.5);
\draw[dashed] (2,0) to (2, 2);
\draw[dashed] (4,0) to (4, 2);
\draw[dashed] (5,0) to (5, 1.5);

\fill(1,0) circle (0.07) node[below]{$c'$};
\fill(1,1.5) circle (0.07) node[left]{$c$};

\fill(2,0) circle (0.07) node[below]{$0$};
\fill(2,1) circle (0.07) node[above right]{$0$};
\fill(2,2) circle (0.07) node[above]{$\infty$};

\fill(4,0) circle (0.07) node[below]{$1$};
\fill(4,1) circle (0.07) node[above right]{$1$};
\fill(4,2) circle (0.07) node[above]{$\lambda(c)$};

\fill(5,0) circle (0.07) node[below]{$\infty$};
\fill(5,1.5) circle (0.07) node[right]{$-c$};

\draw[->] (6.5,1.5) to[out=315, in=45] (6.5,0.5);
\node at (7,1) {$f_c$};
\end{tikzpicture}
\end{center}
Each covering map $f_c$ induces an injective morphism
\[
(f_c)_\ast\colon\pi_1\big(\CP^1\setminus \{0,1,\infty,-c,\lambda(c)\}, c\big)\to \pi_1(\CP^1\setminus\{0,1,\infty\}, c')
\]
whose image is an index-$2$ subgroup of the target. Now, let $\rho\colon \pi_1\big(\CP^1\setminus\{0,1,\infty\},c'\big)\to \SL_2\C$ be such that peripheral loops around $\infty$ are mapped to order 2 elements of $\SL_2\C$ (equivalently, to traceless elements). We don't prescribe any particular order for peripheral loops around $0$ and $1$; let's simply say that there are mapped to elements of trace $t_0\in\C$, respectively $t_1\in\C$. The pullback representation $f_c^\ast\rho=\rho\circ (f_c)_\ast$ is a priori a representation of $\pi_1\big(\CP^1\setminus \{0,1,\infty,-c,\lambda(c)\}, c\big)$ into $\SL_2\C$. However, since $f_c$ ramifies over $\infty$ and the monodromy of $\rho$ around $\infty$ has order 2 by assumption, the monodromy of $f_c^\ast\rho$ around $-c$ is trivial. In other words, the pullback representation $f_c^\ast\rho$ is really a representation of a 4-punctured sphere:
\[
f^\ast_c\rho\colon \pi_1\big(\CP^1\setminus \{0,1,\infty,\lambda(c)\}\big)\to\SL_2\C.
\]
By construction, $f_c^\ast\rho$ maps peripheral loops around $0$ and $\infty$ to elements of trace $t_0$ and peripheral loops around $1$ and $\lambda(c)$ to elements of trace $t_1$. This means that the conjugacy class of $f_c^\ast\rho$ lives inside a relative character variety where two pairs of peripheral conjugacy classes were picked to be equal. We'll explain in Lemma~\ref{lem:pullback-orbits-are-finite} why the mapping class group orbit of $[f_c^\ast\rho]$ is finite. Actually, we'll encounter this exact pullback representation again later in Section~\ref{sec:finite-orbit-n=4} where it'll correspond to a finite mapping class group orbit ``of Type~\Romannum{2}''.
\ensuremath{\lozenge}
\end{ex}

We can try to generalize Example~\ref{ex:pullback} as follows. Say we have a family of ramified covering $f_c\colon\CP^1 \to \CP^1$ where $c$ belongs to some curve $C$ that plays the role of the parameter space of our family. We assume that the coverings have finite degree $d$. The covering $f_c$ induces an injective morphism $(f_c)_\ast\colon\pi_1\big(\CP^1\setminus f_c^{-1}(\{0,1,\infty\})\big)\to \pi_1(\CP^1\setminus\{0,1,\infty\})$ whose image is an index-$d$ subgroup of the target. The pullback of a representation $\rho\colon \pi_1\big(\CP^1\setminus \{0,1,\infty\}\big)\to\SL_2\C$ is the representation
\[
f_c^\ast\rho=\rho\circ (f_c)_\ast\colon \pi_1\big(\CP^1\setminus f_c^{-1}(\{0,1,\infty\})\big)\to \SL_2\C.
\]
We denote by $n_p\in \{2,3,\ldots\}\cup\{\infty\}$ the orders of the peripheral monodromies of $\rho$ around $p\in\{0,1,\infty\}$. To simplify the situation, we assume that for every $c\in C$, there exist $n-3$ distinct points $\lambda_1,\ldots,\lambda_{n-3}$ in $\CP^1\setminus\{0,1,\infty\}$ with the following property: for every $p\in\{0,1,\infty\}$, any point $x\in f_c^{-1}(p)\setminus \{0,1,\infty,\lambda_1,\ldots,\lambda_{n-3}\}$ is ramified to an order divisible by $n_p$. This implies that the pullback representation $f_c^\ast\rho$ has trivial holonomy around each such $x$ and can therefore be seen as a representation of the $n$-punctured sphere $\CP^1\setminus \{0,1,\infty,\lambda_1,\ldots,\lambda_{n-3}\}$. We'll write
\[
f_c^\ast\rho\colon \pi_1\big(\CP^1\setminus \{0,1,\infty,\lambda_1,\ldots,\lambda_{n-3}\}\big)\to \SL_2\C.
\]
\begin{defn}\label{def:pullback-orbits}
Every representation of an $n$-punctured sphere into $\SL_2\C$ that can be constructed in this way is called a \emph{pullback representation}. The mapping class group orbit of its conjugacy class is accordingly called a \emph{pullback orbit}.
\end{defn}

\cm{It follows from Doran's and Diarra's work that pullback orbits are always finite~\cite{doran, diarra-pull-back}. A more modern version of that statement has also been established by Landesman--Litt in~\cite[Proposition~2.1.3]{LL}. We include a rapid summary of the argument, following Landesman's notes~\cite[Section 7]{landesman-notes}, for the sake of completeness.}

\begin{lem}\label{lem:pullback-orbits-are-finite}
Pullback orbits are finite.
\end{lem}
\begin{proof}[Sketch of proof]
\sam{The family of ramified coverings $(f_c)$ can be seen as a fiber bundle map $M\rightarrow M$ where $M$ is a $\CP^1$-bundle over a complex curve $C$  which is our parameter
	space. The points $0,1,\infty,\lambda_1,\ldots,\lambda_{n-3}$ are then reinterpreted as sections $C\to M$ which are pointwise distinct.
	Let $E_\lambda=M\setminus\{0(C),1(C),\infty(C),\lambda_1(C)\ldots,\lambda_{n-3}(C)\}$. A representation $\rho\colon\pi_1(\CP^1\setminus\{0,1,\infty\})\to\SL_2\C$
	now pulls back to a representation $\widetilde\rho\colon\pi_1 E_\lambda\to\SL_2\C$, which fiberwise restricts to a representation
	$\pi_1(\CP^1\setminus\{0,1,\infty,\lambda_1(c),\ldots,\lambda_{n-3}(c)\})\to\SL_2\C$. As we're about to see, all these representations give rise to finite mapping
	class group orbits.}

\sam{Since the section $\lambda=(\lambda_1,\ldots,\lambda_{n-3})$ avoids $\{0,1,\infty\}$, it induces a map $C\to\mathcal{M}_{0,n}$, which is rational
	and nonconstant (otherwise the maps $(f_c)$ were just the same map all along). Hence, it's dominant, i.e.~it's a finite covering over a Zariski dense open subset of
	$\mathcal{M}_{0,n}$ (as in Example~\ref{ex:pullback}). This assumption ensures that the induced map $\pi_1 C\to\pi_1\mathcal{M}_{0,n}$ identifies $\pi_1 C$ with
	a finite index subgroup of $\pi_1\mathcal{M}_{0,n}$.}

\sam{Now, $E_\lambda$ being a fiber bundle over $C$ and $\mathcal{M}_{0,n+1}$ being a fiber bundle over $\mathcal{M}_{0,n}$, both with fibers the $n$-punctured sphere,
	the induced short exact sequences at the level of fundamental groups are related by the morphisms mentioned above, and we obtain the following diagram in which each
	row is a short exact sequence.}

\begin{center}
\begin{tikzcd}
1\arrow[r]&
	\pi_1\big(\CP^1\setminus \{0,1,\infty,\lambda_1(c),\ldots,\lambda_{n-3}(c)\}\big)\arrow[r]\arrow[d,leftrightarrow]&
\pi_1E_\lambda \arrow[r]\arrow[d]&
\pi_1C \arrow[r]\arrow[d]&1\\
	1\arrow[r]&\pi_1\big(\CP^1\setminus \{0,1,\infty,\lambda_1(c),\ldots,\lambda_{n-3}(c)\}\big)\arrow[r]&\pi_1\mathcal{M}_{0,n+1}\arrow[r]&
\pi_1\mathcal{M}_{0,n}\arrow[r]&1
\end{tikzcd}
\end{center}
	From both short exact sequences, we obtain morphisms $\pi_1C\to \Out\big(\pi_1\big(\CP^1\setminus \{0,1,\infty,\lambda(c)\}\big)\big)$ and
	$\pi_1 \mathcal{M}_{0,n}\to \Out\big(\pi_1\big(\CP^1\setminus \{0,1,\infty,\lambda(c)\}\big)\big)$.
	They are obtained by lifting an element to the middle group and letting it act by conjugation
	on~$\pi_1\big(\CP^1\setminus \{0,1,\infty,\lambda(c)\}\big)$. Observe that the second row is really just the Birman exact sequence for mapping
	class groups (see e.g.~\cite[Theorem~4.6]{mcg-primer}).
	This means that the morphism $\pi_1 \mathcal{M}_{0,n}\to \Out\big(\pi_1\big(\CP^1\setminus \{0,1,\infty,\lambda(c)\}\big)\big)$ coincides with
	the morphism $\PMod(\CP^1\setminus \{0,1,\infty,\lambda(c)\})\to \Out\big(\pi_1\big(\CP^1\setminus \{0,1,\infty,\lambda(c)\}\big)\big)$
	from the Dehn--Nielsen--Baer Theorem which is used to define the mapping class group action on the character variety (Section~\ref{sec:character-varieties}).

	We can now conclude the proof that the mapping class group orbits coming from the pullback representations $\pi_1\big(\CP^1\setminus \{0,1,\infty,\lambda(c)\}\big)\to\SL_2\C$ are
	finite. All we need to do is to show that they are stabilized by $\pi_1 C$, because $\pi_1 C$ can be identified with a finite index subgroup of $\PMod(\CP^1\setminus \{0,1,\infty,\lambda(c)\})$. Since pullback representations are all restrictions of $\widetilde\rho\colon \pi_1 E_\lambda\to\SL_2\C$,
	an element $\gamma\in\pi_1 C$ acts on a representation $\pi_1\big(\CP^1\setminus \{0,1,\infty,\lambda(c)\}\big)\to\SL_2\C$ by conjugating it by an element of $\SL_2\C$ of
	the form $\widetilde\rho(\widetilde \gamma)$, where $\widetilde\gamma\in \pi_1E_\lambda$ is any lift of $\gamma$. So, $\gamma$ acts trivially, proving our claim.
\end{proof}

Pullback orbits have been classified by Diarra in~\cite{diarra-pull-back}, generalizing Doran's work in the 4-punctured case~\cite{doran}. Diarra shows that the argument boils down
to classifying finite degree ramified coverings of a 3-punctured sphere (see also~\cite[Section~2.3.4]{litt} for an explanation of algebraic geometric flavor).
Diarra's motivation was to find non-elementary algebraic solutions to Garnier systems of isomonodromy differential equations which generalize the Painlevé~\Romannum{6} equation
by allowing more poles. He gives a complete list of solutions that can be obtained by pulling back a Fuchsian equation, following a method developed by Doran~\cite{doran} and
Kitaev~\cite{kitaev} for Painlevé~\Romannum{6}, and exploited by Boalch~\cite{boalch-2} to find the ''elliptic 237 solutions'' (see the discussion after Theorem~\ref{thm:LT-classification}). The monodromy of every algebraic solution built by Diarra using the Kitaev method gives rise to a pullback orbits in the sense of Definition~\ref{def:pullback-orbits}. It follows from Diarra's classification that there are no finite orbits of pullback type for spheres with seven punctures or more~\cite[Théorème~5.1]{diarra-pull-back}, confirming Tykhyy's Conjecture (Conjecture~\ref{conj:tykhyy}) in this case.

\subsection{Totally non-hyperbolic representations}\label{sec:totally-non-hyperbolic}
\cm{A surface group representation $\rho\colon\pi_1\Sigma\to\SL_2\C$ is \emph{totally non-hyperbolic} if every simple closed curve on $\Sigma$ is mapped to either an elliptic element, a parabolic element, or to $\pm \id$. In other words, the trace of the image of every simple closed curve is real and contained in $[-2,2]$. When all the elliptic images of simple closed curves have finite order, those representations give finite mapping class group orbits.}

\begin{lem}\label{lem:rational-angles-imply-finite-orbit}
	If $\rho\colon\pi_1\Sigma\to\SL_2\C$ is an irreducible representation (Definition~\ref{defn:reductive-and-irreducible-representations}) of a punctured sphere $\Sigma$ that maps every simple closed curve on $\Sigma$ to \sam{an element whose trace belongs to $2\cos(\pi\Q)$}, then the mapping class group \cm{orbit} of $[\rho]$ is finite.
\end{lem}
\begin{proof}
The conjugacy class of an irreducible representation $\rho\colon\pi_1\Sigma\to\SL_2\C$ is completely determined by its \sam{trace function} $\chi\colon \pi_1\Sigma\to \C$ defined
by $\chi(\gamma)=\trace(\rho(\gamma))$. It turns out that the values of $\chi$ on an explicit finite set of elements in $\pi_1\Sigma$ also determines $[\rho]$.
This is a consequence of the work of Procesi~\cite{procesi}. Since we're dealing with matrices of rank 2, it's possible to choose this finite set of elements in $\pi_1\Sigma$ to
consist of \emph{simple} closed curves only, as explained for instance by Goldman-Xia in~\cm{\cite[Proposition~2.2]{goldman-xia}}. We fix an arbitrary such finite
set~$\mathcal C\subset \pi_1\Sigma$ of fundamental group elements representing simple closed curves. 

Let now $\rho\colon\pi_1\Sigma\to\SL_2\C$ be an irreducible representation that maps every simple closed curve to a parabolic element or to an elliptic element of finite order. The restriction of the \cm{trace function} $\chi$ of $\rho$ to simple closed curves therefore takes values in $[-2,2]\cap 2\cos(\pi\Q)$. By Selberg's Lemma, there exists an index-$N$ subgroup of $\rho(\pi_1\Sigma)$ that's torsion-free. This means that all finite order elements in the image of $\rho$ have order at most $N$. Their traces therefore belong to a finite subset $\mathcal S$ of $[-2,2]\cap 2\cos(\pi\Q)$. Since there are only finitely many functions $\mathcal C\to \mathcal S\cup \{\pm 2\}$ and the restriction of a \cm{trace function} to $\mathcal C$ determines the conjugacy class of the associated irreducible representation completely, the mapping class group orbit of $[\rho]$ is finite.
\end{proof}

\begin{rem}
\cm{Lemma~\ref{lem:rational-angles-imply-finite-orbit} remains true if we replace $\Sigma$ by any other surface (it's not particular to punctured spheres) and the same proof applies. We will, however, only make use of it in the context of punctured spheres.}
\end{rem}

\subsection{Arbitrarily long finite orbits}\label{sec:arbitrarily-long-orbits}
\sam{When the peripheral monodromies contain at least one parabolic conjugacy class, finite mapping class group orbits can be of arbitrarily long length.}
When $\Sigma$ is a 4-punctured sphere, arbitrarily long orbits can be found by studying \emph{Cayley orbits} (this name comes from Lisovyy--Tykhyy's nomenclature which we'll introduce in Section~\ref{sec:LT-classification}). Those orbits were first studied by Picard~\cite{picard} and Fuchs~\cite{fuchs}. Here's an example that's also Zariski dense. Consider the following representation $\rho_\theta\colon\pi_1\Sigma\to\SL_2\R$ of the 4-punctured sphere $\Sigma$ with parabolic peripheral monodromy. The representation is defined using a geometric presentation of $\pi_1\Sigma$ with generators $c_1,c_2,c_3,c_4$:
\begin{equation*}
\begin{array}{ll}
\rho_\theta(c_1)=\begin{pmatrix}
1 & 2\cos(\theta/2)-2\\
0 & 1
\end{pmatrix},&\quad   \rho_\theta(c_2)=\begin{pmatrix}
2 & -1 \\
1 & 0
\end{pmatrix},\\
\rho_\theta(c_3)=\begin{pmatrix}
0 & 1\\
-1 & -2
\end{pmatrix}, &\quad  \rho_\theta(c_4)=\begin{pmatrix}
1 & -2\cos(\theta/2)-2 \\
0 & 1
\end{pmatrix}.
\end{array}
\end{equation*}
A rapid computation shows that $\rho_\theta$ is indeed a representation of $\pi_1\Sigma$. Observe that $\rho_\theta(c_1)$ and $\rho_\theta(c_4)$ both fix $\infty$ as a boundary
point of the upper-half plane, whereas $\rho_\theta(c_2)$ fixes $1$ and $\rho_\theta(c_3)$ fixes $-1$. \cm{When $\theta=0$ or $\theta=2\pi$, then $\rho_\theta$ maps either $c_1$ or $c_4$ to $\id$ implying that $[\rho_\theta]$ is a fixed point of
the mapping class group action by Lemma~\ref{lem:fixed-point-implies-abelian} (more precisely, it's a finite orbit ``of Type~\Romannum{1}'' in Lisovyy--Tykhyy's nomenclature).} 
When $\theta\in (0,2\pi)$, the simple closed curve represented by the fundamental group element $c_1c_2$ has elliptic image:
\[
\rho_\theta(c_1c_2)=\begin{pmatrix}
2\cos(\theta/2) & -1\\
1 & 0
\end{pmatrix}.
\]
It fixes the point $e^{i\theta/2}$ in the upper-half plane. In other words, we can represent $[\rho_\theta]$ by a triangle chain consisting of two partially ideal triangles. The first one has vertices $(\infty, 1, e^{i\theta/2})$ and the second one has vertices $(e^{i\theta/2}, -1, \infty)$.
\begin{center}
\begin{tikzpicture}[scale =.8]
\draw[-] (0,0) to (6,0);

\draw (1,0) arc(180:0:2);
\draw[-] (1,0) to (1,6);
\draw[-] (5,0) to (5,6);
\draw[-] (4, 1.732) to (4, 6);

\draw[apricot] (3.7,1.9) arc(170:90:.3);
\node[apricot] at (3.5,2.3) {$\theta/2$};

\fill(1,0) circle (0.07) node[below]{$-1$};
\fill(5,0) circle (0.07) node[below]{$1$};
\fill(4, 1.732) circle (0.07) node[below left]{$e^{i\theta/2}$};
\fill(3,6) circle (0.07) node[below]{$\infty$};
\end{tikzpicture}
\end{center}

\begin{lem}
Assuming that $\theta\in (0,2\pi)$, the mapping class group orbit of $[\rho_\theta]$ is finite if and only if $\theta$ is a rational multiple of $\pi$. Moreover, when this is the case, the length of the orbit is at least equal to the denominator of \cm{$\theta/2\pi$} in irreducible form. 
\end{lem}
\begin{proof}
	\sam{Historically, Cayley orbits have been described in terms of  algebraic isomonodromic leaves, where they admit an explicit description in terms of hyperelliptic functions, as described by Picard~\cite{picard}.}
	
	\sam{Here, we opt for argument of arithmetic flavour. When $\theta\in (0,2\pi)$, the conjugacy classes of the representations $\rho_\theta$ belong to the same bounded component of the corresponding relative $\SL_2\R$-character variety.\footnote{This component is topologically a sphere minus four points which coincides with the bounded component of the regular locus of the surface $X^2+Y^2+Z^2+XYZ=4$ in $\R^3$. We'll make this correspondence more explicit in Section~\ref{sec:fricke-relation}, see also the statement about Cayley orbits in Theorem~\ref{thm:LT-classification}.}
	In particular, it means that $\rho_\theta$ is totally non-hyperbolic---it sends every simple closed curve on $\Sigma$ to an element whose trace is in $[-2,2]$.
	If $\theta$ is not a rational multiple of $\pi$, then the Dehn twist along the simple closed curve $c_1c_2$ would produce infinitely many orbit points when iterated on $[\rho_\theta]$.
	If $\theta$ is a rational multiple of $\pi$, then the representation $\rho_\theta$ is valued in $\SL_2(\Z[2\cos(\theta/2)])$. A quick look at the images of the generators
	ensures that the Galois conjugates of $\rho_\theta$ are of the form $\rho_{\theta'}$ for other values of $\theta'$. In particular they are also
	in the same bounded component of the same relative $\SL_2\R$-character variety. The graph of all the Galois conjugates of $\rho_\theta$ corresponds to a representation $\tilde\rho\colon\pi_1\Sigma\to(\SL_2\R)^n$ for some $n$, whose image is discrete by the Borel--Harish-Chandra Theorem. Moreover, $\tilde\rho$ sends every simple closed curve on $\Sigma$ to either an $n$-tuple of parabolic elements,
	or an $n$-tuple of non-trivial elliptic elements. In the latter case, discreteness implies that the elliptic elements are of finite order.
	This proves that $\rho_\theta$ maps every simple closed curve on $\Sigma$ to an element whose trace is in $2\cos(\pi\Q)$.
	We may now apply Lemma~\ref{lem:rational-angles-imply-finite-orbit}, and obtain that the mapping class group orbit of $[\rho_\theta]$ is finite, as desired. The length of the orbit is at least equal to the denominator of $\theta/2\pi$ in irreducible form because this corresponds to the order of the Dehn twists along the curve $c_1c_2$ when iterated on $[\rho_\theta]$.}
\end{proof}

\subsection{Representations with non-Zariski dense image}\label{sec:non-zariski-dense-orbits}
Not all the representations $\pi_1\Sigma\to\SL_2\C$ have Zariski dense image. Among those with non-Zariski dense image, some may give rise to finite mapping class group orbits. We already encountered the examples of representations with finite image in Section~\ref{sec:finite-image}. There are two other kinds of algebraic subgroups of $\SL_2\C$. According to Sit~\cite{sit}, if the Zariski closure of the image of a representation $\pi_1\Sigma\to\SL_2\C$ is a proper subgroup of $\SL_2\C$, then it is contained, maybe after conjugation, in one of the following:
\begin{enumerate}
\item The subgroup of upper triangular matrices.
\item The subgroup of diagonal and anti-diagonal matrices (sometimes also called the \emph{infinite dihedral group}).
\item A finite subgroup \cm{(see Section~\ref{sec:finite-image})}.
\end{enumerate}

\subsubsection{Upper triangular subgroup}\label{sec:upper-triangular}
Let $\rho\colon\pi_1\Sigma\to\SL_2\C$ be representation valued in the upper triangular subgroup of $\SL_2\C$. We may decompose $\rho$ into a \emph{linear part} $\lambda\colon\pi_1\Sigma\to\C^\times$ and a function $z\colon\pi_1\Sigma\to \C$ by writing the image of $\gamma\in\pi_1\Sigma$ as
\[
\begin{pmatrix}
\lambda(\gamma) & \lambda(\gamma)^{-1}z(\gamma)\\
0 & \lambda(\gamma)^{-1}
\end{pmatrix}.
\]
The relation $\rho(\gamma_1\gamma_2)=\rho(\gamma_1)\rho(\gamma_2)$ implies that the linear part $\lambda\colon\pi_1\Sigma\to\C^\times$ is a group homomoprhism and that the function $z\colon\pi_1\Sigma\to \C$ is a cocycle in the sense that $z(\gamma_1\gamma_2)=z(\gamma_1)+\lambda(\gamma_1)^2z(\gamma_2)$. As such, $z$ defines a cohomology class in $H^1(\pi_1\Sigma, \C_{\lambda^2})$---the first cohomology of the group $\pi_1\Sigma$ with coefficients in the $\pi_1\Sigma$-module $\C_{\lambda^2}$---where the module structure is given by $\gamma.t=\lambda(\gamma)^2t$ for $\gamma\in\pi_1\Sigma$ and $t\in\C$. It's possible to conjugate $\rho$ by an upper triangular matrix so that $\rho$ preserves \cm{its} upper triangular form and its linear part $\lambda$ remains the same. The cocycle $z$ however changes. For instance, if we conjugate $\rho$ by the matrix
\[
\begin{pmatrix}
a & a^{-1}b\\
0 & a^{-1}
\end{pmatrix},
\]
then $z$ changes to \cm{$a^2z+b-b\lambda^2$}. This shows that $z$ is only well-defined up to adding a coboundary and multiplying it by a non-zero complex number. In other words, conjugacy classes of representations with image in the upper triangular subgroup of $\SL_2\C$ form a bundle over $\Hom(\pi_1\Sigma,\C^\times)$ where the fibre over $\lambda$ is $\mathbb{P}H^1(\pi_1\Sigma, \C_{\lambda^2})$---the projectivization of $H^1(\pi_1\Sigma, \C_{\lambda^2})$ as a complex vector space.

Finite mapping class group orbits of conjugacy classes of representations with image in the upper triangular subgroup of $\SL_2\C$ have been classified and explicitely listed by Cousin--Moussard~\cite{cousin-moussard}. Among several other results, they proved the following statement.
\begin{thm}[{\cite[Theorem~2.3.4]{cousin-moussard}}]\label{thm:cousin-moussard}
\cm{Let $\rho\colon\pi_1\Sigma\to\SL_2\C$ be a representation of an $n$-punctured sphere with image in the upper triangular subgroup of $\SL_2\C$. Assume that the linear part of $\rho$ is never equal to $1$ on a family of geometric generators of $\pi_1\Sigma$. If the mapping class group orbit of $[\rho]$ is finite, then $n\leq 6$.}
\end{thm}
\cm{Theorem~\ref{thm:cousin-moussard} is a confirmation of Tykhyy's Conjecture (Conjecture~\ref{conj:tykhyy}) in the special case of representations with image in the upper triangular subgroup of $\SL_2\C$.}

\begin{rem}
	\sam{When studying mapping class group orbits of representations in the upper triangular subgroup of $\SL_2\C$, it is crucial to work in the topological character variety $\Rep(\Sigma,\SL_2\C)$ (Definition~\ref{defn:character-varieties}) rather than with the algebraic GIT quotient. The reason is that, in the GIT quotient, these non-reductive
	representations are identified with their semi-simplifications, which are abelian, hence lead to fixed points of the mapping class group action. In other words, if $\rho$ is upper triangular, then the topological closures of every point in the orbit of $[\rho]$ inside~$\Rep(\Sigma,\SL_2\C)$ intersect in a common point: the conjugacy class of the abelian semi-simplification.}
\end{rem}
\subsubsection{Infinite dihedral subgroup}\label{sec:infinite-dihedral}
\cm{Finite orbits coming from representations in the infinite dihedral subgroup of $\SL_2\C$ have been completely classified by Tykhyy~\cite{tykhyy}. We formulate the conclusions of the classification as Lemmas~\ref{lem:infinite-dihedral-two-anti-diagonal} and~\ref{lem:infinite-dihedral-at-least-four-anti-diagonal}, and give alternative proofs for the sake of completeness.}

Let $\rho\colon\pi_1\Sigma\to\SL_2\C$ be a representation of an $n$-punctured sphere $\Sigma$ whose image consists of diagonal and anti-diagonal matrices. We denote by $\iota(\rho)$ the number of peripheral loops of $\Sigma$ mapped to anti-diagonal matrices by $\rho$. The number $\iota(\rho)$ is well-defined because diagonal matrices remain diagonal after conjugation by an anti-diagonal matrix and vice versa. If $\iota(\rho)=0$, then $\rho$ is an abelian representation with image in the subgroup of diagonal matrices. We're mostly interested in the case where $\iota(\rho)\geq 1$. Actually, in order for $\rho$ to be a representation, $\iota(\rho)$ must be even. 

\begin{lem}\label{lem:infinite-dihedral-two-anti-diagonal}
When $\iota(\rho)=2$, the mapping class group orbit of $[\rho]$ is finite and consists of at most $2^{n-3}$ points.  If we further assume that no peripheral loop of $\Sigma$ is mapped to $\pm\id$ by $\rho$, then the mapping class group orbit of $[\rho]$ has length $2^{n-3}$.
\end{lem}
\begin{proof}
We pick a geometric presentation of $\pi_1\Sigma$ with generators $c_1,\ldots,c_n$ which we choose such that $\rho(c_1)$ and $\rho(c_n)$ are anti-diagonal, while $\rho(c_2),\ldots,\rho(c_{n-1})$ are diagonal. Recall that every point in the mapping class group \cm{orbit} of $[\rho]$ is the conjugacy class of a representation $\rho'$ with the same image as $\rho$. This representation $\rho'$ will then also verify that $\rho'(c_1)$ and $\rho'(c_n)$ are anti-diagonal, and $\rho'(c_2),\ldots,\rho'(c_{n-1})$ are diagonal.

It's convenient to work with the generating family of $\PMod(\Sigma)$ described in Lemma~\ref{lem:GW-mcg-generators} from Appendix~\ref{apx:generators-of-pmod}. It consists of all the Dehn twists $\tau_{i,j}$ along the curves represented by the fundamental group elements $c_i\cdots c_j$ where $1\leq i<j\leq n-1$ and $(i,j)\neq (1,n-1)$. It's useful at this point to remember the explicit description of \cm{those} Dehn twits action \cm{provided in~\eqref{eq:Dehn-twist-action-general}.} For instance, when $2\leq i<j\leq n-1$, the matrix $\rho(c_i\cdots c_j)$ is diagonal and $\tau_{i,j}.[\rho]$ is the conjugacy class of the representation
\[
c_k\mapsto\begin{cases}
\rho(c_i\cdots c_j)\rho(c_k)\rho(c_i\cdots c_j)^{-1} & \text{if } k=1,n,\\
\rho(c_k) & \text{if } \textnormal{else}.
\end{cases}
\]
This representation is itself conjugate to $\rho$, proving that $\tau_{i,j}.[\rho]=[\rho]$. The same argument also shows that $\tau_{i,j}.[\rho']=[\rho']$ for every $[\rho']$ in the mapping class group orbit of $[\rho]$. Now, if $1=i<j\leq n-2$, then $\rho(c_1\cdots c_j)$ is anti-diagonal. This implies that $\rho(c_1\cdots c_j)^2=-\id$ and thus $\tau_{1,j}^2.[\rho]=[\rho]$. Note that the simple closed curves corresponding to the fundamental group elements $c_1\cdots c_j$ with $2\leq j\leq n-2$ are disjoint, which implies that the Dehn twists $\tau_{1,2},\ldots,\tau_{1,n-2}$ commute. We just proved that any point $[\rho']$ in the mapping class group orbit of $[\rho]$ is the image of $[\rho]$ by an element of $\PMod(\Sigma)$ which is a word in $\tau_{1,2},\ldots,\tau_{1,n-2}$ where each generator appears at most once. This shows that the mapping class group \cm{orbit} of $[\rho]$ consists of at most $2^{n-3}$ points.

It remains to prove that the orbit of $[\rho]$ consists of exactly $2^{n-3}$ points when no peripheral curve is mapped to $\pm\id$. By the argument above, it suffices to prove that no non-trivial word in the Dehn twists $\tau_{1,2},\ldots,\tau_{1,n-2}$ fixes $[\rho]$. Note that the intersection of centralizers $Z(\rho(c_1))\cap Z(\rho(c_2))$ is equal to $\{\pm\id\}$---the center of $\SL_2\C$---because $\rho(c_2)$ is a non-trivial diagonal matrix. So, if we would have $\tau_{1,j_1}\cdots\tau_{1,j_k}.[\rho]=[\rho]$ for some $2\leq j_1<\cdots <j_k\leq n-2$, then it would imply $\rho(c_1\cdots c_{j_1})\in Z(\rho(c_{j_1+1}))$. This is however impossible because $\rho(c_{j_1+1})$ is a non-trivial diagonal matrix by assumption and $\rho(c_1\cdots c_{j_1})$ is anti-diagonal.
\end{proof}

\begin{lem}\label{lem:infinite-dihedral-at-least-four-anti-diagonal}
When $\iota(\rho)\geq 4$, if the mapping class group orbit of $[\rho]$ is finite, then the image of $\rho$ is finite.
\end{lem}
\begin{proof}
As in the proof of Lemma~\ref{lem:infinite-dihedral-two-anti-diagonal}, we pick a geometric presentation of $\pi_1\Sigma$ with generators $c_1,\ldots,c_n$ such that $\rho(c_1),\ldots, \rho(c_{2k})$ are anti-diagonal and $\rho(c_{2k+1}),\ldots,\rho(c_n)$ are diagonal (with $2k=\iota(\rho)$).

We first prove that for every pair of distinct generators $c_i,c_j\in\{c_1,\ldots,c_{2k}\}$, the diagonal matrix $\rho(c_ic_j)$ has finite order. It's the case when $\rho(c_ic_j)=\pm \id$. If $\rho(c_ic_j)\neq \pm\id$, then the intersection of centralizers $Z(\rho(c_i))\cap Z(\rho(c_j))$ is equal to $\{\pm\id\}$. Since we're assuming that $\iota(\rho)\geq 4$, there is another generator $c_l\in\{c_1,\ldots,c_{2k}\}\setminus\{c_i,c_j\}$ with anti-diagonal image and we have $\rho(c_ic_j)\notin Z(\rho(c_l))$ because $\rho(c_ic_j)$ is a non-trivial diagonal matrix. This means that the Dehn twist along the curve $c_ic_j$ doesn't fix $[\rho]$. So, in order for the orbit of $[\rho]$ to be finite, $\rho(c_ic_j)$ must have finite order.

We continue by proving that all the diagonal matrices $\rho(c_{2k+1}),\ldots,\rho(c_n)$ have finite order too. If this wouldn't be the case, say some matrix $\rho(c_l)\in \{\rho(c_{2k+1}),\ldots,\rho(c_n)\}$ has infinite order, then the matrix $\rho(c_1c_2c_l)$ would be diagonal and of infinite order because $\rho(c_1c_2)$ has finite order by what we just proved. The Dehn twist along the curve $c_1c_2c_l$ would then act non-trivially on $[\rho]$ and produce infinitely many orbit points.

Now, let $\Gamma$ denote the image of $\rho$ and $\Lambda$ be the subgroup of $\Gamma$ generated by the elements $\{\rho(c_ic_j):1\leq i\neq j\leq 2k\}\cup \{\rho(c_{2k+1}),\ldots,\rho(c_n)\}$. The group $\Lambda$ is generated by diagonal matrices and is therefore abelian. Moreover, we just proved that all the generators of $\Lambda$ have finite order. This shows that $\Lambda$ is a finite group. The generators of $\Gamma$ are $\rho(c_1),\ldots,\rho(c_n)$ and are also all of finite order (anti-diagonal matrices have order 4). This implies that $\Lambda$ is a finite index subgroup of $\Gamma$, proving that $\Gamma$ is a finite group.
\end{proof}

\section{Finite mapping class orbits for 4-punctured spheres}\label{sec:4-punctured}

\subsection{Overview}
This section is a recap of the classification of finite mapping class group orbits for 4-punctured spheres. Before stating the classification in Section~\ref{sec:LT-classification}, we start by recalling the vocabulary in which the list of orbits is traditionally expressed in Section~\ref{sec:fricke-relation}. We then explain in Section~\ref{sec:finite-orbit-n=4} how to extract the list of peripheral angle vectors $\alpha$ satisfying the angle condition~\eqref{eq:angle-condition} and for which the DT component $\RepDT{\alpha}$ contains a finite mapping class group orbit. The resulting list can be found in Table~\ref{tab:finite-mcg-orbits-n=4} in Appendix~\ref{app:tables}.

\sam{For the reader's convenience, we compiled most of the notation introduced and used throughout Chapter~\ref{sec:4-punctured} in the following table (see also Notation~\ref{notation:4-punctured-sphere}).
\begin{table}[ht]
    \centering
    \begin{tblr}{width=\linewidth, colspec={X[1,l]|X[2,l]}}
		Notation & Meaning\\
		\hline\hline
		$\Sigma$ & a $4$-punctured sphere\\
		\hline
		$\mathcal{P}$ & the set of punctures on $\Sigma$\\
		\hline
		$\alpha\in(0,2\pi)^{\mathcal{P}}$ & a peripheral angle vector satisfying $\sum_{p\in\mathcal{P}}\alpha_p>6\pi$\\
		\hline
		$\RepDT{\alpha}$ & the DT component corresponding to $\alpha$ (isomorphic to $\CP^1$)\\
		\hline
		$c_1,c_2,c_3,c_4$ & a set of geometric generators of $\pi_1\Sigma$, i.e.~
		$\pi_1\Sigma=\langle c_1,c_2,c_3,c_4:c_1c_2c_3c_4=1\rangle$\\
		\hline
		$b,d$ & the fundamental group elements $b=(c_1c_2)^{-1}=c_3c_4$ and $d=(c_2c_3)^{-1}=c_4c_1$\\
		\hline
		$\tau_b,\tau_d$ & the Dehn twists along the simple closed curves represented by $b$, respectively by $d$\\
		\hline
		$\beta, \delta\colon\RepDT{\alpha}\rightarrow (0,2\pi)$ & the maps which associate to $[\rho]$ the rotation angle of $\rho(b)$, respectively $\rho(d)$\\
		\hline
		$\mathcal{B},\mathcal{D}$ & the pants decompositions of $\Sigma$ given by $b$, respectively by $d$ \\
		\hline
		$(\beta,\gamma)$ & the action-angle coordinates on $\RepDT{\alpha}$ defined by $\mathcal{B}$ (Section~\ref{sec:DT-representations})\\
		\hline
		$C_1,C_2,C_3,C_4$ & the exterior vertices of a $\mathcal{B}$-triangle chain or a $\mathcal{D}$-triangle chain\\
		\hline
		$B, D$ & the shared vertex of a $\mathcal{B}$-triangle chain, respectively of a $\mathcal{D}$-triangle chain\\
		\hline
		north/south poles & the two singular points of $\beta$ on $\RepDT{\alpha}$, characterized by $\beta=\alpha_3+\alpha_4-2\pi$, respectively $\beta=4\pi-\alpha_1-\alpha_2$\\
		\hline
		$\mathbf{t}\in \C^{\mathcal{P}}$ & a quadruple of peripheral traces, related to $\alpha$ by $\mathbf{t}_p=\pm 2\cos(\alpha_p/2)$\\
		\hline
		$\mathbf{t}=(a,b,c,d)$ & alternative notation for the quadruple of peripheral traces $\mathbf{t}$\\
		\hline
		$(A,B,C,D)$ & the Fricke coefficients (Definition~\ref{defn:Fricke-coefficients}) associated to a quadruple of traces and defined by $A=ab+cd$, $B=bc+ad$, $C=ac+bd$, and $D=4-a^2-b^2-c^2-d^2-abcd$\\
\end{tblr}
\end{table}
}

\subsection{Fricke relation}\label{sec:fricke-relation}
\cm{Let $\Sigma$ be a 4-punctured sphere and $\mathcal{P}$ be the set of punctures.} What makes the case of the 4-punctured sphere so special is the explicit algebraic description of relative character varieties of representations into $\SL_2\C$. To achieve this correspondence, we define the relative $\SL_2\C$-character varieties associated to a quadruple of traces $\mathbf{t}\in \C^{\mathcal{P}}$ as the algebraic \cm{GIT} quotient of the space of all representations that map peripheral loops around each puncture $p\in \mathcal{P}$ to elements of $\SL_2 \C$ with trace $\mathbf{t}_p$.\footnote{\cm{This definition of relative $\SL_2\C$-character varieties in terms of traces is only slightly more general than the one given in Definition~\ref{defn:relative-character-variety}. The point is that every trace $\neq \pm 2$ corresponds to a unique conjugacy class in $\SL_2\C$, while there are two conjugacy classes with trace $2$ (the trivial one and a parabolic one), and two with trace $-2$.}} We'll denote it by $\Rep_\mathbf{t}(\Sigma,\SL_2\C)$. We could alternatively define $\Rep_\mathbf{t}(\Sigma,\SL_2\C)$ as the space of conjugacy classes of reductive representations \cm{(Definition~\ref{defn:reductive-and-irreducible-representations})} with the identical peripheral behaviour.

\begin{rem}\label{rem:relation-between-traces-and-angles}
\cm{Earlier in Section~\ref{sec:character-varieties}, for peripheral angle vectors $\alpha\in (0,2\pi)^\mathcal{P}$, we introduced $\alpha$-relative character varieties as the space of conjugacy classes of representations $\rho\colon\pi_1\Sigma\to\psl$ mapping peripheral loops around the punctures $p\in \mathcal{P}$ to elliptic elements of $\psl$ with rotation angle $\alpha_p$. Since $\pi_1\Sigma$ is a free group, $\rho$ can be lifted to a representation $\overline{\rho}\colon\pi_1\Sigma\to\SL_2\R\subset \SL_2\C$. If $\overline{\rho}$ is reductive, then the conjugacy class of $\overline{\rho}$ belongs to $\Rep_\mathbf{t}(\Sigma,\SL_2\C)$ for some $\mathbf{t}$. The relation between $\mathbf{t}$ and $\alpha$ is given by $\mathbf{t}_p=\pm 2\cos(\alpha_p/2)$, where the signs depend on the choice of lift of $\rho$.}
\end{rem}

\cm{For simplicity and to stick with the traditional notation,} we'll fix an auxiliary ordering of the punctures of $\Sigma$ and write $\mathbf{t}=(a,b,c,d)$. The relative character variety $\Rep_\mathbf{t}(\Sigma,\SL_2\C)$ is isomorphic to the affine algebraic variety
\[
S_{(A,B,C,D)}=\{(X,Y,Z)\in\C^3 : X^2+Y^2+Z^2+XYZ-AX-BY-CZ-D=0\},
\]
where the numbers $(A,B,C,D)$ are complex coefficients that depend on $(a,b,c,d)$ via
\[
A=ab+cd, \quad B=bc+ad,\quad C=ac+bd, \quad D=4-a^2-b^2-c^2-d^2-abcd.
\]
The isomorphism from $\Rep_\mathbf{t}(\Sigma,\SL_2\C)$ to $S_{(A,B,C,D)}$ maps a conjugacy class of representations $[\rho]$ to the traces of $\rho$ evaluated at three different simple closed curves on $\Sigma$. A detailed account of the correspondence can be found in~\cite[Section~1.1]{cantat-loray}. 
\begin{defn}\label{defn:Fricke-coefficients}
We'll refer to the quadruple $(A,B,C,D)$ as the \emph{Fricke coefficients} and we'll call the polynomial relation that defines $S_{(A,B,C,D)}$ the \emph{Fricke relation}.
\end{defn}
The action of the pure mapping class group of $\Sigma$ on $S_{(A,B,C,D)}$ coming from the action on $\Rep_\mathbf{t}(\Sigma,\SL_2\C)$ can be made explicit in the coordinates $(X,Y,Z)$, see~\cite[Section~1.2]{cantat-loray} for precise formulae.

Different quadruples of traces may lead to the same Fricke coefficients. This is reminiscent \cm{of} the fact that different relative character varieties may be associated with the same Fricke relation. Cantat--Loray proved that there are at most twenty-four quadruples of traces corresponding to the same Fricke coefficients and that they can all be obtained from each other by explicit transformations~\cite[Lemma 2.7]{cantat-loray}. Some of them are \emph{Okamoto transformations}. Okamoto transformations are explained in detail in~\cite[Section 2.5]{cantat-loray} and will play an important role later in Section~\ref{sec:finite-orbit-n=4}. Cantat--Loray also proved that if $(A,B,C,D)$ are real numbers and if the smooth locus of the real points of $S_{(A,B,C,D)}$ has a bounded component, then all the corresponding twenty-four quadruples of traces consist of real numbers. Depending on the quadruple of traces, the representations in the bounded component are conjugate inside $\SL_2 \R$ or $\SU(2)$~\cite[Theorem B]{cantat-loray}. It's possible to use Okamoto transformations to switch from a quadruple of traces corresponding to representations into $\SU(2)$ to a quadruple of traces corresponding to representations into $\SL_2 \R$.

\subsection{The classification}\label{sec:LT-classification}
Building \cm{upon} the work of Boalch and others (as related in Section~\ref{sec:history-and-motivations}), Lisovyy--Tykhyy completed the classification of all the quadruples of Fricke coefficients $(A,B,C,D)$ for which $S_{(A,B,C,D)}$ contains a finite mapping class group orbit. Each quadruple of Fricke coefficients corresponds to a \emph{type} of finite orbits, which is really just an equivalence class of finite orbits up to Okamoto transformations. This means that a relative $\SL_2\C$-character variety contains a finite mapping class group orbit if and only if its defining quadruple of traces corresponds to a quadruple of Fricke coefficients in the following list.

\begin{thm}[\cite{LT}]\label{thm:LT-classification}
A finite mapping class group orbit in a relative $\SL_2\C$-character variety of a 4-punctured sphere is of one of the following types.
\begin{itemize}
    \item \emph{(Type \Romannum{1})} These orbits are the fixed points of the mapping class group action. \cm{They correspond to abelian representations or to representations that map a peripheral loop to $\pm\id$ by Lemma~\ref{lem:fixed-point-implies-abelian}.} They actually coincide with singular or isolated points of the relative character variety by a result of~\cite{iwasaki-katsunori-masa-hiko}.
    \item \emph{(Types \Romannum{2}--\Romannum{4})} These finite orbits come into 1-parameter or 2-parameter \cm{complex} families. Their respective length are $2$, $3$, and $4$. The corresponding Fricke coefficients are listed in~\cite[Lemma 39]{LT}.
    \item \emph{(Exceptional orbits)} There are $45$ exceptional finite orbits with a length ranging from $5$ to $72$. They are listed, along with the corresponding Fricke coefficients, in~\cite[Table 4]{LT}.
    \item \emph{(Cayley orbits)} These orbits occur in the case where the quadruple of Fricke coefficients is $(A,B,C,D)=(0,0,0,4)$. The corresponding quadruples of traces are $(0,0,0,0)$ and permutations of $(\pm 2,\pm 2,\pm 2, \mp 2)$.
\end{itemize}
\end{thm}

We already encountered several examples of finite mapping class group for 4-punctured spheres. For instance, the representation $\rho_\theta$ from Section~\ref{sec:arbitrarily-long-orbits} gives rise to a finite orbit of Cayley type when $\theta\in (0,2\pi)$ and of Type~\Romannum{1} when $\theta\in\{0,2\pi\}$. The triangle chain from Example~\ref{ex:example-Dehn-twist-action} corresponds to one of the exceptional orbits which we'll discuss in details in Lemma~\ref{lem:beta_i-for-obits-of-type-8} below. The pullback orbit from Example~\ref{ex:pullback} is a finite orbit of Type~\Romannum{2} which we'll analyse in Lemma~\ref{lem:beta_i-for-orbits-of-type-II}.

A careful analysis of the finite orbits of Types~\Romannum{2}--\Romannum{4} and the 45 exceptional orbits shows that the corresponding representations are always conjugate inside $\SL_2\R$ or $\SU(2)$---the real forms of $\SL_2\C$. \cm{A given orbit type---in other words, an Okamoto equivalence class of finite orbits---always contains finite orbits made of conjugacy classes of representations into $\SL_2\R$ and into $\SU(2)$. The fact that the two eventualities always occur is a consequence of~\cite[Theorem~B]{cantat-loray}.} 
\begin{fact}\label{fact:finite-SL2-orbit-in-compact-component}
The conjugacy class of a representation into $\SL_2\R$ that belongs to an exceptional orbit or to a finite orbit of Types~\Romannum{2}--\Romannum{4} necessarily belongs to the compact component of the corresponding relative $\SL_2\R$-character variety. 
\end{fact}
\begin{proof}
\cm{The above discussion shows that a finite orbit $\mathcal{O}$ coming from $\SL_2\R$ representations which is either exceptional or of Type~\Romannum{2}--\Romannum{4} is equivalent to a finite orbit $\mathcal{O}'$ coming from $\SU(2)$ representations under a certain Okamoto transformation. Okamoto transformations not only map finite orbits to finite orbits, they map entire connected components of relative character varieties homeomorphically to each other. Since $\SU(2)$ is compact, the component containing $\mathcal{O}'$ is compact, and hence so is the component containing $\mathcal{O}$.}
\end{proof}

Lisovyy--Tykhyy only \cm{exhibit} one orbit of each type in~\cite[p.155-162]{LT}. Cousin worked out all their Okamoto transformations in~\cite[p.49]{cousin-thesis}. Cousin also determined which \cm{ones} are pullback orbits (Definition~\ref{def:pullback-orbits}). It turns out that at least one orbit of each exceptional type is a pullback (see~\cite[Théorème~2.5.1]{cousin-thesis} for a precise statement). Among the 45 exceptional orbit types, 41 are the equivalence classes under Okamoto transformations of a representation with \emph{finite} image in $\SU(2)$ (1 tetrahedral, 7 octahedral, and 33 icosahedral). The remaining 4 orbit types are those numbered 8, 32, 33, and 34 by Lisovyy--Tykhyy. They correspond, respectively, to the \emph{Klein solution} from~\cite{boalch-3} and the three \emph{elliptic 237 solutions} discovered by Boalch~\cite{boalch-2} and Kitaev~\cite{kitaev-2}.\footnote{This sub-classification of the 45 exceptional orbit types motivates Boalch's nomenclature:
$$
T06,\quad O07-O13, \quad I20-I52,\quad 237 \text{ (three times)},\quad K. 
$$}
It's quite remarkable \cm{that} all their Okamoto transforms have infinite image, even those with image in $\SU(2)$. Finite orbits of Type~8 have one more particularity: they're all of pullback type (the only other exceptional orbit type with this property is Type~30).\footnote{Finite orbits of Types~8 and~30 are also characterized by the property that Okamoto transformations coincide with Galois conjugation. This is a consequence of Cousin's list of Galois equivalence classes~\cite[p.51]{cousin-thesis}.} We'll encounter finite orbits of Types~8 and~33 again later in Section~\ref{sec:finite-orbits-of-type-8-33}.

Let us rapidly say a few words about Lisovyy--Tykhyy's proof of Theorem~\ref{thm:LT-classification}. The first step in their proof is a careful study of the possible orders of a Dehn twist when iterated on a point in a finite orbit. They work out explicit formulae for the traces of simple closed curves after iterations by a Dehn twist (\cite[Lemma 13]{LT}). In order to infer bounds on the lengths of finite orbits, a technical lemma is required. It lists all rational solutions to equations of the type $\sum_i \cos(2\pi x_i)=0$ for a certain number of variables $x_i$ (\cite[Section~2.4]{LT}). The last step uses a computer to run an algorithm that searches for all the finite orbits while making use of the bounds on their lengths (\cite[Section~2.6]{LT}).

\subsection{Finite orbits inside DT components}\label{sec:finite-orbit-n=4}
For the purpose of this work, we are interested in identifying all the angle vectors $\alpha\in (0,2\pi)^\mathcal{P}$ satisfying $\Sigma_{p\in\mathcal{P}}\alpha_p>6\pi$ and such that the DT component $\RepDT{\alpha}$ of the relative character variety $\Rep_\alpha(\Sigma,\psl)$ \cm{(Definition~\ref{defn:relative-character-variety})} contains a finite mapping class group orbit. Recall \cm{from Fact~\ref{fact:finite-SL2-orbit-in-compact-component}} that when $\Sigma$ is a 4-punctured sphere, if $\Rep_\alpha(\Sigma,\psl)$ contains a finite orbit, it automatically lies in $\RepDT{\alpha}$.

\cm{We'll work with an auxiliary ordering of the punctures of $\Sigma$, as we did in Section~\ref{sec:fricke-relation}, and write $\alpha=(\alpha_a,\alpha_b,\alpha_c,\alpha_d)$. Recall that angles and traces are related by $t=\pm 2\cos(\alpha_t/2)$ for $t\in\{a,b,c,d\}$.} We first observe that such an orbit is never of Type~\Romannum{1} since the assumption $\alpha_a+\alpha_b+\alpha_c+\alpha_d>6\pi$ implies that the relative character variety $\Rep_\alpha(\Sigma,\psl)$ is smooth and has no isolated points \cm{by Fact~\ref{fact:rep_alpha-smooth}}. Similarly, Cayley orbits can't occur since, in order to satisfy $(A,B,C,D)=(0,0,0,4)$ while picking angles $\alpha\in(0,2\pi)^4$, one has to take $\alpha=(\pi,\pi,\pi,\pi)$ \cm{(see Remark~\ref{rem:relation-between-traces-and-angles})}, which violates the inequality. So, the only candidates are the finite orbits of Type \Romannum{2}--\Romannum{4} or the exceptional ones. In order to determine all the angle vectors $\alpha$ for which $\RepDT{\alpha}$ admits finite orbits, we apply the following routine (a version of which was employed by Cousin in~\cite[Section~2.2]{cousin-thesis}).

\begin{enumerate}
    \item  For each of the Types \Romannum{2}--\Romannum{4} and for each exceptional orbit type, we first determine all the quadruples of traces whose corresponding relative $\SL_2 \C$-character varieties contains a finite orbit of the given type. The information we get from~\cite{LT} is the list of all quadruples of Fricke coefficients \cm{(Definition~\ref{defn:Fricke-coefficients})}, along with one corresponding quadruple of traces. In order to compute the twenty-four quadruples of traces (\cite[Lemma 2.7]{cantat-loray}) corresponding to the same quadruple of Fircke coefficients, we proceed as follows. \cm{Following the notation of~\cite{LT}, we start with the quadruple of traces $(a,b,c,d)=(2\cos(\pi\theta_a), 2\cos(\pi\theta_b), 2\cos(\pi\theta_c), 2\cos(\pi\theta_d))$, where $\bm{\theta}=(\theta_a,\theta_b,\theta_c,\theta_d)$ is provided in~\cite[p.~156--162]{LT}.} We apply the following two Okamoto transformations to the vector $\bm{\theta}$ in order to obtain two new vectors and thus two new quadruple of traces:
    \begin{align*}
    \Ok(\bm{\theta})&=\frac{1}{2}\begin{pmatrix}
    1 & -1 & -1 & -1\\
    -1 & 1 & -1 & -1\\
    -1 & -1 & 1 & -1\\
    -1 & -1 & -1 & 1
    \end{pmatrix}\begin{pmatrix}
    \theta_a\\ \theta_b\\ \theta_c\\ \theta_d    
    \end{pmatrix}+\begin{pmatrix}
    1 \\ 1 \\ 1 \\ 1   
    \end{pmatrix},\\
    \widetilde{\Ok}(\bm{\theta})&=\frac{1}{2}\begin{pmatrix}
    1 & -1 & -1 & 1\\
    -1 & 1 & -1 & 1\\
    -1 & -1 & 1 & 1\\
    1 & 1 & 1 & 1
    \end{pmatrix}\begin{pmatrix}
    \theta_a\\ \theta_b\\ \theta_c\\ \theta_d    
    \end{pmatrix}.
    \end{align*}
    As explained in~\cite[Section~2.5]{cantat-loray}, all the possible twenty-four quadruples of traces can now be obtained from the three we just described by permuting the entries and switching all four signs. Since we are ultimately interested in relative $\psl$-character varieties, permuting traces and changing their sign will only result into a permutation of the entries of the peripheral angle vector. \cm{So, for each $\bm{\theta}$, we'll only consider the three quadruples of traces obtained from $\bm{\theta}$, $\Ok(\bm{\theta})$, and $\widetilde{\Ok}(\bm{\theta})$.}
    
    \cm{Among all the quadruples of traces obtained in this way, we only keep those whose entries are real numbers in the interval $(-2,2)$. This is because DT representations have non-trivial elliptic peripheral monodromies by definition (see Definition~\ref{defn:DT-representation}).}
    \item For each remaining quadruple of traces, we apply Benedetto--Goldman's criterion\footnote{The original criterion appears in~\cite{benedetto-goldman} as Proposition 6.1. However, there seems to be a computational mistake at the bottom of the first column, on page 102, leading to a wrong expression for $r^2-p-q$. A corrected version of the criterion was stated as Theorem 3.12 in~\cite{cantat-loray}.} to discriminate between representations into $\SU(2)$ and $\SL_2 \R$. The criterion says that the compact component of the real points of the relative $\SL_2\C$-character variety determined by the quadruple of traces $(a,b,c,d)$ consists of $\SL_2\R$ representations if
    \[
    2(a^2+b^2+c^2+d^2)-abcd-16>\sqrt{(4-a^2)(4-b^2)(4-c^2)(4-d^2)}.
    \]
    If the inequality is violated, then the compact component consists of representations into $\SU(2)$. We only keep the quadruples of traces corresponding to $\SL_2\R$ representations. 
    \item \cm{At this point, we are left with a bunch of quadruples of traces $\mathbf{t}=(a,b,c,d)$ that haven't been eliminated yet. There's a first subtlety in the fact that a trace in $(-2,2)$ doesn't determine a unique conjugacy class in $\SL_2\R$, but two. So, each remaining $\mathbf{t}$ is the peripheral trace vector of multiple relative $\SL_2\R$-character varieties (in the sense of Definition~\ref{defn:relative-character-variety}) whose compact components carry a finite orbit. All these compact components are mapped via the projection $\SL_2 \R\to \psl$ to DT components $\RepDT{\alpha}$ (Section~\ref{sec:definition-DT}). Our immediate goal is to determine all the peripheral angle vectors $\alpha$ that arise by this construction from all remaining $\mathbf{t}$.}
    
\cm{We proceed as follows. Let $\mathbf{t}=(a,b,c,d)$ be one of the remaining quadruple of traces. For every trace $t\in\{a,b,c,d\}\subset (-2,2)$, there are two angles $\overline{\theta}_t\in (0,4\pi)$ such that $2\cos(\overline{\theta}_t/2)=t$. If $\overline{\theta}_t$ denotes one of them, then $4\pi-\overline{\theta}_t$ is the other one. Assume that $\overline{\theta}_t < 4\pi-\overline{\theta}_t$, or equivalently that $\overline{\theta}_t<2\pi$. Under this assumption, the angles $\overline{\theta}_t$ and $2\pi-\overline{\theta}_t$ both lie in $(0,2\pi)$. Let $\overline{\alpha}_t=\max\{\overline{\theta}_t, 2\pi-\overline{\theta}_t\}$ and $\overline{\alpha}=(\overline{\alpha}_a, \overline{\alpha}_b, \overline{\alpha}_c, \overline{\alpha}_d)$. Since the quadruple $(a,b,c,d)$ satisfies Benedetto--Goldman's criterion, it holds that $\overline{\alpha}_a+\overline{\alpha}_b+\overline{\alpha}_c+\overline{\alpha}_d>6\pi$. Assume that $\overline{\alpha}_a\geq \overline{\alpha}_b\geq \overline{\alpha}_c\geq \overline{\alpha}_d$. There is at most one other combination $\overline{\alpha}'=(\overline{\alpha}'_a, \overline{\alpha}'_b, \overline{\alpha}'_c, \overline{\alpha}'_d)$ with $\overline{\alpha}'_t\in\{\overline{\theta}_t, 2\pi-\overline{\theta}_t\}$ which also satisfies $\overline{\alpha}'_a+\overline{\alpha}'_b+\overline{\alpha}'_c+\overline{\alpha}'_d>6\pi$. It is given by $\overline{\alpha}'=(\overline{\alpha}_a, \overline{\alpha}_b, \overline{\alpha}_c, 2\pi-\overline{\alpha}_d)$. 
The angle vector $\alpha$ that we are looking for is either $\overline{\alpha}$ or $\overline{\alpha}'$; only one of the two is the peripheral angle vector of $\psl$ representations that can be lifted to $\SL_2\R$ with peripheral traces $\mathbf{t}=(a,b,c,d)$ (the other one correspond to representations with peripheral traces $(-a,b,c,d)$ or any other quadruple with an odd number of minus signs). So, in order to decide between $\overline{\alpha}$ and $\overline{\alpha}'$, we apply the following criterion. If $\alpha_d=\pi$, then $\overline{\alpha}=\overline{\alpha}'$ and there is no decision to make. If $\alpha_d\neq \pi$, then we look at the sign of $abcd\neq 0$. If $abcd<0$, then we take $\alpha=\overline{\alpha}$. If $abcd>0$, then we take instead $\alpha=\overline{\alpha}'$. This leads to the right angle vector by Lemma~\ref{lem:sign-product-of-traces-n=4}.}
\end{enumerate}

\cm{We first apply the above routine to finite orbits of Type~\Romannum{2}--\Romannum{4} for which, respectively,
\[
\bm{\theta}=(\theta_a,\theta_b,\theta_b,1-\theta_a),\quad (2\theta,\theta,\theta,2/3),\quad (\theta,\theta,\theta,1/2),
\]
where $\theta,\theta_a,\theta_b$ are complex parameters. Noting that $2\cos(\pi z)\in\R\cap (-2,2)$ if and only if $z\in \R\setminus \Z$, we may assume that $\theta,\theta_a,\theta_b$ are real numbers.} Step (1) is a simple computation. For each quadruple of traces obtained in (1), we need to determine the range of the parameters $\theta,\theta_a,\theta_a$ for which the quadruple of traces satisfies Benedetto--Goldman's criterion. This amounts to solving an inequality. We confirmed our computations with a computer algebra system (namely Mathematica). The resulting angle vectors $\alpha\in (0,2\pi)^4$ are:
\begin{itemize}
\item \cm{For orbits of Type~\Romannum{2}, all angle vectors $\alpha$ whose entries are a permutation of $\{\theta_1,\theta_1,\theta_2,\theta_2\}$ for some angles $\theta_i\in(0,2\pi)$ satisfying $\theta_1+\theta_2>3\pi$.}
\item \cm{Orbits of Types~\Romannum{3} give rise to the angle vectors $\alpha$ whose entries are a permutation of $\{4\pi/3, 2\theta-2\pi, \theta, \theta\}$ for an angle $\theta\in (5\pi/3,2\pi)$.} 
\item \cm{The routine produces two families of angle vectors $\alpha$ for orbits of Type~\Romannum{4} whose entries are permutations of $\{\pi, \theta, \theta, \theta\}$ and $\{\theta, \theta, \theta, 3\theta-4\pi\}$, where $\theta\in (5\pi/3,2\pi)$ in both cases. To avoid confusion, we'll refer to orbits of the former type by \emph{Type~\Romannum{4} orbits} and of the latter type by \emph{Type~\Romannum{4}$^\ast$ orbits}. }
\end{itemize}
The case of the 45 exceptional orbit types is more straightforward because there are no abstract parameters involved. We again verified our computations with a computer (the code is available in our GitHub repository).\footnote{Interestingly enough, the exceptional orbit types that don't correspond to any finite orbit in a DT component (2, 3, 5, 9, 16, 17, 21, 28, 29, 31, 35, 36, and 42) are precisely those that were ``rayé(s)'' by Cousin~\cite[p.48-49]{cousin-thesis} because they lead to same monodromy up permuting traces and changing their signs.} The resulting list of angle vectors can be found in Table~\ref{tab:finite-mcg-orbits-n=4} in Appendix~\ref{app:tables} (the list is coherent with those obtained in~\cite{cousin-thesis, tykhyy}).

\subsection{Some particular orbits in detail}\label{sec:some-finite-orbits}
In order to make the list of Table~\ref{tab:finite-mcg-orbits-n=4} from Appendix~\ref{app:tables} more concrete, we describe some finite orbits in more details by computing the action-angle coordinates from Section~\ref{sec:DT-representations} of every orbit point. We don't plan to carry out this analysis in every case, but only for some particular orbits that will be relevant later. We take the opportunity to describe the image of the representations in these finite orbits. \cm{Throughout the section, we'll use the following notation.}

\begin{notation}\label{notation:4-punctured-sphere}
\cm{We'll work with a fixed geometric presentation of $\pi_1\Sigma$ given by
\[
\pi_1\Sigma=\langle c_1,c_2,c_3,c_4: c_1c_2c_3c_4=1\rangle.
\]
The fundamental group element 
\[
b=(c_1c_2)^{-1}
\]
represents a simple closed curve on $\Sigma$ and defines a pants decomposition $\mathcal{B}$ of $\Sigma$. The $\mathcal{B}$-triangle chain of a point in a DT component $\RepDT{\alpha}$ consists of two triangles with vertices $(C_1,C_2,B)$ and $(B,C_3,C_4)$. The associated action-angle coordinates defined in Section~\ref{sec:action-angle-coordinates} are denoted by $(\beta,\gamma)$. (We're not adding indices to $b$ and $(\beta,\gamma)$ since $\mathcal{B}$ consists of a unique curve here.)}

\cm{We remind the reader that the angle coordinate $\gamma$ is well defined for all but two points of $\RepDT{\alpha}$. These two points are the two singular points of $\beta$ (where the maximum and minimum of $\beta$ are achieved). We'll refer to the one with $\beta=\alpha_3+\alpha_4-2\pi$ as the \emph{north pole} of $\RepDT{\alpha}$ and to the one with $\beta=4\pi-\alpha_1-\alpha_2$ as the \emph{south pole}. Equivalently, the singular points of $\beta$ are those whose $\mathcal{B}$-triangle chain is singular (Definition~\ref{def:regular-singular-triangle-chains}), i.e.~the points for which $C_1=C_2=B$ (south pole) or $B=C_3=C_4$ (north pole). They are also the only two fixed points of the Dehn twist $\tau_b$ by Fact~\ref{lem:fixed-points-Dehn-twists}.}

\cm{We'll also need the fundamental group element 
\[
d=(c_{2}c_{3})^{-1}
\]
which represents the simple closed curve on $\Sigma$ illustrated below.}
\begin{center}
\vspace{2mm}
\begin{tikzpicture}[scale=1.1, decoration={
    markings,
    mark=at position 0.7 with {\arrow{>}}}]
  \draw[postaction={decorate}] (2,-.5) arc(270:90: .25 and .5) node[midway, left]{$c_1$};
  \draw[smoked] (2,.5) arc(90:-90: .25 and .5);
  \draw[apricot, postaction={decorate}] (4,.5) arc(90:270: .25 and .5) node[midway, below left]{$b$};
  \draw[lightapricot] (4,.5) arc(90:-90: .25 and .5);
  \draw[postaction={decorate}] (6,.5) arc(90:270: .25 and .5) node[midway, left]{$c_4$};
  \draw (6,.5) arc(90:-90: .25 and .5);

  \draw[postaction={decorate}, mauve] (5.65, .65) arc(0:-180: 1.65 and .5) node[near end, below]{$d$};
  \draw[lightmauve] (2.35, .65) arc(180:0: 1.65 and -.25);
  
  \draw (2.5,1) arc(180:0: .5 and .25)node[midway, above]{$c_2$};
  \draw[postaction={decorate}] (2.5,1) arc(-180:0: .5 and .25);
  \draw (4.5,1) arc(180:0: .5 and .25)node[midway, above]{$c_3$};
  \draw[postaction={decorate}] (4.5,1) arc(-180:0: .5 and .25);
  
  \draw (2,.5) to[out=0,in=-90] (2.5,1);
  \draw (3.5,1) to[out=-90,in=180] (4,.5);
  \draw (2,-.5) to[out=0,in=180] (4,-.5);
  
  \draw (4,.5) to[out=0,in=-90] (4.5,1);
  \draw (5.5,1) to[out=-90,in=180] (6,.5);
  \draw (4,-.5) to[out=0,in=180] (6,-.5);
\end{tikzpicture}
\vspace{2mm}
\end{center}
\cm{The Dehn twists along the curves $b$ and $d$ are $\tau_b$ and $\tau_d$. Note that those Dehn twists were denoted by $\tau_{1,2}$ and $\tau_{2,3}$ in Section~\ref{sec:action-of-Dehn-twists}, but we favour a more concise notation here. The Dehn twists $\tau_b$ and $\tau_d$ together generate $\PMod(\SpherePk)$. This is a consequence of Lemma~\ref{lem:mcg-generators} from Appendix~\ref{apx:generators-of-pmod}. We'll also denote by 
\[
\delta\colon\RepDT{\alpha}\to (0,2\pi)
\]
the function that maps $[\rho]$ to the rotation angle (Definition~\ref{defn:rotation-angle}) of $\rho(d)$.}

\cm{The simple closed curve represented by $d$ defines a pants decomposition of $\SpherePk$ which we'll denote by $\mathcal{D}$. This is the pants decomposition $\mathcal{B}[2]$ in the notation of Section~\ref{sec:action-of-Dehn-twists} and the cyclically permuted generators $(c_2,c_3,c_4,c_1)$ are compatible with the pants decomposition $\mathcal{D}$ in the sense of Definition~\ref{defn:compatible-generators}. The $\mathcal{D}$-triangle chain of a point in $\RepDT{\alpha}$ consists of two triangles with vertices $(C_2,C_3,D)$ and $(D,C_4,C_1)$. Recall from Section~\ref{sec:action-of-Dehn-twists} and Example~\ref{ex:example-Dehn-twist-action} as well that we can juggle between the $\mathcal{B}$-triangle chain and the $\mathcal{D}$-triangle chain of an orbit point in $\RepDT{\alpha}$ to compute its images under $\tau_b$ and $\tau_d$.}
\ensuremath{\lozenge}
\end{notation}

\cm{Before diving into computations,} it's also useful to observe that there can be at most two points in $\RepDT{\alpha}$ on which $\beta$ and $\delta$ agree.

\begin{fact}\label{fact:beta-delta-given-at-most-two-points}
For every pair $(x,y)$ of real numbers, the intersection $\beta^{-1}(x)\cap \delta^{-1}(y)\subset \RepDT{\alpha}$ consists of at most two points.
\end{fact}
\begin{proof}
\cm{All the representations whose conjugacy class lie inside $\RepDT{\alpha}$ can be lifted to representations into $\SL_2\R\subset \SL_2\C$. Those lifts are irreducible (Fact~\ref{fact:rep_alpha-smooth}), hence reductive (Definition~\ref{sec:character-varieties}).} Their conjugacy classes are therefore determined by the traces of their evaluations on three simple closed curves of $\Sigma$ as we explained in Section~\ref{sec:fricke-relation}. The values of $\beta$ and $\delta$ determines two of the traces. The possible values for the third trace are obtained by solving a second degree polynomial equation (the Fricke relation). So, there are at most two different values for the third trace which proves the claim.
\end{proof}

\subsubsection{Orbits of Type~\Romannum{2}}
Let's start with finite orbits of Type~\Romannum{2}. They are the only finite orbits of length 2 inside DT components and their discovery is usually attributed to Hitchin~\cite{hitchin-1}. The peripheral angles for an orbit of Type~\Romannum{2} are $\{\theta_1,\theta_1,\theta_2,\theta_2\}$ for two angles $\theta_i\in(0,2\pi)$ satisfying $\theta_1+\theta_2>3\pi$.

\begin{lem}\label{lem:beta_i-for-orbits-of-type-II}
The action-angle coordinates of orbit points in a finite orbit of Type~\Romannum{2} are provided in the table below. \cm{Each line corresponds to an angle vector $\alpha$ obtained by ordering the entries of $\{\theta_1,\theta_1,\theta_2,\theta_2\}$ in one of three possible ways. When $\theta_1\neq \theta_2$, the DT component $\RepDT{\alpha}$ contains a unique finite orbit of length $2$. If $\theta_1=\theta_2$, then $\RepDT{\alpha}$ contains all three finite orbits of length $2$ whose coordinates are provided in the table below.} The image of a representation associated to a finite orbit of Type~\Romannum{2} is a rotation triangle group\footnote{As defined in Section~\ref{sec:triangle-groups}.} $D(2,\overline\theta_1,\overline\theta_2)\subset \psl$, where $\overline\theta_i=(1-\theta_i/2\pi)^{-1}$.
\begin{table}[ht]
    \centering
    \begin{tblr}{c|c|c}
        $\alpha$ & $\beta$ & $\gamma$ \\
        \hline\hline
        \SetCell[r=2]{c} $(\theta_1,\theta_1,\theta_2,\theta_2)$ & $4\pi-2\theta_1$ & south pole \\
        \hline
        & $2\theta_2-2\pi$ & north pole \\
        \hline
        $(\theta_1,\theta_2,\theta_1,\theta_2)$ & $\pi$ & $\{\frac{\pi}{2} ,\frac{3\pi}{2}\}$ \\
        \hline
		$(\theta_1,\theta_2,\theta_2,\theta_1)$ & $\pi$ & $\{0, \pi\}$ \\
    \end{tblr}
\end{table}
\end{lem}
\begin{proof}
We start with the case where the vector $\alpha$ of peripheral angles is $(\theta_1,\theta_1,\theta_2,\theta_2)$. The Dehn twist $\tau_b$ has two fixed points: the north pole and the south pole of $\RepDT{\alpha}$. The the north pole is characterized by $\beta=2\theta_2-2\pi$ and its $\mathcal{B}$-triangle chain consists of a single triangle with vertices $(C_1,C_2,C_3=C_4)$. The interior angles are, respectively, $(\pi-\theta_1/2, \pi-\theta_1/2, 2\pi-\theta_2)$. The south pole has $\beta=4\pi-2\theta_1$. Its $\mathcal{B}$-triangle chain consists of a single triangle with vertices $(C_1=C_2,C_3,C_4)$; and with respective interior angles $(2\pi-\theta_1, \pi-\theta_2/2, \pi-\theta_2/2)$. 
\begin{center}
\resizebox{10cm}{!}{
\begin{tikzpicture}[font=\sffamily,decoration={
    markings,
    mark=at position 1 with {\arrow{>}}}]
    
\node[anchor=south west,inner sep=0] at (0,0) {\includegraphics[width=3cm]{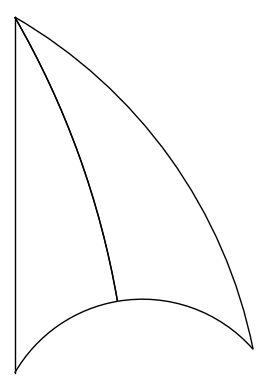}};
\node[anchor=south west,inner sep=0] at (7,0) {\includegraphics[height=5cm]{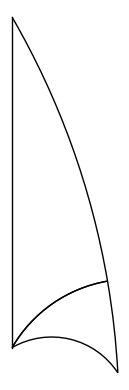}};

\begin{scope}
\fill (0.18,0.22) circle (0.07) node[left]{$C_2$};
\fill (0.18,4.2) circle (0.07) node[above]{$B=C_3=C_4$};
\fill (2.84,0.48) circle (0.07) node[right]{$C_1$};
\fill (1.35,1) circle (0.07) node[below]{$D$};
\end{scope}

\draw[mauve, thick] (.45,.55) arc (30:90:.3) node[at end, left]{\small $\pi-\theta_1/2$};
\draw[mauve, thick] (2.72,1) arc (105:163:.3) node[at start, right]{\small $\pi-\theta_1/2$};
\draw[apricot, thick] (0.18,3.5) arc (270:334:.6) node[at start, left]{\small $2\pi-\theta_2$};

\node[anchor=south west,inner sep=0] at (.7,-.3) {\textnormal{north pole}};

\begin{scope}
\fill (7.175,0.57) circle (0.07) node[left]{$C_1=C_2=B$};
\fill (7.175,4.8) circle (0.07) node[left]{$C_3$};
\fill (8.5,0.25) circle (0.07) node[below right]{$C_4$};
\fill (8.4,1.4) circle (0.07) node[right]{$D$};
\end{scope}

\draw[mauve, thick] (8.25,.5) arc (150:90:.3) node[at end, right]{\small $\pi-\theta_2/2$};
\draw[mauve, thick] (7.175,4.3) arc (270:310:.3) node[at end, right]{\small $\pi-\theta_2/2$};
\draw[apricot, thick] (7.6,0.68) arc (5:90:.45) node[at end, left]{\small $2\pi-\theta_1$};

\node[anchor=south west,inner sep=0] at (6.8,-.3) {\textnormal{south pole}};
\end{tikzpicture}
}
\end{center}

We claim that $\tau_d$ maps the north pole to the south pole, showing that the two poles form a finite orbit of length 2. To see that, we apply the procedure described in Section~\ref{sec:action-of-Dehn-twists}. We start by describing the $\mathcal{D}$-triangle chain of the north pole. The angle bisector of the triangle $(C_1,C_2,C_3=C_4)$ at the vertex $C_3=C_4$ cuts the geodesic segment $C_1C_2$ perpendicularly at a point $D$. The triangle with vertices $(C_2,C_3,D)$ has respective interior angles $(\pi-\theta_1/2, \pi-\theta_2/2,\pi/2)$. The triangle with vertices $(D,C_4,C_1)$ has respective interior angles $(\pi/2, \pi-\theta_2/2, \pi-\theta_1/2)$. The two triangles $(C_2,C_3,D)$ and $(D,C_4,C_1)$ form a chain of two triangles (chained at $D$) that has the same defining geometric features as the $\mathcal{D}$-triangle of the north pole \cm{(Lemma~\ref{lem:properties-triangle-chains})}. The two therefore coincide and we deduce that the north pole has $\delta=\pi$, where $\delta$ is the angle function associated to $d$. We can now compute the image of the north pole by $\tau_d$. The $\mathcal{D}$-triangle chain of the image is obtained from the $\mathcal{D}$-triangle chain of the north pole by rotating the triangle $(D,C_4,C_1)$ anti-clockwise around the vertex $D$ by an angle $\delta=\pi$. It's easy to see that the resulting triangle chain has $C_1=C_2$ and coincides with the $\mathcal{B}$-triangle chain of the south pole.

From the shape of the $\mathcal{D}$-triangle chain of the north pole, it's immediate that the image  of any representation in the conjugacy class of the north pole is a rotation triangle group $D(2,\overline\theta_1,\overline\theta_2)$.

We just described a finite mapping class group orbit of length 2 inside $\RepDT{\alpha}$ when $\alpha=(\theta_1,\theta_1,\theta_2,\theta_2)$. We now prove that $\RepDT{\alpha}$ does not contain any other orbit of length 2 if $\theta_1\neq \theta_2$. Assume for the sake of contradiction that it does. We denote the two hypothetical orbit points by $[\rho_1]$ and $[\rho_2]$. Since we excluded the orbit made of the two poles, we can assume that the $\mathcal{B}$-triangle chains of $[\rho_1]$ and $[\rho_2]$ are both regular \cm{(Definition~\ref{def:regular-singular-triangle-chains})}. In particular, $\tau_b$ does not fix any of these two points by \cm{Lemma~\ref{lem:fixed-points-Dehn-twists}} and therefore permutes them. We deduce that $\tau_b^2$ fixes both orbit points which implies that both $[\rho_1]$ and $[\rho_2]$ have $\beta=\pi$. It remains to determine the angle coordinates of $[\rho_1]$ and $[\rho_2]$. The two triangles in the $\mathcal{B}$-chain of $[\rho_1]$ have vertices $(C_1,C_2,B)$, respectively $(B,C_3,C_4)$. Both triangles are isosceles, with a right angle at $B$. The other two angles are equal to $\pi-\theta_1/2$, respectively $\pi-\theta_2/2$. Since we're assuming that $\theta_1\neq \theta_2$, we always have $C_2\neq C_3$ and $C_4\neq C_1$, no matter what the value of the angle coordinate $\gamma$ of $[\rho_1]$ is. This means, by \cm{Lemma~\ref{lem:fixed-points-Dehn-twists}}, that $\tau_d$ does not fix $[\rho_1]$ and therefore maps it to $[\rho_2]$, exactly like $\tau_b$. So, both $[\rho_1]$ and $[\rho_2]$ have $\delta=\pi$. This means that the $\mathcal{D}$-triangle chain of $[\rho_1]$ consists of two right triangles with vertices $(C_2,C_3,D)$, respectively $(D,C_4,C_1)$.

\begin{center}
\resizebox{5cm}{!}{
\begin{tikzpicture}[font=\sffamily,decoration={
    markings,
    mark=at position 1 with {\arrow{>}}}]
    
\node[anchor=south west,inner sep=0] at (0,0) {\includegraphics[width=3cm]{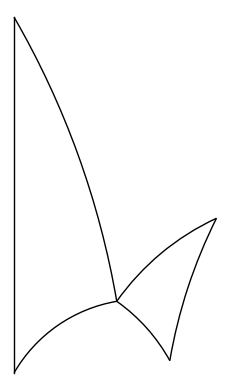}};

\begin{scope}
\fill (0.18,0.22) circle (0.07) node[left]{$C_1$};
\fill (0.18,4.9) circle (0.07) node[above]{$C_2$};
\fill (2.8,2.2) circle (0.07) node[above right]{$C_3$};
\fill (2.2,.4) circle (0.07) node[right]{$C_4$};
\fill (1.52,1.15) circle (0.07) node[below]{$B$};
\end{scope}

\draw[mauve, thick] (.45,.55) arc (33:90:.3) node[at end, left]{\small $\pi-\theta_1/2$};
\draw[mauve, thick] (0.18,4.2) arc (270:334:.3) node[at start, left]{\small $\pi-\theta_1/2$};
\draw[mauve, thick] (1.95,.8) arc (140:60:.3) node[at end, right]{\small $\pi-\theta_2/2$};
\draw[mauve, thick] (2.7,1.9) arc (250:200:.3) node[at start, right]{\small $\pi-\theta_2/2$};

\draw[sky, postaction={decorate}, thick] (2,1.7) arc (45:105:.6) node[near start, above]{\textnormal{?}};
\end{tikzpicture}
}
\end{center}

It's now time for some hyperbolic trigonometry. The hyperbolic law of cosines applied to the triangle $(C_1,C_2,B)$ and $(B,C_3,C_4)$ gives\footnote{Unless otherwise mentioned, we'll use $XY$ to denote the hyperbolic distance between two points $X$ and $Y$ in the hyperbolic plane.}
\[
\cosh(BC_2)=\frac{-\cos(\theta_1/2)}{\sin(\theta_1/2)} \quad\text{and}\quad \cosh(BC_3)=\frac{-\cos(\theta_2/2)}{\sin(\theta_2/2)}.
\]
Similarly, in the triangle $(C_2,C_3,D)$, the hyperbolic law of cosines gives
\[
\cosh(C_2C_3)=\frac{\cos(\theta_1/2)\cos(\theta_2/2)}{\sin(\theta_1/2)\sin(\theta_2/2)}.
\]
We conclude that $\cosh(C_2C_3)=\cosh(BC_2)\cosh(BC_3)$ proving that the triangle $(B, C_2, C_3)$ has a right angle at $B$. This implies that $[\rho_1]$ has angle coordinate $\gamma=\pi/2$ or $\gamma=3\pi/2$. Furthermore, $D$ must belong, for symmetry reasons, to the common perpendicular bisector of the geodesic segments $C_1C_2$ and $C_3C_4$, which goes through $B$. Since the rotation of angle $\pi$ around $D$ maps $C_1$ to a point on the geodesic line through $C_1$ and $C_3$, $D$ must also belong to that geodesic line. This means that $B=D$ which is a contradiction since we're assuming $\theta_1\neq \theta_2$. So, we conclude that $\RepDT{\alpha}$ contains a unique finite mapping class group orbit of length 2 (made of north and south poles) when $\alpha=(\theta_1,\theta_1,\theta_2,\theta_2)$ and $\theta_1\neq \theta_2$.

If instead $\alpha=(\theta_1,\theta_2,\theta_2,\theta_1)$ and $\theta_1\neq \theta_2$, then we can use the previous case to conclude that the $\mathcal{D}$-triangle chain of any point in a finite orbit of length 2 inside $\RepDT{\alpha}$ is singular. It's not hard to see that the corresponding $\mathcal{B}$-triangle chains have $\beta=\pi$ and $\gamma\in\{0,\pi\}$.

Similarly, when  $\alpha=(\theta_1,\theta_2,\theta_1,\theta_2)$, then we can flip the orientation of the second triangle in any $\mathcal{B}$-chain to obtain a valid $\mathcal{B}$-triangle chain for the case $\alpha=(\theta_1,\theta_2,\theta_2,\theta_1)$. This observation highlights a bijective correspondence between finite orbits proving that $\RepDT{\alpha}$ also contains a unique finite orbit of length 2 with action-angle coordinates $\beta=\pi$ and $\gamma\in\{\pi/2,3\pi/2\}$ when $\alpha=(\theta_1,\theta_2,\theta_1,\theta_2)$ and $\theta_1\neq \theta_2$.
\end{proof}

\subsubsection{Orbits of Type~\Romannum{3}} We move on to orbits of Type~\Romannum{3}. These orbits were first discovered by Dubrovin~\cite{dubrovin} and have length 3. The peripheral angles are $\{4\pi/3, \theta, \theta, 2\theta-2\pi\}$ for some angle $\theta\in(5\pi/3,2\pi)$. We state the analogous statement to Lemma~\ref{lem:beta_i-for-orbits-of-type-II} for finite orbits of Type~\Romannum{3}.

\begin{lem}\label{lem:beta_i-for-orbits-of-type-III}
For each permutation of the peripheral angles, the action-angle coordinates of every orbit point in a finite orbit of Type~\Romannum{3} are provided in the following table. There is always a unique orbit of length 3 in the corresponding DT component. The image of a representation associated to a finite orbit of Type~\Romannum{3} is a rotation triangle group $D(2,3,\overline{\theta})$, where $\overline{\theta}=(1-\theta/2\pi)^{-1}$.

\begin{table}[ht]
    \centering
    \begin{tblr}{c|c|c}
        $\alpha$ & $\beta$ & $\gamma$ \\
        \hline\hline
        $(4\pi/3, 2\theta-2\pi, \theta, \theta)$ & \SetCell[r=2]{c} $4\pi/3$ & $\{0 ,\frac{2\pi}{3}, \frac{4\pi}{3}\}$ \\
        \hline
		$(2\theta-2\pi, 4\pi/3, \theta, \theta)$ & & $\{\frac{\pi}{3}, \pi ,\frac{5\pi}{3}\}$ \\
		\hline
		$(\theta, \theta, 4\pi/3, 2\theta-2\pi)$ & \SetCell[r=2]{c} $2\pi/3$ & $\{\frac{\pi}{3}, \pi ,\frac{5\pi}{3}\}$ \\
        \hline
		$(\theta, \theta, 2\theta-2\pi, 4\pi/3)$ & & $\{0 ,\frac{2\pi}{3}, \frac{4\pi}{3}\}$ \\
		\hline
		\SetCell[r=2]{c} $(4\pi/3, \theta, 2\theta-2\pi, \theta), (\theta, 4\pi/3, \theta, 2\theta-2\pi)$ & $\pi$ & $\{0, \pi\}$ \\
        \hline
         & $3\theta-4\pi$ & north pole \\
        \hline
        \SetCell[r=2]{c} $(2\theta-2\pi,\theta, 4\pi/3,\theta), (\theta, 2\theta-2\pi, \theta, 4\pi/3)$ & $\pi$ & $\{0, \pi\}$ \\
        \hline
         & $6\pi-3\theta$ & south pole \\
        \hline
		\SetCell[r=2]{c}$(4\pi/3, \theta, \theta, 2\theta-2\pi), (\theta, 4\pi/3, 2\theta-2\pi, \theta)$ & $\pi$ & $\{\frac{\pi}{2}, \frac{3\pi}{2}\}$ \\
		\hline
         & $3\theta-4\pi$ & north pole
         \\
         \hline
         \SetCell[r=2]{c} $(\theta,2\theta-2\pi, 4\pi/3,\theta), (2\theta-2\pi,\theta,\theta, 4\pi/3)$ & $\pi$ & $\{\frac{\pi}{2}, \frac{3\pi}{2}\}$ \\
		\hline
         & $6\pi-3\theta$ & south pole
         \\
    \end{tblr}
\end{table}
\end{lem}
\begin{proof}
We'll only do the proof in the case where $\alpha=(4\pi/3, 2\theta-2\pi, \theta, \theta)$; the other cases can be treated with similar arguments. 

First, we prove that the north and south poles do not belong to a finite orbit of length 3 when $\alpha=(4\pi/3, 2\theta-2\pi, \theta, \theta)$. The north pole is the point of $\RepDT{\alpha}$ whose $\mathcal{B}$-triangle chain is singular and consists of a single triangle with vertices $(C_1,C_2,C_3=C_4)$. It has interior angles $(\pi/3, 2\pi-\theta, 2\pi-\theta)$. The angle bisector at vertex $C_3=C_4$ cuts the geodesic segment $C_1C_2$ at a point $D$. Similarly, as in the proof of Lemma~\ref{lem:beta_i-for-orbits-of-type-II}, we can analyse the geometric features of the triangles $(C_2, C_3, D)$ and $(D,C_4, C_1)$ and observe that the two triangles chained at $D$ are the $\mathcal{D}$-triangle chain of the north pole. The angle $\angle C_2DC_3$ is equal to $\pi-\delta/2$, where $\delta$ is the angle function of the curve $d$ evaluated at the north pole.
\begin{center}
\begin{tikzpicture}[font=\sffamily,decoration={
    markings,
    mark=at position 1 with {\arrow{>}}}]
    
\node[anchor=south west,inner sep=0] at (0,0) {\includegraphics[width=3cm]{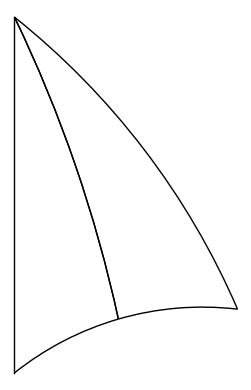}};

\begin{scope}
\fill (0.18,0.22) circle (0.07) node[left]{$C_2$};
\fill (0.18,4.45) circle (0.07) node[above]{$B=C_3=C_4$};
\fill (2.83,1) circle (0.07) node[right]{$C_1$};
\fill (1.4,.85) circle (0.07) node[above right]{$D$};
\end{scope}

\draw[mauve, thick] (.45,.445) arc (33:90:.3) node[at end, left]{\small $2\pi-\theta$};
\draw[apricot, thick] (0.18,4) arc (270:330:.4) node[at start, left]{\small $2\pi-\theta$};
\draw[thick] (1.1,.75) arc (200:100:.3) node[near start, below right]{\small $\pi-\delta/2$};
\draw[mauve, thick] (2.5,1) arc (180:112:.3) node[near end, left]{\small $\pi/3$};
\end{tikzpicture}
\end{center}
When we apply the hyperbolic law of cosines to the two triangles whose vertices are~$(D, C_2, C_3=C_4)$ and~$(C_1, C_2, C_3=C_4)$, we find
\[
\cosh(C_2C_3)=\frac{\cos(\delta/2)+\cos(\theta)\cos(\theta/2)}{\sin(\theta)\sin(\theta/2)}\quad\text{and}\quad \cosh(C_2C_3)=\frac{1/2+\cos(\theta)^2}{\sin(\theta)^2}.
\]
After simplifications, this leads to
\begin{equation}\label{eq:delta-orbits-of-type-3}
\cos(\delta/2)=\frac{1-2\cos(\theta)}{4\cos(\theta/2)}.
\end{equation}
Since the $\mathcal{D}$-triangle chain of the north pole is regular, it is not fixed by $\tau_d$. If the north pole were to belong to a finite orbit of length 3, then the order of $\tau_d$ when applied to the north pole would be $2$ or $3$. This would mean that $\delta\in\{2\pi/3, \pi, 4\pi/3\}$. We can check that these three values, once plugged into~\eqref{eq:delta-orbits-of-type-3}, would not lead to any solution for $\theta$ with $5\pi/3<\theta<2\pi$. We conclude that the north pole does not belong to any finite orbit of Type~\Romannum{3}. The same statement holds for the south pole and can be established by similar arguments.

We can assume from now on that every point in a finite orbit of length 3 has a regular $\mathcal{B}$-triangle chain. Let's fix such a point $[\rho]$. We'll prove that the action coordinate of $[\rho]$ is $\beta=4\pi/3$ and that its action coordinate $\gamma$ belongs to $\{0,2\pi/3,4\pi/3\}$. Since the $\mathcal{B}$-triangle chain of $[\rho]$ is regular, it is not fixed by $\tau_b$ and thus $\beta\in\{2\pi/3,\pi,4\pi/3\}$ (depending on whether $\tau_b$ has order 2 or 3 when applied to $[\rho]$). One of the inequalities from~\eqref{eq:inequalities-polytope-beta} reads $\beta\geq 4\pi-\alpha_1-\alpha_2$. Recalling that we have $\alpha_1=4\pi/3$ and $\alpha_2=2\theta-2\pi$, we conclude that $\beta\geq 14\pi/3-2\theta$. Using that $\theta<2\pi$, we obtain $\beta>2\pi/3$. Now, if $\beta=\pi$, then this would mean that $\tau_{b}$ permutes two of three orbit points in the mapping class group orbit of $[\rho]$ and fixes the third point. This is impossible since we already excluded both poles. So, we conclude that $\beta=4\pi/3$.

It remains to prove that the angle coordinate $\gamma$ belongs to $\{0,2\pi/3,4\pi/3\}$. Since $\beta=4\pi/3$, the Dehn twist $\tau_{b}$ applied to $[\rho]$ has order 3. The $\mathcal{B}$-triangle chain of $[\rho]$ consists of two triangles with angles $(\pi/3, 2\pi-\theta, \pi/3)$ and $(2\pi/3, \pi-\theta/2, \pi-\theta/2)$. 
\begin{center}
\resizebox{7cm}{!}{
\begin{tikzpicture}[font=\sffamily,decoration={
    markings,
    mark=at position 1 with {\arrow{>}}}]
    
\node[anchor=south west,inner sep=0] at (0,0) {\includegraphics[width=5cm]{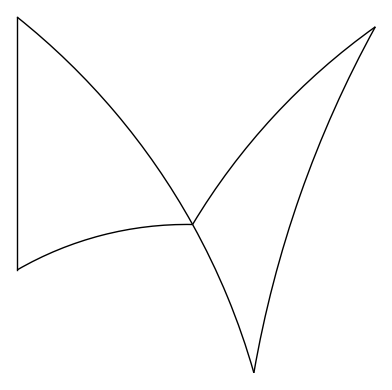}};

\begin{scope}
\fill (0.2,1.55) circle (0.07) node[below left]{$C_1$};
\fill (0.25,4.75) circle (0.07) node[above]{$C_2$};
\fill (4.75,4.6) circle (0.07) node[right]{$C_3$};
\fill (3.25,0.25) circle (0.07) node[right]{$C_4$};
\fill (2.45,2.15) circle (0.07) node[right]{$B$};
\end{scope}

\draw[mauve, thick] (.5,1.7) arc (20:90:.3) node[at end, left]{\small $\pi/3$};
\draw[mauve, thick] (0.25,4.3) arc (270:328:.4) node[at start, left]{\small $2\pi-\theta$};
\draw[apricot, thick] (2.1,2.12) arc (180:112:.3) node[near end, left]{\small $\pi/3$};
\draw[mauve, thick] (3.15,0.6) arc (120:85:.3) node[near end, right]{\small $\pi-\theta/2$};
\draw[mauve, thick] (4.55,4.15) arc (250:205:.3) node[near start, right]{\small $\pi-\theta/2$};

\draw[sky, postaction={decorate}, thick] (2.75,2.6) arc (45:112:.6) node[near start, above]{\textnormal{?}};
\end{tikzpicture}
}
\end{center}
Similar trigonometric computations as before show that $BC_1\neq BC_4$. This means that every point in the orbit of $[\rho]$ has $C_1\neq C_4$. In particular, any point in the orbit of $[\rho]$ that is fixed by $\tau_d$ must have $C_2=C_3$ by \cm{Lemma~\ref{lem:fixed-points-Dehn-twists}} and thus $\gamma=0$. Now, let's prove that $\tau_{d}$ must indeed fix an orbit point of $[\rho]$. If not, then $\tau_{d}$ applied to $[\rho]$ would have order 3. This would imply that the angle functions $\beta$ and $\delta$ agree on all three orbit points of $[\rho]$, which is impossible by Fact~\ref{fact:beta-delta-given-at-most-two-points}. Therefore, $\tau_{d}$ fixes at least one of the orbit points of $[\rho]$ and this point has $\gamma=0$. The other two orbit points are obtained by applying $\tau_{b}$ and thus have $\gamma=4\pi/3$ and $\gamma=2\pi/3$ as we explained in Section~\ref{sec:action-of-Dehn-twists}. It's not hard to see that they are permuted by $\tau_d$, effectively giving us a finite orbit of length~3.

It remains to determine the image of a representation in a finite orbit of Type~\Romannum{3}. Take the orbit point with $\beta=4\pi/3$ and $\gamma=0$. Its $\mathcal{B}$-triangle chain is made of two triangles, that can be each decomposed into two triangles with interior angles $(\pi/2, \pi/3, \pi-\theta/2)$. These smaller triangles are images of each other by reflections through their sides. This shows that the image of a corresponding representation is a rotation triangle group $D(2,3,\overline{\theta})$.
\begin{center}
\resizebox{9cm}{!}{
\begin{tikzpicture}[font=\sffamily,decoration={
    markings,
    mark=at position 1 with {\arrow{>}}}]
    
\node[anchor=south west,inner sep=0] at (0,0) {\includegraphics[width=6cm]{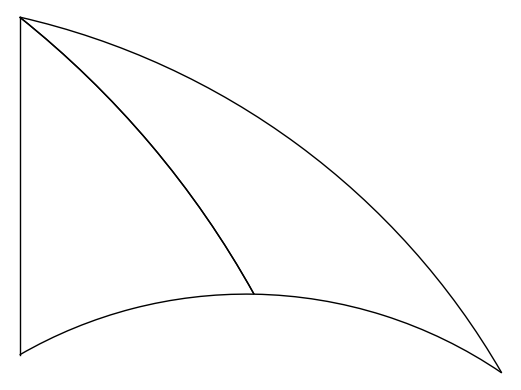}};

\begin{scope}
\fill (0.25,0.4) circle (0.07) node[below left]{$C_1$};
\fill (0.24,4.3) circle (0.07) node[above]{$C_2=C_3$};
\fill (5.75,0.25) circle (0.07) node[right]{$C_4$};
\fill (2.95,1.1) circle (0.07) node[below]{$B$};
\end{scope}

\draw[mauve, thick] (.5,.55) arc (20:90:.3) node[at end, left]{\small $\pi/3$};
\draw[mauve, thick] (0.25,3.8) arc (270:328:.4) node[at start, left]{\small $2\pi-\theta$};
\draw[apricot, thick] (2.6,1.1) arc (180:112:.3) node[near end, left]{\small $\pi/3$};
\draw[mauve, thick] (0.7,3.9) arc (315:350:.45) node[near end, above right]{\small $\pi-\theta/2$};
\draw[mauve, thick] (5.1,0.6) arc (160:110:.45) node[at end, right]{\small $\pi-\theta/2$};

\draw[dashed] (0.24,4.3) to[bend left = 8] (1.5,.95);
\draw[dashed] (2.95,1.1) to[bend left = 8] (3.83,2.53);
\end{tikzpicture}
}
\end{center}
\end{proof}

\subsubsection{Orbits of Types~\Romannum{4} and~\Romannum{4}$^\ast$} We continue with a similar analysis for orbits of Types~\Romannum{4} and~\Romannum{4}$^\ast$. Both types of finite orbits have length 4 and were discovered by Dubrovin~\cite{dubrovin}. The peripheral angles are $\{\pi, \theta, \theta, \theta\}$ for orbits of Type~\Romannum{4} and $\{\theta,\theta,\theta,3\theta-4\pi\}$ for orbits of Type~\Romannum{4}$^\ast$, with $\theta\in (5\pi/3,2\pi)$ for both. Here's the analogous statement to Lemmas~\ref{lem:beta_i-for-orbits-of-type-II} and~\ref{lem:beta_i-for-orbits-of-type-III}.

\begin{lem}\label{lem:beta_i-for-orbits-of-type-IV}
The following table provides the action-angle coordinates of every orbit point in a finite orbit of Type~\Romannum{4}, respectively of Type~\Romannum{4}$^\ast$, for every permutation of the peripheral angles. In every case, there is a unique orbit of length 4 in the corresponding DT component. The image of a representation associated to a finite orbit of Type~\Romannum{4} or of Type~\Romannum{4}$^\ast$ is a rotation triangle group $D(2,3,\overline{\theta})$, where $\overline{\theta}=(1-\theta/2\pi)^{-1}$.
\begin{table}[ht]
    \centering
    \begin{tblr}{c|c|c}
        $\alpha$ & $\beta$ & $\gamma$ \\
        \hline\hline
        \SetCell[r=2]{c} $(\theta,\theta,\theta,\pi)$ & $2\pi/3$ & $\{0 ,\frac{2\pi}{3}, \frac{4\pi}{3}\}$ \\
        \hline
        & $4\pi-2\theta$ & south pole \\
        \hline
		\SetCell[r=2]{c} $(\theta,\theta,\pi,\theta)$ & $2\pi/3$ & $\{\frac{\pi}{3} ,\pi, \frac{5\pi}{3}\}$ \\
		\hline
		 & $4\pi-2\theta$ & south pole \\
		\hline
		\SetCell[r=2]{c} $(\pi,\theta,\theta,\theta)$ & $4\pi/3$ & $\{0 ,\frac{2\pi}{3}, \frac{4\pi}{3}\}$ \\
        \hline
         & $2\theta-2\pi$ & north pole \\
        \hline
		\SetCell[r=2]{c} $(\theta,\pi,\theta,\theta)$ & $4\pi/3$ & $\{\frac{\pi}{3} ,\pi, \frac{5\pi}{3}\}$ \\
		\hline
		 & $2\theta-2\pi$ & north pole \\
    \end{tblr}
    \begin{tblr}{c|c|c}
        $\alpha$ & $\beta$ & $\gamma$ \\
        \hline\hline
		\SetCell[r=2]{c} $(\theta,\theta,\theta,3\theta-4\pi)$ & $2\pi/3$ & $\{\frac{\pi}{3} ,\pi, \frac{5\pi}{3}\}$ \\
		\hline
		 & $4\theta-6\pi$ & north pole \\
		\hline
		\SetCell[r=2]{c} $(\theta,\theta,3\theta-4\pi,\theta)$ & $2\pi/3$ & $\{0 ,\frac{2\pi}{3}, \frac{4\pi}{3}\}$ \\
        \hline
        & $4\theta-6\pi$ & north pole \\
        \hline
		\SetCell[r=2]{c} $(\theta,3\theta-4\pi,\theta,\theta)$ & $4\pi/3$ & $\{0 ,\frac{2\pi}{3}, \frac{4\pi}{3}\}$ \\
        \hline
         & $8\pi-4\theta$ & south pole \\
        \hline
		\SetCell[r=2]{c} $(3\theta-4\pi,\theta,\theta,\theta)$ & $4\pi/3$ & $\{\frac{\pi}{3} ,\pi, \frac{5\pi}{3}\}$ \\
		\hline
		 & $8\pi-4\theta$ & south pole \\
    \end{tblr}
\end{table}
\end{lem}
\begin{proof}
We'll give a detailed proof when $\alpha$ is either $(\theta,\theta,\theta,\pi)$ or $(\theta,\theta,\theta,3\theta-4\pi)$; the other cases can be treated similarly. We start with the following general observation. No matter the ordering of the peripheral angles and whether we're looking at a finite orbit of Type~\Romannum{4} or of Type~\Romannum{4}$^\ast$, there will always be an orbit point $[\rho]$ whose $\mathcal{B}$-triangle chain is regular since the orbit has length 4. The action coordinate $\beta$ of such a $[\rho]$ belongs to $\{\pi/2,2\pi/3,\pi,4\pi/3,3\pi/2\}$ depending on whether $\tau_b$ has order 2, 3, or 4 when applied to $[\rho]$. 

Let's start by showing that $\beta\notin \{\pi/2,3\pi/2\}$. Assume for the sake of contradiction that $\beta\in \{\pi/2,3\pi/2\}$. This would mean that $\tau_{b}$ has order 4 when applied to $[\rho]$. In other words, iterating $\tau_{b}$ on $[\rho]$ would give the whole mapping class group orbit of $[\rho]$, showing that $\beta$ is constant along the whole orbit of $[\rho]$. This implies that the angle function $\delta$ associated to the curve $d$ can take at most twice the same value along the orbit of $[\rho]$ by Fact~\ref{fact:beta-delta-given-at-most-two-points}. So, $\tau_{d}$ could either permute two pairs of points in the orbit of $[\rho]$, or it could fix two points and permutes the other two. Each point in a pair that is permuted by $\tau_{d}$ has $\delta=\pi$. So, in the first case, we would have $\delta=\pi$ on all four points of the orbit, which is impossible by Fact~\ref{fact:beta-delta-given-at-most-two-points}. In the second case, $\tau_{d}$ has two fixed points in the orbit of $[\rho]$. Their $\mathcal{B}$-triangle chains must have $C_{2}=C_{3}$, respectively $C_1=C_4$. They can also be mapped to each others and to $[\rho]$ by an iteration of $\tau_b$, implying that the $\mathcal{B}$-triangle chain of $[\rho]$ satisfies $BC_2=BC_3$ and $BC_1=BC_4$. However, since the angle vectors corresponding to orbits of Type~\Romannum{4} or~\Romannum{4}$^\ast$ have three identical entries, one of the two triangles in the $\mathcal{B}$-triangle chain of $[\rho]$ is isosceles at $B$, implying that actually $BC_1=BC_2=BC_3=BC_4$. This means that the two triangles in the $\mathcal{B}$-chain of $[\rho]$ are isosceles. This is impossible since the angle vector never consists of four identical entries. We've just proved that $\beta\neq \pi/2, 3\pi/2$.

We now specialize to the case $\alpha=(\theta,\theta,\theta,\pi)$. One of the inequalities~\eqref{eq:inequalities-polytope-beta} applied at the point $[\rho]$ reads $\beta\leq \theta-\pi$. Using $\theta<2\pi$, we deduce that $\beta<\pi$ and thus $\beta=2\pi/3$. This means that $\tau_b$ has order 3 when applied to $[\rho]$. Since $\delta$ cannot take the same value on all three points in the $\tau_b$-orbit of $[\rho]$, it must fix one of them. It's not hard to see through trigonometric computations that the $\mathcal{B}$-triangle chain of that point cannot have $C_1=C_4$, and therefore must have $C_2=C_3$ by \cm{Lemma~\ref{lem:fixed-points-Dehn-twists}} and thus $\gamma=0$. The other two points in the $\tau_b$-orbit of $[\rho]$ therefore have $\gamma=2\pi/3$ and $\gamma=4\pi/3$. Finally, we can apply the procedure of Section~\ref{sec:action-of-Dehn-twists} to see that the image by $\tau_d$ of the point $(\beta,\gamma)=(2\pi/3,2\pi/3)$ is the south pole. Moreover, its image by $\tau_d^2$ is the point $(\beta,\gamma)=(2\pi/3,4\pi/3)$ and its image by $\tau_d^3$ is itself. In conclusion, when $\alpha=(\theta,\theta,\theta,\pi)$, there is a unique finite orbit of Type~\Romannum{4} given by the south pole and the three points with $\beta=2\pi/3$ and $\gamma\in\{0,2\pi/3,4\pi/3\}$.

Let's finally consider the case $\alpha=(\theta,\theta,\theta,3\theta-4\pi)$. We already explained why the action coordinate $\beta$ of a point $[\rho]$ with a regular $\mathcal{B}$-triangle chain cannot be equal to $\pi/2$ or $3\pi/2$. We also want to exclude $\pi$. Assume for the sake of contradiction that $\beta=\pi$. This means that $[\rho]$ has period 2 under $\tau_{b}$. If the other two points in the mapping class group orbit of $[\rho]$ were also permuted by $\tau_{b}$, then all four points would have $\beta=\pi$ and we would obtain a contradiction to Fact~\ref{fact:beta-delta-given-at-most-two-points} similarly as in the case $\beta\in\{\pi/2,3\pi/2\}$. Now, let's assume that the other two points are fixed by $\tau_{b}$ and thus coincide with the two poles of $\RepDT{\alpha}$. The $\mathcal{B}$-triangle chain of the north pole consists of a single triangle with vertices $(C_1,C_2,C_3=C_4)$; and respective interior angles $(\pi-\theta/2,\pi-\theta/2, 4\pi-2\theta)$. Let $D$ be the point on the geodesic segment $C_1C_2$ such that $\angle DC_3C_2=\pi-\theta/2$ and $\angle DC_4C_1=3\pi-3\theta/2$. As in the proof of Lemma~\ref{lem:beta_i-for-orbits-of-type-II}, we observe that the $\mathcal{D}$-triangle chain of the north pole is given by the triangles $(C_2,C_3,D)$ and $(D,C_4,C_1)$. So, $\angle C_2DC_3=\pi-\delta/2$, where $\delta$ is the value of angle function of the curve $d$ applied to the north pole. 
\begin{center}
\begin{tikzpicture}[font=\sffamily]
    
\node[anchor=south west,inner sep=0] at (0,0) {\includegraphics[width=6cm]{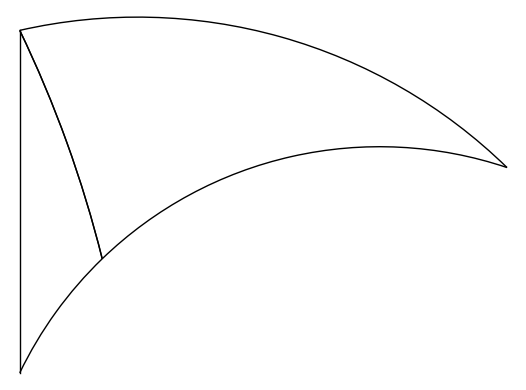}};
\node[anchor=south west,inner sep=0] at (8,0) {\includegraphics[height=5.5cm]{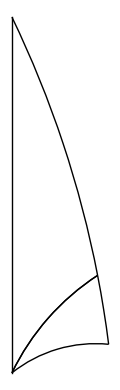}};

\begin{scope}
\fill (0.25,.2) circle (0.07) node[left]{$C_2$};
\fill (0.22,4.15) circle (0.07) node[above left]{$B=C_3=C_4$};
\fill (5.8,2.55) circle (0.07) node[right]{$C_1$};
\fill (1.15,1.5) circle (0.07) node[right]{$D$};
\end{scope}

\draw[mauve, thick] (.5,.7) arc (60:90:.45) node[at end, left]{\small $\pi-\theta/2$};
\draw[apricot, thick] (0.25,3.6) arc (270:362:.6) node[at start, left]{\small $4\pi-2\theta$};
\draw[thick] (.85,1.2) arc (230:98:.4) node[at start, below right]{\small $\pi-\delta/2$};
\draw[mauve, thick] (5.2,2.7) arc (170:118:.3) node[at start, above left]{\small $\pi-\theta/2$};

\node[anchor=south west,inner sep=0] at (.7,-.6) {\textnormal{north pole}};

\begin{scope}
\fill (8.2,.25) circle (0.07) node[below left]{$C_1=C_2=B$};
\fill (8.2,5.25) circle (0.07) node[above left]{$C_3$};
\fill (9.55,.65) circle (0.07) node[below right]{$C_4$};
\fill (9.4,1.65) circle (0.07) node[right]{$D$};
\end{scope}

\draw[mauve, thick] (8.2,4.8) arc (270:300:.3) node[at end, right]{\small $\pi-\theta/2$};
\draw[apricot, thick] (8.5,.45) arc (30:90:.4) node[at end, left]{\small $2\pi-\theta$};
\draw[thick] (9.33,1.9) arc (100:215:.3) node[at start, above right]{\small $\pi-\delta/2$};
\draw[mauve, thick] (9.3,.65) arc (180:105:.3) node[at end, right]{\small $3\pi-3\theta/2$};

\node[anchor=south west,inner sep=0] at (8.2,-.6) {\textnormal{south pole}};
\end{tikzpicture}
\end{center}
A trigonometric computation, similar to the one that lead to~\eqref{eq:delta-orbits-of-type-3} where we computed $\cosh(C_2C_3)$ in two different ways, shows that $\cos(\delta/2)=1/2$ and thus $\delta=2\pi/3$. So, the $\tau_d$-orbit of the north pole has length 3 and therefore contains the south pole (because $\tau_d$ doesn't fix the south pole). This means that the south pole also has $\delta=2\pi/3$. When we conduct the same trigonometric computations as we did for the north pole but this time for the south pole, we obtain $\cos(\delta/2)=1-\cos(\theta)$ and thus $\cos(\theta)=1/2$ because $\delta=2\pi/3$ at the south pole too. This leads to $\theta\in\{\pi/3, 5\pi/3\}$ which is impossible because we're assuming $\theta>5\pi/3$. We conclude that $\beta\neq \pi$.

There are only two possibilities left for $\beta$, namely $\beta=2\pi/3$ or $\beta=4\pi/3$. In both cases, this means that $\tau_b$ has order 3 when applied to $[\rho]$, which implies that the mapping class group orbit of $[\rho]$ contains one of the two poles of $\RepDT{\alpha}$ since it has length $4$. As above, the order of $\tau_d$ when applied to the pole in the mapping class group orbit of $[\rho]$ is necessarily 3, giving $\delta=2\pi/3$ or $\delta=4\pi/3$ at this pole. We already computed that $\delta$ is equal to $2\pi/3$ on the north pole and satisfies $\cos(\delta/2)=1-\cos(\theta)$ on the south pole. Since we're assuming $5\pi/3<\theta<2\pi$, it's impossible to have $\delta\in\{2\pi/3,4\pi/3\}$ at the south pole. We conclude that the pole in the mapping class group orbit of $[\rho]$ is the north pole. We can compute its successive images under $\tau_d$ using the procedure described in Section~\ref{sec:action-of-Dehn-twists} and we obtain the points $(\beta,\gamma)=(2\pi/3, \pi/3)$ and $(\beta,\gamma)=(2\pi/3,5\pi/3)$ before finding the north pole again after three iterations. The last orbit point is necessarily $(\beta,\gamma)=(2\pi/3,\pi)$ (which is fixed by $\tau_d$) and its $\mathcal{B}$-triangle chain is the following.
\begin{center}
\begin{tikzpicture}[font=\sffamily,decoration={
    markings,
    mark=at position 1 with {\arrow{>}}}]
    
\node[anchor=south west,inner sep=0] at (0,0) {\includegraphics[width=7cm]{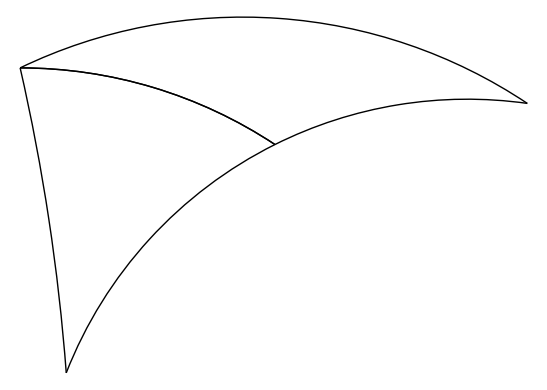}};

\begin{scope}
\fill (.85,.25) circle (0.07) node[left]{$C_3$};
\fill (0.25,4.1) circle (0.07) node[above left]{$C_4=C_1$};
\fill (6.75,3.7) circle (0.07) node[right]{$C_2$};
\fill (3.5,3.15) circle (0.07) node[below]{$B$};
\end{scope}

\draw[mauve, thick] (.8,.9) arc (95:55:.45) node[at start, left]{\small $\pi-\theta/2$};
\draw[mauve, thick] (0.4,3.6) arc (275:358:.5) node[at start, left]{\small $3\pi-3\theta/2$};
\draw[mauve, thick] (1,4.1) arc (-2:30:.6) node[at end, above]{\small $\pi-\theta/2$};
\draw[apricot, thick] (3,2.85) arc (220:147:.5) node[midway, left]{\small $\pi/3$};
\draw[mauve, thick] (6.3,3.7) arc (180:152:.5) node[at start, below]{\small $\pi-\theta/2$};

\draw[dashed] (3.5,3.15) to[bend left = 4] (3.55,4.8);
\draw[dashed] (0.25,4.1) to[bend left = 10] (2.5,2.5);
\draw[dashed] (0.25,4.1) to[bend left = 8] (1.8,1.8);
\draw[dashed] (1.8,1.8) to[bend right = 12] (.7,1.8);
\end{tikzpicture}
\end{center}
The two triangles decompose into two, respectively four, copies of a smaller triangle with angles $(\pi/2, \pi/3, \pi-\theta/2)$, each of them obtained from the others by reflections through their sides. This shows that the image of any representation in the conjugacy class is a rotation triangle group $D(2,3,\overline{\theta})$. 
\end{proof}

\subsubsection{Orbits of Types~8 and~33}\label{sec:finite-orbits-of-type-8-33}
We also compute the action-angle coordinates of all points in finite orbits of Types~8 and~33.
The peripheral angles for a Type~8 orbit are $\{10\pi/7, 12\pi/7, 12\pi/7, 12\pi/7\}$ and for an orbit of Type~33 they are $\{4\pi/3, 12\pi/7, 12\pi/7, 12\pi/7\}$. The Type~8 orbit corresponds to the "Klein solution" from~\cite{boalch-3} and has length 7. The Type~33 orbit originates from~\cite{kitaev-2} and has length 18. 

\begin{lem}\label{lem:beta_i-for-obits-of-type-33}
The action-angle coordinates of orbit points in a finite orbit of Type~33 are provided by the following table, according to the ordering of the peripheral angles. Sometimes, the action coordinate is given by a number $\gamma_0=0.58697\ldots$ which is defined by the formula
\[
\tan(\gamma_0)=\frac{\sqrt{3}(1+\cos(4\pi/7))\tan(\pi/7)}{\sin(4\pi/7)}.
\]
In every case, the finite orbit is unique in the corresponding DT component and has length 18. The image of a representation associated to a finite orbit of Type~33 is the triangle group $D(2,3,7)$, and is thus discrete.
\begin{table}[ht]
    \centering
    \begin{tblr}{c|c|c}
        $\alpha$ & $\beta$ & $\gamma$ \\
        \hline\hline
        \SetCell[r=4]{c} $\left(\frac{4\pi}{3},\frac{12\pi}{7},\frac{12\pi}{7},\frac{12\pi}{7}\right)$ & $\pi$ & $\{\frac{\pi}{4},\frac{5\pi}{4}\}\cup \{\frac{3\pi}{4},\frac{7\pi}{4}\}$\\
        \hline
        & $\frac{8\pi}{7}$ & $\{\frac{2k\pi}{7}:k=0,\ldots,6\}$\\
        \hline
        & $\frac{4\pi}{3}$ & $\{\frac{\pi}{3}\pm\gamma_0, \pi\pm\gamma_0, \frac{5\pi}{3}\pm\gamma_0\}$\\
        \hline
        & $\frac{10\pi}{7}$ & north pole\\
        \hline
        \SetCell[r=4]{c} $\left(\frac{12\pi}{7},\frac{4\pi}{3},\frac{12\pi}{7},\frac{12\pi}{7}\right)$ & $\pi$ & $\{\frac{\pi}{4},\frac{5\pi}{4}\}\cup \{\frac{3\pi}{4},\frac{7\pi}{4}\}$\\
        \hline
        & $\frac{8\pi}{7}$ & $\{\frac{\pi}{7}+\frac{2k\pi}{7}:k=0,\ldots,6\}$\\
        \hline
        & $\frac{4\pi}{3}$ & $\{\pm\gamma_0,\pm\gamma_0 +\frac{2\pi}{3},\pm\gamma_0+\frac{4\pi}{3}\}$\\
        \hline
        & $\frac{10\pi}{7}$ & north pole\\
        \hline
        \SetCell[r=4]{c} $\left(\frac{12\pi}{7},\frac{12\pi}{7},\frac{4\pi}{3},\frac{12\pi}{7}\right)$ & $\pi$ & $\{\frac{\pi}{4},\frac{5\pi}{4}\}\cup \{\frac{3\pi}{4},\frac{7\pi}{4}\}$\\
        \hline
        & $\frac{6\pi}{7}$ & $\{\frac{\pi}{7}+\frac{2k\pi}{7}:k=0,\ldots,6\}$\\
        \hline
        & $\frac{2\pi}{3}$ & $\{\pm\gamma_0,\pm\gamma_0+\frac{2\pi}{3},\pm\gamma_0+\frac{4\pi}{3}\}$\\
        \hline
        & $\frac{4\pi}{7}$ & south pole\\
        \hline
        \SetCell[r=4]{c} $\left(\frac{12\pi}{7},\frac{12\pi}{7},\frac{12\pi}{7},\frac{4\pi}{3}\right)$ & $\pi$ & $\{\frac{\pi}{4},\frac{5\pi}{4}\}\cup \{\frac{3\pi}{4},\frac{7\pi}{4}\}$\\
        \hline
        & $\frac{6\pi}{7}$ & $\{\frac{2k\pi}{7}:k=0,\ldots,6\}$\\
        \hline
        & $\frac{2\pi}{3}$ & $\{\frac{\pi}{3}\pm\gamma_0, \pi\pm\gamma_0, \frac{5\pi}{3}\pm\gamma_0\}$\\
        \hline
        & $\frac{4\pi}{7}$ & south pole\\
    \end{tblr}
\end{table}
\end{lem}
\begin{proof}
Unlike in the proofs of Lemmas~\ref{lem:beta_i-for-orbits-of-type-II}--\ref{lem:beta_i-for-orbits-of-type-IV}, we won't reprove the uniqueness of the finite orbit in the corresponding DT component, but simply rely on~\cite{LT} for this. In other words, we'll simply show the existence of a finite orbit whose action-angle coordinates are the ones provided in the table. We only consider the case where $\alpha=(4\pi/3, 12\pi/7, 12\pi/7, 12\pi/7)$; the other cases can be treated similarly.

It turns out that the north pole is part of the desired finite orbit. We'll compute its orbit and observe that it has length 18. The north pole is fixed by the Dehn twist $\tau_b$. Let's compute its image by the Dehn twist $\tau_d$ following the routine of Section~\ref{sec:action-of-Dehn-twists}. For simplicity, we'll write $[\rho_{NP}]$ for the north pole, and $\tau_d.[\rho_{NP}]$ for its image by $\tau_d$. The $\mathcal{B}$-triangle chain of $[\rho_{NP}]$ consists of a single triangle with vertices $(C_1,C_2,C_3=C_4)$ and respective angles $(\pi/3, \pi/7, 2\pi/7)$. This triangle fits well on a $(2,3,7)$-tessellation of the hyperbolic plane, similarly to the triangle chain of Example~\ref{ex:discrete-DT}.
\begin{center}
\begin{tikzpicture}[font=\sffamily,decoration={markings, mark=at position 1 with {\arrow{>}}}]
\node[anchor=south west,inner sep=0] at (0,0) {\includegraphics[width=8cm]{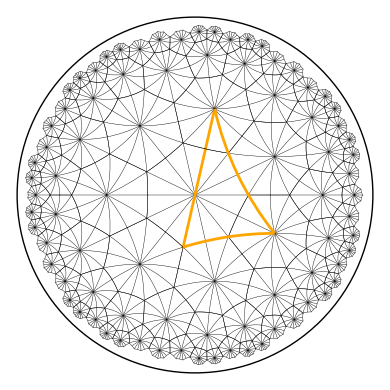}};

\fill (3.75,2.95) circle (0.07) node[left]{$C_1$};
\fill (4.4,5.8) circle (0.07) node[left]{$C_2$};
\fill (5.6,3.25) circle (0.07) node[below]{$C_3=C_4$};
\end{tikzpicture}
\end{center}
The angle bisector at the vertex $C_3=C_4$ cuts the side $C_1C_2$ at $D$. The triangle $(C_2,C_3,D)$ and $(D,C_4,C_1)$ form the $\mathcal{D}$-triangle chain of $[\rho_{NP}]$. We conclude that $\delta=6\pi/7$ at $[\rho_{NP}]$.
\begin{center}
\begin{tikzpicture}[font=\sffamily,decoration={markings, mark=at position 1 with {\arrow{>}}}]
\node[anchor=south west,inner sep=0] at (0,0) {\includegraphics[width=8cm]{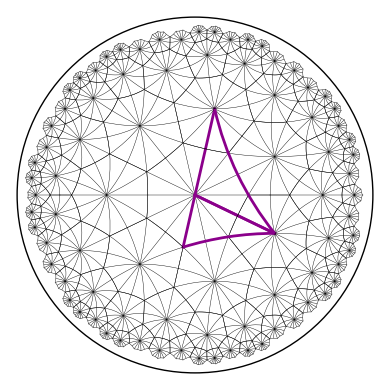}};

\fill (3.75,2.95) circle (0.07) node[left]{$C_1$};
\fill (4.4,5.8) circle (0.07) node[left]{$C_2$};
\fill (5.6,3.25) circle (0.07) node[below]{$C_3=C_4$};
\fill (4,4) circle (0.07) node[left]{$D$};
\end{tikzpicture}
\end{center}
If we rotate the triangle $(D,C_4,C_1)$ anti-clockwise around $D$ by angle $\delta=6\pi/7$, we obtain a new triangle with vertices $(D,C_4',C_1')$. The two triangles $(C_2,C_3,D)$ and $(D,C_4',C_1')$ form the $\mathcal{D}$-triangle chain of $\tau_d.[\rho_{NP}]$.
\begin{center}
\begin{tikzpicture}[font=\sffamily,decoration={markings, mark=at position 1 with {\arrow{>}}}]
\node[anchor=south west,inner sep=0] at (0,0) {\includegraphics[width=8cm]{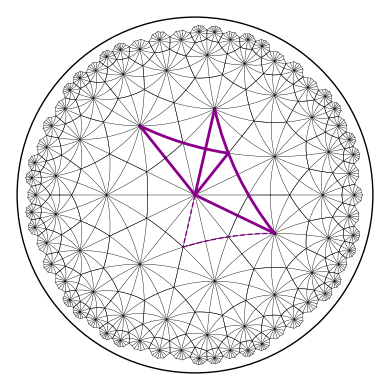}};

\fill (2.9,5.4) circle (0.07) node[left]{$C_4'$};
\fill (4.4,5.8) circle (0.07) node[left]{$C_2$};
\fill (5.6,3.25) circle (0.07) node[below]{$C_3$};
\fill (4.7,4.85) circle (0.07) node[right]{$C_1'$};
\fill (4,4) circle (0.07) node[left]{$D$};
\end{tikzpicture}
\end{center}
We let $B'$ be the unique point such that the triangle $(C_1',C_2,B')$ is clockwise oriented and has angles $\pi/3$ at $C_1'$ and $\pi/7$ at $C_2$. We observe that the angle at $B'$ is a right angle and that the triangle $(B,C_3,C_4')$ has angles $\pi/2$, $\pi/7$, and $\pi/7$. The two triangles $(C_1',C_2,B')$ and $(B',C_3,C_4')$ form the $\mathcal{B}$-triangle chain of $\tau_d.[\rho_{NP}]$.
\begin{center}
\begin{tikzpicture}[font=\sffamily,decoration={markings, mark=at position 1 with {\arrow{>}}}]
\node[anchor=south west,inner sep=0] at (0,0) {\includegraphics[width=8cm]{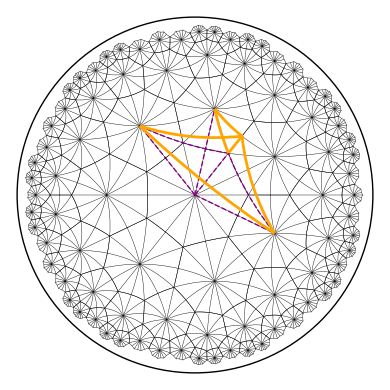}};

\fill (2.9,5.4) circle (0.07) node[left]{$C_4'$};
\fill (4.4,5.8) circle (0.07) node[left]{$C_2$};
\fill (5.6,3.25) circle (0.07) node[below]{$C_3$};
\fill (4.7,4.85) circle (0.07) node[right]{\small{$C_1'$}};
\fill (4.95,5.2) circle (0.07) node[right]{$B'$};
\end{tikzpicture}
\end{center}
We just proved that the action coordinate of $\tau_d.[\rho_{NP}]$ is $\beta=\pi$ and its angle coordinate $\gamma$ is given by the angle between the rays $B'C_3$ and $B'C_2$. It's not hard to see that $\gamma=5\pi/4$ by running some trigonometric computations. In conclusion, the Dehn twist $\tau_d$ maps the north pole to the points with action-angle coordinates $(\beta,\gamma)=(\pi,5\pi/4)$. The other orbit points can obtained in similar fashion. Since all the above triangle chains fit well on the $(2,3,7)$-tessellation \cm{(in the sense of Example~\ref{ex:discrete-DT})}, it's immediate that image of a representation associated to a finite orbit of Type~33 is a rotation triangle group $D(2,3,7)$.
\end{proof}

\begin{lem}\label{lem:beta_i-for-obits-of-type-8}
Depending on the ordering of the peripheral angles, the action-angle coordinates of every orbit point in a finite orbit of Type~8 are provided by the following table. In every case, the finite orbit is unique in the corresponding DT component and has length 7. The image of a representation associated to a finite orbit of Type~8 is the triangle group $D(2,3,7)$. This means that the orbit is made of discrete representations.
\begin{table}[ht]
    \centering
    \begin{tblr}{c|c|c}
        $\alpha$ & $\beta$ & $\gamma$ \\
        \hline\hline
        \SetCell[r=2]{c} $\left(\frac{10\pi}{7},\frac{12\pi}{7},\frac{12\pi}{7},\frac{12\pi}{7}\right)$ & $\pi$ & $\{\frac{\pi}{4},\frac{5\pi}{4}\}\cup \{\frac{3\pi}{4},\frac{7\pi}{4}\}$\\
        \hline
        & $\frac{4\pi}{3}$ & $\{0,\frac{2\pi}{3},\frac{4\pi}{3}\}$\\
        \hline
        \SetCell[r=2]{c} $\left(\frac{12\pi}{7},\frac{10\pi}{7},\frac{12\pi}{7},\frac{12\pi}{7}\right)$ & $\pi$ & $\{\frac{\pi}{4},\frac{5\pi}{4}\}\cup \{\frac{3\pi}{4},\frac{7\pi}{4}\}$\\
        \hline
        & $\frac{4\pi}{3}$ & $\{\frac{\pi}{3},\pi,\frac{5\pi}{3}\}$\\
        \hline
        \SetCell[r=2]{c} $\left(\frac{12\pi}{7},\frac{12\pi}{7},\frac{10\pi}{7},\frac{12\pi}{7}\right)$ & $\pi$ & $\{\frac{\pi}{4},\frac{5\pi}{4}\}\cup \{\frac{3\pi}{4},\frac{7\pi}{4}\}$\\
        \hline
        & $\frac{2\pi}{3}$ & $\{\frac{\pi}{3},\pi,\frac{5\pi}{3}\}$\\
        \hline
        \SetCell[r=2]{c} $\left(\frac{12\pi}{7},\frac{12\pi}{7},\frac{12\pi}{7},\frac{10\pi}{7}\right)$ & $\pi$ & $\{\frac{\pi}{4},\frac{5\pi}{4}\}\cup \{\frac{3\pi}{4},\frac{7\pi}{4}\}$\\
        \hline
        & $\frac{2\pi}{3}$ & $\{0,\frac{2\pi}{3},\frac{4\pi}{3}\}$\\
    \end{tblr}
\end{table}
\end{lem}
\begin{proof}
As in the proof of Lemma~\ref{lem:beta_i-for-obits-of-type-33}, we can use our methods to show the existence of a finite orbit of length 7 and rely on~\cite{LT} for the uniqueness statement. We omit the details as they are similar to the proof of Lemma~\ref{lem:beta_i-for-obits-of-type-33}.
\end{proof}

\subsubsection{Orbits of Types~26 and~27}
It'll be useful for the upcoming arguments that we have a better understanding of finite orbits of Types~26 and~27.  Finite orbits of Type~26 have peripheral angles $\{8\pi/5, 26\pi/15, 26\pi/15, 26\pi/15\}$ and those of Type~27 have $\{6\pi/5, 28\pi/15, 28\pi/15, 28\pi/15\}$. They both have length 15 and were originally discovered by Boalch~\cite{boalch-1}. We won't compute the action-angle coordinates of every orbit point, but rather list the possible action coordinates. That'll be enough for our work.

\begin{lem}\label{lem:beta_i-for-obits-of-type-26-or-27}
The possible action coordinates for orbit points in a finite orbit of Type~26 or~27 are provided by the following tables, according to the permutation of the peripheral angles. 
\begin{table}[ht]
    \centering
    \begin{tblr}{c|c}
        $\alpha$ & $\beta$ \\
        \hline\hline
        $\left(\frac{8\pi}{5} ,\frac{26\pi}{15}, \frac{26\pi}{15}, \frac{26\pi}{15}\right)$ & \SetCell[r=2]{c} $\{\frac{4\pi}{5}, \pi, \frac{4\pi}{3}\}$ \\
        \hline
        $\left(\frac{26\pi}{15},\frac{8\pi}{5} , \frac{26\pi}{15}, \frac{26\pi}{15}\right)$ &  \\
        \hline
        $\left(\frac{26\pi}{15}, \frac{26\pi}{15},\frac{8\pi}{5} , \frac{26\pi}{15}\right)$ & \SetCell[r=2]{c} $\{\frac{2\pi}{3}, \pi, \frac{6\pi}{5}\}$ \\
        \hline
        $\left(\frac{26\pi}{15}, \frac{26\pi}{15}, \frac{26\pi}{15}, \frac{8\pi}{5} ,\right)$ &  \\
    \end{tblr}
    \begin{tblr}{c|c}
        $\alpha$ & $\beta$ \\
        \hline\hline
        $\left(\frac{6\pi}{5}, \frac{28\pi}{15}, \frac{28\pi}{15} , \frac{28\pi}{15}\right)$ & \SetCell[r=2]{c} $\{\pi, \frac{4\pi}{3}, \frac{8\pi}{5}\}$ \\
        \hline
        $\left(\frac{28\pi}{15}, \frac{6\pi}{5}, \frac{28\pi}{15} , \frac{28\pi}{15}\right)$ &  \\
        \hline
        $\left(\frac{28\pi}{15}, \frac{28\pi}{15}, \frac{6\pi}{5}, \frac{28\pi}{15}\right)$ & \SetCell[r=2]{c} $\{\frac{2\pi}{5}, \frac{2\pi}{3}, \pi\}$ \\
        \hline
        $\left(\frac{28\pi}{15}, \frac{28\pi}{15}, \frac{28\pi}{15}, \frac{6\pi}{5}\right)$ &  \\
    \end{tblr}
\end{table}
\end{lem}
\begin{proof}
Our approach is similar to that of Lemmas~\ref{lem:beta_i-for-obits-of-type-33} and~\ref{lem:beta_i-for-obits-of-type-8}. We know from~\cite{LT} that the finite orbits are unique in the corresponding relative character variety. This means that it's enough to identify one point and then compute all orbit points using the method of Section~\ref{sec:action-of-Dehn-twists}. This will give us all possible values for the action coordinates of orbit points. Again, details are omitted and can be filled in using similar arguments as in the proof of Lemma~\ref{lem:beta_i-for-obits-of-type-33}.
\end{proof}

\section{Tykhyy's conjecture for DT representations}\label{chap:tykhyy}

\subsection{Overview}\label{sec:overview-tykhyy}
We are now ready to prove that DT components don't admit any finite mapping class group orbit if the underlying $n$-punctured sphere $\Sigma$ has $7$ punctures or more. Our argument works as follows. Starting from some $[\rho]$ inside a finite orbit, we'll obtain increasingly more restrictive conditions on the action coordinates $\beta_1,\ldots,\beta_{n-3}$ of $[\rho]$ until we're able to deduce that $n\leq 6$, confirming Tykhyy's Conjecture (Conjecture~\ref{conj:tykhyy}) for DT representations (Theorem~\ref{thm:no-finite-orbit-for-n-geq-7}). We start by restricting the values of $\beta_1,\ldots,\beta_{n-3}$ to a finite set (Section~\ref{sec:restrecting-beta_i-to-a-finite-set}) which already implies $n\leq 18$ (Remark~\ref{rem:n_leq_18}). By carrying out a more precise analysis of some finite orbits from Table~\ref{tab:finite-mcg-orbits-n=4} in Appendix~\ref{app:tables}, we'll be able to bring the upper bound on $n$ down to $6$ (Section~\ref{sec:thykyy-s-conjecture}).

\cm{For the reader's convenience, we provide a table with the notation that we'll introduce or use throughout Chapter~\ref{chap:tykhyy}.
\begin{table}[ht]
    \centering
    \begin{tblr}{width=\linewidth, colspec={X[1,l]|X[2,l]}}
		Notation & Meaning\\
		\hline\hline
		$\Sigma$ & an $n$-punctured sphere\\
		\hline
		$c_1,\ldots,c_n$ & a set of geometric generators of $\pi_1\Sigma$, i.e.~
		$\pi_1\Sigma=\langle c_1,\ldots,c_n:c_1\cdots c_n=1\rangle$\\
		\hline
		$\alpha\in(0,2\pi)^n$ & a peripheral angle vector satisfying $\sum_{i=1}^n\alpha_i>2\pi(n-1)$\\
		\hline
		$\RepDT{\alpha}$ & the DT component of $\Sigma$ and $\alpha$ (isomorphic to $\CP^{n-3}$)\\
		\hline
		$b_i$ & the fundamental group elements $b_i=(c_1\cdots c_{i+1})^{-1}$ (note that $b_0=c_1^{-1}$ and $b_{n-2}=c_n$)\\
		\hline
		$\tau_{b_1},\ldots, \tau_{b_{n-3}}$ & the Dehn twists along the simple closed curves represented by $b_1,\ldots,b_{n-3}$\\
		\hline
		$\beta_i\colon\RepDT{\alpha}\rightarrow (0,2\pi)$ & the maps which associate to $[\rho]$ the rotation angle of $\rho(b_i)$, also called action coordinates (note that $\beta_0=2\pi-\alpha_1$ and $\beta_{n-2}=\alpha_n$)\\
		\hline
		$\mathcal{B}$ & the pants decompositions of $\Sigma$ given by the curves $b_1,\ldots, b_{n-3}$\\
		\hline
		$(\beta_1,\ldots,\beta_{n-3},\gamma_1,\ldots,\gamma_{n-3})$ & the action-angle coordinates on $\RepDT{\alpha}$ defined by $\mathcal{B}$ (Section~\ref{sec:DT-representations})\\
        \hline 
        $C_1,\ldots,C_{n}$ and $B_1,\ldots,B_{n-3}$ & the exterior and shared vertices of a $\mathcal{B}$-triangle chain\\
        \hline 
        $\Sigma^{(i)}$ & the sub-sphere of $\Sigma$ with the four boundary curves $(b_{i-1}^{-1},c_{i+1},c_{i+2},b_{i+1})$\\
        \hline 
        $\alpha^{(i)}$ & the angle vector $(2\pi-\beta_{i-1},\alpha_{i+1},\alpha_{i+2},\beta_{i+1})$\\
        \hline 
        $\RepDT{\alpha^{(i)}}$ & the DT component of the sphere $\Sigma^{(i)}$ and $\alpha^{(i)}$\\
\end{tblr}
\end{table}
}

\subsection{Restricting the values of action coordinates}\label{sec:restrecting-beta_i-to-a-finite-set}
In order to restrict the range of admissible values for the action coordinates of points inside a finite orbit, we start with a simple observation.
\begin{lem}\label{lem:beta_i-strictly-increasing}
For any chained pants decomposition of an $n$-punctured sphere $\Sigma$ and any compatible geometric presentation of $\pi_1\Sigma$ (Section~\ref{sec:triangle-chains}), the associated action coordinates $\beta_1,\ldots, \beta_{n-3}$ satisfy $\beta_1<\cdots <\beta_{n-3}$ everywhere on the DT component $\RepDT{\alpha}$.
\end{lem}
\begin{proof}
This statement is a weaker version of the system of inequalities~\eqref{eq:inequalities-polytope-beta} defining the polytope of admissible values for $\beta_1,\ldots,\beta_{n-3}$. Some of the polytope inequalities read $\beta_{i+1}-\beta_{i}\geq 2\pi-\alpha_{i+2}$ for $i=1,\ldots, n-4$. Since $\alpha_{i+2}<2\pi$ by assumption, we obtain $\beta_i<\beta_{i+1}$.
\end{proof}

Let $\rho\colon\pi_1\Sigma\to \psl$ be a DT representation of an $n$-punctured sphere $\Sigma$ and assume that the mapping class group orbit of $[\rho]$ in $\RepDT{\alpha}$ is finite. It's crucial for several steps of the argument that we work with a geometric presentation of $\pi_1\Sigma$ and a compatible pants decomposition $\mathcal{B}$ of $\Sigma$ such that the triangle chain associated to $[\rho]$ is regular \cm{(Definition~\ref{def:regular-singular-triangle-chains})}. Such a presentation, with a compatible pants decomposition, always exists by Lemma~\ref{lem:existence-regular-triangle-chain}. We label the generators $c_1,\ldots, c_n$ in such a way that the pants curves of $\mathcal{B}$ are given by the fundamental group elements $b_i=(c_1\cdots c_{i+1})^{-1}$ for $i=1,\ldots, n-3$. We denote the associated action coordinates by $\beta_1,\ldots, \beta_{n-3}\colon\RepDT{\alpha}\to (0,2\pi)$. We'll explain how to find increasingly stronger restrictions on the possible values that the action coordinates $\beta_1,\ldots,\beta_{n-3}$ of $[\rho]$ can take. We start by showing that the action coordinates must be rational multiples of $\pi$.

\begin{lem}\label{lem:finite-mcg-orbit-implies-beta_i-rational}
If the mapping class group orbit of $[\rho]$ is finite and if the $\mathcal{B}$-triangle chain of $[\rho]$ is regular, then the action coordinates $\beta_1,\ldots,\beta_{n-3}$ of $[\rho]$ are all rational multiples of $\pi$.
\end{lem}
\begin{proof}
Since we're assuming that the triangle chain of $[\rho]$ is regular, the Dehn twists $\tau_{b_i}$ along the curves $b_i$ act non-trivially on $[\rho]$ \cm{(Lemma~\ref{lem:fixed-points-Dehn-twists})}. More precisely, they each act by changing the angle coordinate $\gamma_i$ of $[\rho]$ to $\gamma_i-\beta_i$ \cm{(Lemma~\ref{lem:Dehn-twist-action-on-action-angle-coordinates})}. So, if $\beta_i$ was not a rational multiple of $\pi$, iterating $\tau_{b_i}$ on $[\rho]$ would produce infinitely many orbit points.
\end{proof}

In order to constrain the possible values of the action coordinates of $[\rho]$ even more, we'll apply the complete knowledge we have of finite mapping class group orbits for 4-punctured spheres from~\cite{LT}. To make this possible, we consider restrictions of $\rho$ to smaller spheres. More precisely, let $\Sigma^{(i)}$ be the 4-punctured sphere defined as the sub-surface of $\Sigma$ with boundary loops $b_{i-1}^{-1}$, $c_{i+1}$, $c_{i+2}$, and $b_{i+1}$.
\begin{center}
\vspace{2mm}
\begin{tikzpicture}[scale=1.1, decoration={
    markings,
    mark=at position 0.6 with {\arrow{>}}}]
  \draw[postaction={decorate}, black!40] (0,-.5) arc(-90:-270: .25 and .5);
  \draw[black!40] (0,.5) arc(90:-90: .25 and .5);
  \draw[postaction={decorate}] (2,.5) arc(90:270: .25 and .5) node[midway, left]{$b_{i-1}$};
  \draw (2,.5) arc(90:-90: .25 and .5);
  \draw[apricot, postaction={decorate}] (4,.5) arc(90:270: .25 and .5) node[midway, left]{$b_i$};
  \draw[lightapricot] (4,.5) arc(90:-90: .25 and .5);
  \draw[postaction={decorate}] (6,.5) arc(90:270: .25 and .5) node[midway, left]{$b_{i+1}$};
  \draw (6,.5) arc(90:-90: .25 and .5);
  \draw[postaction={decorate}, black!40] (8,.5) arc(90:270: .25 and .5);
  \draw[black!40] (8,.5) arc(90:-90: .25 and .5);
  
  \draw[black!40](.5,1) arc(180:0: .5 and .25)node[midway, above]{$c_{i}$};
  \draw[postaction={decorate, black!40}, black!40] (.5,1) arc(-180:0: .5 and .25);
  \draw (2.5,1) arc(180:0: .5 and .25)node[midway, above]{$c_{i+1}$};
  \draw[postaction={decorate}] (2.5,1) arc(-180:0: .5 and .25);
  \draw (4.5,1) arc(180:0: .5 and .25)node[midway, above]{$c_{i+2}$};
  \draw[postaction={decorate}] (4.5,1) arc(-180:0: .5 and .25);
  \draw[black!40] (6.5,1) arc(180:0: .5 and .25)node[midway, above]{$c_{i+3}$};
  \draw[postaction={decorate}, black!40] (6.5,1) arc(-180:0: .5 and .25);
   
  \draw[black!40] (0,.5) to[out=0,in=-90] (.5,1);
  \draw[black!40] (1.5,1) to[out=-90,in=180] (2,.5);
  \draw[black!40] (0,-.5) to[out=0,in=180] (2,-.5);
  
  \draw (2,.5) to[out=0,in=-90] (2.5,1);
  \draw (3.5,1) to[out=-90,in=180] (4,.5);
  \draw (2,-.5) to[out=0,in=180] (4,-.5);
  
  \draw (4,.5) to[out=0,in=-90] (4.5,1);
  \draw (5.5,1) to[out=-90,in=180] (6,.5);
  \draw (4,-.5) to[out=0,in=180] (6,-.5);
  
  \draw[black!40] (6,.5) to[out=0,in=-90] (6.5,1);
  \draw[black!40] (7.5,1) to[out=-90,in=180] (8,.5);
  \draw[black!40] (6,-.5) to[out=0,in=180] (8,-.5);
\end{tikzpicture}
\vspace{2mm}
\end{center}
We follow the convention that $b_0=c_1^{-1}$ and $b_{n-2}=c_n$. We'll write $\rho\vert_{\Sigma^{(i)}}$ for the restriction of $\rho$ to $\pi_1\Sigma^{(i)}$. It was already observed in~\cite{deroin-tholozan} that the restrictions $\rho\vert_{\Sigma^{(i)}}$ are again a DT representations. The conujugacy class of $\rho\vert_{\Sigma^{(i)}}$ lives in the $\alpha^{(i)}$-relative character variety of $\Sigma^{(i)}$ given by the angle vector $\alpha^{(i)}=(2\pi-\beta_{i-1},\alpha_{i+1},\alpha_{i+2},\beta_{i+1})$. A convenient geometric presentation of $\pi_1\Sigma^{(i)}$ is given by the fundamental group elements $b_{i-1}^{-1}$, $c_{i+1}$, $c_{i+2}$, and $b_{i+1}$. The fundamental group element $b_i=c_{i+1}^{-1}b_{i-1}$ defines a compatible pants decomposition of $\Sigma^{(i)}$. The corresponding triangle chain of $[\rho\vert_{\Sigma^{(i)}}]$ is regular because it's simply the sub-chain of the $\mathcal{B}$-triangle chain of $[\rho]$ given by the $(i-1)$th and the $i$th triangles. Since we're assuming that $[\rho]$ has a finite orbit for the action of $\PMod (\Sigma)$ on $\RepDT{\alpha}$, the orbit of $[\rho\vert_{\Sigma^{(i)}}]$ in $\RepDT{\alpha^{(i)}}$ for the action of $\PMod(\Sigma^{(i)})$ is also finite. The group $\PMod(\Sigma^{(i)})$ can be realized as the subgroup of $\PMod (\Sigma)$ consisting of isotopy classes of homeomorphisms that are supported inside $\Sigma^{(i)}\subset \Sigma$. Since $\Sigma^{(i)}$ was chosen to be a 4-punctured sphere and $[\rho\vert_{\Sigma^{(i)}}]$ belongs to a finite orbit inside the DT component $\RepDT{\alpha^{(i)}}$ of $\Sigma^{(i)}$, the angle vector $\alpha^{(i)}$ must appear as one of the entries in Table~\ref{tab:finite-mcg-orbits-n=4}. We'll say that the restriction $[\rho\vert_{\Sigma^{(i)}}]$ is \emph{of Type~X}, with $X\in\{\text{\Romannum{2}}, \text{\Romannum{3}}, \text{\Romannum{4}}, \text{\Romannum{4}}^\ast, 1,\ldots, 45\}$, if $\alpha^{(i)}$ appears in the line corresponding to the solution labelled $X$ in Table~\ref{tab:finite-mcg-orbits-n=4}.

In order to use this inductive scheme efficiently, it's best to first restrict the possible values for the action coordinates of $[\rho]$ further. By studying the restriction $[\rho\vert_{\Sigma^{(i)}}]$, we'll obtain constrains on the possible value of $\beta_i$. We'll do so by considering a helpful tool which is a variation of a classical invariant called the \emph{trace field} of a representation, described for instance in~\cite[Chapter 3]{reid-book}. Our variation is defined as follows. 

\begin{defn}\label{def:non-peripheral-trace-field}
To a representation $\phi\colon\pi_1\Sigma\to \psl$, we associate the field 
\[
\Q(\phi)=\Q\big(\pm\trace(\phi(\gamma)):\gamma\in\pi_1\Sigma \text{ is a non-peripheral simple loop}\big).
\]
Here, \emph{non-peripheral} fundamental group elements are the ones that are not freely homotopic to a puncture of $\Sigma$. We call it the \emph{non-peripheral trace field} of $\phi$. Note that since $\phi$ takes values in $\psl$, the trace of an element in the image of $\phi$ is only defined up to a sign, which is not a source of concerns for the definition of $\Q(\phi)$. Clearly, the non-peripheral trace field is a conjugacy invariant and it's constant along mapping class group orbits in the character variety.
\end{defn}

The non-peripheral trace field differs from the classical trace field of a representation which is defined as the extension of $\Q$ by the traces of every element in $\pi_1\Sigma$, including non-simple loops and loops around punctures. Our next task is to find a set of generators for non-peripheral trace fields when $\Sigma$ is a 4-punctured sphere.

\begin{lem}\label{lem:generators-non-peripheral-trace-fields}
Let $\phi\colon\pi_1\Sigma\to \psl$ be a representation of a 4-punctured sphere $\Sigma$. Denote by $(A,B,C,D)$ the Fricke coefficients \cm{(Definition~\ref{defn:Fricke-coefficients})} parametrizing the relative character variety to which $[\phi]$ belongs, in the sense of Section~\ref{sec:fricke-relation}. The non-peripheral trace field $\Q(\phi)$ is generated over $\Q$ by the traces of $\phi$ along three non-peripheral simple closed curves and by the numbers $A,B,C$. In other words, if $c_1,c_2,c_3,c_4$ are geometric generators of $\pi_1\Sigma$, then
\[
\Q(\phi)=\Q(\pm\trace(\phi(c_1c_2)), \pm\trace(\phi(c_2c_3)), \pm\trace(\phi(c_1c_3)), A, B, C).
\]
\end{lem}
\begin{proof}
Up to inverses and conjugates, any non-peripheral simple loop $\gamma\in\pi_1\Sigma$ can be mapped to either $c_1c_2$, $c_2c_3$, or $c_1c_3$ by an element of $\PMod(\Sigma)$. As explained for instance in~\cite[Section 2.2]{cantat-loray}, the action of $\PMod(\Sigma)$ can be described explicitly and it's possible to write down a polynomial expression for $\pm\trace(\phi(\gamma))$ in the variables $\pm\trace(\phi(c_1c_2)), \pm\trace(\phi(c_1c_3)), \pm\trace(\phi(c_2c_3)), A, B, C$.

Conversely, it's possible to write $A,B,C$ as rational functions in the traces of the three simple curves $c_1c_2$, $c_2c_3$, and $c_1c_3$, and their images under the Dehn twists along the same curves.
\end{proof}

\begin{rem}
Table~\ref{tab:finite-mcg-orbits-n=4} from Appendix~\ref{app:tables} provides many examples of representations for which the non-peripheral trace field is a proper sub-field of the classical trace field.    
\end{rem}

It turns out that the non-peripheral trace fields of exceptional finite orbits in the case of 4-punctured spheres are fairly restricted. This leads to the following finer version of Lemma~\ref{lem:finite-mcg-orbit-implies-beta_i-rational}.

\begin{lem}\label{lem:finite-orbit-implies-beta_i-in-some-finite-list}
If the mapping class group orbit of $[\rho]$ is finite and if the $\mathcal{B}$-triangle chain of $[\rho]$ is regular, then the action coordinates $\beta_1,\ldots,\beta_{n-3}$ of $[\rho]$ belong to
\[
\left\{\frac{2\pi}{3},\pi,\frac{4\pi}{3}\right\}\cup\left\{\frac{\pi}{2},\frac{3\pi}{2}\right\}\cup\left\{\frac{2k\pi}{5}:1\leq k\leq 4\right\}\cup\left\{\frac{2k\pi}{7}:1\leq k\leq 6\right\}.
\]
\end{lem}
\begin{proof}
We start by computing the non-peripheral trace fields $K$ of each exceptional finite orbit. We do so by finding the values of the six generators of Lemma~\ref{lem:generators-non-peripheral-trace-fields}. This data can be found in~\cite[Table~4]{LT}. The Fricke coefficients are explicitly given in the third column of~\cite[Table~4]{LT} and the traces of the three simple closed curves can be computed by taking the cosines of the three angles provided in the fourth column of~\cite[Table~4]{LT}. The resulting field for each exceptional orbit can be found in the third column of Table~\ref{tab:finite-mcg-orbits-n=4} from Appendix~\ref{app:tables}. We observe that only the following fields appear:
\[
K\in\left\{\Q,\Q\big(\sqrt{2}\big),\Q\big(\sqrt{5}\big),\Q\big(\cos(\pi/7)\big)\right\}.
\]

Lemma~\ref{lem:finite-mcg-orbit-implies-beta_i-rational} tells us that the action coordinates $\beta_1,\ldots,\beta_{n-3}$ are rational mutiples of $\pi$. Now, if the restriction $[\rho\vert_{\Sigma^{(i)}}]$ corresponds to a finite orbit with non-peripheral trace field $K$, then $2\cos(\beta_i/2)\in K$. Our next task therefore consists in finding all the angles $\beta_i\in (0,2\pi)\cap \pi\Q$ such that $2\cos(\beta_i/2)\in K$ for every $K$ above. In the case $K=\Q$, this problem is solved by Niven's Theorem and the only possible values for $\beta_i$ are $2\pi/3$, $\pi$, $4\pi/3$. In the other cases, we use the algorithm of~\cite{niven-generalization} that we described in Section~\ref{sec:number-fields}. According to Table~\ref{tab:rational-angles}, when $K=\Q(\cos(\pi/7))$, the possible new values for $\beta_i$ are $2k\pi/7$ with $k=1,\ldots, 6$. Similarly, we obtain the new values $\pi/2$, $3\pi/2$ for $K=\Q(\sqrt{2})$ and $2k\pi/5$ with $k=1,\ldots, 4$ for $K=\Q(\sqrt{5})$. In conclusion, if the restriction $[\rho\vert_{\Sigma^{(i)}}]$ is of Type~1--45, then $\beta_i$ takes one of the anticipated values.

It remains to consider the case where $[\rho\vert_{\Sigma^{(i)}}]$ is of Type~\Romannum{2}--\Romannum{4} or~\Romannum{4}$^\ast$. We'll use our complete knowledge of the action-angle coordinates of orbit points in those cases (Section~\ref{sec:some-finite-orbits}). Since we're assuming that the triangle chain of $[\rho]$ is regular, the triangle chain of $[\rho\vert_{\Sigma^{(i)}}]$ is also regular. Lemmas~\ref{lem:beta_i-for-orbits-of-type-II}--\ref{lem:beta_i-for-orbits-of-type-IV} thus imply that $\beta_i\in\{2\pi/3,\pi, 4\pi/3\}$. We conclude that the restrictions $[\rho\vert_{\Sigma^{(i)}}]$ that are not of exceptional type don't contribute any new possible values of $\beta_i$.
\end{proof}

\begin{rem}\label{rem:n_leq_18} 
Lemma~\ref{lem:finite-orbit-implies-beta_i-in-some-finite-list} already implies that there are no finite mapping class group orbits inside DT components of a sphere with $19$ punctures or more. This is because there are $15$ possible values for each action coordinate $\beta_i$ by Lemma~\ref{lem:finite-orbit-implies-beta_i-in-some-finite-list} and they must form an increasing sequence by Lemma~\ref{lem:beta_i-strictly-increasing}. So, the longest increasing sequence of action coordinates has length 15 showing that the triangle chain of $[\rho]$ can consist of at most 16 triangles. By studying more precisely which finite orbits from Table~\ref{tab:finite-mcg-orbits-n=4} are actually possible for the restrictions $[\rho\vert_{\Sigma^{(i)}}]$, we'll be able to bring the bound on the number of punctures from 19 down to 7 and prove Tykhyy's Conjecture for DT representations (Theorem~\ref{thm:no-finite-orbit-for-n-geq-7}).
\end{rem}

\subsection{Bringing down the bound on the number of punctures}\label{sec:thykyy-s-conjecture}
We achieved in Lemma~\ref{lem:finite-orbit-implies-beta_i-in-some-finite-list} our goal to restrict the possible values for the action coordinates to a reasonably small and explicit finite set. We are now ready to prove that regular triangle chains associated to DT representations whose conjugacy classes belong to finite mapping class group orbits consist of at most $4$ triangles. Here's our strategy.

Our first task is go through the list of angle vectors $\alpha$ in Table~\ref{tab:finite-mcg-orbits-n=4} and eliminate those for which $\alpha$ does not contain any of the angles listed in Lemma~\ref{lem:finite-orbit-implies-beta_i-in-some-finite-list}. These angle vectors cannot be the angle vector corresponding to any of the $n-3$ restrictions $[\rho\vert_{\Sigma^{(i)}}]$ if $\rho$ is a representation of an $n$-punctured sphere with $n\geq 5$. After elimination, the only remaining angle vectors are those of the following types
\[
\{\text{\Romannum{2}},\text{\Romannum{3}},\text{\Romannum{4}},\text{\Romannum{4}}^\ast,1, 8, 11, 12, 15, 26, 27, 33, 39\}.
\]
We further divide them into two sub-families $\mathcal{F}_1$ and $\mathcal{F}_2$. Those in $\mathcal{F}_1$ are the ones for which $\alpha$ contains only one value among the ones of Lemma~\ref{lem:finite-orbit-implies-beta_i-in-some-finite-list}, and the others are in $\mathcal{F}_2$. A rapid check gives
\[
\mathcal{F}_1=\{1,11,12,15,26,27,39\},\quad \mathcal{F}_2=\{\text{\Romannum{2}},\text{\Romannum{3}},\text{\Romannum{4}},\text{\Romannum{4}}^\ast, 8, 33\}.
\]
The finite orbits from $\mathcal{F}_1$ correspond to triangle chains that must necessarily be at the beginning or at the end of a longer chain (they can only be continued on one side). They tell us what are the possible values for the first two and last two action coordinates $\beta_1<\beta_2$ and $\beta_{n-4}<\beta_{n-3}$. The orbits in $\mathcal{F}_2$ can be anywhere in the chain, as they can a priori be continued at both extremities. When they lie at the start or at the end of the chain, they contribute to new possible values for $\beta_1<\beta_2$ and $\beta_{n-4}<\beta_{n-3}$. When they instead lie in the middle of the chain, they tell us what partial sequences $\beta_{i-1}<\beta_i<\beta_{i+1}$ are possible. Here's what we obtain when we go through these few steps.

\begin{lem}\label{lem:possible-partial-sequences-beta-angles}
Assume that the mapping class group orbit of $[\rho]$ is finite and that the $\mathcal{B}$-triangle chain of $[\rho]$ is regular. The action coordinates of $[\rho]$ form an increasing sequence $\beta_1<\ldots<\beta_{n-3}$ by Lemma~\ref{lem:beta_i-strictly-increasing}. The sequence is made of the following $n-3$ partial sequences
\[
(\beta_1<\beta_2), \ldots, (\beta_{i-1}<\beta_i<\beta_{i+1}), \ldots, (\beta_{n-4}<\beta_{n-3}),
\]
where $i$ ranges from $2$ to $n-4$. Each partial sequence $(\beta_{i-1}<\beta_i<\beta_{i+1})$ must be one of those listed in Table~\ref{tab:partial-sequences-beta_i-1-beta_i-beta_i+1}. The partial sequences $(\beta_1<\beta_2)$ and $(\beta_{n-4}<\beta_{n-3})$ must be among the ones listed in Table~\ref{tab:partial-sequences-beta_1-beta_2-beta_n-4-beta-n-3}.
\begin{table}[ht]
    \centering
    \begin{tblr}{c|c|c|c}
        \makecell{Lisovyy--Tykhyy's\\ numbering} & $\beta_{i-1}$ & $\beta_i$ & $\beta_{i+1}$\\
        \hline\hline
        Sol. \Romannum{2} & $2\pi- \theta_1$, $2\pi- \theta_2$ & $\pi$ & $\theta_1$, $\theta_2$\\
        \hline
        \SetCell[r=5]{c} Sol. \Romannum{3} & $\frac{2\pi}{3}, 4\pi-2\theta$ & $\pi,\frac{4\pi}{3}$ & $\theta$ \\
        \hline
        & $\frac{2\pi}{3}$ & $\pi$ & $2\theta-2\pi$ \\
        \hline
        & $2\pi-\theta$ & $\frac{2\pi}{3}, \pi$ & $ 2\theta-2\pi, \frac{4\pi}{3}$ \\
        \hline
        & $2\pi-\theta$ & $\pi$ & $\theta$ \\
        \hline
        & $4\pi-2\theta$ & $\pi$ & $\frac{4\pi}{3}$ \\
        \hline
        \SetCell[r=2]{c} Sol. \Romannum{4} & $\pi$ & $\frac{4\pi}{3}$ & $\theta$ \\
        \hline
        & $2\pi-\theta$ & $\frac{2\pi}{3}$ & $\pi$ \\
        \hline
        Sol. \Romannum{4} \& \Romannum{4}$^\ast$ & $2\pi- \theta$ & $\frac{2\pi}{3}, \frac{4\pi}{3}$ & $\theta$ \\
        \hline
        \SetCell[r=2]{c} Sol. \Romannum{4}$^\ast$ & $2\pi- \theta$ & $\frac{2\pi}{3}$ & $3\theta-4\pi$ \\
        \hline
        & $6\pi- 3\theta$ & $\frac{4\pi}{3}$ & $\theta$ \\
        \hline\hline
        Sol. 1 \& 15 & $\frac{2\pi}{5}$ & $\frac{2\pi}{3},\pi,\frac{4\pi}{3}$ & $\frac{8\pi}{5}$ \\
        \hline
        \SetCell[r=3]{c} Sol. 8 & $\frac{2\pi}{7}$ & $\frac{2\pi}{3},\pi,\frac{4\pi}{3}$ & $\frac{12\pi}{7}$ \\
        \hline
        & $\frac{2\pi}{7}$ & $\frac{2\pi}{3},\pi$ & $\frac{10\pi}{7}$ \\
        \hline
        & $\frac{4\pi}{7}$ & $\pi,\frac{4\pi}{3}$ & $\frac{12\pi}{7}$ \\
        \hline
        \SetCell[r=3]{c} Sol. 33 & $\frac{2\pi}{7}$ & $\frac{2\pi}{3},\frac{6\pi}{7},\pi$ & $\frac{4\pi}{3}$ \\
        \hline
        & $\frac{2\pi}{3}$ & $\pi, \frac{8\pi}{7},\frac{4\pi}{3}$ & $\frac{12\pi}{7}$ \\
        \hline
        & $\frac{2\pi}{7}$ & $\frac{2\pi}{3},\frac{6\pi}{7},\pi,\frac{8\pi}{7},\frac{4\pi}{3}$ & $\frac{12\pi}{7}$ \\ 
        \end{tblr}
        \caption{All the possible partial sequences $(\beta_{i-1}<\beta_i<\beta_{i+1})$ when $2\leq i\leq n-4$.}
        \label{tab:partial-sequences-beta_i-1-beta_i-beta_i+1}
\end{table}
\begin{table}[ht]
    \centering
    \begin{tblr}{c|c|c||c|c}
        \makecell{Lisovyy--Tykhyy's\\ numbering} & $\beta_{1}$ & $\beta_2$ & $\beta_{n-4}$ & $\beta_{n-3}$\\
        \hline\hline
        Sol. \Romannum{2} & $\pi$ & $\theta$ & $2\pi-\theta$ & $\pi$\\
        \hline 
        \SetCell[r=2]{c} Sol. \Romannum{3} & $\pi,\frac{4\pi}{3}$ & $\theta$ & $2\pi-\theta$ & $\frac{2\pi}{3},\pi$\\
        \hline
        & $\frac{2\pi}{3},\pi$ & $2\theta-2\pi, \frac{4\pi}{3}$ & $4\pi-2\theta, \frac{2\pi}{3}$ & $\pi,\frac{4\pi}{3}$\\
        \hline
        Sol. \Romannum{4} \& \Romannum{4}$^\ast$ & $\frac{2\pi}{3},\frac{4\pi}{3}$ & $\theta$ & $2\pi-\theta$ & $\frac{2\pi}{3},\frac{4\pi}{3}$\\
        \hline
        Sol. \Romannum{4} & $\frac{2\pi}{3}$ & $\pi$ & $\pi$ & $\frac{4\pi}{3}$\\
        \hline
        Sol. \Romannum{4}$^\ast$ & $\frac{2\pi}{3}$ & $3\theta-4\pi$ & $6\pi-3\theta$ & $\frac{4\pi}{3}$\\
        \hline\hline
        Sol. 1 \& 15 & $\frac{2\pi}{3},\pi,\frac{4\pi}{3}$ & $\frac{8\pi}{5}$ & $\frac{2\pi}{5}$ & $\frac{2\pi}{3},\pi,\frac{4\pi}{3}$\\
        \hline
        \SetCell[r=2]{c} Sol. 8 & $\frac{2\pi}{3},\pi,\frac{4\pi}{3}$ & $\frac{12\pi}{7}$ & $\frac{2\pi}{7}$ & $\frac{2\pi}{3},\pi,\frac{4\pi}{3}$\\
        \hline
        & $\frac{2\pi}{3},\pi$ & $\frac{10\pi}{7}$ & $\frac{4\pi}{7}$ & $\pi,\frac{4\pi}{3}$\\
        \hline
        \SetCell[r=2]{c} Sol. 33 & $\frac{2\pi}{3},\frac{6\pi}{7},\pi,$ & $\frac{4\pi}{3}$ & $\frac{2\pi}{3}$ & $\pi,\frac{8\pi}{7},\frac{4\pi}{3}$\\
        \hline
        & $\frac{2\pi}{3},\frac{6\pi}{7},\pi,\frac{8\pi}{7},\frac{4\pi}{3}$ & $\frac{12\pi}{7}$ & $\frac{2\pi}{7}$ & $\frac{2\pi}{3},\frac{6\pi}{7},\pi,\frac{8\pi}{7},\frac{4\pi}{3}$\\
        \hline
        Sol. 11 \& 12 & $\frac{2\pi}{5}, \frac{2\pi}{3},\frac{4\pi}{5}, \pi,\frac{6\pi}{5}, \frac{4\pi}{3}$ & $\frac{3\pi}{2}$ & $\frac{\pi}{2}$ &  $\frac{2\pi}{3}, \frac{4\pi}{5}, \pi,\frac{6\pi}{5}, \frac{4\pi}{3},\frac{8\pi}{5}$\\
        \hline 
        Sol. 26 & $\frac{2\pi}{3},\pi,\frac{6\pi}{5}$ & $\frac{8\pi}{5}$ & $\frac{2\pi}{5}$ &  $\frac{4\pi}{5}, \pi, \frac{4\pi}{3}$\\
        \hline
        Sol. 27 & $\frac{2\pi}{5}, \frac{2\pi}{3}, \pi$ & $\frac{6\pi}{5}$ & $\frac{4\pi}{5}$ &  $\pi, \frac{4\pi}{3}, \frac{8\pi}{5}$\\
        \hline
         Sol. 39 & $\frac{2\pi}{5}, \frac{2\pi}{3},\frac{4\pi}{5}, \pi,\frac{6\pi}{5}, \frac{4\pi}{3}$ & $\frac{3\pi}{2}$ & $\frac{\pi}{2}$ &  $\frac{2\pi}{3}, \frac{4\pi}{5}, \pi,\frac{6\pi}{5}, \frac{4\pi}{3},\frac{8\pi}{5}$\\      
    \end{tblr}
    \caption{All the possible partial sequences $(\beta_1<\beta_2)$ and $(\beta_{n-4}<\beta_{n-3})$.}
    \label{tab:partial-sequences-beta_1-beta_2-beta_n-4-beta-n-3}
\end{table}
\end{lem}
\begin{proof}
Let's start with exceptional orbits. For instance, the angle vector corresponding to a restriction of Type~26 is $\{8\pi/5, 26\pi/15, 26\pi/15, 26\pi/15\}$. It contains just one value from the list of Lemma~\ref{lem:finite-orbit-implies-beta_i-in-some-finite-list}, namely $8\pi/5$. This means that the triangle chain can only be at the start or at the end of a longer chain. If it's at the start, then $\beta_2=8\pi/5$, and if it's at the end then $2\pi-\beta_{n-4}=8\pi/5$ giving $\beta_{n-4}=2\pi/5$. We can infer the possible values for $\beta_1$ and $\beta_{n-3}$ from Lemma~\ref{lem:beta_i-for-obits-of-type-26-or-27}. It says that $\beta_1\in\{2\pi/3, \pi, 6\pi/5\}$ if the triangle chain is at the start of the longer chain and $\beta_{n-3}\in\{4\pi/5, \pi, 4\pi/3\}$ if it's at the end. These are the expected values listed in Table~\ref{tab:partial-sequences-beta_1-beta_2-beta_n-4-beta-n-3}. Type~27 is treated similarly using Lemma~\ref{lem:beta_i-for-obits-of-type-26-or-27}.

We move to restrictions of Type~39. Their angle vector $\{3\pi/2, 11\pi/6, 11\pi/6, 11\pi/6\}$ contains only one value from the list of Lemma~\ref{lem:finite-orbit-implies-beta_i-in-some-finite-list}, namely $3\pi/2$. As before, we get $\beta_2=3\pi/2$ if the chain is at the start of the longer chain of $\beta_{n-4}=\pi/2$ if it's at the end. The values of $\beta_1$ and $\beta_{n-3}$, must be such that $2\cos(\beta_1/2)$, respectively $2\cos(\beta_{n-3}/2)$, belongs to the non-peripheral trace field of Solution~39 which turns out to be $\Q(\sqrt{5})$. Using Table~\ref{tab:rational-angles}, we conclude that $\beta_1$ and $\beta_{n-3}$ belong to $\{2\pi/3,\pi,4\pi/3\}\cup\{2k\pi/5:1\leq k\leq 4\}$. Since $\beta_1<\beta_2$ and $\beta_{n-4}<\beta_{n-3}$ by Lemma~\ref{lem:beta_i-strictly-increasing}, the possible values for $\beta_1$ are $\{2\pi/5, 2\pi/3, 4\pi/5, \pi, 6\pi/5, 4\pi/3\}$. Similarly, the possible values for $\beta_{n-3}$ are $\{2\pi/3, 4\pi/5,\pi, 6\pi/5, 4\pi/3, 8\pi/5\}$. This leaves us with the possibilities listed in Table~\ref{tab:partial-sequences-beta_1-beta_2-beta_n-4-beta-n-3}. Restrictions of Types~11,~12, and~27 are treated similarly.

Let's now consider restrictions of Type~8. Their angle vector $\{10\pi/7, 12\pi/7,12\pi/7,12\pi/7\}$ contains four values from the list of Lemma~\ref{lem:finite-orbit-implies-beta_i-in-some-finite-list}. So, we are dealing with a shorter triangle chain that can be anywhere in the longer chain. Assume first that the shorter chain is not at the start or at the end of the longer chain. The possible values for the pair $(2\pi-\beta_{i-1},\beta_{i+1})$ are $\{(10\pi/7,12\pi/7), (12\pi/7, 10\pi/7), (12\pi/7, 12\pi/7)\}$ which implies that the pair $(\beta_{i-1},\beta_{i+1})$ is one among $\{(4\pi/7,12\pi/7), (2\pi/7, 10\pi/7), (2\pi/7, 12\pi/7)\}$. As above, $2\cos(\beta_i/2)$ belongs to the non-peripheral trace field of Solution~8 which is $\Q$. This means that $\beta_i\in\{2\pi/3,\pi,4\pi/3\}$ by Niven's Theorem. We can eliminate some combinations by using Lemma~\ref{lem:beta_i-for-obits-of-type-8}, leaving us with the possibilities listed in Table~\ref{tab:partial-sequences-beta_i-1-beta_i-beta_i+1}. Assume now that the shorter chain is at the start of the longer chain. This means that $\beta_2\in\{10\pi/7,12\pi/7\}$ and $\beta_1$ is any value among $\{2\pi/3,\pi,4\pi/3\}$ that's smaller than $\beta_2$ and satisfies the conclusions of Lemma~\ref{lem:beta_i-for-obits-of-type-8}. This gives the possibilities listed in Table~\ref{tab:partial-sequences-beta_1-beta_2-beta_n-4-beta-n-3}. The case where the shorter chain is at the end of the longer chain is treated similarly. Restrictions of Type ~33 are treated similarly using Lemma~\ref{lem:beta_i-for-obits-of-type-33} instead of Lemma~\ref{lem:beta_i-for-obits-of-type-8}; restrictions of Types~1 and~15 can also be dealt with in a similar fashion.

We can argue about orbits of Type~\Romannum{2}--\Romannum{4} and of Type~\Romannum{4}$^\ast$ in the same way. Lemmas~\ref{lem:beta_i-for-orbits-of-type-II}--\ref{lem:beta_i-for-orbits-of-type-IV} give a complete description of the possible values of $\beta_i$ given the ordering of the peripheral angles in $\alpha^{(i)}$.
\end{proof}

Now that we have established in Lemma~\ref{lem:possible-partial-sequences-beta-angles} what partial sequences of action coordinates are possible, it's time to prove that there are no finite mapping class group orbits in a DT component $\RepDT{\alpha}$ if the underlying sphere has $7$ punctures or more.

\begin{thm}\label{thm:no-finite-orbit-for-n-geq-7}
Assume that mapping class group action on $\RepDT{\alpha}$ admits a finite orbit. Then the underlying sphere $\Sigma$ has at most $6$ punctures.
\end{thm}
\begin{proof}
Let $[\rho]\in \RepDT{\alpha}$ be a point with a finite mapping class group orbit. We choose a pants decomposition $\mathcal{B}$ of $\Sigma$ and an assorted geometric presentation of $\pi_1\Sigma$ for which the triangle chain of $[\rho]$ is regular. This is always possible by Lemma~\ref{lem:existence-regular-triangle-chain}. We write $n$ for the number of punctures on $\Sigma$. Let $\beta_1,\ldots,\beta_{n-3}\in (0,2\pi)$ denote the action coordinates of $[\rho]$. By Lemma~\ref{lem:beta_i-strictly-increasing}, it holds that $\beta_1<\cdots<\beta_{n-3}$. The sequence must satisfy the conclusion of Lemma~\ref{lem:finite-orbit-implies-beta_i-in-some-finite-list}, meaning that each of the $n-3$ partial sequence of $\beta_1<\cdots<\beta_{n-3}$ must be among those listed in Tables~\ref{tab:partial-sequences-beta_i-1-beta_i-beta_i+1} and~\ref{tab:partial-sequences-beta_1-beta_2-beta_n-4-beta-n-3}. We'll prove that the longest sequence $\beta_1<\cdots<\beta_{n-3}$ we can build with that requirement has $n=6$. Assume for the sake of contradiction that such a sequence exists for some $n\geq 7$.

Let's start with the partial sequence $(\beta_1<\beta_2)$. The possible values for $\beta_1<\beta_2$ are provided in Table~\ref{tab:partial-sequences-beta_1-beta_2-beta_n-4-beta-n-3}. Remembering that $\theta >\pi $ for orbits of Type~\Romannum{2} and $\theta> 5\pi/3$ for orbits of Types~\Romannum{3}--\Romannum{4} and~\Romannum{4}$^\ast$, we observe that $\beta_2\geq \pi$, with equality if and only if we are looking at an orbit of Type~\Romannum{4} and $\beta_1=2\pi/3$. Since we're assuming that $n\geq 7$, the partial sequence $(\beta_2<\beta_3<\beta_4)$ must appear in Table~\ref{tab:partial-sequences-beta_i-1-beta_i-beta_i+1} with $i=3$. When we go over the possible values for $\beta_2$ in Table~\ref{tab:partial-sequences-beta_i-1-beta_i-beta_i+1}, we observe that $\beta_2\leq \pi$ with equality if and only if we are again looking at an orbit of Type~\Romannum{4} and $\beta_3=4\pi/3$. The next partial sequence is either $(\beta_3<\beta_4<\beta_5)$ if $n>7$, or $(\beta_3<\beta_4)$ if $n=7$. In the first case, we obtain $\beta_3\leq \pi$ in the same way as we obtained $\beta_2\leq \pi$ previously, which contradicts $\beta_3=4\pi/3$. In the second case, we can also infer that $\beta_3\leq \pi$ by a similar analysis, leading to the same contradiction. This shows that any sequence $\beta_1<\cdots<\beta_{n-3}$ that satisfies the conclusion of Lemma~\ref{lem:possible-partial-sequences-beta-angles} has $n\leq 6$, proving that $\Sigma$ has at most $6$ punctures.
\end{proof}

\section{Classifying finite mapping class group orbits}\label{sec:classification}

\subsection{Overview}
We just proved in Theorem~\ref{thm:no-finite-orbit-for-n-geq-7} that there are no finite orbits coming from DT representations of an $n$-punctured sphere whenever $n\geq 7$. It's now time to list all the possible finite mapping class group orbits when $n\leq 6$. The case $n=4$ was already studied extensively in Chapter~\ref{sec:4-punctured} and the list of finite orbits can be found in Table~\ref{tab:finite-mcg-orbits-n=4} in Appendix~\ref{app:tables}. The classification of finite orbits for $n=5$ was completed by Tykhyy~\cite{tykhyy}. We remind the reader about it in Section~\ref{sec:classification-n=5} and explain how to recover it with our methods in the case of DT representations. The new contribution of our paper is the classification of finite orbits coming from DT representations when $n=6$, which we explain in Section~\ref{sec:classification-n=6}. 

\cm{Throughout Chapter~\ref{sec:classification}, we'll use the same notation as in Chapter~\ref{chap:tykhyy}, so the reader may consult the table from Section~\ref{sec:overview-tykhyy} for help.}

\subsection{Classification of finite orbits for \texorpdfstring{$n=6$}{n=6}}\label{sec:classification-n=6}
Let's start with the classification of finite orbits for 6-punctured spheres as it was previously unknown. Diarra already proved that there is a unique finite orbit of pullback type (Definition~\ref{def:pullback-orbits}) when $n=6$~\cite{diarra-pull-back}. This finite orbit lives inside a DT component $\RepDT{\alpha}$ where all the entries of $\alpha$ are equal. As we'll see in Theorem~\ref{thm:angle-vector-alpha-with-finite-orbits-n=6} below, it turns out that there are no other DT components that carry finite orbits when $n=6$. In other words, if the entries of $\alpha$ are not all equal, then the associated DT component doesn't contain any finite orbit. Furthermore, we'll prove in Theorem~\ref{thm:existence-finite-orbit-n=6} that Diarra's orbit is the only finite orbit contained in $\RepDT{\alpha}$ when all the entries of $\alpha$ are equal.

\begin{thm}\label{thm:angle-vector-alpha-with-finite-orbits-n=6}
For a 6-punctured sphere, if $\RepDT{\alpha}$ contains a finite mapping class group orbit, then $\alpha=(\theta,\theta, \theta, \theta, \theta, \theta)$ with $\theta\in (5\pi/3,2\pi)$. 
\end{thm}
\begin{proof}
Let's pick $[\rho]\in\RepDT{\alpha}$ inside a finite mapping class group orbit. We start by choosing a pants decomposition $\mathcal{B}$ of $\Sigma$ for which the triangle chain of $[\rho]$ is regular using Lemma~\ref{lem:existence-regular-triangle-chain}, as we did in the proof of Theorem~\ref{thm:no-finite-orbit-for-n-geq-7}. Since we're working in the context of a 6-punctured sphere, $\mathcal{B}$ consists of three pants curves. The associated angle coordinates of $[\rho]$ are $\beta_1<\beta_2<\beta_3$ (recall that they form an increasing sequence by Lemma~\ref{lem:beta_i-strictly-increasing}). We count three partial sequences: $(\beta_1<\beta_2)$, $(\beta_1<\beta_2<\beta_3)$, and $(\beta_2<\beta_3)$.

We carry out the same analysis about the partial sequences $(\beta_1<\beta_2)$ and $(\beta_2<\beta_3)$ as we did in the proof of Theorem~\ref{thm:no-finite-orbit-for-n-geq-7} to conclude that $\beta_2=\pi$, along with $\beta_1=2\pi/3$ and $\beta_3=4\pi/3$. Moreover, the restrictions $[\rho\vert_{\Sigma^{(1)}}]$ and $[\rho\vert_{\Sigma^{(3)}}]$ are of Type~\Romannum{4}, with peripheral angles $\alpha^{(1)}=(\theta,\theta,\theta,\pi)$ and $\alpha^{(3)}=(\pi,\theta',\theta',\theta')$. The partial sequence $(\beta_1<\beta_2<\beta_3)$ turns out to be $(2\pi/3<\pi<4\pi/3)$ which is possible only when the restriction $[\rho\vert_{\Sigma^{(2)}}]$ is of Type~\Romannum{2} with $\alpha^{(2)}=(4\pi/3,\theta'',\theta'',4\pi/3)$. This forces $\theta=\theta''=\theta'$, proving that $\alpha$ consists of identical entries.
\end{proof}

Our next result shows that if $\alpha=(\theta,\theta,\theta,\theta,\theta,\theta)$, then $\RepDT{\alpha}$ contains a unique finite orbit. We start by fixing some notation. First, pick an auxiliary geometric presentation of $\pi_1\Sigma$. It induces an ordering of the punctures and of the entries of $\alpha$. We write $\mathcal{B}$ for the standard pants decomposition associated to the geometric presentation we just fixed. Since $n=6$, the DT component $\RepDT{\alpha}$ is isomorphic to $\CP^3$. There are four points inside $\RepDT{\alpha}$ whose $\mathcal{B}$-triangle chain is degenerate to a single point; they correspond to the pre-images of the vertices of the moment polytope described by~\eqref{eq:inequalities-polytope-beta}.

\begin{thm}\label{thm:existence-finite-orbit-n=6}
When $\alpha=(\theta,\theta,\theta,\theta,\theta,\theta)$ with $\theta\in (5\pi/3,2\pi)$, the DT component $\RepDT{\alpha}$ contains a unique finite mapping class group orbit: Diarra's pullback orbit~\cite{diarra-pull-back}. It consists of 40 orbit points and the image of any representation associated to the finite orbit is the rotation triangle group\footnote{Defined in Section~\ref{sec:triangle-groups}.} $D(2,3,\overline{\theta})$, where $\overline{\theta}=(1-\theta/2\pi)^{-1}$. The action-angle coordinates of every orbit point are provided in Table~\ref{tab:UFO-orbit} in Appendix~\ref{app:tables}.
\begin{table}[ht]
    \centering
    \begin{tblr}{c|c|c|c|c}
        & $\alpha$ & Orbit length & Image & Coordinates \\
        \hline\hline
        \makecell{jester's hat\\ orbit} & $(\theta,\theta,\theta,\theta,\theta,\theta)$ & 40 & $D(2,3,\overline{\theta})$ & Table~\ref{tab:UFO-orbit} \\
    \end{tblr}
\end{table}
\end{thm}
\begin{proof}
We start with the existence statement. Consider the ``jester's hat'' triangle chain, which is
parametrized by $(\beta_1,\beta_2,\beta_3)=(2\pi/3,\pi,4\pi/3)$ and $(\gamma_1,\gamma_2,\gamma_3)=(2\pi/3,0,2\pi/3)$. It consists of four triangles $(C_1,C_2,B_1)$, $(B_1,C_3,B_2)$, $(B_2,C_4,B_3)$, and $(B_3,C_5,C_6)$ with $C_3=C_4$ and $C_1=C_6$. The angles at the exterior vertices $C_1,\ldots, C_6$ are all equal to $\pi-\theta/2$.
\begin{center}
\begin{tikzpicture}[font=\sffamily] 
\node[anchor=south west,inner sep=0] at (0.2,0.2) {\includegraphics[width=5.6cm]{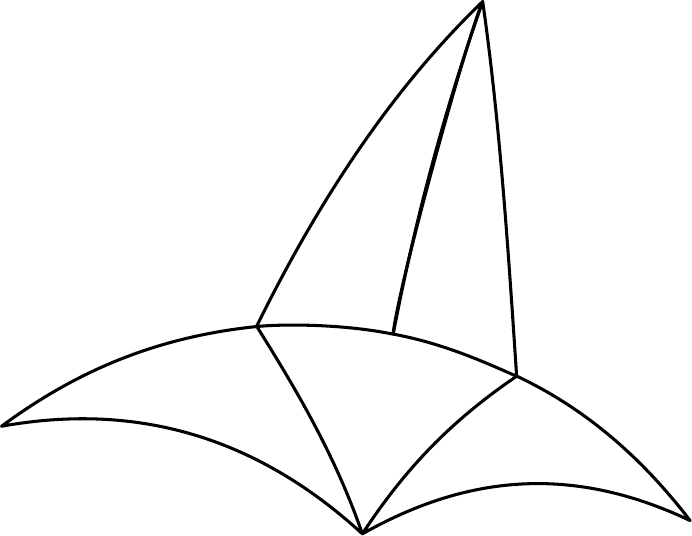}};

\begin{scope}
\fill (3.15,.22) circle (0.07) node[below]{$C_1=C_6$};
\fill (0.22,1.1) circle (0.07) node[left]{$C_2$};
\fill (2.33,1.9) circle (0.07) node[above left]{$B_1$};
\fill (4.1,4.5) circle (0.07) node[above]{$C_3=C_4$};
\fill (4.4,1.52) circle (0.07) node[above right]{$B_3$};
\fill (3.41,1.82) circle (0.07) node[below]{$B_2$};
\fill (5.8,0.35) circle (0.07) node[right]{$C_5$};
\end{scope}
\end{tikzpicture}
\end{center}
It parametrizes some point $[\rho]\in\RepDT{\alpha}$. The image of $\rho$ is the subgroup of $\psl$ generated by rotations of angle $\theta$ around each of the exterior vertices. We would like to see that the image of $\rho$ is conjugate to the rotation triangle group $D(2,3,\theta')$ which can be realized as the rotation triangle group of the triangle $(B_1,C_3,B_2)$. It's not hard to see that $C_1$ (hence $C_6$ too) is the reflection of $C_3$ through the geodesic line $B_1B_2$. Similarly, $C_2$ and $C_5$ are the reflections of $C_1$ through the geodesic lines $C_3B_1$, respectively $C_3B_3$. This shows that the image of $\rho$ is contained in $D(2,3,\theta')$. Conversely, the image of $\rho$ also contains the rotation of angle $2\pi/3$ around $B_1$ and the rotation of angle $\pi$ around $B_2$. So, it's actually equal to $D(2,3,\theta')$.

Now, we would like to prove that the mapping class group orbit of $[\rho]$---the \emph{jester's hat orbit}---is finite and has length 40. Finiteness is immediate because, by construction, $\rho$ is the pullback of the triangle group $D(2,3,\theta')$ and so the orbit of $[\rho]$ is finite by Lemma~\ref{lem:pullback-orbits-are-finite}. One could alternatively compute the action-angle coordinates of all the points in the jester's hat orbit using the methods of Section~\ref{sec:action-of-Dehn-twists} and verify that we do get a total of 40 points. To avoid long computations, we used a computer to run the algorithm described in Appendix~\ref{apx:algo-orbit}. The resulting coordinates can be found in Table~\ref{tab:UFO-orbit}.

Let's now prove the uniqueness statement. Pick a finite orbit inside $\RepDT{\alpha}$. If it contains an orbit point $[\rho]$ whose $\mathcal{B}$-triangle chain is regular, then we know from the proof of Theorem~\ref{thm:angle-vector-alpha-with-finite-orbits-n=6} that $[\rho]$ has action coordinates $(\beta_1,\beta_2,\beta_3)=(2\pi/3, \pi, 4\pi/3)$. Moreover, the combination of solutions has to be \Romannum{4}--\Romannum{2}--\Romannum{4}. When studying the restrictions of $[\rho]$ to smaller spheres more carefully, we'll be able to determine its action coordinates $(\gamma_1,\gamma_2,\gamma_3)$ as follows.
\begin{itemize}
\item First, the restriction $[\rho\vert_{\Sigma^{(1)}}]$ is of Type~\Romannum{4} and has $\alpha^{(1)}=(\theta, \theta, \theta, \pi)$. So, Lemma~\ref{lem:beta_i-for-orbits-of-type-IV} tells us that $\gamma_1\in\{0,2\pi/3,4\pi/3\}$. These three values of $\gamma_1$ can be exchanged by applying $\tau_{b_1}$ once or twice to $[\rho]$. This shows the existence of an orbit point with $\gamma_1=2\pi/3$. 

\item Next, we observe that the restriction $[\rho\vert_{\Sigma^{(2)}}]$ is of Type~\Romannum{2} and has $\alpha^{(2)}=(4\pi/3, \theta, \theta, 4\pi/3)$. Lemma~\ref{lem:beta_i-for-orbits-of-type-II} implies that $\gamma_2\in\{0,\pi\}$. Up to applying $\tau_{b_2}$ to $[\rho]$ once, this shows the existence an orbit point with $\gamma_2=0$.

\item Finally, the restriction $[\rho\vert_{\Sigma^{(3)}}]$ being of Type~\Romannum{4} with $\alpha^{(3)}=(\pi, \theta, \theta, \theta)$, we learn from Lemma~\ref{lem:beta_i-for-orbits-of-type-IV} that $\gamma_2\in\{0,2\pi/3, 4\pi/3\}$. Similarly, we can apply $\tau_{b_3}$ to $[\rho]$ once or twice to obtain an orbit point with $\gamma_3=2\pi/3$.
\end{itemize}
Since the Dehn twists $\tau_{b_1}$, $\tau_{b_2}$, and $\tau_{b_3}$ commute, we conclude from the above that there exists a point in the orbit of $[\rho]$ whose $\mathcal{B}$-triangle chain is the jester's hat triangle chain, showing that the orbit of $[\rho]$ coincides with the jester's hat orbit. 

We just proved that there is a unique finite orbit that contains an orbit point whose $\mathcal{B}$-triangle chain is regular. Next, we prove that there are no finite orbits in $\RepDT{\alpha}$ for which every orbit point has a singular $\mathcal{B}$-triangle chain. We'll assume by contradiction that such a finite orbit exists and prove that we can always construct an orbit point whose $\mathcal{B}$-triangle chain is regular, contradicting our assumption. We start with an arbitrary point $[\phi]$ belonging to such a finite orbit. Lemma~\ref{lem:existence-regular-triangle-chain} implies the existence of an other pants decomposition $\mathcal{B}'$ of $\Sigma$ such that the $\mathcal{B}'$-triangle chain of $[\phi]$ is regular. By the above, this means that the mapping class group orbit of $[\phi]$ is the jester's hat orbit for the pants decomposition $\mathcal{B}'$. In particular, it has length 40. Since $\RepDT{\alpha}$ contains only four points with a $\mathcal{B}$-triangle chain that consists of a single triangle, we can assume (up to replacing $[\phi]$ by another orbit point) that the $\mathcal{B}$-triangle chain of $[\phi]$ has one or two degenerate triangles. If $(\beta_1,\beta_2,\beta_3)$ denote the action coordinates of $[\phi]$ computed from $\mathcal{B}$, then the inequalities~\eqref{eq:inequalities-polytope-beta} read
\begin{equation*}
\begin{cases}
     \beta_1 \geq  4\pi-2\theta,  \\
     \beta_{i+1}-\beta_i \geq 2\pi-\theta, \quad i=1,2,\\
     \beta_{3} \leq 2\theta-2\pi. 
\end{cases}
\end{equation*}
We claim that this system of inequalities imply $\beta_i\in (2\pi-\theta,\theta)$ for $i=1,2,3$. For instance, we can infer that $\beta_1>2\pi-\theta$ from the first inequality by using that $\theta<2\pi$. Similarly, we deduce from $\beta_{2}-\beta_1 \geq 2\pi-\theta$ that $\beta_1<\theta$ by using that $\beta_2<2\pi$. We obtain $\beta_2, \beta_3\in (2\pi-\theta,\theta)$ by analogous arguments.

First, let's consider the case where the $\mathcal{B}$-triangle chain of $[\phi]$ contains a single degenerate triangle. The idea is to find a Dehn twist $\tau$ that maps $[\phi]$ to a point with a regular $\mathcal{B}$-triangle chain, providing the desired contradiction. Let $i\in\{1,2,3,4\}$ be the number such that the $i$th triangle in the chain of $[\phi]$ is degenerate. If $i\in\{1,2,3\}$, then we can take $\tau$ to be $\tau_{i+1,i+2}$ (following the notation of Section~\ref{sec:action-of-Dehn-twists}), and if $i=4$, then we instead take $\tau$ to be $\tau_{4,5}$. We'll explain why this choice works in details for the case $i=1$. Namely, we'll prove that when $i=1$, $\tau:=\tau_{2,3}$ sends $[\phi]$ to a point with a regular $\mathcal{B}$-triangle chain. The first non-degenerate triangle in the $\mathcal{B}$-chain of $[\phi]$ has vertices $(C_1=C_2,C_3,B_2)$, with respective angles $(2\pi-\theta,\pi-\theta/2,\pi-\beta_2/2)$. Let $D:=D_{2,3}$ be the intersection of the angle bisector at $C_1=C_2$ with the geodesic segment $C_3B_2$. We let $\delta:=\delta_{2,3}$ be the number such that $\angle C_2DC_3=\pi-\delta/2$.
\begin{center}
\begin{tikzpicture}[font=\sffamily,decoration={
    markings,
    mark=at position 1 with {\arrow{>}}}]
    
\node[anchor=south west,inner sep=0] at (0,0) {\includegraphics[width=6cm]{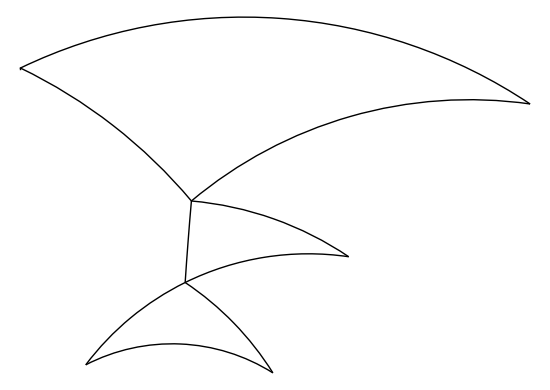}};

\begin{scope}
\fill (.2,3.51) circle (0.07) node[above left]{$C_1=C_2$};
\fill (5.8,3.15) circle (0.07) node[right]{$C_3$};
\fill (2.1,2.1) circle (0.07) node[left]{$B_2$};
\fill (2.05,1.2) circle (0.07) node[above left]{$B_3$};
\fill (3.8,1.45) circle (0.07) node[right]{$C_4$};
\fill (3,.215) circle (0.07) node[right]{$C_5$};
\fill (.95,.23) circle (0.07) node[left]{$C_6$};
\end{scope}

\fill (4,3.05) circle (0.07) node[below]{$D$};
\draw[dashed] (.2,3.51) to[bend left = 15] (4,3.05);

\draw[mauve, thick] (.8,3.2) arc (-40:15:.4) node[midway, right]{\small $\pi-\theta/2$};
\draw[mauve, thick] (1.1,3.6) arc (0:24:.6) node[at end, above]{\small $\pi-\theta/2$};
\draw[mauve, thick] (5.3,3.18) arc (180:160:.6) node[at end, above right]{\small $\pi-\theta/2$};
\draw[thick] (4.5,3.14) arc (5:163:.5) node[at start, below right]{\small $\pi-\delta/2$};
\draw[apricot, thick] (2.5,2.4) arc (30:133:.5) node[midway, above]{\small $\pi-\beta_2/2$};
\end{tikzpicture}
\end{center}
The Dehn twist $\tau$ acts on the chain of $[\phi]$ by rotating the triangle $(C_2, C_3,D)$ anti-clockwise around $D$ by an angle $\delta$ and leaves the other vertices fixed. The resulting $\mathcal{B}$-triangle chain has $C_1\neq C_2$ since $0<\delta<2\pi$. It is thus regular as long as $\tau$ doesn't send $C_3$ to $B_2$, which only happens if $\delta=\pi$ and $DC_3=DB_2$. The latter is equivalent to $\pi-\theta/2=\pi-\beta_2/2$ because $D$ lies on the angle bisector at $C_1=C_2$. This, however, cannot hold because we observed above that $\beta_2<\theta$. The other cases where $i$ belongs to $\{2,3,4\}$ can be treated similarly.

It remains to consider the case where the $\mathcal{B}$-triangle chain of $[\phi]$ has two degenerate triangles. The strategy here is to first find a Dehn twist $\tau$ such that the $\mathcal{B}$-triangle chain of $\tau.[\phi]$ has only one degenerate triangle. We can then apply the argument above to find a second Dehn twist $\tau'$ such that the $\mathcal{B}$-triangle chain of $(\tau'\tau).[\phi]$ has no degenerate triangles any more, contradicting our assumption on $[\phi]$. We find the Dehn twist $\tau$ as follows. Let $i<j\in\{1,2,3,4\}$ denote the two numbers such that the $i$th and $j$th triangles in the $\mathcal{B}$-triangle chain of $[\phi]$ are degenerate. If $j\in\{2,3\}$, then we take $\tau$ to be $\tau_{j+1,j+2}$ as in the previous case. If $j=4$ and $i\in\{1,2\}$, then we pick $\tau$ to be $\tau_{4,5}$ instead. Finally, if $j=4$ and $i=3$, then we take $\tau$ to be $\tau_{3,4}$. An analogous argument as before shows that when $\tau$ is picked that way, then the $\mathcal{B}$-triangle chain of $\tau.[\phi]$ has only one degenerate triangle. This finishes the proof of the uniqueness statement.
\end{proof} 

\begin{rem}
The jester's hat orbit is of pullback type (Definition~\ref{def:pullback-orbits}) for any value of $\theta>5\pi/3$, see~\cite[Table~4]{diarra-pull-back}. It's interesting to observe that the jester's hat orbit consists of discrete representations if and only if $\overline{\theta}$ is equal to some integer $k\geq 7$ by Theorem~\ref{thm:felikson}. This is equivalent to $\theta$ being of the form $2\pi\frac{k-1}{k}$ for some integer $k\geq 7$. For instance, when $\theta=12\pi/7$ ($k=7$), the jester's hat orbit consists of discrete representations whose image is a rotation triangle group $D(2,3,7)$.
\end{rem}

\subsection{Classification of finite orbits for \texorpdfstring{$n=5$}{n=5}}\label{sec:classification-n=5}
We explain how our methods also provide a classification of finite orbits in the case of a 5-punctured sphere. Recall that this classification was already achieved by Tykhyy using a different (computer-aided) approach. Here's what Tykhyy proved for DT representations\footnote{The description of the images of the representations associated to the finite orbits in Theorem~\ref{thm:angle-vector-alpha-with-finite-orbits-n=5} in terms of rotation triangle groups is not mentionned in~\cite{tykhyy}, but can immediately be deduced from the triangle chain descriptions of orbit points provided in Sections~\ref{sec:proof-thm-finite-orbits-n=5-1}--\ref{sec:proof-thm-finite-orbits-n=5-3}.}, see~\cite[Table~9]{tykhyy}. 

\begin{thm}[\cite{tykhyy}]\label{thm:angle-vector-alpha-with-finite-orbits-n=5}
In the case of a 5-punctured sphere, if $\RepDT{\alpha}$ contains a finite mapping class group orbit, then $\alpha$ is one of the following (unordered) angle vectors:
\[
\left\{\frac{4\pi}{3}, \theta, \theta, \theta, \theta\right\},\quad \left\{2\theta-2\pi, \theta, \theta, \theta, \theta\right\}, \quad \left\{\frac{12\pi}{7}, \frac{12\pi}{7}, \frac{12\pi}{7}, \frac{12\pi}{7}, \frac{12\pi}{7}\right\}.
\]
In every case, the parameter $\theta$ satisfies $\theta\in (5\pi/3,2\pi)$. Moreover, the finite orbit is always unique in the corresponding DT component and its length is provided below. The images of representations associated to the finite orbits are rotation triangle groups which can be found in the following table, where we use the notation $\overline{\theta}=(1-\theta/2\pi)^{-1}$.
\begin{table}[ht]
    \centering
    \begin{tblr}{c|c|c|c|c}
        & $\alpha$ & Orbit length & Image & Coordinates \\
        \hline\hline
        \makecell{hang-glider\\orbit} & $\{\frac{4\pi}{3},\theta,\theta,\theta,\theta\}$ & 9 & $D(2,3,\overline{\theta})$ & Table~\ref{tab:hang-glider-orbit} \\
        \hline
        \makecell{sand clock\\orbit} & $\{2\theta-2\pi,\theta,\theta,\theta,\theta\}$ & 12 & $D(2,3,\overline{\theta})$ & Table~\ref{tab:sand-clock-orbit}  \\
        \hline
        bat orbit & $\{\frac{12\pi}{7},\frac{12\pi}{7},\frac{12\pi}{7},\frac{12\pi}{7},\frac{12\pi}{7}\}$ & 105 & $D(2,3,7)$ & Table~\ref{tab:bat-orbit} \\
    \end{tblr}
\end{table}
\end{thm}

\begin{rem}
Every finite orbit in Theorem~\ref{thm:angle-vector-alpha-with-finite-orbits-n=5} is of pullback type (Definition~\ref{def:pullback-orbits}), see~\cite[Table~2]{diarra-pull-back}. They are respectively labelled 3, 4, and 7 in~\cite{tykhyy}. We also observe that the orbit lengths we find coincide with those found by Tykhyy~\cite[Table~9]{tykhyy}.
\end{rem}


Our next goal is to explain how our methods can be used to re-prove Theorem~\ref{thm:angle-vector-alpha-with-finite-orbits-n=5}. Our arguments unfold over the next sections (Sections~\ref{sec:plan-of-proof-thm-n=5}--\ref{sec:third-bad-angle-vector}).

\subsubsection{Plan of proof for Theorem~\ref{thm:angle-vector-alpha-with-finite-orbits-n=5}}\label{sec:plan-of-proof-thm-n=5}
In order to prove Theorem~\ref{thm:angle-vector-alpha-with-finite-orbits-n=5}, we start by fixing some $[\rho]\in\RepDT{\alpha}$ within a finite mapping class group orbit. We then apply Lemma~\ref{lem:existence-regular-triangle-chain} in order to find a pants decomposition $\mathcal{B}$ of $\Sigma$ and a compatible geometric presentation of $\pi_1\Sigma$ such that the $\mathcal{B}$-triangle chain of $[\rho]$ is regular. We'll write $\beta_1<\beta_2$ for the action coordinates of $[\rho]$ computed from $\mathcal{B}$. 

The first step in the proof consists in restricting the possible values of the pair $(\beta_1,\beta_2)$. Our argument, as before, studies the possible partial sequences provided by Table~\ref{tab:partial-sequences-beta_1-beta_2-beta_n-4-beta-n-3}. Since $[\rho]$ has two restrictions $[\rho\vert_{\Sigma^{(1)}}]$ and $[\rho\vert_{\Sigma^{(2)}}]$, the partial sequence $(\beta_1<\beta_2)$ must appear both as $(\beta_1<\beta_2)$ and as $(\beta_{n-4}<\beta_{n-3})$ (with $n=5$) in Table~\ref{tab:partial-sequences-beta_1-beta_2-beta_n-4-beta-n-3}. It turns out that only finitely many angles appear simultaneously in the second and fourth, respectively third and fifth, columns of Table~\ref{tab:partial-sequences-beta_1-beta_2-beta_n-4-beta-n-3}, restricting the possible values of $\beta_1$ and $\beta_2$ to a finite set (Section~\ref{sec:shortlist-beta_1-n=5}). We'll then make this analysis more precise by considering each value of $\beta_1$ individually and computing all possible corresponding values for $\beta_2$ (Sections~\ref{sec:beta_1=pi}--\ref{sec:beta_1=2pi/5}). This leads to a finite list of candidate angle vectors $\alpha$ (Section~\ref{sec:shortlist-angle-vector-alpha}). Next, we'll eliminate all angles vectors $\alpha$ obtained in the previous step for which $\RepDT{\alpha}$ doesn't contain a finite mapping class group orbit (Sections~\ref{sec:first-bad-angle-vector}--\ref{sec:third-bad-angle-vector}). Finally, we'll compute the action-angle coordinates of every orbit point, as well as their images, for the three angle vectors of Theorem~\ref{thm:angle-vector-alpha-with-finite-orbits-n=5} and prove that the finite orbit is unique in all three cases (Sections~\ref{sec:proof-thm-finite-orbits-n=5-1}--\ref{sec:proof-thm-finite-orbits-n=5-3}).

\subsubsection{Shortlisting the values of $\beta_1$ and $\beta_2$}\label{sec:shortlist-beta_1-n=5}

\begin{lem}\label{lem:beta_1-in-case-n=5}
The action coordinates of $[\rho]$ satisfy 
\[
\beta_1\in \left\{\frac{2\pi}{5}, \frac{2\pi}{3}, \frac{4\pi}{5}, \frac{6\pi}{7}, \pi\right\},\quad \beta_2\in \left\{\pi, \frac{8\pi}{7}, \frac{6\pi}{5}, \frac{4\pi}{3},\frac{8\pi}{5}\right\}.
\]    
\end{lem}
\begin{proof}
Since $[\rho]$ has two restrictions $[\rho\vert_{\Sigma^{(1)}}]$ and $[\rho\vert_{\Sigma^{(2)}}]$, the partial sequence $(\beta_1<\beta_2)$ must appear in Table~\ref{tab:partial-sequences-beta_1-beta_2-beta_n-4-beta-n-3} both as $(\beta_1<\beta_2)$ (second and third columns) and as $(\beta_{n-4}<\beta_{n-3})$ (fourth and fifth columns, with $n=5$). The possible values for $\beta_1$ are those that appear both in the second and the fourth columns of Table~\ref{tab:partial-sequences-beta_1-beta_2-beta_n-4-beta-n-3}. Studying the values in the fourth column, we first observe that $\beta_1\leq \pi$ (as we did in the proof of Theorem~\ref{thm:no-finite-orbit-for-n-geq-7}). Now, when we study the values in the second column of Table~\ref{tab:partial-sequences-beta_1-beta_2-beta_n-4-beta-n-3} and use that $\beta_1\leq \pi$, we obtain desired list of values for $\beta_1$. Similarly, the third column of Table~\ref{tab:partial-sequences-beta_1-beta_2-beta_n-4-beta-n-3} contains the possible values for $\beta_2$ and gives $\beta_2\geq \pi$. The fifth column then provides the expected values for $\beta_2$.
\end{proof}

\subsubsection{$\beta_1=\pi$}\label{sec:beta_1=pi}
When studying the fourth column of Table~\ref{tab:partial-sequences-beta_1-beta_2-beta_n-4-beta-n-3} that gives the possible values for $\beta_1$, we see that the value $\beta_1=\pi$ is only possible when $\beta_2=4\pi/3$ and the finite orbit of $[\rho\vert_{\Sigma^{(2)}}]$ is of Type~\Romannum{4}. This gives $\alpha^{(2)}=(\pi,\theta,\theta,\theta)$ with the angle $\pi$ playing the role of $2\pi-\beta_1$. In order to find out the type of the finite orbit of $[\rho\vert_{\Sigma^{(1)}}]$, we look for the pair $(\pi,4\pi/3)$ the second and third column of Table~\ref{tab:partial-sequences-beta_1-beta_2-beta_n-4-beta-n-3}. There are three matching lines telling us that the finite orbit of $[\rho\vert_{\Sigma^{(1)}}]$ can be of Type~\Romannum{2} with one angle equal to $4\pi/3$, of Type~\Romannum{3}, or of Type~33. Let us consider these three cases individually.

\begin{itemize}
\item In the first case, since $\beta_1=\pi$, Lemma~\ref{lem:beta_i-for-orbits-of-type-II} implies that $\alpha^{(1)}=(4\pi/3,\theta',\theta',4\pi/3)$ or $\alpha^{(1)}=(\theta',4\pi/3,\theta',4\pi/3)$, with the last angle $4\pi/3$ being the angle $\beta_2$. Since $\alpha^{(1)}_3=\alpha_3=\alpha^{(2)}_2$, we obtain $\alpha_3=\theta=\theta'$. We conclude that $\alpha$ is either the angle vector $(4\pi/3,\theta,\theta,\theta,\theta)$ or the same angle vector but with the first two entries permuted. 
\item Now, if the finite orbit of $[\rho\vert_{\Sigma^{(1)}}]$ is of Type~\Romannum{3} instead, then $\alpha^{(1)}$ is one of the three angle vectors $(2\theta'-2\pi,\theta',\theta',4\pi/3)$, $(\theta',2\theta'-2\pi,\theta',4\pi/3)$, or $(\theta',\theta',2\theta'-2\pi,4\pi/3)$, where $4\pi/3$ is the angle $\beta_2$. Arguing that $\alpha_3$ belongs to both $\alpha^{(1)}$ and $\alpha^{(2)}$, we have either $\alpha_3=\theta'=\theta$ and $\alpha=(2\theta-2\pi,\theta,\theta,\theta,\theta)$ (or the same angle vector but with the first two entries permuted), or $\alpha_3=2\theta'-2\pi=\theta$ and $\alpha=(\theta',\theta',2\theta'-2\pi,2\theta'-2\pi,2\theta'-2\pi)$. The latter case occurs when $\alpha^{(1)}=(\theta',\theta',2\theta'-2\pi, 4\pi/3)$ which does not allow for $\beta_1=\pi$ according to Lemma~\ref{lem:beta_i-for-orbits-of-type-III}. So, we can actually eliminate the angle vector $(\theta',\theta',2\theta'-2\pi,2\theta'-2\pi,2\theta'-2\pi)$ right away.
\item Finally, if the finite orbit of $[\rho\vert_{\Sigma^{(1)}}]$ is of Type~33, then $\alpha^{(1)}=(12\pi/7,12\pi/7,12\pi/7,4\pi/3)$ with $4\pi/3$ being the angle $\beta_2$. Now, since $\alpha_3$ belongs to both $\alpha^{(1)}$ and $\alpha^{(2)}$, we must have $\theta=12\pi/7$ giving $\alpha=(12\pi/7,12\pi/7,12\pi/7,12\pi/7, 12\pi/7)$.
\end{itemize}

     \begin{table}[ht]
         \centering
         \begin{tblr}{c|c|c|c|c}
             $\beta_1<\beta_2$ & \makecell{Combination\\ of solutions} &  $\alpha^{(1)}$ & $\alpha^{(2)}$ & $\alpha$ \\
             \hline\hline
             \SetCell[r=3]{c} $\left(\pi, \frac{4\pi}{3}\right)$ & 33--\Romannum{4} & $\left(\frac{12\pi}{7},\frac{12\pi}{7},\frac{12\pi}{7}, \frac{4\pi}{3}\right)$ & $\left(\pi,\frac{12\pi}{7},\frac{12\pi}{7}, \frac{12\pi}{7}\right)$ & $\left(\frac{12\pi}{7}, \frac{12\pi}{7}, \frac{12\pi}{7},\frac{12\pi}{7},\frac{12\pi}{7}\right)$\\
             \cline{2-5}
             & \Romannum{2}--\Romannum{4} & $\left(\frac{4\pi}{3},\theta,\theta, \frac{4\pi}{3}\right)$ & $\left(\pi,\theta,\theta, \theta\right)$ & $\left(\frac{4\pi}{3}, \theta, \theta,\theta,\theta \right)$\\
             \cline{2-5}
             & \Romannum{3}--\Romannum{4} & $\left(2\theta-2\pi,\theta,\theta, \frac{4\pi}{3}\right)$ & $\left(\pi,\theta,\theta, \theta\right)$ & $\left(2\theta-2\pi, \theta, \theta,\theta,\theta \right)$\\
         \end{tblr}
     \end{table}

\subsubsection{$\beta_1=6\pi/7$}\label{sec:beta_1=6pi/7}
According to the fourth column of Table~\ref{tab:partial-sequences-beta_1-beta_2-beta_n-4-beta-n-3} which gives the values for $\beta_1$, the value $\beta_1=6\pi/7$ is only possible when $\beta_2=\pi$ and the finite orbit of $[\rho\vert_{\Sigma^{(2)}}]$ is of Type~\Romannum{2} with one angle equal to $8\pi/7$, or when $\beta_2=4\pi/3$ and the finite orbit of $[\rho\vert_{\Sigma^{(2)}}]$ is of Type~\Romannum{4}$^\ast$ with $\theta=12\pi/7$. However, when we compare with the third column of Table~\ref{tab:partial-sequences-beta_1-beta_2-beta_n-4-beta-n-3} that gives the possible values for $\beta_2$, we see that $\beta_2=\pi$ is impossible when $\beta_1=6\pi/7$. This means that only $\beta_2=4\pi/3$ is possible when $\beta_1=6\pi/7$. The only combination of solutions is 33--\Romannum{4}$^\ast$ with $\alpha^{(1)}=(12\pi/7,12\pi/7,12\pi/7,4\pi/3)$ and $\alpha^{(2)}=(8\pi/7,12\pi/7,12\pi/7,12\pi/7)$. This gives $\alpha=(12\pi/7,12\pi/7,12\pi/7,12\pi/7,12\pi/7)$.

     \begin{table}[ht]
         \centering
         \begin{tblr}{c|c|c|c|c}
             $\beta_1<\beta_2$ & \makecell{Combination\\ of solutions} &  $\alpha^{(1)}$ & $\alpha^{(2)}$ & $\alpha$ \\
             \hline\hline
             $\left(\frac{6\pi}{7}, \frac{4\pi}{3}\right)$ & 33--\Romannum{4}$^\ast$ & $\left(\frac{12\pi}{7},\frac{12\pi}{7},\frac{12\pi}{7}, \frac{4\pi}{3}\right)$ & $\left(\frac{8\pi}{7},\frac{12\pi}{7},\frac{12\pi}{7}, \frac{12\pi}{7}\right)$ & $\left(\frac{12\pi}{7}, \frac{12\pi}{7}, \frac{12\pi}{7},\frac{12\pi}{7},\frac{12\pi}{7}\right)$\\
         \end{tblr}
     \end{table}

\subsubsection{$\beta_1=4\pi/5$}\label{sec:beta_1=4pi/5}
We proceed as above. We first look at the third column of Table~\ref{tab:partial-sequences-beta_1-beta_2-beta_n-4-beta-n-3} giving the possible values for $\beta_2$, we conclude that $\beta_2=3\pi/2$ when $\beta_1=4\pi/5$. However, the value $\beta_2=3\pi/2$ never appears in the fifth column of Table~\ref{tab:partial-sequences-beta_1-beta_2-beta_n-4-beta-n-3} giving the possible values for $\beta_{n-4}$. We conclude that $\beta_1=4\pi/5$ is actually not possible.

\subsubsection{$\beta_1=2\pi/3$}\label{sec:beta_1=2pi/3}
This case is lengthier because $\beta_2$ can take three different values when $\beta_1=2\pi/3$ as we're about to see. We proceed as in the other cases and study the last two columns of Table~\ref{tab:partial-sequences-beta_1-beta_2-beta_n-4-beta-n-3}. Using $\beta_1=2\pi/3$ and remembering that $\beta_2\geq \pi$ from Lemma~\ref{lem:beta_1-in-case-n=5}, we obtain the following list
\[
\beta_2\in\left\{\pi, \frac{8\pi}{7}, \frac{4\pi}{3}\right\}.
\]
We treat each of these sub-cases individually.

\begin{itemize}
    \item  $\beta_2=\pi$. In this case, the finite orbit of $[\rho\vert_{\Sigma^{(1)}}]$ is necessarily of Type~\Romannum{4} giving $\alpha^{(1)}=(\theta,\theta,\theta,\pi)$ where $\pi$ plays the role of $\beta_2$. The finite orbit of $[\rho\vert_{\Sigma^{(2)}}]$ can be of Type~\Romannum{2}, Type~\Romannum{3}, or of Type~33. 
\begin{itemize}
\item If it is of Type~\Romannum{2}, then $\alpha^{(2)}$ is either $(4\pi/3,\theta',\theta',4\pi/3)$ or $(4\pi/3,\theta',4\pi/3,\theta')$, with the first $4\pi/3$ being the angle $2\pi-\beta_1$. Arguing on $\alpha_3$ as we did before and using that $\theta>5\pi/3$, we conclude that $\theta=\theta'$ and $\alpha$ is the angle vector $(\theta,\theta,\theta,\theta,4\pi/3)$ or the one obtained by permuting the last two entries. 
\item If the finite orbit of $[\rho\vert_{\Sigma^{(2)}}]$ is instead of Type~\Romannum{3}, then $\alpha^{(2)}$ is either $(4\pi/3, \theta', \theta',2\theta'-2\pi)$ or $(4\pi/3, \theta', 2\theta'-2\pi, \theta')$, since the angle vector $(4\pi/3,2\theta'-2\pi,\theta', \theta')$ does not allow for $\beta_2=\pi$ by Lemma~\ref{lem:beta_i-for-orbits-of-type-III}. We argue on $\alpha_3$ to deduce that $\alpha$ is the angle vector $\alpha=(\theta,\theta,\theta,\theta,2\theta-2\pi)$ or the one obtained by permuting the last two entries.
\item Finally, if the finite orbit of $[\rho\vert_{\Sigma^{(2)}}]$ is of Type~33, then the tuple $\alpha^{(2)}$ is given by $(4\pi/3, 12\pi/7, 12\pi/7, 12\pi/7)$. This implies $\theta=12\pi/7$ and the tuple $\alpha$ is given
	by $(12\pi/7,12\pi/7, 12\pi/7,12\pi/7, 12\pi/7)$.
\end{itemize}

    \item $\beta_2=8\pi/7$. By studying Table~\ref{tab:partial-sequences-beta_1-beta_2-beta_n-4-beta-n-3}, one would observe that the pair $\beta_1=2\pi/3$ and $\beta_2=8\pi/7$ is only possible with the combination of solutions \Romannum{4}$^\ast$--33, with $\theta=12\pi/7$ for the solution of Type~\Romannum{4}$^\ast$. This means that $\alpha^{(1)}=(12\pi/7, 12\pi/7, 12\pi/7,8\pi/7)$ and $\alpha^{(2)}=(4\pi/3,12\pi/7, 12\pi/7, 12\pi/7)$. We get $\alpha=(12\pi/7,12\pi/7, 12\pi/7,12\pi/7, 12\pi/7)$.

    \item $\beta_2=4\pi/3$. This case is the richest. By studying Table~\ref{tab:partial-sequences-beta_1-beta_2-beta_n-4-beta-n-3} again, one would observe that the pair $\beta_1=2\pi/3$ and $\beta_2=4\pi/3$ is possible when the finite orbits of $[\rho\vert_{\Sigma^{(1)}}]$ and $[\rho\vert_{\Sigma^{(2)}}]$ are of Type~\Romannum{3}, Type~\Romannum{4}$^\ast$ with $\theta=16\pi/9$, or of Type~33. This represents nine possible combinations of solutions. 
\begin{itemize}
    \item Let's start with the case where the finite orbit of $[\rho\vert_{\Sigma^{(1)}}]$ is of Type~33. This means that $\alpha^{(1)}=(12\pi/7,12\pi/7,12\pi/7,4\pi/3)$. When the finite orbit of $[\rho\vert_{\Sigma^{(2)}}]$ is of Type~33, we have $\alpha^{(2)}=(4\pi/3,12\pi/7,12\pi/7,12\pi/7)$. This leads to $\alpha=(12\pi/7, 12\pi/7, 12\pi/7, 12 \pi/7, 12\pi/7)$. If instead the finite orbit of $[\rho\vert_{\Sigma^{(2)}}]$ is of Type~\Romannum{3}, then $\alpha^{(2)}=(4\pi/3,2\theta-2\pi,\theta, \theta)$ because it's the only angle vector that starts with $4\pi/3$ and allows for $\beta_2=4\pi/3$ according to Lemma~\ref{lem:beta_i-for-orbits-of-type-III}. Arguing on the value of $\alpha_3$, we conclude that either $\theta=13\pi/7$ and $\alpha=(12\pi/7, 12\pi/7, 12\pi/7, 13\pi/7, 13\pi/7)$. Finally, if the finite orbit of $[\rho\vert_{\Sigma^{(2)}}]$ is of Type~\Romannum{4}$^\ast$ with $\theta=16\pi/9$, then $\alpha^{(2)}=(4\pi/3, 16\pi/9, 16\pi/9, 16\pi/9)$. This case can be eliminated as there is no choice for the value of $\alpha_3$.

    \item Let's now move to the case where the finite orbit of $[\rho\vert_{\Sigma^{(1)}}]$ is of Type~\Romannum{3}. In that case, $\alpha^{(1)}=(\theta, \theta,2\theta-2\pi, 4\pi/3)$ since, again, it's the only angle vector that ends with $4\pi/3$ and allows for $\beta_1=2\pi/3$ according to Lemma~\ref{lem:beta_i-for-orbits-of-type-III}. If the finite orbit of $[\rho\vert_{\Sigma^{(2)}}]$ is of Type~33, then we obtained $\alpha^{(2)}=(4\pi/3,12\pi/7,12\pi/7,12\pi/7)$ and $\alpha=(13\pi/7, 13\pi/7, 12\pi/7, 12\pi/7, 12\pi/7)$ as above. If the finite orbit of $[\rho\vert_{\Sigma^{(2)}}]$ is of Type~\Romannum{3} too, then $\alpha^{(2)}=(4\pi/3, 2\theta'-2\pi, \theta', \theta')$. Arguing on the value of $\alpha_3$, we see that $\theta=\theta'$ and $\alpha=(\theta, \theta, 2\theta-2\pi, \theta, \theta)$. If the finite orbit of $[\rho\vert_{\Sigma^{(2)}}]$ is of Type~\Romannum{4}$^\ast$ with $\theta'=16\pi/9$ instead, then $\alpha^{(2)}=(4\pi/3, 16\pi/9, 16\pi/9, 16\pi/9)$. When looking for the value of $\alpha_3$, we obtain $\theta=17\pi/9$ and $\alpha=(17\pi/9, 17\pi/9, 16\pi/9, 16\pi/9, 16\pi/9)$.

    \item Finally, we consider the case where the finite orbit of $[\rho\vert_{\Sigma^{(1)}}]$ is of Type~\Romannum{4}$^\ast$ with $\theta=16\pi/9$, meaning that $\alpha^{(1)}=(16\pi/9, 16\pi/9, 16\pi/9,4\pi/3)$. As before, the finite orbit of $[\rho\vert_{\Sigma^{(2)}}]$ can't be of Type~33. If it is instead of Type~\Romannum{3}, then we obtain $\alpha^{(2)}=(4\pi/3, 16\pi/9, 17\pi/9, 17\pi/9)$ and $\alpha=(16\pi/9, 16\pi/9, 16\pi/9, 17\pi/9, 17\pi/9)$. Lastly, if the finite orbit of $[\rho\vert_{\Sigma^{(2)}}]$ is of Type~\Romannum{4}$^\ast$, then we obtain that $\alpha^{(2)}$ is $(4\pi/3, 16\pi/9, 16\pi/9, 16\pi/9)$ and $\alpha=(16\pi/9, 16\pi/9, 16\pi/9, 16\pi/9, 16\pi/9)$.
    \end{itemize}

     \begin{table}[ht]
         \centering
         \begin{tblr}{c|c|c|c|c}
             $\beta_1<\beta_2$ & \makecell{Combination\\ of solutions} &  $\alpha^{(1)}$ & $\alpha^{(2)}$ & $\alpha$ \\
             \hline\hline
             \SetCell[r=3]{c} $\left(\frac{2\pi}{3}, \pi\right)$ & \Romannum{4}--33 & $\left(\frac{12\pi}{7},\frac{12\pi}{7},\frac{12\pi}{7}, \pi\right)$ & $\left(\frac{4\pi}{3},\frac{12\pi}{7},\frac{12\pi}{7}, \frac{12\pi}{7}\right)$ & $\left(\frac{12\pi}{7}, \frac{12\pi}{7}, \frac{12\pi}{7},\frac{12\pi}{7},\frac{12\pi}{7}\right)$\\
             \cline{2-5}
             & \Romannum{4}--\Romannum{2} & \SetCell[r=2]{c} $\left(\theta,\theta,\theta, \pi\right)$ & $\left(\frac{4\pi}{3},\theta,\theta, \frac{4\pi}{3}\right)$ & $\left(\theta, \theta,\theta,\theta,\frac{4\pi}{3} \right)$\\
             \cline{2}\cline{4-5}
             & \Romannum{4}--\Romannum{3}  &  & $\left(\frac{4\pi}{3},\theta,\theta, 2\theta-2\pi\right)$ & $\left(\theta, \theta,\theta,\theta,2\theta-2\pi \right)$\\
             \hline
    
            $\left(\frac{2\pi}{3}, \frac{8\pi}{7}\right)$ & \Romannum{4}$^\ast$--33 & $\left(\frac{12\pi}{7},\frac{12\pi}{7},\frac{12\pi}{7}, \frac{8\pi}{7}\right)$ & $\left(\frac{4\pi}{3},\frac{12\pi}{7},\frac{12\pi}{7}, \frac{12\pi}{7}\right)$ & $\left(\frac{12\pi}{7}, \frac{12\pi}{7}, \frac{12\pi}{7},\frac{12\pi}{7},\frac{12\pi}{7}\right)$\\
             \hline
            
             \SetCell[r=7]{c} $\left(\frac{2\pi}{3}, \frac{4\pi}{3}\right)$ & 33--33 & \SetCell[r=2]{c} $\left(\frac{12\pi}{7},\frac{12\pi}{7},\frac{12\pi}{7}, \frac{4\pi}{3}\right)$ & $\left(\frac{4\pi}{3},\frac{12\pi}{7},\frac{12\pi}{7}, \frac{12\pi}{7}\right)$ & $\left(\frac{12\pi}{7}, \frac{12\pi}{7}, \frac{12\pi}{7},\frac{12\pi}{7},\frac{12\pi}{7}\right)$\\
            \cline{2}\cline{4-5}
              & 33--\Romannum{3} & & $\left(\frac{4\pi}{3},\frac{12\pi}{7},\frac{13\pi}{7}, \frac{13\pi}{7}\right)$ & $\left(\frac{12\pi}{7}, \frac{12\pi}{7}, \frac{12\pi}{7},\frac{13\pi}{7},\frac{13\pi}{7}\right)$\\
             \cline{2-5}
             & \Romannum{3}--33 & $\left(\frac{13\pi}{7},\frac{13\pi}{7},\frac{12\pi}{7}, \frac{4\pi}{3}\right)$ & $\left(\frac{4\pi}{3},\frac{12\pi}{7},\frac{12\pi}{7}, \frac{12\pi}{7}\right)$ & $\left(\frac{13\pi}{7}, \frac{13\pi}{7}, \frac{12\pi}{7},\frac{12\pi}{7},\frac{12\pi}{7}\right)$\\
             \cline{2-5}
             & \Romannum{3}--\Romannum{3} & $\left(\theta,\theta,2\theta-2\pi, \frac{4\pi}{3}\right)$ & $\left(\frac{4\pi}{3}, 2\theta-2\pi,\theta,\theta\right)$ & $\left(\theta, \theta, 2\theta-2\pi,\theta,\theta\right)$\\
              \cline{2-5}
             & \Romannum{3}--\Romannum{4}$^\ast$ & $\left(\frac{17\pi}{9},\frac{17\pi}{9},\frac{16\pi}{9}, \frac{4\pi}{3}\right)$ & $\left(\frac{4\pi}{3},\frac{16\pi}{9},\frac{16\pi}{9}, \frac{16\pi}{9}\right)$ & $\left(\frac{17\pi}{9}, \frac{17\pi}{9}, \frac{16\pi}{9},\frac{16\pi}{9},\frac{16\pi}{9}\right)$\\
             \cline{2-5}
             & \Romannum{4}$^\ast$--\Romannum{3} & $\left(\frac{16\pi}{9},\frac{16\pi}{9},\frac{16\pi}{9}, \frac{4\pi}{3}\right)$ & $\left(\frac{4\pi}{3},\frac{16\pi}{9},\frac{17\pi}{9}, \frac{17\pi}{9}\right)$ & $\left(\frac{16\pi}{9}, \frac{16\pi}{9}, \frac{16\pi}{9},\frac{17\pi}{9},\frac{17\pi}{9}\right)$\\
             \cline{2-5}
             & \Romannum{4}$^\ast$--\Romannum{4}$^\ast$ & $\left(\frac{16\pi}{9},\frac{16\pi}{9},\frac{16\pi}{9}, \frac{4\pi}{3}\right)$ & $\left(\frac{4\pi}{3},\frac{16\pi}{9},\frac{16\pi}{9}, \frac{16\pi}{9}\right)$ & $\left(\frac{16\pi}{9}, \frac{16\pi}{9}, \frac{16\pi}{9},\frac{16\pi}{9},\frac{16\pi}{9}\right)$\\
         \end{tblr}
     \end{table}
 \end{itemize}

\subsubsection{$\beta_1=2\pi/5$}\label{sec:beta_1=2pi/5}
We start by studying the third column of Table~\ref{tab:partial-sequences-beta_1-beta_2-beta_n-4-beta-n-3} and we observe that if $\beta_1=2\pi/5$, then $\beta_2\in\{3\pi/2, 6\pi/5\}$. Since $3\pi/2$ doesn't appear in the last column of Table~\ref{tab:partial-sequences-beta_1-beta_2-beta_n-4-beta-n-3}, we are left with the case $\beta_2=6\pi/5$. The combination of solutions must be 27--26. This case can be eliminated since there will be no choice for the value of $\alpha_3$. So, $\beta_1=2\pi/5$ is impossible.

\subsubsection{Resulting shortlist of angle vectors $\alpha$}\label{sec:shortlist-angle-vector-alpha}

A finite list of candidate angle vectors $\alpha$ for which $\RepDT{\alpha}$ may contain a finite mapping class group orbit have been computed in Sections~\ref{sec:beta_1=pi}--\ref{sec:beta_1=2pi/5}. Up to permutation of the entries, the candidates for $\alpha$ are
\begin{itemize}
\item $(12\pi/7, 12\pi/7, 12\pi/7, 12\pi/7, 12\pi/7)$
\item $(4\pi/3, \theta, \theta, \theta, \theta)$
\item $(2\theta-2\pi, \theta, \theta, \theta, \theta)$
\item $(12\pi/7, 12\pi/7, 12\pi/7,13\pi/7,13\pi/7)$
\item $(16\pi/9, 16\pi/9, 16\pi/9,17\pi/9,17\pi/9)$
\item $(16\pi/9, 16\pi/9, 16\pi/9,16\pi/9,16\pi/9)$
\end{itemize}
In order to finish the proof of Theorem~\ref{thm:angle-vector-alpha-with-finite-orbits-n=5}, it remains to first prove that the last three angle vectors don't correspond to any DT component supporting a finite orbit and then exhibit a finite orbit inside the DT component corresponding to the other three angle vectors (recall that we don't show the uniqueness of those orbits).

\subsubsection{Getting rid of $\alpha=(12\pi/7, 12\pi/7, 12\pi/7,13\pi/7,13\pi/7)$}\label{sec:first-bad-angle-vector}

When we look back in the Sections~\ref{sec:beta_1=pi} to \ref{sec:beta_1=2pi/5}, we observe that all orbit points with a regular triangle chain in a finite mapping class group orbit inside $\RepDT{\alpha}$
where $\alpha=(12\pi/7, 12\pi/7, 12\pi/7,13\pi/7,13\pi/7)$ must have $(\beta_1,\beta_2)=(2\pi/3,4\pi/3)$ with the combination of solutions 33--\Romannum{3}. However, we learned in Lemma~\ref{lem:beta_i-for-obits-of-type-33} that since the orbit of
the restriction $[\rho\vert_{\Sigma^{(1)}}]$ is of Type~33 with $\alpha^{(1)}$ being the tuple$(12\pi/7, 12\pi/7, 12\pi/7, 4\pi/3)$. It also contains orbit points with $\beta_1=4\pi/7$ and $\beta_1=\pi$ which is not the case here.
This shows that $\RepDT{\alpha}$ doesn't contain any finite orbit when $\alpha=(12\pi/7, 12\pi/7, 12\pi/7,13\pi/7,13\pi/7)$. The argument for $\alpha=(13\pi/7, 13\pi/7, 12\pi/7,12\pi/7,12\pi/7)$ is analogous.

\subsubsection{Getting rid of $\alpha=(16\pi/9, 16\pi/9, 16\pi/9,17\pi/9,17\pi/9)$}\label{sec:second-bad-angle-vector}

Getting rid of the angle vector $\alpha=(16\pi/9, 16\pi/9, 16\pi/9,17\pi/9,17\pi/9)$ requires more work than the previous case. Our strategy is to find an interior simple closed curve whose image is an irrational rotation. This will force the orbit to be infinite (Lemma~\ref{lem:finite-mcg-orbit-implies-beta_i-rational}).

Assume for the sake of contradiction that $\RepDT{\alpha}$ contains a finite orbit and let's study the triangle chains of points in it. We learned from Sections~\ref{sec:beta_1=pi}--\ref{sec:beta_1=2pi/5} that a point $[\rho]$ in a finite orbit inside $\RepDT{\alpha}$ whose triangle chain is regular has action coordinates $(\beta_1, \beta_2)=(2\pi/3, 4\pi/3)$. Moreover, the combination of solutions is \Romannum{4}$^\ast$-\Romannum{3}, with $\alpha^{(1)}=(16\pi/9, 16\pi/9, 16\pi/9, 4\pi/3)$ and $\alpha^{(2)}=(4\pi/3, 16\pi/9, 17\pi/9, 17\pi/9)$. A possible value for the pair $(\gamma_1,\gamma_2)$ of angle coordinates for $[\rho]$ is $(\pi, 0)$ according to Lemmas~\ref{lem:beta_i-for-orbits-of-type-III} and~\ref{lem:beta_i-for-orbits-of-type-IV}, giving the following triangle chain.

\begin{center}
\begin{tikzpicture}[font=\sffamily,decoration={
    markings,
    mark=at position 1 with {\arrow{>}}}]
    
\node[anchor=south west,inner sep=0] at (0,0) {\includegraphics[width=9cm]{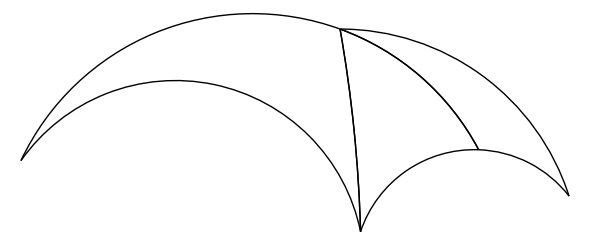}};

\begin{scope}
\fill (5.5,.22) circle (0.07) node[below]{$C_3=C_4$};
\fill (.3,1.3) circle (0.07) node[left]{$C_5$};
\fill (5.2,3.3) circle (0.07) node[above right]{$C_1=B_2$};
\fill (8.7,.8) circle (0.07) node[right]{$C_2$};
\fill (7.3,1.45) circle (0.07) node[below]{$B_1$};
\end{scope}

\fill (3.65,2.35) circle (0.07) node[below]{$X$};
\draw[dashed] (5.2,3.3) to[bend right=18] (3.65,2.35);

\draw[mauve, thick] (6.4,3.1) arc (-10:-40:.7) node[at end, right]{\small $\pi/9$};
\draw[mauve, thick] (8.1,1.3) arc (150:122:.6) node[midway, above left]{\small $\pi/9$};
\draw[mauve, thick] (5.9,.9) arc (50:95:.6) node[near start, above]{\small $\pi/9$};
\draw[mauve, thick] (5.45,1.5) arc (95:132:.7) node[at start, above left]{\tiny $\pi/18$};
\draw[thick] (4,2.2) arc (-30:60:.3) node[midway, right]{\small $\pi-\vartheta/2$};
\draw[apricot, thick] (5.3,2.7) arc (275:340:.54) node[near end, below]{\small $\pi/3$};
\end{tikzpicture}
\end{center}

The geodesic lines through $C_1$ and $C_2$, and through $C_4$ and $C_5$, intersect at $X$. The triangle with vertices $(C_2,C_3,X)$ has interior angle $\pi/9$ at $C_2$ and $\pi/9+\pi/18$ at $C_3=C_4$. This means that $X$ is the fixed point of the elliptic element $\rho((c_2c_3c_4)^{-1})$. Moreover, if $\vartheta$ denote the rotation angle of $\rho((c_2c_3c_4)^{-1})$, then $\angle C_3XC_2=\pi-\vartheta/2$. Now, the hyperbolic law of cosines applied to the triangle $(C_3=C_4,C_1,C_2)$ gives
\[
\cosh C_2C_3=\frac{\cos(4\pi/9)+\cos(\pi/9)^2}{\sin(\pi/9)^2}.
\]
When applied to the triangle $(C_2,X,C_3=C_4)$, the hyperbolic law of cosines further gives
\[
\cosh C_2C_3=\frac{-\cos(\vartheta/2)+\cos(\pi/9)\cos(\pi/6)}{\sin(\pi/9)\sin(\pi/6)}.
\]
After equalling the two expressions for $\cosh C_2C_3$ and simplifying, we obtain
\[
2\cos(\vartheta/2)=\frac{-1}{2\sin(\pi/9)}=\frac{-1}{2\cos(7\pi/18)}.
\]
This implies that $2\cos(\vartheta/2)\in\Q(\cos(\pi/18))$. According to Table~\ref{tab:rational-angles}, if $\vartheta$ was a rational multiple of $\pi$, then  $-1/(2\cos(7\pi/18))$ would belong to $\{2\cos(k\pi/18):k=0,\ldots,18\}$. We can check that this is not the case and conclude that $\vartheta$ is not a rational multiple of $\pi$. The Dehn twist along the curve $(c_2c_3c_4)^{-1}$ doesn't fix $[\rho]$ by \cm{Lemma~\ref{lem:fixed-points-Dehn-twists}}, showing that the mapping class group orbit of $[\rho]$ is infinite and providing the desired contradiction.

\subsubsection{Getting rid of $\alpha=(16\pi/9, 16\pi/9, 16\pi/9,16\pi/9,16\pi/9)$}\label{sec:third-bad-angle-vector}
The argument is similar to the previous case. A point $[\rho]$ in a finite mapping class group orbit inside $\RepDT{\alpha}$ has $(\beta_1,\beta_2)=(2\pi/3, 4\pi/3)$. The combination of solutions is \Romannum{4}$^\ast$-\Romannum{4}$^\ast$, with $\alpha^{(1)}=(16\pi/9, 16\pi/9, 16\pi/9, 4\pi/3)$ and $\alpha^{(2)}=(4\pi/3, 16\pi/9, 16\pi/9, 16\pi/9)$. Lemma~\ref{lem:beta_i-for-orbits-of-type-IV} says that the action coordinates of $[\rho]$ can be taken to be $(\gamma_1,\gamma_2)=(\pi,\pi)$, giving the following triangle chain.
\begin{center}
\begin{tikzpicture}[font=\sffamily,decoration={
    markings,
    mark=at position 1 with {\arrow{>}}}]
    
\node[anchor=south west,inner sep=0] at (0,0) {\includegraphics[width=9cm]{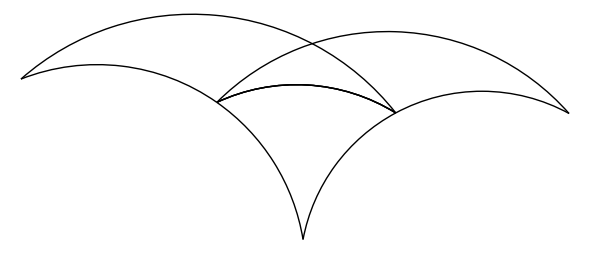}};

\begin{scope}
\fill (4.63,.22) circle (0.07) node[below]{$C_3$};
\fill (3.32,2.3) circle (0.07) node[below left]{$C_1=B_2$};
\fill (6.05,2.15) circle (0.07) node[below right]{$B_1=C_5$};
\fill (8.67,2.2) circle (0.07) node[right]{$C_2$};
\fill (0.37,2.67) circle (0.07) node[below]{$C_4$};
\end{scope}

\fill (4.77,3.22) circle (0.07) node[above]{$X$};

\draw[mauve, thick] (5.5,2.4) arc (160:130:.5) node[at end, right]{\small $\pi/9$};
\draw[mauve, thick] (3.8,2.5) arc (20:40:.5) node[at end, left]{\small $\pi/9$};
\draw[mauve, thick] (4.45,.9) arc (110:63:.5) node[midway, above]{\small $\pi/9$};
\draw[apricot, thick] (3.7,2.45) arc (20:-40:.4) node[midway, right]{\small $\pi/3$};
\draw[apricot, thick] (5.6,2.35) arc (170:230:.4) node[midway, left]{\small $\pi/3$};
\draw[thick] (4.6,3.15) arc (210:320:.2) node[midway, below]{\tiny $\pi-\vartheta/2$};
\end{tikzpicture}
\end{center}
The geodesic lines $C_4C_5$ and $C_1C_2$ intersect at $X$. The point $X$ coincides with the fixed points of $\rho((c_5c_1)^{-1})$. The rotation angle $\vartheta$ of $\rho((c_5c_1)^{-1})$ is also given by the relation $\angle C_5XC_1=\pi-\vartheta/2$. The hyperbolic law of cosines applied to the triangle $(C_1,C_3,C_5)$ gives
\[
\cosh C_1C_5=\frac{\cos(\pi/9)+\cos(\pi/3)^2}{\sin(\pi/3)^2}=\frac{4\cos(\pi/9)+1}{3}.
\]
If we apply it instead to the triangle $(C_1,X,C_5)$, we obtain
\[
\cosh C_1C_5=\frac{-\cos(\vartheta/2)+\cos(\pi/9)^2}{\sin(\pi/9)^2}.
\]
When combined, these two relations become
\[
2\cos(\vartheta/2)=\frac{1}{3}(3-2\cos(\pi/9)+4\cos(2\pi/9)).
\]
This implies that $2\cos(\vartheta/2)\in\Q(\cos(\pi/9))$. We can consult Table~\ref{tab:rational-angles} one more time to deduce that if $\vartheta$ was a rational multiple of $\pi$, then $(3-2\cos(\pi/9)+4\cos(2\pi/9))/3$ would belong to $\{0, 2\cos(\frac{k\pi}{9}):k=0,\ldots,9\}$. We can again check that this is not the case and conclude that $\vartheta$ is not a rational multiple of $\pi$. Since the Dehn twist along the curve $(c_5c_1)^{-1}$ doesn't fix $[\rho]$ by \cm{Lemma~\ref{lem:fixed-points-Dehn-twists}}, it shows that the mapping class group orbit of $[\rho]$ is infinite and gives the desired contradiction.

\subsubsection{The finite orbit for $\alpha=(4\pi/3,\theta,\theta,\theta,\theta)$}\label{sec:proof-thm-finite-orbits-n=5-1}
We start by considering the point $[\rho]$ of $\RepDT{\alpha}$ with action coordinates $(\beta_1,\beta_2)=(\pi,4\pi/3)$ and angle coordinates $(\gamma_1,\gamma_2)=(\pi,0)$. It's $\mathcal{B}$-triangle chain has the following ``hang-glider'' shape. It consists of three triangles with vertices $(C_1,C_2,B_1)$, $(B_1,C_3,B_2)$, and $(B_2,C_4,C_5)$ arranged as below.
\begin{center}
\begin{tikzpicture}[font=\sffamily,decoration={
    markings,
    mark=at position 1 with {\arrow{>}}}]
    
\node[anchor=south west,inner sep=0] at (.2,.25) {\includegraphics[width=3.58cm]{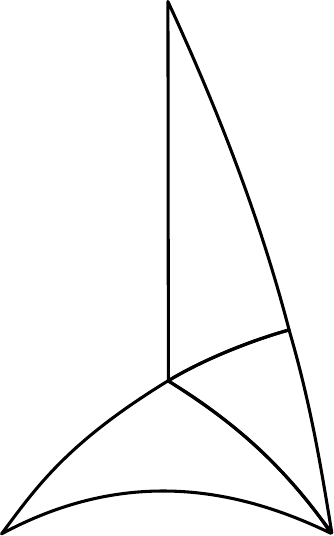}};

\begin{scope}
\fill (2.02,1.95) circle (0.07) node[above left]{$C_1=B_2$};
\fill (2.02,5.95) circle (0.07) node[below left]{$C_2$};
\fill (3.33,2.45) circle (0.07) node[right]{$B_1$};
\fill (3.75,.3) circle (0.07) node[below]{$C_3=C_4$};
\fill (0.25,.3) circle (0.07) node[below]{$C_5$};
\end{scope}
\end{tikzpicture}
\end{center}
We can see from the configuration of the triangle chain that the image of $\rho$ is conjugate to a rotation triangle group $D(2,3,\overline{\theta})$. Recall that the image of $\rho$ is generated by the rotation of angle $2\pi/3$ around $C_1$ and the rotations of angle $\theta$ around the other exterior vertices. The group $D(2,3,\overline{\theta})$ can be realized as the rotation triangle group of the triangle $(B_1,C_3,B_2)$. Since $C_1$ is the reflection of $C_3$ through the geodesic line $B_1B_2$, and $C_5$ is the reflection of $C_2$, this shows that image of $\rho$ is contained in $D(2,3,\overline{\theta})$. Conversely, the image of $\rho$ contains the rotation of angle $\pi$ around $B_i$ and the rotation of angle $2\pi/3$ around $B_2$, so it's actually equal to $D(2,3,\overline{\theta})$.

We now claim that the mapping class group orbit of $[\rho]$---the \emph{hang-glider orbit}---is finite and has length 9. As in the proof of Theorem~\ref{thm:existence-finite-orbit-n=6}, finiteness is immediate from Lemma~\ref{lem:pullback-orbits-are-finite} because $\rho$ is a pullback representation by construction. We rely on the algorithm of Appendix~\ref{apx:algo-orbit} to determine the coordinates of every point in the hang-glider orbit which can be found in Table~\ref{tab:hang-glider-orbit} in Appendix~\ref{app:tables}. There are exactly 9 orbit points.

It remains to prove that the hang-glider orbit is the unique finite orbit in $\RepDT{\alpha}$. Following the same strategy as in the proof of the uniqueness statement in Theorem~\ref{thm:existence-finite-orbit-n=6}, it's enough to prove that any finite orbit contains a point whose $\mathcal{B}$-triangle chain is regular. Moreover, any finite orbit is the hang-glider orbit for some (maybe different) pants decomposition of $\Sigma$ and, in particular, a finite orbit is always of length 9. So, if all the points of a finite orbit have a singular $\mathcal{B}$-triangle chain, then there exists an orbit point $[\phi]$ with only one degenerate triangle (actually, there are at least six such orbit points). Arguing as in the proof of Theorem~\ref{thm:existence-finite-orbit-n=6}, we can see that if the third triangle in the chain of $[\phi]$ is degenerate, then the $\mathcal{B}$-triangle chain of $\tau_{3,4}.[\phi]$ is either regular or has $B_1=C_3=B_2$, meaning that its second triangle is degenerate. Similarly, if the second triangle in the chain of $[\phi]$ is degenerate, then the $\mathcal{B}$-triangle chain of $\tau_{2,3}.[\phi]$ is either regular or has $C_1=C_2=B_1$, meaning that its first triangle is degenerate. In conclusion, if all the points in the mapping class group orbit of $[\phi]$ have a singular $\mathcal{B}$-triangle chain, then there always exists one of them, say $[\phi']$, for which only the first triangle is degenerate. Now, in order for the $\mathcal{B}$-triangle chain of $\tau_{2,3}.[\phi']$ to be singular, $[\phi']$ must have $\beta_2=\alpha_3$ by the same geometric arguments as in the proof of Theorem~\ref{thm:existence-finite-orbit-n=6}. 

The angle $\alpha_3$ is either equal to $\theta$ or to $4\pi/3$. We can see that $\alpha_3\neq \theta$. Assume for the sake of contradiction that $\alpha_3=\theta$ and note that at least one of the two angles $\alpha_4$ or $\alpha_5$ is equal to $\theta$ too. However, since $\beta_2$ satisfies $\beta_2\leq \alpha_4+\alpha_5-2\pi$ by~\eqref{eq:inequalities-polytope-beta}, we obtain a contradiction. So, we must have $\alpha_3=4\pi/3$. In that case, $[\phi']$ can be seen as a point in the DT component of the sub-sphere $\Sigma^{(2)}$ of $\Sigma$ (following the notation of Section~\ref{sec:restrecting-beta_i-to-a-finite-set}) with peripheral angles $(2\theta-2\pi,4\pi/3,\theta,\theta)$. The $\PMod(\Sigma^{(2)})$-orbit of $[\phi']$ in the DT component of $\Sigma^{(2)}$ is of course finite and thus of Type~\Romannum{3}. Using Lemma~\ref{lem:beta_i-for-orbits-of-type-III}, we conclude that the angle coordinate $\gamma_2$ of $[\phi']$ belongs to $\{\pi/3,\pi, 5\pi/3\}$. In other words, as a point of the DT component of the full sphere $\Sigma$ again, $[\phi']$ has action-angle coordinate $(\beta_1,\beta_2)=(4\pi-2\theta,4\pi/3)$ and $\gamma_2\in \{\pi/3,\pi, 5\pi/3\}$. It turns out that all three points of $\RepDT{\alpha}$ with coordinates $(\beta_1,\beta_2)=(4\pi-2\theta,4\pi/3)$ and $\gamma_2\in \{\pi/3,\pi, 5\pi/3\}$ belong to the hang-glider orbit associated to the pants decomposition $\mathcal{B}$, and thus so does $[\phi']$. This finishes the proof that the hang-glider orbit is the unique finite orbit in $\RepDT{\alpha}$.

\subsubsection{The finite orbit for $\alpha=(2\theta-2\pi,\theta,\theta,\theta,\theta)$}\label{sec:proof-thm-finite-orbits-n=5-2}
We deal with this case as we did in Section~\ref{sec:proof-thm-finite-orbits-n=5-1}. Consider the point $[\rho]$ with coordinates $(\beta_1,\beta_2)=(\pi,4\pi/3)$ and $(\gamma_1,\gamma_2)=(\pi/2, 0)$. Its $\mathcal{B}$-triangle chain is shaped as a ``sand clock'' and is made of three triangles $(C_1,C_2,B_1)$, $(B_1,C_3,B_2)$, and $(B_2,C_3,C_4)$ with $C_3=C_4$.
\begin{center}
\begin{tikzpicture}[font=\sffamily,decoration={
    markings,
    mark=at position 1 with {\arrow{>}}}]
    
\node[anchor=south west,inner sep=0] at (0.25,0.2) {\includegraphics[width=6.55cm]{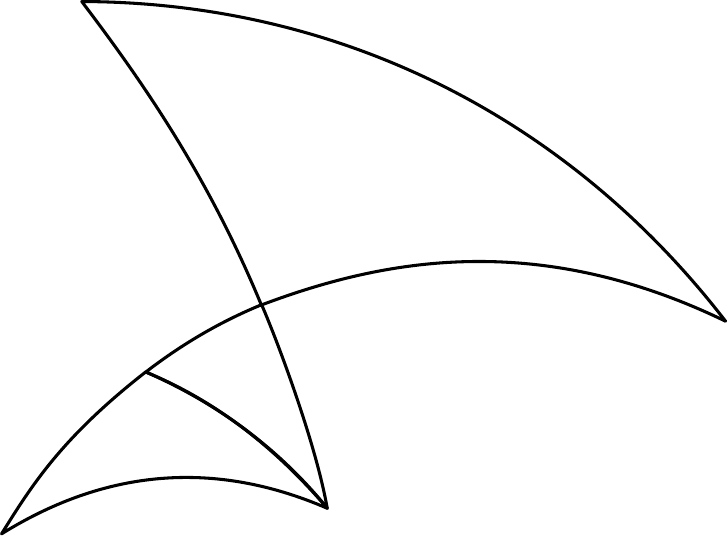}};

\begin{scope}
\fill (1,5) circle (0.07) node[left]{$C_1$};
\fill (6.75,2.15) circle (0.07) node[right]{$C_2$};
\fill (2.6,2.3) circle (0.07) node[above left]{$B_1$};
\fill (3.2,.5) circle (0.07) node[below]{$C_3=C_4$};
\fill (1.6,1.7) circle (0.07) node[above left]{$B_2$};
\fill (0.3,.25) circle (0.07) node[below]{$C_5$};
\end{scope}
\end{tikzpicture}
\end{center}
Similar geometric considerations as in Section~\ref{sec:proof-thm-finite-orbits-n=5-1} show that the image of $\rho$ is conjugate to a rotation triangle group $D(2,3,\overline{\theta})$. The orbit is finite by Lemma~\ref{lem:pullback-orbits-are-finite} and consists of 12 points. When we compute the coordinates of every other point in the orbit of $[\rho]$---the \emph{sand clock orbit}---using the algorithm of Appendix~\ref{apx:algo-orbit}, we indeed observe that the orbit consists of 12 points. All their coordinates can be found in Table~\ref{tab:sand-clock-orbit} in Appendix~\ref{app:tables}. An analogous argument as in Section~\ref{sec:proof-thm-finite-orbits-n=5-1} shows that the sand clock orbit is the unique finite orbit in this DT component.

\subsubsection{The finite orbit for $\alpha=(12\pi/7,12\pi/7,12\pi/7,12\pi/7,12\pi/7)$}\label{sec:proof-thm-finite-orbits-n=5-3}
This is the last angle vector to consider. The $\mathcal{B}$-triangle chain of the point $[\rho]$ with coordinates $(\beta_1,\beta_2)=(2\pi/3,8\pi/7)$ and $(\gamma_1,\gamma_2)=(\pi/3, 4\pi/7)$ has the shape of a ``bat''. It consists of three triangles $(C_1,C_2,B_1)$, $(B_1,C_3,B_2)$, and $(B_2,C_3,C_4)$.
\begin{center}
\begin{tikzpicture}[font=\sffamily,decoration={
    markings,
    mark=at position 1 with {\arrow{>}}}]
    
\node[anchor=south west,inner sep=0] at (0.25,0.25) {\includegraphics[width=6.5cm]{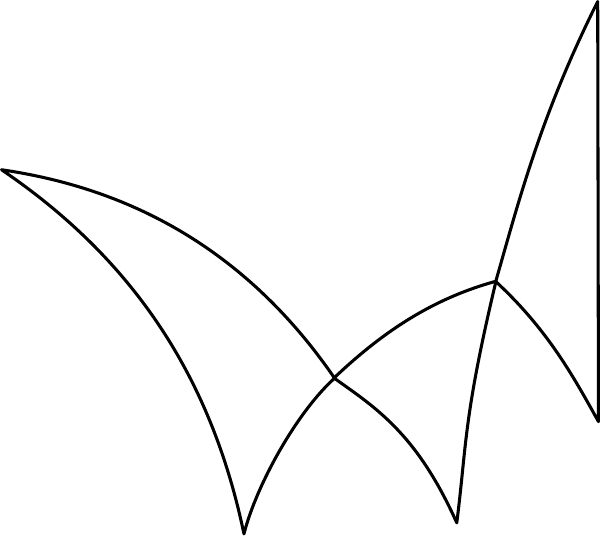}};

\begin{scope}
\fill (6.75,6) circle (0.07) node[right]{$C_1$};
\fill (6.75,1.5) circle (0.07) node[right]{$C_2$};
\fill (5.65,3) circle (0.07) node[above left]{$B_1$};
\fill (5.22,.4) circle (0.07) node[below]{$C_3$};
\fill (3.85,2) circle (0.07) node[left]{$B_2$};
\fill (2.92,.3) circle (0.07) node[below]{$C_4$};
\fill (0.27,4.23) circle (0.07) node[left]{$C_5$};
\end{scope}
\end{tikzpicture}
\end{center}
The bat triangle chain fits well onto a tessellation of the hyperbolic plane by triangles with interior angles $(\pi/2,\pi/3,\pi/7)$ \cm{(in the sense of Example~\ref{ex:discrete-DT})}. This observation tells us that the image of $\rho$ is conjugate to the rotation triangle group $D(2,3,7)$, hence discrete. The orbit of $[\rho]$ is therefore finite by Corollary~\ref{cor:DT-discrete-implies-finite-orbit} (note that $\rho$ is also a pullback representation).

We compute the coordinates of every other point in the orbit of $[\rho]$---the \emph{bat orbit}---using the algorithm of Appendix~\ref{apx:algo-orbit} as we did in Sections~\ref{sec:proof-thm-finite-orbits-n=5-1} and~\ref{sec:proof-thm-finite-orbits-n=5-2}. Doing so, we can confirm that the bat orbit is made of 105 points. Their coordinates are provided in Table~\ref{tab:bat-orbit} in Appendix~\ref{app:tables}. An analogous argument as in Section~\ref{sec:proof-thm-finite-orbits-n=5-1} shows that the bat orbit is the unique finite orbit in its DT component.

\section{Proof of Tykhyy's Conjecture}\label{chap:tykhyy-conjecture-general}
\subsection{Overview}
In this section we explain how Theorems~\ref{thm:angle-vector-alpha-with-finite-orbits-n=6} \&~\ref{thm:existence-finite-orbit-n=6}, and Theorem~\ref{thm:no-finite-orbit-for-n-geq-7} cover the last remaining open cases in Tykhyy's Conjecture (Conjecture~\ref{conj:tykhyy}) and complete its proof. Section~\ref{sec:statement-conjecture} contains the statement of the conjecture. A recap on existing partial achievements, along with a reduction of the proof to the case of DT representations are presented in Section~\ref{sec:corlette-simpson}. The proof of Tykhyy's Conjecture can be found in Section~\ref{sec:proof-tykhyy-conjecture}.

\subsection{Tykhyy's Conjecture}\label{sec:statement-conjecture}
Tykhyy conjectured that every representation of a sphere $\Sigma$ with $n\geq 3$ punctures whose conjugacy class lies in a \emph{finite} mapping class group orbit belong to certain list that can be found in~\cite[Section~11]{tykhyy}. We reformulate Tykhyy's Conjecture by regrouping the representations into different sub-families and by often mentioning only one representative of each conjugacy class. For simplicity, we fix a geometric presentation of $\pi_1\Sigma$ with generators $c_1,\ldots,c_n$.

\cm{Before stating the conjecture, there are three operations that trivially produce new finite orbits from a given one. Let $\rho\colon\pi_1\Sigma\to\SL_2\C$ be a representation of a punctured sphere whose conjugacy class lies in a finite orbit.}
\begin{itemize}
\item \cm{\emph{Adding trivial peripheral loops.} If we add punctures on $\Sigma$ to create a new punctured sphere $\Sigma'$ and define $\rho'\colon\pi_1\Sigma'\to\SL_2\C$ from $\rho$ by mapping the new peripheral loops to $\pm \id$, changing the sign of the image of an old peripheral loop if necessary to make sure $\rho'$ is a representation, then the conjugacy class of $\rho'$ lies in a finite orbit.}
\item \cm{\emph{Reduction.} Conversely, by multiplying two consecutive peripheral monodromies, then one turns $\rho$ into a representation of a sphere with one puncture less than $\Sigma$. The conjugacy class of this newly built representation also lies in a finite mapping class group orbit.}
\item \cm{\emph{Reordering punctures.} Since we are studying the action of the \emph{pure} mapping class group (Definition~\ref{defn:pure-mcg}), acting on a representation $\pi_1\Sigma\to\SL_2\C$ by pre-composition by an automorphism of $\pi_1\Sigma$ that permutes punctures, like the automorphisms $\sigma_1,\ldots,\sigma_{n-1}$ given by
\[
\sigma_j(c_i)=\begin{cases}
c_i & \text{if } i\neq j, j+1,\\
c_jc_{j+1}c_j^{-1} & \text{if } i=j,\\
c_j & \text{if } i=j+1,\\
\end{cases}
\]
typically produces different finite orbits.}

\end{itemize}

\cm{We're now ready to state the conjecture. We consider the following eight kinds of representations that form \emph{Tykhyy's list}. The first kinds are those whose image is \emph{not} Zariski dense in $\SL_2\C$. There are examples for every $n\geq 3$.}
\begin{enumerate}
\item Representations with finite image, \cm{which is either cyclic or} contained in the dihedral, tetrahedral, octahedral, or icosahedral finite subgroups of $\SL_2\C$ (Section~\ref{sec:finite-image}).
\item Representations with \emph{infinite} image contained in the subgroup of diagonal and anti-diagonal matrices (also known as the infinite dihedral subgroup) of $\SL_2\C$. These representations \cm{are either diagonal or} map exactly two of the generators $c_1,\ldots,c_n$ to anti-diagonal matrices (Section~\ref{sec:infinite-dihedral}).
\item Representations with \cm{\emph{infinite}} image contained in the subgroup of upper triangular matrices (Section~\ref{sec:upper-triangular}).
\end{enumerate}
\cm{Next are the representations with Zariski dense images and non-trivial peripheral monodromies.}
\begin{enumerate}
\setcounter{enumi}{3}
\item For $n=6$: a 1-parameter family of representations $\rho_\theta\colon\pi_1\Sigma\to\SL_2\C$ \cm{parametrized by $\theta\in\C$} and given by
\begin{equation*}
\begin{array}{lll}
\rho_\theta(c_1)=\begin{pmatrix}
e^{i\theta/2} & 0\\
1 & e^{-i\theta/2}
\end{pmatrix}, &\quad 
\rho_\theta(c_2)=\begin{pmatrix}
e^{-i\theta/2} & -1\\
0 & e^{i\theta/2}
\end{pmatrix}, &\quad
\rho_\theta(c_3)=\begin{pmatrix}
-e^{-i\theta/2} & 1\\
0 & -e^{i\theta/2}
\end{pmatrix}, \\
\rho_\theta(c_4)=\begin{pmatrix}
e^{i\theta/2} & 0\\
1 & e^{-i\theta/2}
\end{pmatrix}, &\quad 
\rho_\theta(c_5)=\begin{pmatrix}
e^{-i\theta/2} & -1\\
0 & e^{i\theta/2}
\end{pmatrix}, &\quad
\rho_\theta(c_6)=\begin{pmatrix}
e^{-i\theta/2} & -1\\
0 & e^{i\theta/2}
\end{pmatrix},
\end{array}
\end{equation*}
\cm{and all the representations obtained from $\rho_\theta$ by reordering the punctures.}
\item For $n=5$: first, representations obtained by reduction of the representation $\rho_\theta$ from~(4) \cite[Equations~(17) \& (18)]{tykhyy}. There are also four exceptional representations: one with purely parabolic peripheral monodromies~\cite[Equation~(20)]{tykhyy} and the other three with purely elliptic peripheral monodromies where eigenvalues are seventh roots of unity~\cite[Equation~(19)]{tykhyy} (one is valued in $\SL_2\R$, the other two are valued in $\SU(2)$).
\item For $n=4$: all the representations in the list of Boalch and Lisovyy--Tykhyy~\cite{boalch, LT}.
\item For $n=3$: any representation.
\end{enumerate}
\cm{Finally, Tykhyy's list contains all the representations with Zariski dense images that map some peripheral loops to $\pm\id$.}
\begin{enumerate}
\setcounter{enumi}{7}
\item For any $n\geq 3$: representations that map some of the generators $c_1,\ldots,c_n$ to $\pm\id$ and can be reduced to one of the representations of types~(4)--(7). In particular, when $n\geq 7$, the only Zariski dense representations in Tykhyy's list are those that send at least $n-6$ generators to~$\pm\id$ and can be reduced to one of the representations listed previously.
\end{enumerate} 

\begin{conj}[\cite{tykhyy}]\label{conj:tykhyy}
Let $\Sigma$ be a sphere with $n\geq 3$ punctures and $\rho\colon\pi_1\Sigma\to\SL_2\C$ be a representation. If the conjugacy class $[\rho]$ belongs to a finite mapping class group orbit, then orbit of $[\rho]$ contains the conjugacy class of a representation from~\textnormal{(1)--(8)}.
\end{conj}

\begin{rem}\label{rem:complete-list-of-orbits}
\cm{We point out that Conjecture~\ref{conj:tykhyy} only gives an explicit list of all possible finite orbits in the case of Zariski dense representations. Complete classifications for representations into the subgroup of upper-triangular matrices and into the infinite dihedral subgroup of $\SL_2\C$ have been achieved already, as related in Section~\ref{sec:non-zariski-dense-orbits}. On the other hand, it would require more effort to classify all the conjugacy classes of representations with finite image according to their mapping class group orbits. As related in Section~\ref{sec:finite-image}, this task already required substantial work in the case of 4-punctured spheres. Our methods cannot be directly used to address that question.}
\end{rem}

The finite orbits that we encountered in Chapter~\ref{sec:classification} fit well in the scope of Conjecture~\ref{conj:tykhyy}. The jester's hat orbit from Theorem~\ref{thm:existence-finite-orbit-n=6} contains the conjugacy class of the representation $\rho_\theta$ of type~(4) \cm{when $\theta\in(5\pi/3,2\pi)$, as we'll explain after Table~\ref{tab:UFO-orbit} in Appendix~\ref{app:tables} where we also provide an explicit conjugate of $\rho_\theta$ with values in $\SL_2\R$.} The hang-glider orbit (Section~\ref{sec:proof-thm-finite-orbits-n=5-1}) corresponds to a reduction of the representation $\rho_\theta$ from~(4) obtained by multiplying the first two generators, whereas the sand clock orbit (Section~\ref{sec:proof-thm-finite-orbits-n=5-2}) corresponds to a reduction of $\rho_\theta$ where the second and third generators are multiplied together, \cm{both in the case $\theta\in(5\pi/3,2\pi)$.} The bat orbit (Section~\ref{sec:proof-thm-finite-orbits-n=5-3}) is one of the four exceptional orbits of type~(5).

Some cases of Conjecture~\ref{conj:tykhyy} have already been treated. We mentioned the case of representations with finite image in Section~\ref{sec:finite-image}, as well as the other representations with non-Zariski dense image in Section~\ref{sec:non-zariski-dense-orbits}, including the work of Cousin--Moussard~\cite{cousin-moussard} that covers the upper triangular case. As we already cited a number of times, a complete classification is available for $n=4$~\cite{LT} \cm{and a computer-aided classification exists for $n=5$~\cite{tykhyy}.} At this point, what remains to do in order to prove Conjecture~\ref{conj:tykhyy} is to identify all the Zariski dense representations $\rho\colon\pi_1\Sigma\to\SL_2\C$ that give rise to finite mapping class group orbits when $\Sigma$ is a sphere with $6$ punctures or more. There have been two major contributions towards this goal: Diarra classified the ones of pullback type (Definition~\ref{def:pullback-orbits})~\cite{diarra-pull-back} and Lam--Landesman--Litt those with at least one peripheral monodromy of infinite order~\cite{lll}. What remains is the case of Zariski dense representations $\rho\colon\pi_1\Sigma\to\SL_2\C$ with all peripheral monodromies of finite order. As we're about to see (Proposition~\ref{prop:Zd+non-pullback-implies-DT}), these representations will always be Galois conjugate to DT representations. This observation and our good understanding of DT representations will be enough to conclude the proof of Conjecture~\ref{conj:tykhyy} in Theorem~\ref{thm:tykhyy}.

As a corollary of Theorem~\ref{thm:tykhyy}, we obtain an alternative formulation of Conjecture~\ref{conj:tykhyy} for Zariski dense representations in terms of Katz's middle convolution (which generalizes Okamoto transformations). It was suggested to us by Litt and comes as an answer to~\cite[Question~2.4.1]{litt}. Before stating it, let's first recall what Lam--Landesman--Litt proved. They showed that a Zariski dense representation $\rho\colon\pi_1\Sigma\to\SL_2\C$ with at least one peripheral monodromy of infinite order gives rise to a finite orbit if and only if $\rho$ is a pullback representation or if it's obtained via Katz's middle convolution from a representation with image in a finite complex reflection group~\cite[Corollary~1.1.7]{lll}. We also know from our discussion in Section~\ref{thm:LT-classification} that this dichotomy fails for some of the exceptional orbits in the case of 4-punctured spheres: all Okamoto transforms of the Klein solution (Type~8) and the three elliptic 237 solutions (Types~32---34) have infinite image and some are not of pullback type. According to Vayalinkal~\cite{amal}, they're also not obtainable via Katz's middle convolution from finite complex reflection groups of higher rank.

\begin{cor}\label{cor:litt-version}
Let $\Sigma$ be a sphere with $n\geq 4$ punctures. If $\rho\colon\pi_1\Sigma\to\SL_2\C$ is a Zariski dense representation such that the mapping class group orbit of $[\rho]$ is finite, then $\rho$ is of one of the following three types:
\begin{enumerate}
\item $\rho$ is a pullback representation.
\item $\rho$ is obtained via Katz's middle convolution from a representation with image in a finite complex reflection group.
\item $\rho$ is obtained via Katz's middle convolution from a pullback representation.
\end{enumerate}
\end{cor}
\begin{proof}
The statement holds for 4-punctured spheres by~\cite{LT} and when $\rho$ has a peripheral monodromy of infinite order by~\cite{lll}. We'll see in the proof of Theorem~\ref{thm:tykhyy} that when $n\geq 5$ and $\rho$ only has finite order peripheral monodoromy, then $\rho$ is Galois conjugate to a DT representation of pullback type. Since Galois conjugation and applying the pullback construction of Section~\ref{sec:pullback-orbits} are two commuting operations, this finishes the proof.
\end{proof}

\subsection{The Corlette--Simpson alternative}\label{sec:corlette-simpson}
A particularly useful result that helps classifying finite mapping class group orbits is a theorem by Corlette--Simpson~\cite{corlette-simpson}. It's about Zariski dense
representations of the fundamental group of a quasi-projective variety $X$ into $\SL_2\C$ with \emph{quasi-unipotent monodromy} at infinity, meaning that all eigenvalues are roots of unity.
It's been extended by Loray--Pereira--Touzet~\cite{loraypereiratouzet} to omit
the quasi-unipotent monodromy assumption. We formulate both results using a different vocabulary than what can be found in the original papers. We would like to thank Tholozan for explaining us these alternative formulations and the subsequent arguments involving Galois conjugations and variations of Hodge structures, which originate from Lam--Landesman--Litt's paper~\cite{lll}.
\begin{thm}[\cite{corlette-simpson}]\label{thm:corlette-simpson}
	Let $X$ be a quasi-projective manifold and
	let $\rho\colon\pi_1 X\rightarrow\SL_2\C$ be a Zariski dense representation
	with quasi-unipotent monodromy at infinity. Then one of the following holds.
	\begin{itemize}
	\item The representation $\rho$ projectively \emph{factorizes through an orbisurface}, meaning that $\rho$ factorizes through a morphism $f\colon X\rightarrow Y$ to an orbisurface $Y$.  
	
	\item \cm{The representation $\rho$ is \emph{integral and preserves a Hermitian form}, meaning that, up to conjugating $\rho$, it's valued in $\SU(h,\mathcal{O}_L)$,
	where $L\subset \C$ is a totally imaginary extension of a totally real number field, $\mathcal{O}_L$ is its ring of integers, and $h$ is a Hermitian form on $L^2$.}
	In that case, $\rho$ is rigid and the associated local system supports a variation of Hodge structures.
	\end{itemize}
\end{thm}
Loray--Pereira--Touzet showed that without the quasi-unipotent monodromy assumption,
the only possibility is to factorize through an orbisurface.
\begin{thm}[\cite{loraypereiratouzet}]\label{thm:loray-pereira-touzet}
	Let $\rho\colon\pi_1 X\rightarrow\SL_2\C$ be a Zariski dense representation
	which is not quasi-unipotent at infinity. Then $\rho$ projectively
	factors through an orbisurface.
\end{thm}
\cm{Note that the second case in Theorem~\ref{thm:corlette-simpson} arises in two flavours according to the signature of $h$ seen as a Hermitian form on $\C^2$.
If $h$ is positive (or negative) definite, then $\SU(h)\cong \SU(2)$ is compact and so, up to conjugating $\rho$, its image is contained in $\SU(2)$.
Otherwise, if $h$ is of signature $(1,1)$, then $\SU(h)\cong \SU(1,1)$ and, up to conjugating $\rho$, its image is contained in $\SL_2\R$.}

The relation between finite mapping class group orbits and Theorems~\ref{thm:corlette-simpson} \&~\ref{thm:loray-pereira-touzet} is quite well-known. It relies on the following observation.  Let $\mathcal{M}_{0,n}$ denote the moduli space of $n$-punctured spheres. As we already mentioned in the proof of Lemma~\ref{lem:pullback-orbits-are-finite}, the pure mapping class group of a sphere $\Sigma$ with $n$ punctures is isomorphic to the fundamental group $\pi_1\mathcal{M}_{0,n}$. Now, if $\rho\colon\pi_1\Sigma \rightarrow\SL_2\C$ is a representation of an $n$-punctured sphere $\Sigma$ whose conjugacy class $[\rho]$ has finite mapping class group orbit, then the stabilizer $\Gamma\subset \pi_1\mathcal{M}_{0,n}\cong\PMod(\Sigma)$ of $[\rho]$ is a subgroup of finite index. Its pre-image $\widetilde \Gamma\subset \pi_1\mathcal{M}_{0,n+1}$ under the morphism $\pi_1\mathcal{M}_{0,n+1}\to\pi_1\mathcal{M}_{0,n}$ from the Birman exact sequence (see e.g.~\cite[Theorem~4.6]{mcg-primer}) has finite index as well. We can realize $\pi_1\mathcal{M}_{0,n+1}$ as a subgroup of $\Aut^\star(\pi_1\Sigma)$ (using the notation introduced in Section~\ref{sec:character-varieties}) from the following commutative diagram, in which each row is exact.
\begin{center}
\begin{tikzcd}
1\arrow[r]&
\pi_1\Sigma\arrow[r]\arrow[d,leftrightarrow]&
\Aut^\star(\pi_1\Sigma) \arrow[r]&
\Out^\star(\pi_1\Sigma) \arrow[r]&1\\
1\arrow[r]&\pi_1\Sigma\arrow[r]&\pi_1\mathcal{M}_{0,n+1}\arrow[r]\arrow[hook, u]&
\pi_1\mathcal{M}_{0,n}\arrow[r]\arrow[hook, u]&1
\end{tikzcd}
\end{center}
This produces a representation $\tilde\rho\colon\widetilde\Gamma\to\mathrm{PSL}_2\C$ defined by sending $f\in \widetilde\Gamma$ to the unique element $g\in \mathrm{PSL}_2\C$ such that $\rho\circ f=g\rho g^{-1}$, where we think of $f$ as an automorphism of $\pi_1\Sigma$. The element $g$ exists because the projection of $f$ inside $\Gamma\subset \pi_1\mathcal{M}_{0,n}$ fixes $[\rho]$, and $g$ is unique because $\mathrm{PSL}_2\C$ has trivial center. The representation $\tilde\rho\colon\widetilde\Gamma\to\mathrm{PSL}_2\C$ extends $\rho$ in the sense that the image of the morphism $\pi_1\Sigma\to \pi_1\mathcal{M}_{0,n+1}$ from the Birman exact sequence is a subgroup of $\widetilde\Gamma$ and the composition $\pi_1\Sigma\to\widetilde\Gamma\to \mathrm{PSL}_2\C$ is the projectivization of $\rho$. In other words, the following diagram commutes.
\begin{center}
\begin{tikzcd}
\pi_1\Sigma\arrow[hook, d]\arrow[r, "\rho"] & \SL_2\C\arrow[d]\\
\widetilde\Gamma\arrow[r, "\tilde\rho"]& \mathrm{PSL}_2\C
\end{tikzcd}
\end{center}
This is easy to see. If $f\in\widetilde\Gamma$ is the image of $\gamma\in\pi_1\Sigma$, then $\rho\circ f(x)=\rho(\gamma x \gamma^{-1})=\rho(\gamma)\rho(x)\rho(\gamma)^{-1}$.

\cm{The group $\Gamma$ is isomorphic to the fundamental group of some (possibly ramified) finite order covering $Y$ of $\mathcal{M}_{0,n}$. Furthermore, $\widetilde{\Gamma}$ is the fundamental group of the pullback $X$ of $\mathcal{M}_{0,n+1}$ by $Y\to \mathcal{M}_{0,n}$, which turns $X$ into an $n$-punctured sphere bundle over $Y$.} 
This procedure can be described precisely in more geometric terms and it's possible to lift $\tilde\rho\colon \widetilde\Gamma\cong\pi_1 X\to \mathrm{PSL}_2\C$ to a linear representation $\pi_1 X\to\SL_2\C$ to which Theorems~\ref{thm:corlette-simpson} \&~\ref{thm:loray-pereira-touzet} apply. We refer the reader to~\cite[Section~1.1]{lll} for a precise statement and to~\cite[Chapter~2]{LL} for details and proofs. 

\begin{lem}\label{lem:Zd+non-pullback-implies-values-in-integers}
If $\rho\colon\pi_1\Sigma\rightarrow\SL_2\C$ is a Zariski dense representation such that the mapping class group \cm{orbit} of $[\rho]$ is finite and \emph{not} of pullback type, then \cm{up to conjugating $\rho$}, it is valued in the ring of integers of a number field $L$ and the associated local system supports a variation of Hodge structures for any choice of complex structure on $\Sigma$.
\end{lem}
\begin{proof}
The procedure described above produces a Zariski dense representation $\tilde\rho\colon\pi_1 X\to\SL_2\C$ which extends $\rho$ and is of two possible kind according to Theorems~\ref{thm:corlette-simpson} \&~\ref{thm:loray-pereira-touzet}. Either it factorizes through an orbisurface and the orbit of $[\rho]$ is actually of pullback type, or it has quasi-unipotent monodromy at infinity and, \cm{up to conjugation}, it's valued in $\SU(h,\mathcal{O}_L)$ for some number field $L$ and Hermitian metric $h$. \cm{In the latter case, the local system associated to $\widetilde\rho\colon\pi_1 X\to\SL_2\C$ supports a variation of Hodge structures. Restricting to the fibres of the $n$-punctured sphere bundle $X\to Y$, we conclude that the local system associated to $\rho\colon\pi_1\Sigma\rightarrow\SL_2\C$ supports a variation of Hodge structures for any choice of complex structures on $\Sigma$.}
\end{proof}

\begin{defn}\label{defn:universal-CVHS}
\cm{When the local system associated to $\rho\colon\pi_1\Sigma\to\SL_2\C$ supports a variation of Hodge structures for \emph{any} choice of complex structure on $\Sigma$, or in short when $\rho$ comes from a variation of Hodge structures for any choice of complex structure on $\Sigma$, then we say that $\rho$ is a \emph{universal variation of Hodge structures}.}
\end{defn}
The universal variations of Hodge structures from Lemma~\ref{lem:Zd+non-pullback-implies-values-in-integers} come in two flavours, depending
on the signature of $h$. If $h$ is of signature $(2,0)$ or $(0,2)$, then $\SU(h)\cong \SU(2)$ is compact and any representation valued in $\SU(h)$ is a universal variation of Hodge structures (the corresponding Higgs fields always vanish).
If $h$ is of signature $(1,1)$, then $\SU(h)\cong\SL_2\R$ and $\rho$ is a universal variation of
Hodge structures if and only if for every complex structure on $\Sigma$, there exists a holomorphic or anti-holomorphic $\rho$-equivariant map
$\widetilde\Sigma\to \HH$ \cm{which is non-constant}. As shown by Deroin--Tholozan, DT representations \cm{(Definition~\ref{defn:DT-representation})} are universal variations of Hodge structures \cite[Theorem~5]{deroin-tholozan}. \cm{Moreover, a representation $\rho\colon\pi_1\Sigma\to\psl$ with non-trivial elliptic peripheral monodromies which is a universal variation of Hodge structures is either a DT representation or an abelian representation~\cite[discussion after Theorem~5]{deroin-tholozan}. We'll come back to this correspondence in Proposition~\ref{prop:universal-CVHS-implies-DT}.}

\subsubsection{The Galois action}
When a representation $\rho\colon\pi_1\Sigma\to\SL_2\C$ is valued in $\SU(h,\mathcal{O}_L)$ for some number field $L$ (as it's the case in the second alternative in Theorem~\ref{thm:corlette-simpson}), there is an action of the Galois group $\mathrm{Gal}$ of $L$ over $\Q$ on $\rho$ by Galois conjugation. In other words, given $\sigma\in\rm{Gal}$, there is a representation $\rho_\sigma\colon\pi_1\Sigma\to\SL_2\C$ defined by $\rho_\sigma(\gamma)=\sigma(\rho(\gamma))$.
\begin{lem}\label{lem:Galois-conjugate-bijection-orbits}
The mapping class group orbit of $[\rho]$ is in bijection with the mapping class group orbit of  $[\rho_\sigma]$.
\end{lem}
\begin{proof}
The Galois group acts by post-composition while the mapping class
group acts by pre-composition. Both actions thus commute and the mapping class group orbits are identified via the action of $\sigma$.
\end{proof}
The Galois conjugate $\rho_\sigma$ is valued in $\SU(h^\sigma,\mathcal{O}_L)$ for some Hermitian metric $h^\sigma$. For a fixed $\rho$, if the Hermitian metrics $h^\sigma$ have signature $(0,2)$ or $(2,0)$ for every $\sigma\in\rm{Gal}$, then $\rho$ has finite image by the Borel--Harish-Chandra Theorem. If some Hermitian metric $h^\sigma$ has signature $(1,1)$, then $\rho_\sigma$ is valued in $\SU(1,1)\cong\SL_2\R$.

\begin{lem}\label{lem:Zd+non-pullback-implies-values-in-SL2R}
If $\rho\colon\pi_1\Sigma\rightarrow\SL_2\C$ is a Zariski dense representation such that the mapping class group \cm{orbit} of $[\rho]$ is finite and \emph{not} of pullback type, then some Galois conjugate $\rho_\sigma$ of $\rho$ is conjugate to a representation $\rho'\colon\pi_1\Sigma\to\SL_2\R$ which is a universal variation of Hodge structures.
\end{lem}
\begin{proof}
By Lemma~\ref{lem:Zd+non-pullback-implies-values-in-integers}, $\rho$ is valued in some $\SU(h,\mathcal{O}_L)$ because it's not of pullback type and it's a universal variation of Hodge structures. Since $\rho$ is Zariski dense, it has infinite image. So, there exists $\sigma\in \mathrm{Gal}$ such that $h^\sigma$ is of signature $(1,1)$. We can now conjugate $\rho_\sigma$ in order to turn it into a representation $\rho'$ valued in $\SL_2\R$. Being a universal variation of Hodge structures is preserved by Galois conjugation, so $\rho'$ is also a universal variation of Hodge structures.
\end{proof}

\subsubsection{Reducing to DT representations}\label{sec:uVHS-implies-DT}
\cm{The next step consists in proving that among the conjugacy classes of representations $\pi_1\Sigma\to\psl$ with non-trivial elliptic peripheral monodromies, the universal variations of Hodge structures (Definition~\ref{defn:universal-CVHS} are the ones in compact components of relative character varieties. As we explained in Section~\ref{sec:DT-representations}, all compact components in relative character varieties with elliptic peripheral monodromies are either DT components or isolated points which are conjugacy classes of abelian representations (see Definition~\ref{defn:DT-representation} and the discussion beforehand). In particular, Zariski dense representations $\pi_1\Sigma\to\psl$ which are universal variations of Hodge structures must be DT representations. This statement was as already observed by Deroin--Tholozan~\cite[discussion after Theorem~5]{deroin-tholozan}. We give a brief argument here for the sake of completeness.}
\begin{prop}\label{prop:universal-CVHS-implies-DT}
	\sam{If $\rho\colon\pi_1\Sigma\to\psl$ is a universal complex variation of Hodge structures (Definition~\ref{defn:universal-CVHS}) with non-trivial elliptic peripheral monodromies, then it's either valued in a conjugate of $\mathrm{PSO}(2)$ or it's a DT representation.}
\end{prop}
\begin{proof}
\sam{By Definition~\ref{defn:universal-CVHS}, for any choice of complex structure on $\Sigma$, there is a holomorphic or anti-holomoprhic map $f\colon\widetilde \Sigma\to\HH$ which is $\rho$-equivariant. If we assume that $\rho$ is not valued in a conjugate of $\mathrm{PSO}(2)$, then $f$ is non-constant. When $f$ is non-constant, we may pullback the hyperbolic metric on $\HH$ via $f$ to give $\Sigma$ a branched hyperbolic structure with holonomy $\rho$ and cone singularities at the punctures of $\Sigma$.} 

	\sam{Let $d$ denote the number of branch points, counted with multiplicity, of those structures.
	Since we obtain such a branched hyperbolic structure for every choice of complex structure on $\Sigma$ and all have the same holonomy $\rho$, we obtain a first inequality:
	\begin{equation}\label{eq:first-inequality}
		n-3=\mathrm{dim}_{\C}\Teich(\Sigma)\leq d,
	\end{equation}
	where $\Teich(\Sigma)$ is the Teichmüller space of $\Sigma$.
	On the other hand, the Gauss--Bonnet inequality gives
	\begin{equation}\label{eq:second-inequality}
		2-n+d=\chi(\Sigma)+d<0,
	\end{equation}
	where $\chi(\Sigma)=2-n$ is the Euler characteristic of $\Sigma$.
	Combining these two inequalities forces $d=n-3$. In particular, this means that $\rho$ is a DT representation by~\cite[Lemma~4.2]{deroin-tholozan}.}
\end{proof}

\begin{rem}
\sam{Note that we don't need to assume that $\Sigma$ is an $n$-punctured sphere in Proposition~\ref{prop:universal-CVHS-implies-DT}. Indeed, if $\Sigma$ is a general surface of genus $g\geq 0$ with $n$ punctures, then~\eqref{eq:first-inequality} becomes $3g+n-3\leq d$ and~\eqref{eq:second-inequality} becomes $2-2g-n+d<0$. The two inequalities force $g=0$ and $d=n-3$, and the same conclusion applies.}
\end{rem}

We're reaching the point where one should pay attention to the distinction between $\SL_2\R$ and $\psl$, similarly as we did in Section~\ref{sec:finite-orbit-n=4}. Recall that we defined DT representations in Definition~\ref{defn:DT-representation} to be $\psl$-valued representations. We'll write $\iota\colon\SL_2\R\to\psl$ for the quotient map.

\begin{prop}\label{prop:Zd+non-pullback-implies-DT}
Let $\rho\colon\pi_1\Sigma\rightarrow\SL_2\C$ be a Zariski dense representation with only finite order peripheral monodromies and no peripheral loop mapped to $\pm\id$. If the mapping class group orbit of $[\rho]$ is finite and \emph{not} of pullback type, then some Galois conjugate of $\rho$ is conjugate to a representation $\rho'\colon \pi_1\Sigma\to\SL_2\R$ for which $\iota\circ\rho'$ is a DT representation.
\end{prop}
\begin{proof}
Lemma~\ref{lem:Zd+non-pullback-implies-values-in-SL2R} says that some Galois conjugate of $\rho$ is conjugate to a representation $\rho'$ valued in $\SL_2\R$ which is a universal variation of Hodge structures. Because of our assumptions on $\rho$, the peripheral monodromies of $\rho'$ have finite order and are non-trivial. This means that the peripheral monodromies of $\iota\circ\rho'$ are all elliptic. Since $\rho'$ is a universal variation of Hodge structures and has Zariski dense image, it follows from Proposition~\cm{\ref{prop:universal-CVHS-implies-DT}} that $\iota\circ\rho'$ is a DT representation.
\end{proof}

\subsection{Proof of the conjecture}\label{sec:proof-tykhyy-conjecture}
\begin{thm}\label{thm:tykhyy}
Tykhyy's Conjecture (Conjecture~\ref{conj:tykhyy}) is true.
\end{thm}
\begin{proof}
Let's start by fixing a representation $\rho\colon\pi_1\Sigma\to\SL_2\C$ whose conjugacy class $[\rho]$ belongs to a finite mapping class group orbit. As usual, $\Sigma$ denotes a sphere with $n\geq 3$ punctures. We have to prove that $\rho$ belongs to one of the families (1)--(8) from Section~\ref{sec:statement-conjecture}. We already explained in Section~\ref{sec:statement-conjecture} that it's enough to consider the eventuality where $\rho$ is Zariski dense, with only finite order peripheral monodromies, and where the orbit of $[\rho]$ is not of pullback type because the complementary cases have been classified already \cite{cousin-moussard, lll, diarra-pull-back}. \cm{We may also assume that $n\geq 5$ because of the complete classification in the $n=4$ case \cite{LT}.} We can exclude representations of type (8) if we suppose that no peripheral loops of $\Sigma$ is mapped by $\rho$ to $\pm \id$.

With all these extra assumptions on $\rho$, it now satisfies the hypotheses of Lemma~\ref{lem:Zd+non-pullback-implies-values-in-SL2R} and Proposition~\ref{prop:Zd+non-pullback-implies-DT}. So, there exists a representation $\rho'\colon\pi_1\Sigma\to\SL_2\R$ which is conjugate to a Galois conjugate of $\rho$ and becomes a DT representation once post-composed with $\iota\colon\SL_2\R\to\psl$. As we're considering $\SL_2\C$-conjugations (and not only $\SL_2\R$-conjugations), we may assume that the sum of the peripheral angles of $\iota\circ\rho'\colon\pi_1\Sigma\to\psl$ is larger than $2\pi(n-1)$. Theorem~\ref{thm:no-finite-orbit-for-n-geq-7} now applies and forces \cm{$n=5$ or $n=6$}. 

\cm{When $n=6$,} Theorems~\ref{thm:angle-vector-alpha-with-finite-orbits-n=6} \&~\ref{thm:existence-finite-orbit-n=6} further imply that $[\iota\circ\rho']$ belongs to a jester's hat orbit with common elliptic peripheral monodromy given by an angle $\theta\in (5\pi/3,2\pi)$. As we explain after Table~\ref{tab:UFO-orbit} in Appendix~\ref{app:tables}, the orbit point in the jester's hat orbit with action-angle coordinates $(\beta_1,\beta_2,\beta_3)=(2\pi/3,\pi,2\theta-2\pi)$ and $(\gamma_1,\gamma_2)=(0,3\pi/2)$ is the conjugacy class of the representation $\iota \circ\rho_\theta$ of type (4) described in Conjecture~\ref{conj:tykhyy}. In conclusion, the mapping class group orbits of $[\rho_\theta]$ and $[\rho']$ coincide. Since Galois conjugates of $\rho_\theta$ are of the form $\rho_\theta'$ for some other $\theta'$, we conclude that the mapping class group orbit of $[\rho]$ is also of type (4), as desired.

\cm{When $n=5$, Theorem~\ref{thm:angle-vector-alpha-with-finite-orbits-n=5} implies that $[\iota\circ\rho']$ belongs to either a hang-glider orbit, a sand clock orbit, or a bat orbit. The first two correspond to reductions of the jester hat orbit in the sense of (5) in Conjecture~\ref{conj:tykhyy}. The bat orbit corresponds to the exceptional orbit from (5) with elliptic peripheral monodromy and values in $\SL_2\R$, as we explained after Conjecture~\ref{conj:tykhyy}. Like in the case $n=6$, we conclude that the the mapping class group orbit of $[\rho]$ is of type (5), as desired.}
\end{proof}

\bibliographystyle{amsalpha}
\bibliography{references.bib}
\addresseshere

\newpage
\appendix
\setcounter{section}{0}
\renewcommand\sectionname{Appendix}
\addcontentsline{toc}{section}{Appendix}
\renewcommand{\thesection}{\Alph{section}}

\section{A generating family for \texorpdfstring{$\PMod(\Sigma)$}{PMod(S)}}\label{apx:generators-of-pmod}
The pure mapping class group $\PMod(\Sigma)$ of a punctured sphere $\Sigma$ was introduced in \cm{Definition~\ref{defn:pure-mcg}}. As explained in~\cite[Theorem~4.9]{mcg-primer} for instance, the group $\PMod(\Sigma)$ is finitely generated and generators can be taken to be Dehn twists. Explicit generating families can be obtained from presentations of braid groups which closely relate to mapping class groups of punctured spheres. Working with Artin's presentation of braid groups (described for instance in~\cite[p.~251]{mcg-primer}), Ghaswala--Winarski worked out a presentation of $\PMod(\Sigma)$~\cite[Lemma~4.1]{mcg-generators}. It's given as follows. The $n$-punctured sphere $\Sigma$ is homeomorphic to $\C$ minus the real points $1,\ldots, n-1$. For each pair of points $1\leq i<j\leq n-1$, pick a simple closed curve $s_{i,j}$ in $\C$ that loops clockwise around $i$, passes below all the points $i+1,\ldots, j-1$ and loops clockwise around $j$ before closing up in the lower half-plane. The Dehn twist along the curve $s_{i,j}$ is denoted $\sigma_{i,j}\in\PMod(\Sigma)$.
\begin{center}
\begin{tikzpicture}
\fill (1,0) circle (0.07) node[below]{$1$}; 
\fill (3,0) circle (0.07) node[below]{$i$}; 
\fill (5,0) circle (0.07) node[below]{$j$}; 
\fill (7,0) circle (0.07) node[below]{$n-1$}; 

\draw (2,0) node{$\ldots$};
\draw (4,0) node{$\ldots$};
\draw (6,0) node{$\ldots$};

\draw[mauve, thick] (2.5,0) arc(180:0:.5);
\draw[mauve, thick] (4.5,0) arc(180:0:.5);
\draw[-, mauve, thick] (3.5,0) to[out=270, in=270] (4.5,0);
\draw[-, mauve, thick] (2.5,0) to[out=270, in=270] (5.5,0);

\draw[mauve] (4,-1.2) node{$s_{i,j}$};
\end{tikzpicture}
\end{center}

\begin{lem}[\cite{mcg-generators}]\label{lem:GW-mcg-generators}
The group $\PMod(\Sigma)$ is generated by $\{\sigma_{i,j}:1\leq i<j\leq n-1\}$ and the relations are
\begin{equation*}
\begin{cases}
[\sigma_{i,j},\sigma_{k,l}]=1, \quad\forall i<j<k<l,\\
[\sigma_{i,l},\sigma_{j,k}]=1, \quad\forall i<j<k<l,\\
\sigma_{i,k}\sigma_{j,k}\sigma_{i,j}=\sigma_{j,k}\sigma_{i,j}\sigma_{i,k}=\sigma_{i,j}\sigma_{i,k}\sigma_{j,k}, \quad\forall i<j<k,\\
[\sigma_{k,l}\sigma_{i,k}\sigma_{k,l}^{-1},\sigma_{j,l}]=1, \quad\forall i<j<k<l,\\
(\sigma_{1,2}\sigma_{1,3}\cdots\sigma_{1,n-1})\cdots(\sigma_{n-3,n-2}\sigma_{n-3,n-1})\sigma_{n-2,n-1}=1.
\end{cases}
\end{equation*} 
\end{lem}

It's interesting to observe that the first four relations are trivial in the abelianization of $\PMod(\Sigma)$ and that the last relation can be used to get rid of one generator. In other words, the abelianization of $\PMod(\Sigma)$ is the free group on $\binom{n-1}{2}-1$ generators. This shows that the minimal number of generators of $\PMod(\Sigma)$ is $\binom{n-1}{2}-1$ and these can be taken to be Dehn twists. 

Another remarkable observation is that $\PMod(\Sigma)$ is \emph{positively} generated by the Dehn twists $\sigma_{i,j}$. This means that any element of $\PMod(\Sigma)$ can be written as a word in the $\sigma_{i,j}$ where all the exponents are non-negative (no inverses are needed). This is because the last relation can also be used to express the inverse of any of the $\sigma_{i,j}$ as a positive product of the other generators.

For the reasons explained in Section~\ref{sec:action-of-Dehn-twists}, we prefer to work with another generating family. For every $1\leq i<j\leq n-1$, we pick a simple closed curve $t_{i,j}$ which loops around all the points $i,\ldots, j$ in the simplest fashion. The Dehn twist about the curve $t_{i,j}$ is denoted $\tau_{i,j}\in\PMod(\Sigma)$.
\begin{center}
\begin{tikzpicture}
\fill (1,0) circle (0.07) node[below]{$1$}; 
\fill (3,0) circle (0.07) node[below]{$i$}; 
\fill (5,0) circle (0.07) node[below]{$j$}; 
\fill (7,0) circle (0.07) node[below]{$n-1$}; 

\draw (2,0) node{$\ldots$};
\draw (4,0) node{$\ldots$};
\draw (6,0) node{$\ldots$};

\draw[mauve, thick] (2.5,0) arc(180:90:.5);
\draw[mauve, thick] (5,.5) arc(90:0:.5);
\draw[-, mauve, thick] (3,.5) to (5,.5);
\draw[-, mauve, thick] (2.5,0) to[out=270, in=270] (5.5,0);

\draw[mauve] (4,-1.2) node{$t_{i,j}$};
\end{tikzpicture}
\end{center}
Note that since we're working inside the pure mapping class group of $\Sigma$, the Dehn twist $\tau_{1,n-1}$ is trivial.

\begin{lem}\label{lem:mcg-generators}
The group $\PMod(\Sigma)$ is also generated by $\{\tau_{i,j}:1\leq i<j\leq n-1, (i,j)\neq (1,n-1)\}$.
\end{lem}
\begin{proof}
We'll show how to write each of the generators $\sigma_{i,j}$ from Lemma~\ref{lem:GW-mcg-generators} as a word in the Dehn twists $\tau_{k,l}$. We'll proceed by induction on $j-i\geq 1$. First, for the base case $j-i=1$, observe that $\sigma_{i,i+1}=\tau_{i,i+1}$ by definition. Now, assume that we can write each $\sigma_{i,j}$ as a word in the twists $\tau_{k,l}$ whenever $j-i\leq d$. Here, $d$ is some integer with $1\leq d\leq n-3$. If $i<j$ are two indices with $|i-j|=d+1$, then we'll use the identity
\begin{equation}\label{eq:relation-two-types-of-Dehn-twists}
\tau_{i,j}=(\sigma_{i,i+1}\sigma_{i,i+2}\cdots \sigma_{i,j})\cdots(\sigma_{j-2,j-1}\sigma_{j-2,j})\sigma_{j-1,j}.
\end{equation}
The relation~\eqref{eq:relation-two-types-of-Dehn-twists} can be seen at the level of curves directly, as explained in~\cite[p.~250]{mcg-primer}. All the $\sigma_{k,l}$ appearing on the right-hand side of~\eqref{eq:relation-two-types-of-Dehn-twists} have $l-k\leq d$ except for $\sigma_{i,j}$. This expresses $\sigma_{i,j}$ as a product of the twists $\tau_{k,l}$ by the induction hypothesis. 
\end{proof}

We point out that the generating family of Lemma~\ref{lem:mcg-generators} is minimal because it consists of precisely $\binom{n-1}{2}-1$ Dehn twists.

\section{Computing orbit points}\label{apx:algo-orbit}
In Section~\ref{sec:action-of-Dehn-twists}, we presented a method to geometrically compute the image of a triangle chain under a particular kind of Dehn twists. Recall that triangle chains are parametrized by the action-angle coordinates introduced in Section~\ref{sec:action-angle-coordinates}. We just explained in Lemma~\ref{lem:mcg-generators} that this family of Dehn twists generate the pure mapping class group of $\Sigma$. We now describe a procedure to compute the coordinates of all orbit points within a finite mapping class group orbit (Section~\ref{sec:algorithm-to-compute-orbit-points}). The input of the algorithm is any orbit point. It will always terminate if the orbit of the input point is finite. The idea behind it is quite naive and was certainly used in numerous other occasions.

The algorithm, however, involves a lot of computations like those conducted in Example~\ref{ex:example-Dehn-twist-action} which can be quite time consuming, especially when \cm{run} by a human being. It might therefore be useful to be able to approximate the coordinates of every orbit point first (typically by using a computer). In Section~\ref{sec:application-to-mcg}, we describe an alternative to the procedure of Section~\ref{sec:action-of-Dehn-twists} to compute images of orbit points by Dehn twists, which is easier to implement on a computer (the code is available on our GitHub repository). 

\subsection{The abstract algorithm}\label{sec:algorithm-to-compute-orbit-points}
Assume we're in the setting of a group $G$ acting on a space $X$. We wish to compute all the orbit points in the $G$-orbit $\mathcal{O}$ of some element $x\in X$ in the case where we expect $\mathcal{O}$ to be finite. \cm{We also assume that $G$ is generated by a finite collection of elements $\Gamma=\{g_1,\ldots, g_k\}\subset G$ in the sense that every element of $G$ is a product of finitely many $g_i$ and $g_i^{-1}$.} The idea is to construct the orbit points in $\mathcal O$ recursively by successively computing the images by all the generators in $\Gamma$ of previously constructed orbit points. Formally, The procedure has the following steps.

\begin{enumerate}
\setcounter{enumi}{-1}
\item Let $\mathcal O_0=\{x\}$ and $n_0=1$ (representing the cardinality of $\mathcal O_0$).
\item Initiate the step by starting with an empty set $\mathcal O_1$. Going over all $g\in \Gamma$ one after another, compute $gx$. If $gx$ doesn't belong yet to $\mathcal O_0\cup \mathcal O_1$, add it to $\mathcal O_1$. End the step by computing the cardinality $n_1$ of $\mathcal O_1$.
\item[($i+1$)] Assume that we constructed the set $\mathcal O_i$ of cardinality $n_i$ during the previous step. Start with an empty set $\mathcal O_{i+1}$. Going over all $g\in \Gamma$ and all $y\in \mathcal O_i$, compute $gy$. If $gy$ doesn't belong yet to $\mathcal O_0\cup \mathcal O_1\cup\ldots\cup \mathcal O_i$, add it to $\mathcal O_{i+1}$. Once this is done, end the step by computing the cardinality $n_{i+1}$ of $\mathcal O_{i+1}$.
\end{enumerate}

If, after performing the step $i+1$, we observe that $n_{i+1}=0$, the algorithm ends. 

\begin{fact}\label{fact:finite-orbit-when-algo-terminates}
\cm{If the algorithm ends after finitely many steps, it means that the orbit of $x$ was indeed finite. When this happens after the step $i+1$, the orbit $\mathcal O$ of $x$ is equal to $\mathcal O_0\cup \mathcal O_1\cup\cdots\cup \mathcal O_i$ and it has cardinality $1+n_1+\cdots+n_i$.}
\end{fact}
\begin{proof}
\cm{Let $\Gamma^+$ be the set of all positive products of the generators $g_i\in\Gamma$. By assumption, the set $\{gx:g\in \Gamma^+\}$ is finite. In particular, every $g_i\in\Gamma$ has finite order when iterated on $gx$ for all $g\in \Gamma^+$. So, any negative power $g_i^{-a}(gx)$ is actually equal to some positive power $g_i^b(gx)$. In particular, if $h\in G$ is an arbitrary element, then $hx=gx$ for some $g\in \Gamma^+$. This shows that the orbit of $x$ is finite.
}
\end{proof}

\subsection{Application to mapping class group orbits}\label{sec:application-to-mcg}
In our case, $G=\PMod(\Sigma)$ is the pure mapping class group of $\Sigma$ and $X=\RepDT{\alpha}$ is a DT component. The generating family $\Gamma$ of $\PMod(\Sigma)$ that we want to use is the one of Lemma~\ref{lem:mcg-generators}. In order to apply the algorithm of Section~\ref{sec:algorithm-to-compute-orbit-points}, we need a parametrization of the points in $\RepDT{\alpha}$. We'll use the action-angle coordinates described in Section~\ref{sec:action-angle-coordinates} for that. We'll work with a fixed geometric presentation of $\pi_1\Sigma$ with generators $c_1,\ldots,c_n$ and its standard pants decomposition $\mathcal{B}$. Every point in $\RepDT{\alpha}$ is then parametrized by $2(n-3)$ numbers $\beta_1,\ldots,\beta_{n-3}$ and $\gamma_1,\ldots,\gamma_{n-3}$.

In order to run the algorithm from Section~\ref{sec:algorithm-to-compute-orbit-points}, we need to be able to compute the coordinates of the point $\tau.[\rho]$ for $\tau\in\Gamma$ and $[\rho]\in \RepDT{\alpha}$ from the coordinates of $[\rho]$. We could do it by applying the routine described in Section~\ref{sec:action-of-Dehn-twists}, but we'll use a slightly modified version of it that is better suited for computer simulations. Here are the several steps we implement.
\begin{enumerate}
\item We use a hyperbolic geometry package on SageMath to \cm{model} the $\mathcal{B}$-triangle chain of $[\rho]$ in the upper half-plane from its action-angle coordinates. The exterior vertices are $C_1,\ldots,C_n$. We can assume that $\rho\colon\pi_1\Sigma\to\psl$ is the representation that sends $c_i$ to the unique elliptic element of $\psl$ that fixes $C_i$ and has rotation angle $\alpha_i$. There are explicit formulae to write $\rho(c_i)$ as a $2\times 2$ real matrix with determinant 1. Namely, if $C_i$ is the point $x+iy$ in the upper half-plane, then $\rho(c_i)$ is the matrix
\begin{equation*}
\pm\begin{pmatrix}
  \cos(\alpha_i/2)-xy^{-1}\sin (\alpha_i/2) & (x^2y^{-1}+y)\sin(\alpha_i/2)  \\
  -y^{-1}\sin(\alpha_i/2) & \cos(\alpha_i/2)+xy^{-1}\sin (\alpha_i/2)
 \end{pmatrix}.
\end{equation*}
\item We picked $\tau$ as an element of the generating family $\Gamma$ of Lemma~\ref{lem:mcg-generators}. This means that $\tau$ is a Dehn twist along a simple closed curve looping around some consecutive punctures $i,\ldots,j$ of $\Sigma$. We can simply think of $\tau$ as the Dehn twist along the simple curve associated to the fundamental group element $c_i\cdots c_j$. Consider the new representation $\rho'\colon\pi_1\Sigma\to\psl$ defined by
\[
c_k\mapsto\begin{cases}
\rho(c_k) & \text{if } k\notin\{i,\ldots,j\}\\
\rho(c_i\cdots c_j)\rho(c_k)\rho(c_i\cdots c_j)^{-1} &\text{else.}
\end{cases}
\]
The conjugacy class of $\rho'$ coincide with $\tau.[\rho]$. 
\item We are now ready to compute the action coordinates of $\tau.[\rho]$. We start with the action coordinates $\beta_1',\ldots, \beta_{n-3}'$ which are given by the angles of rotation of the elliptic elements $\rho'(c_1c_2)^{-1},\ldots, \rho'(c_1\cdots c_{n-2})^{-1}$. They can be computed using the following formula. If $
\pm\begin{pmatrix}
a & b\\
c& d
\end{pmatrix}$
denotes an elliptic element of $\psl$, then its angle of rotation, seen as a number inside $(0,2\pi)$, is given by
\[
\arctan\left(\frac{-c}{\vert c\vert}\cdot\frac{a+d}{(a+d)^2-2}\sqrt{4-(a+d)^2}\right)+\varepsilon,
\]
where
\begin{equation*}
\varepsilon= \left\{\begin{array}{ll}
0 &\text{if } (a+d)^2>2 \quad\text{and}\quad (a+d)\frac{-c}{\vert c\vert}>0,\\
\pi &\text{if } (a+d)^2<2, \\
2\pi &\text{if } (a+d)^2>2 \quad\text{and}\quad (a+d)\frac{-c}{\vert c\vert}<0.\\
\end{array}\right.
\end{equation*}

\item It remains to compute the angle coordinates $\gamma_1,\ldots,\gamma_{n-3}$ of $\tau.[\rho]$. For that, we use SageMath to draw the $\mathcal{B}$-triangle chain of $\rho'$ and measure each $\gamma_i$. We get the triangle chain by computing the coordinates of all its vertices. The exterior vertices $C_1',\ldots, C_n'$ are the fixed points of $\rho'(c_1),\ldots, \rho'(c_n)$ and the shared vertices $B_1,\ldots, B_{n-3}$ are the fixed points of $\rho'(c_1c_2),\ldots, \rho'(c_1\cdots c_{n-2})$. Recall that the fixed point of an elliptic element $
\pm\begin{pmatrix}
a & b\\
c& d
\end{pmatrix}$ in $\psl$ is given by
\[
\frac{a-d}{2c}+i\cdot \frac{\sqrt{4-(a+d)^2}}{2\vert c\vert}.
\]
The angle coordinate $\gamma_i$ is now given by the angle $\angle C_{i+2}B_iC_{i+1}$ which we can measure using SageMath.
\end{enumerate}

We just explained how to compute (using SageMath) the action-angle coordinates $\beta_1',\ldots, \beta_{n-3}'$ and $\gamma_1',\ldots,\gamma_{n-3}'$ of the point $\tau.[\rho]\in\RepDT{\alpha}$ for an arbitrary generator $\tau\in \Gamma$. This is all we need in order to run the algorithm of Section~\ref{sec:algorithm-to-compute-orbit-points} for $G=\PMod(\Sigma)$ and $X=\RepDT{\alpha}$.

\section{Tables}\label{app:tables}

\subsection{Finite orbits in DT components}
Here's the list of all DT components of 4-punctured spheres that contain a finite mapping class group \cm{orbit}. The third column contains all angle vectors $\alpha\in (0,2\pi)^4$ satisfying $\alpha_1+\alpha_2+\alpha_3+\alpha_4>6\pi$ and for which $\RepDT{\alpha}$ contains a finite orbit. The last column indicates the non-peripheral trace field (Definition~\ref{def:non-peripheral-trace-field}) for exceptional orbits.

\begin{table}[H]
    \centering
    \resizebox{!}{9.5cm}{
    \begin{tblr}{c|c|c|c}
        \makecell{Lisovyy--Tykhyy's\\ numbering} & Orbit length & Angle vector $\alpha$ & \makecell{Non-peripheral\\ trace field}\\
        \hline\hline
        Sol. \Romannum{2} & 2 & $\left\{\theta_1, \theta_1, \theta_2, \theta_2\right\}$, \, $\theta_1+\theta_2>3\pi$ & \\
        \hline
        Sol. \Romannum{3} & 3 & $\left\{\frac{4\pi}{3},2\theta-2\pi, \theta, \theta\right\}$,\, $\theta>5\pi/3$ & \\
        \hline
        Sol. \Romannum{4} & 4 & $\left\{\pi, \theta, \theta, \theta\right\}$,\, $\theta>5\pi/3$ & \\
        \hline
        Sol. \Romannum{4}$^\ast$ & 4 & $\left\{\theta, \theta, \theta, 3\theta-4\pi\right\}$,\, $\theta>5\pi/3$ & \\
        \hline\hline
        Sol. 1 & 5 & $\left\{\frac{22\pi}{15}, \frac{8\pi}{5}, \frac{8\pi}{5}, \frac{28\pi}{15}\right\}$ & $\Q$ \\
        \hline
        Sol. 4 & 6 & $\left\{\frac{19\pi}{12}, \frac{19\pi}{12}, \frac{23\pi}{12}, \frac{23\pi}{12}\right\}$ & $\Q(\sqrt{2})$\\
        \hline
        Sol. 6 & 6 & $\left\{\frac{23\pi}{15}, \frac{23\pi}{15}, \frac{5\pi}{3}, \frac{29\pi}{15}\right\}$ & $\Q(\sqrt{5})$\\
        \hline
        Sol. 7 & 6 & $\left\{\frac{17\pi}{15}, \frac{5\pi}{3}, \frac{29\pi}{15}, \frac{29\pi}{15}\right\}$ & $\Q(\sqrt{5})$\\
        \hline
        Sol. 8 & 7 & $\left\{\frac{10\pi}{7}, \frac{12\pi}{7}, \frac{12\pi}{7}, \frac{12\pi}{7}\right\}$ & $\Q$\\
        \hline
        Sol. 10 & 8 & $\left\{\frac{17\pi}{12}, \frac{7\pi}{4}, \frac{7\pi}{4}, \frac{23\pi}{12}\right\}$ & $\Q(\sqrt{2})$ \\
        \hline
        Sol. 11 & 8 & $\left\{\frac{13\pi}{10}, \frac{3\pi}{2}, \frac{19\pi}{10}, \frac{19\pi}{10}\right\}$  & $\Q(\sqrt{5})$ \\
        \hline
        Sol. 12 & 8 & $\left\{\frac{3\pi}{2}, \frac{17\pi}{10}, \frac{17\pi}{10}, \frac{19\pi}{10}\right\}$  & $\Q(\sqrt{5})$ \\
        \hline
        Sol. 13 & 9 & $\left\{\frac{26\pi}{15}, \frac{26\pi}{15}, \frac{26\pi}{15}, \frac{28\pi}{15}\right\}$ & $\Q(\sqrt{5})$\\
        \hline
        Sol. 14 & 9 & $\left\{\frac{14\pi}{15}, \frac{28\pi}{15}, \frac{28\pi}{15}, \frac{28\pi}{15}\right\}$ & $\Q(\sqrt{5})$ \\
        \hline
        Sol. 15 & 10 & $\left\{\frac{8\pi}{5}, \frac{8\pi}{5}, \frac{9\pi}{5}, \frac{9\pi}{5}\right\}$ & $\Q$\\
        \hline
        Sol. 18 & 10 & $\left\{\frac{23\pi}{15}, \frac{23\pi}{15}, \frac{23\pi}{15}, \frac{9\pi}{5}\right\}$ & $\Q(\sqrt{5})$\\
        \hline
        Sol. 19 & 10 & $\left\{\frac{7\pi}{5}, \frac{29\pi}{15}, \frac{29\pi}{15}, \frac{29\pi}{15}\right\}$ & $\Q(\sqrt{5})$\\
        \hline
        Sol. 20 & 12 & $\left\{\frac{11\pi}{6}, \frac{11\pi}{6}, \frac{11\pi}{6}, \frac{11\pi}{6}\right\}$ & $\Q(\sqrt{2})$\\
        \hline
        Sol. 22 & 12 & $\left\{\frac{19\pi}{15}, \frac{9\pi}{5}, \frac{9\pi}{5}, \frac{29\pi}{15}\right\}$ & $\Q(\sqrt{5})$\\
        \hline
        Sol. 23 & 12 & $\left\{\frac{37\pi}{30}, \frac{47\pi}{30}, \frac{11\pi}{6}, \frac{11\pi}{6}\right\}$ & $\Q(\sqrt{5})$\\
        \hline
        Sol. 24 & 12 & $\left\{\frac{49\pi}{30}, \frac{11\pi}{6}, \frac{11\pi}{6}, \frac{59\pi}{30}\right\}$ & $\Q(\sqrt{5})$\\
        \hline
        Sol. 25 & 12 & $\left\{\frac{43\pi}{30}, \frac{49\pi}{30}, \frac{53\pi}{30}, \frac{59\pi}{30}\right\}$ & $\Q(\sqrt{5})$\\
        \hline
        Sol. 26 & 15 & $\left\{\frac{8\pi}{5}, \frac{26\pi}{15}, \frac{26\pi}{15}, \frac{26\pi}{15}\right\}$ & $\Q(\sqrt{5})$ \\
        \hline
        Sol. 27 & 15 & $\left\{\frac{6\pi}{5}, \frac{28\pi}{15}, \frac{28\pi}{15}, \frac{28\pi}{15}\right\}$ & $\Q(\sqrt{5})$\\
        \hline
        Sol. 30 & 16 & $\left\{\frac{7\pi}{4}, \frac{7\pi}{4}, \frac{7\pi}{4}, \frac{7\pi}{4}\right\}$ & $\Q$\\
        \hline
        Sol. 32 & 18 & $\left\{\frac{37\pi}{21}, \frac{37\pi}{21}, \frac{37\pi}{21}, \frac{41\pi}{21}\right\}$ & $\Q(\cos(\pi/7))$\\
        \hline
        Sol. 33 & 18 & $\left\{\frac{4\pi}{3}, \frac{12\pi}{7}, \frac{12\pi}{7}, \frac{12\pi}{7}\right\}$ & $\Q(\cos(\pi/7))$\\
        \hline
        Sol. 34 & 18 & $\left\{\frac{25\pi}{21}, \frac{41\pi}{21}, \frac{41\pi}{21}, \frac{41\pi}{21}\right\}$ & $\Q(\cos(\pi/7))$\\
        \hline
        Sol. 37 & 20 & $\left\{\frac{47\pi}{30}, \frac{53\pi}{30}, \frac{19\pi}{10}, \frac{19\pi}{10}\right\}$ & $\Q(\sqrt{5})$\\
        \hline
        Sol. 38 & 20 & $\left\{\frac{41\pi}{30}, \frac{17\pi}{10}, \frac{17\pi}{10}, \frac{59\pi}{30}\right\}$ & $\Q(\sqrt{5})$\\
        \hline
        Sol. 39 & 24 & $\left\{\frac{3\pi}{2}, \frac{11\pi}{6}, \frac{11\pi}{6}, \frac{11\pi}{6}\right\}$ & $\Q(\sqrt{5})$\\
        \hline
        Sol. 40 & 30 & $\left\{\frac{23\pi}{15}, \frac{23\pi}{15}, \frac{28\pi}{15}, \frac{28\pi}{15}\right\}$ & $\Q(\sqrt{5})$\\
        \hline
        Sol. 41 & 30 & $\left\{\frac{26\pi}{15}, \frac{26\pi}{15}, \frac{29\pi}{15}, \frac{29\pi}{15}\right\}$ & $\Q(\sqrt{5})$\\
        \hline
        Sol. 43 & 40 & $\left\{\frac{17\pi}{10}, \frac{17\pi}{10}, \frac{17\pi}{10}, \frac{17\pi}{10}\right\}$ & $\Q(\sqrt{5})$\\
        \hline
        Sol. 44 & 40 & $\left\{\frac{19\pi}{10}, \frac{19\pi}{10}, \frac{19\pi}{10}, \frac{19\pi}{10}\right\}$ & $\Q(\sqrt{5})$\\
        \hline
        Sol. 45 & 72 & $\left\{\frac{11\pi}{6}, \frac{11\pi}{6}, \frac{11\pi}{6}, \frac{11\pi}{6}\right\}$ & $\Q(\sqrt{5})$\\
    \end{tblr}
    }
    \caption{The list of finite orbits of DT representations for 4-punctured spheres.}
    \label{tab:finite-mcg-orbits-n=4}
\end{table}

\subsection{Finite orbits}
Below are the tables that contain the action-angles coordinates of all orbit points in the finite orbits for 5-punctured and 6-punctured spheres that we studied in Chapter~\ref{sec:classification}. The tables are meant to be read in the following way.

There is one table for each possible tuple of action coordinates. Each cell corresponds to a possible combination of angle coordinates for the action coordinates associated to the table. We start with a ``basepoint'' $[\rho]\in \RepDT{\alpha}$ which is our preferred orbit point. Inside every other cell you'll find a product of Dehn twists, expressed with the notation of Section~\ref{sec:action-of-Dehn-twists}. When a cell contains a product of Dehn twists defining a mapping class $f\in \PMod(\Sigma)$, it means that $f.[\rho]$ is an orbit point whose action-angle coordinates are those corresponding to the cell containing $f$. The action coordinates marked with a $\clubsuit$ correspond to regular triangle chains. When marked with a $\spadesuit$ instead, it means that one triangle is degenerate (and so one of the angle coordinates is irrelevant). Two triangles are degenerate when the action coordinates are marked with a $\heartsuit$ (two action coordinates are irrelevant), and three triangles are degenerate when the action coordinates are marked with a $\diamondsuit$ (three action coordinates are irrelevant).

\subsubsection{The hang-glider orbit}
The parameters are $n=5$ and $\alpha=(4\pi/3,\theta, \theta, \theta,\theta)$ with $\theta>5\pi/3$. The basepoint $[\rho]$ is the hang-glider triangle chain parametrized by $(\beta_1,\beta_2)=(\pi,4\pi/3)$ and $(\gamma_1,\gamma_2)=(\pi,0)$ (illustrated in Section~\ref{sec:proof-thm-finite-orbits-n=5-1}). The other orbit points are the following.
\begin{table}[H]
    \centering
    \begin{tblr}{colspec={c||c|c|c}, row{4}={6ex}}
        \diaghead(-3,2){\hskip.8cm}%
        {$\gamma_1$}{$\gamma_2$} & $0$ & $\frac{2\pi}{3}$ & $\frac{4\pi}{3}$\\
        \hline\hline
        0 & $\tau_{1,2}$ & $\tau_{2,3}\tau_{1,2}$ & $\tau_{1,2}\tau_{1,3}$\\
        \hline
        $\pi$ & $[\rho]$ & $\tau_{2,3}$ & $\tau_{1,3}$\\
        \SetCell[c=4]{c}{$(\beta_1,\beta_2)=(\pi,\frac{4\pi}{3})^\clubsuit$} & & &\\
    \end{tblr}
\end{table}
\vspace{-2.5em}
\begin{table}[H]
	\centering
    \begin{tblr}{colspec={c||c|c}, row{3}={6ex}}
        $\gamma_1$ & $\frac{\pi}{2}$ & $\frac{3\pi}{2}$\\
        \hline\hline
         & $\tau_{1,3}\tau_{3,4}$ & $\tau_{1,3}\tau_{2,4}$ \\
        \SetCell[c=3]{c}{$(\beta_1,\beta_2)=(\pi,2\theta-2\pi)^\spadesuit$} & &\\
    \end{tblr}
    \begin{tblr}{colspec={c}, row{2}={6ex}}
         $\tau_{2,4}$ \\
        $(\beta_1,\beta_2)=(3\theta-4\pi,2\theta-2\pi)^\heartsuit$\\
    \end{tblr}
    \caption{The list of all 9 orbit points in the hang-glider orbit.}
    \label{tab:hang-glider-orbit}
\end{table}

\subsubsection{The sand clock orbit}
The parameters are $n=5$ and $\alpha=(2\theta-2\pi,\theta, \theta, \theta,\theta)$ with $\theta>5\pi/3$. The basepoint $[\rho]$ is the sand clock triangle chain parametrized by $(\beta_1,\beta_2)=(\pi,4\pi/3)$ and $(\gamma_1,\gamma_2)=(\pi/2,0)$ (illustrated in Section~\ref{sec:proof-thm-finite-orbits-n=5-2}). The other orbit points are the following.

\begin{table}[H]
    \centering
    \begin{tblr}{colspec={c||c|c|c}, row{4}={6ex}}
        \diaghead(-3,2){\hskip.8cm}%
        {$\gamma_1$}{$\gamma_2$} & $0$ & $\frac{2\pi}{3}$ & $\frac{4\pi}{3}$\\
        \hline\hline
        $\frac{\pi}{2}$ & $[\rho]$ & $(\tau_{1,3})^2$ & $\tau_{1,3}$\\
        \hline
        $\frac{3\pi}{2}$ & $\tau_{1,2}$ & $(\tau_{2,3})^2$ & $\tau_{1,2}\tau_{1,3}$\\
        \SetCell[c=4]{c}{$(\beta_1,\beta_2)=(\pi,\frac{4\pi}{3})^\clubsuit$} & & &\\
    \end{tblr}
\end{table}
\vspace{-2.5em}
\begin{table}[H]
    \centering
    \begin{tblr}{colspec={c||c|c|c}, row{3}={6ex}}
        $\gamma_2$ & $\frac{\pi}{3}$ & $\pi$ & $\frac{5\pi}{3}$\\
        \hline\hline
         & $\tau_{2,3}$ & $(\tau_{1,3})^2\tau_{2,3}$ & $\tau_{1,3}\tau_{2,3}$ \\
        \SetCell[c=4]{c}{$(\beta_1,\beta_2)=(6\pi-3\theta,\frac{4\pi}{3})^\spadesuit$} & &\\
    \end{tblr}
    \begin{tblr}{colspec={c||c|c}, row{3}={6ex}}
        $\gamma_1$ & $0$ & $\pi$ \\
        \hline\hline
         & $\tau_{1,3}\tau_{3,4}$ & $\tau_{1,2}\tau_{1,3}\tau_{3,4}$ \\
        \SetCell[c=3]{c}{$(\beta_1,\beta_2)=(\pi,2\theta-2\pi)^\spadesuit$} & &\\
    \end{tblr}
\end{table}
\vspace{-2.5em}
\begin{table}[H]
    \centering
    \begin{tblr}{colspec={c}, row{2}={6ex}}
         $\tau_{2,3}\tau_{3,4}$ \\
        $(\beta_1,\beta_2)=(6\pi-3\theta,3\theta-4\pi)^\heartsuit$\\
    \end{tblr}
    \caption{The list of all 12 orbit points in the sand clock orbit.}
    \label{tab:sand-clock-orbit}
\end{table}

\subsubsection{The bat orbit}
The parameters are $n=5$ and $\alpha=(12\pi/7,12\pi/7, 12\pi/7, 12\pi/7,12\pi/7)$ with $\theta>5\pi/3$. The basepoint $[\rho]$ is the bat triangle chain parametrized by $(\beta_1,\beta_2)=(2\pi/3,8\pi/7)$ and $(\gamma_1,\gamma_2)=(\pi/3,4\pi/7)$ (illustrated in Section~\ref{sec:proof-thm-finite-orbits-n=5-3}). The other orbit points are the following.

\begin{table}[H]
    \centering
    \resizebox{\columnwidth}{!}{\begin{tblr}{colspec={c||c|c|c|c|c|c|c}, row{5}={6ex}}
        \diaghead(-3,2){\hskip.8cm}%
        {$\gamma_1$}{$\gamma_2$} & $0$ & $\frac{2\pi}{7}$ & $\frac{4\pi}{7}$ & $\frac{6\pi}{7}$ & $\frac{8\pi}{7}$ & $\frac{10\pi}{7}$ & $\frac{12\pi}{7}$ \\
        \hline\hline
        $\frac{\pi}{3}$ & $(\tau_{1,3})^3$ & $(\tau_{1,2})^2(\tau_{2,3})^2$ & $[\rho]$ & $(\tau_{1,3})^2$ & $\tau_{1,2}\tau_{3,4}(\tau_{2,3})^2$ & $(\tau_{1,2})^2\tau_{3,4}$ & $\tau_{1,3}$ \\
        \hline
        $\pi$ & $\tau_{1,2}(\tau_{1,3})^3$ & $\tau_{3,4}\tau_{2,3}$ & $\tau_{1,2}$ & $\tau_{1,2}(\tau_{1,3})^2$ & $\tau_{3,4}(\tau_{2,3})^2$ & $\tau_{3,4}$ & $\tau_{1,2}\tau_{1,3}$ \\
        \hline
        $\frac{5\pi}{3}$ & $\tau_{1,3}\tau_{2,3}$ & $\tau_{3,4}\tau_{2,3}\tau_{1,2}$ & $(\tau_{1,2})^2$ & $\tau_{2,3}$ & $(\tau_{1,3})^2\tau_{2,3}$ & $\tau_{1,2}\tau_{3,4}$ & $(\tau_{1,2})^2\tau_{1,3}$ \\
        \SetCell[c=8]{c}{$(\beta_1,\beta_2)=(\frac{2\pi}{3},\frac{8\pi}{7})^\clubsuit$} & & & & & & &\\
    \end{tblr}}
\end{table}
\vspace{-2em}
\begin{table}[H]
    \centering
    \resizebox{\columnwidth}{!}{\begin{tblr}{colspec={c||c|c|c|c|c|c|c}, row{5}={6ex}}
        \diaghead(-3,2){\hskip.8cm}%
        {$\gamma_2$}{$\gamma_1$} & $0$ & $\frac{2\pi}{7}$ & $\frac{4\pi}{7}$ & $\frac{6\pi}{7}$ & $\frac{8\pi}{7}$ & $\frac{10\pi}{7}$ & $\frac{12\pi}{7}$ \\
        \hline\hline
        $\frac{\pi}{3}$ & $\tau_{2,4}(\tau_{1,2})^2\tau_{1,3}$ & $\tau_{2,4}\tau_{1,3}$ & $\tau_{1,3}\tau_{2,4}\tau_{1,2}$ & $(\tau_{1,3})^2\tau_{3,4}\tau_{2,3}$ & $\tau_{2,4}\tau_{1,2}\tau_{1,3}$ & $\tau_{2,4}\tau_{3,4}\tau_{1,3}$ & $\tau_{1,3}\tau_{2,4}$ \\
        \hline
        $\pi$ & $\tau_{2,4}(\tau_{1,2})^2$ & $\tau_{2,4}$ & $\tau_{1,2}(\tau_{2,4})^2$ & $(\tau_{2,3})^2\tau_{3,4}\tau_{1,3}$ & $\tau_{2,4}\tau_{1,2}$ & $\tau_{2,4}\tau_{3,4}$ & $\tau_{1,2}(\tau_{2,4})^2\tau_{3,4}$ \\
        \hline
        $\frac{5\pi}{3}$ & $\tau_{1,3}(\tau_{2,3})^2\tau_{3,4}$ & $\tau_{2,4}(\tau_{1,3})^2$ & $(\tau_{1,2})^2\tau_{2,3}\tau_{3,4}$ & $(\tau_{2,3})^2\tau_{3,4}$ & $\tau_{2,4}\tau_{1,2}(\tau_{1,3})^2$ & $\tau_{2,3}\tau_{3,4}\tau_{2,3}$ & $\tau_{1,3}\tau_{2,4}\tau_{1,3}$ \\
        \SetCell[c=8]{c}{$(\beta_1,\beta_2)=(\frac{6\pi}{7},\frac{4\pi}{3})^\clubsuit$} & & & & & & &\\
    \end{tblr}}
\end{table}
\vspace{-2em}
\begin{table}[H]
    \centering
    \resizebox{\columnwidth}{!}{\begin{tblr}{colspec={c||c|c|c|c|c|c}, row{8}={6ex}}
        \diaghead(-3,2){\hskip.8cm}%
        {$\gamma_2$}{$\gamma_1$} & $\frac{\pi}{3}-\gamma_0$ & $\frac{\pi}{3}+\gamma_0$ & $\pi-\gamma_0$ & $\pi+\gamma_0$ & $\frac{5\pi}{3}-\gamma_0$ & $\frac{5\pi}{3}+\gamma_0$\\
        \hline\hline
        $\frac{\pi}{3}-\gamma_0$ & $\tau_{1,3}\tau_{3,4}$ & & $\tau_{1,2}\tau_{1,3}\tau_{3,4}$ & & $\tau_{2,3}\tau_{2,4}\tau_{1,3}$ & \\
        \hline
        $\frac{\pi}{3}+\gamma_0$ & & $\tau_{1,2}\tau_{2,4}$ & & $\tau_{1,2}\tau_{2,4}\tau_{1,2}$ & & $\tau_{2,3}\tau_{2,4}\tau_{1,3}$ \\
        \hline
        $\pi-\gamma_0$ & $(\tau_{1,2})^2\tau_{2,4}\tau_{1,3}$ & & $(\tau_{1,3})^2\tau_{2,3}\tau_{3,4}$ & & $\tau_{2,3}\tau_{2,4}$ \\
        \hline
        $\pi+\gamma_0$ & & $\tau_{2,3}\tau_{3,4}$ & & $(\tau_{1,3})^2\tau_{3,4}$ & & $\tau_{2,3}\tau_{2,4}$ \\
        \hline
        $\frac{5\pi}{3}-\gamma_0$ & $(\tau_{1,2})^2\tau_{2,4}$ & & $(\tau_{1,2})^2\tau_{2,4}\tau_{1,2}$ & & $(\tau_{2,4})^2\tau_{3,4}$ & \\
        \hline
        $\frac{5\pi}{3}+\gamma_0$ & & $\tau_{1,2}\tau_{2,4}\tau_{1,3}$ & & $\tau_{1,2}\tau_{2,4}\tau_{1,2}\tau_{1,3}$ & & $\tau_{2,3}\tau_{2,4}(\tau_{1,3})^2$ \\
        \SetCell[c=7]{c}{$(\beta_1,\beta_2)=(\frac{2\pi}{3},\frac{4\pi}{3})^\clubsuit$} & & & & & &\\
    \end{tblr}}
\end{table}
\vspace{-2em}
\begin{table}[H]
    \centering
    \begin{tblr}{colspec={c||c|c|c|c}, row{6}={6ex}}
        \diaghead(-3,2){\hskip.8cm}%
        {$\gamma_1$}{$\gamma_2$} & $\frac{\pi}{4}$ & $\frac{3\pi}{4}$ & $\frac{5\pi}{4}$ & $\frac {7\pi}{4}$\\
        \hline\hline
        $0$ & $\tau_{1,3}(\tau_{3,4})^2$ & $(\tau_{2,4})^2(\tau_{1,2})^2$ & $\tau_{1,2}\tau_{2,4}\tau_{1,2}\tau_{3,4}$ & $\tau_{3,4}\tau_{2,3}\tau_{1,2}\tau_{3,4}$ \\
        \hline
        $\frac{2\pi}{3}$ & $\tau_{1,2}\tau_{1,3}(\tau_{3,4})^2$ & $(\tau_{2,4})^2$ & $\tau_{2,3}\tau_{2,4}\tau_{1,3}\tau_{3,4}$ & $(\tau_{2,4})^2\tau_{1,3}$ \\
        \hline
        $\frac{4\pi}{3}$ & $\tau_{1,2}\tau_{2,4}\tau_{3,4}\tau_{1,3}$ & $(\tau_{2,4})^2\tau_{1,2}$ & $\tau_{1,2}\tau_{2,4}\tau_{3,4}$ & $\tau_{3,4}\tau_{2,3}\tau_{3,4}$ \\
        \SetCell[c=5]{c}{$(\beta_1,\beta_2)=(\frac{2\pi}{3},\pi)^\clubsuit$} & & & &\\
    \end{tblr}
\end{table}
\vspace{-2em}
\begin{table}[H]
    \centering
    \begin{tblr}{colspec={c||c|c|c|c}, row{6}={6ex}}
        \diaghead(-3,2){\hskip.8cm}%
        {$\gamma_2$}{$\gamma_1$} & $\frac{\pi}{4}$ & $\frac{3\pi}{4}$ & $\frac{5\pi}{4}$ & $\frac {7\pi}{4}$\\
        \hline\hline
        $0$ & $\tau_{1,3}\tau_{2,3}\tau_{2,4}$ & $\tau_{1,2}(\tau_{1,3})^3\tau_{2,4}$ & $\tau_{1,3}\tau_{2,3}\tau_{2,4}\tau_{1,2}$ & $(\tau_{1,3})^3\tau_{2,4}$ \\
        \hline
        $\frac{2\pi}{3}$ & $\tau_{1,2}\tau_{1,3}\tau_{2,4}\tau_{1,2}\tau_{1,3}$ & $\tau_{1,3}\tau_{2,4}\tau_{1,2}\tau_{2,4}$ & $\tau_{1,2}\tau_{1,3}\tau_{2,4}\tau_{1,3}$ & $(\tau_{2,3})^2\tau_{2,4}\tau_{1,3}$ \\
        \hline
        $\frac{4\pi}{3}$ & $\tau_{1,2}\tau_{1,3}\tau_{2,4}\tau_{1,2}$ & $(\tau_{2,3})^2\tau_{2,4}\tau_{1,2}$ & $\tau_{1,2}\tau_{1,3}\tau_{2,4}$ & $(\tau_{2,3})^2\tau_{2,4}$ \\
        \SetCell[c=5]{c}{$(\beta_1,\beta_2)=(\pi,\frac{4\pi}{3})^\clubsuit$} & & & &\\
    \end{tblr}
\end{table}
\vspace{-2em}
\begin{table}[H]
	\centering
    \resizebox{\columnwidth}{!}{\begin{tblr}{colspec={c||c|c|c|c|c|c|c}, row{3}={6ex}}
        $\gamma_2$ & $0$ & $\frac{2\pi}{7}$ & $\frac{4\pi}{7}$ & $\frac {6\pi}{7}$ & $\frac {8\pi}{7}$ & $\frac {10\pi}{7}$ & $\frac {12\pi}{7}$\\
        \hline\hline
         & $\tau_{1,2}\tau_{3,4}\tau_{2,3}\tau_{1,2}$ & $(\tau_{1,2})^2\tau_{2,3}$ & $(\tau_{2,3})^2$ & $(\tau_{1,3})^2(\tau_{2,3})^2$ & $\tau_{1,2}\tau_{3,4}\tau_{2,3}$ & $(\tau_{2,3})^2\tau_{1,2}$ & $\tau_{1,3}(\tau_{2,3})^2$ \\
        \SetCell[c=8]{c}{$(\beta_1,\beta_2)=(\frac{6\pi}{7},\frac{8\pi}{7})^\spadesuit$} & & & &\\
    \end{tblr}}
\end{table}
\vspace{-2em}
\begin{table}[H]
    \centering
    \resizebox{\columnwidth}{!}{\begin{tblr}{colspec={c||c|c|c|c}, row{3}={6ex}}
        $\gamma_2$ & $\frac{\pi}{4}$ & $\frac{3\pi}{4}$ & $\frac{5\pi}{4}$ & $\frac {7\pi}{4}$\\
        \hline\hline
         & $\tau_{2,3}(\tau_{2,4})^2\tau_{1,3}$ & $\tau_{1,2}\tau_{1,3}(\tau_{3,4})^2\tau_{2,3}$ & $\tau_{2,3}(\tau_{2,4})^2$ & $\tau_{2,3}\tau_{2,4}\tau_{1,3}\tau_{3,4}\tau_{2,3}$ \\
        \SetCell[c=5]{c}{$(\beta_1,\beta_2)=(\frac{4\pi}{7},\pi)^\spadesuit$} & & & &\\
    \end{tblr}}
\end{table}
\vspace{-2em}
\begin{table}[H]
    \centering
    \begin{tblr}{colspec={c||c|c|c|c}, row{3}={6ex}}
        $\gamma_1$ & $\frac{\pi}{4}$ & $\frac{3\pi}{4}$ & $\frac{5\pi}{4}$ & $\frac {7\pi}{4}$\\
        \hline\hline
         & $(\tau_{2,3})^2\tau_{2,4}\tau_{1,2}\tau_{3,4}$ & $\tau_{1,2}\tau_{1,3}\tau_{2,4}\tau_{3,4}$ & $(\tau_{2,3})^2\tau_{2,4}\tau_{3,4}$ & $(\tau_{2,4})^2\tau_{1,3}\tau_{2,4}$ \\
        \SetCell[c=5]{c}{$(\beta_1,\beta_2)=(\pi,\frac{10\pi}{7})^\spadesuit$} & & & &\\
    \end{tblr}
\end{table}
\vspace{-2em}
\begin{table}[H]
    \centering
    \begin{tblr}{colspec={c||c|c|c}, row{3}={6ex}}
        $\gamma_1$ & $0$ & $\frac{2\pi}{3}$ & $\frac{4\pi}{3}$\\
        \hline\hline
         & $\tau_{2,4}\tau_{1,2}\tau_{2,4}\tau_{1,2}$ & $\tau_{1,2}\tau_{2,4}\tau_{1,2}\tau_{2,4}$ & $\tau_{2,4}\tau_{1,2}\tau_{2,4}$\\
        \SetCell[c=5]{c}{$(\beta_1,\beta_2)=(\frac{2\pi}{3},\frac{10\pi}{7})^\spadesuit$} & & & &\\
    \end{tblr}
    \begin{tblr}{colspec={c||c|c|c}, row{3}={6ex}}
        $\gamma_2$ & $0$ & $\frac{2\pi}{3}$ & $\frac{4\pi}{3}$\\
        \hline\hline
         & $\tau_{1,2}(\tau_{1,3})^2\tau_{2,4}\tau_{1,3}$ & $\tau_{1,2}(\tau_{1,3})^2\tau_{2,4}$ & $\tau_{1,3}\tau_{2,4}\tau_{1,2}\tau_{2,3}\tau_{2,4}$\\
        \SetCell[c=5]{c}{$(\beta_1,\beta_2)=(\frac{4\pi}{7},\frac{4\pi}{3})^\spadesuit$} & & & &\\
    \end{tblr}
    \caption{The list of all 105 orbit points in the bat orbit.}
    \label{tab:bat-orbit}
\end{table}

\subsubsection{The jester's hat orbit}
The parameters are $n=6$ and $\alpha=(\theta,\theta,\theta,\theta,\theta,\theta)$ with $\theta>5\pi/3$. The basepoint $[\rho]$ is the jester's hat triangle chain parametrized by $(\beta_1,\beta_2,\beta_3)=(2\pi/3,\pi,4\pi/3)$ and $(\gamma_1,\gamma_2,\gamma_3)=(2\pi/3,0,2\pi/3)$ (illustrated in the proof of Theorem~\ref{thm:existence-finite-orbit-n=6}). The other orbit points are the following.

\begin{table}[H]
    \centering
    \vspace{.3cm}
    \begin{tblr}{colspec={c||c|c|c|c|c|c}, row{6}={6ex}}
        \diaghead(-3,2){\hskip.8cm}%
        {$\gamma_1$}{$\gamma_3$} & \SetCell[c=2]{c}{0} & & \SetCell[c=2]{c}{$\frac{2\pi}{3}$} & & \SetCell[c=2]{c}{$\frac{4\pi}{3}$} &\\
        \hline\hline
        0 & $(\tau_{1,2})^2\tau_{1,4}$ & $\tau_{4,5}\tau_{3,4}$ & $(\tau_{1,2})^2$ & $(\tau_{1,2})^2\tau_{1,3}$ & $\tau_{1,2}\tau_{3,5}$ & $\tau_{1,3}\tau_{3,4}$\\
        \hline
        $\frac{2\pi}{3}$ & $\tau_{1,4}$ & $\tau_{1,3}\tau_{1,4}$ & $[\rho]$ & $\tau_{1,3}$ & $\tau_{1,3}\tau_{4,5}$ & $\tau_{4,5}$\\
        \hline
        $\frac{4\pi}{3}$ & $\tau_{1,2}\tau_{1,4}$ & $\tau_{1,2}\tau_{1,3}\tau_{1,4}$ & $\tau_{1,2}$ &  $\tau_{1,2}\tau_{1,3}$ & $\tau_{2,4}$ & $\tau_{1,2}\tau_{4,5}$\\
        \hline\hline
        $\gamma_2$ & 0 & $\pi$ & 0 & $\pi$ & 0 & $\pi$\\
        \SetCell[c=7]{c}{$(\beta_1,\beta_2,\beta_3)=(\frac{2\pi}{3},\pi,\frac{4\pi}{3})^\clubsuit$} & & & & & &\\
    \end{tblr}
\end{table}
\vspace{-2.5em}
\begin{table}[H]
	\centering
    \begin{tblr}{colspec={c||c|c|c}, row{4}={6ex}}
        \diaghead(-3,2){\hskip.8cm}%
        {$\gamma_2$}{$\gamma_3$} & $0$ & $\frac{2\pi}{3}$ & $\frac{4\pi}{3}$\\
        \hline\hline
        $\frac{\pi}{2}$ & $\tau_{1,3}\tau_{2,4}$ & $\tau_{2,3}$ & $\tau_{1,3}\tau_{1,4}\tau_{2,4}$\\
        \hline
        $\frac{3\pi}{2}$ & $\tau_{1,2}\tau_{4,5}\tau_{2,4}$ & $\tau_{1,3}\tau_{2,3}$ & $\tau_{2,3}\tau_{4,5}$\\
        \SetCell[c=4]{c}{$(\beta_1,\beta_2,\beta_3)=(4\pi-2\theta,\pi,\frac{4\pi}{3})^\spadesuit$} & & & \\
    \end{tblr}
    \begin{tblr}{colspec={c||c|c|c}, row{4}={6ex}}
        \diaghead(-3,2){\hskip.8cm}%
        {$\gamma_2$}{$\gamma_1$} & $0$ & $\frac{2\pi}{3}$ & $\frac{4\pi}{3}$\\
        \hline\hline
        $\frac{\pi}{2}$ & $\tau_{1,3}\tau_{3,5}$ & $(\tau_{4,5})^2$ & $\tau_{1,3}\tau_{2,5}$\\
        \hline
        $\frac{3\pi}{2}$ & $\tau_{1,2}\tau_{3,5}\tau_{4,5}$ & $\tau_{1,3}\tau_{2,3}\tau_{2,5}$ & $\tau_{2,4}\tau_{4,5}$\\
        \SetCell[c=4]{c}{$(\beta_1,\beta_2,\beta_3)=(\frac{2\pi}{3},\pi,2\theta-2\pi)^\spadesuit$} & & & \\
    \end{tblr}
\end{table}
\vspace{-2.5em}
\begin{table}[H]
	\centering
    \begin{tblr}{colspec={c||c|c|c}, row{3}={6ex}}
        $\gamma_1$ & $\frac{\pi}{3}$ & $\pi$ & $\frac{5\pi}{3}$\\
        \hline\hline
         & $\tau_{1,3}\tau_{2,4}\tau_{2,5}$ & $\tau_{2,4}\tau_{4,5}\tau_{3,4}$ & $\tau_{1,3}\tau_{1,4}\tau_{3,5}$\\
         \SetCell[c=4]{c}{$(\beta_1,\beta_2,\beta_3)=(\frac{2\pi}{3},3\theta-4\pi,2\theta-2\pi)^\heartsuit$} & & & \\
    \end{tblr}
    \begin{tblr}{colspec={c||c|c|c}, row{3}={6ex}}
        $\gamma_3$ & $\frac{\pi}{3}$ & $\pi$ & $\frac{5\pi}{3}$\\
        \hline\hline
         & $\tau_{1,3}\tau_{2,4}\tau_{3,4}$ & $\tau_{2,3}\tau_{3,4}$ & $\tau_{1,3}\tau_{2,3}\tau_{3,5}$\\
         \SetCell[c=4]{c}{$(\beta_1,\beta_2,\beta_3)=(4\pi-2\theta,6\pi-3\theta,\frac{4\pi}{3})^\heartsuit$} & & & \\
    \end{tblr}
    \begin{tblr}{colspec={c||c|c}, row{3}={6ex}}
        $\gamma_2$ & $0$ & $\pi$ \\
        \hline\hline
         & $\tau_{1,2}\tau_{2,5}$ & $\tau_{1,2}\tau_{2,5}\tau_{1,3}$\\
         \SetCell[c=3]{c}{$(\beta_1,\beta_2,\beta_3)=(4\pi-2\theta,\pi,2\theta-2\pi)^\heartsuit$} & & \\
    \end{tblr}
\end{table}
\vspace{-2.5em}
\begin{table}[H]
	\centering
    \begin{tblr}{colspec={c}, row{2}={6ex}}
        $\tau_{1,2}\tau_{2,5}\tau_{3,5}$ \\
        $(\beta_1,\beta_2,\beta_3)=(4\pi-2\theta,6\pi-3\theta,8\pi-4\theta)^\diamondsuit$ \\
    \end{tblr}
    \begin{tblr}{colspec={c}, row{2}={6ex}}
        $\tau_{1,4}\tau_{2,5}$ \\
        $(\beta_1,\beta_2,\beta_3)=(4\theta-6\pi,3\theta-4\pi,2\theta-2\pi)^\diamondsuit$ \\
    \end{tblr}
    \caption{The list of all 40 orbit points in the jester's hat orbit.}
    \label{tab:UFO-orbit}
\end{table}

One point in the jester's hat orbit played a special role in the proof of Theorem~\ref{thm:tykhyy}: the one with coordinates $(\beta_1,\beta_2,\beta_3)=(2\pi/3,\pi,2\theta-2\pi)$ and $(\gamma_1,\gamma_2)=(0,3\pi/2)$. This point is the conjugacy class of the representation $\rho\colon\pi_1\Sigma\to\psl$ defined by
\[
\rho(c_1)=\rho(c_4)=\pm\begin{pmatrix}
\cos(\theta/2) & \frac{1+\sqrt{-1+2\cos(\theta)}}{2}\\
\frac{-1+\sqrt{-1+2\cos(\theta)}}{2} & \cos(\theta/2)
\end{pmatrix}
\]
and
\[
\rho(c_2)=\rho(c_3)=\rho(c_5)=\rho(c_6)=\pm\begin{pmatrix}
\cos(\theta/2) & \frac{1-\sqrt{-1+2\cos(\theta)}}{2}\\
\frac{-1-\sqrt{-1+2\cos(\theta)}}{2} & \cos(\theta/2)
\end{pmatrix}.
\]
The representation $\rho$ is conjugate to the representation $\rho_\theta$ of type (4) in Conjecture~\ref{conj:tykhyy}.

\end{document}